\documentclass[12pt]{article}%
\usepackage{graphicx}
\usepackage{amsmath}%
\usepackage{amsfonts}%
\usepackage{amssymb}
\newtheorem{theorem}{Theorem}[section]

\newtheorem{claim}[theorem]{Claim}

\newtheorem{corollary}[theorem]{Corollary}

\newtheorem{definition}[theorem]{Definition}
\newtheorem{example}[theorem]{Example}

\newtheorem{lemma}[theorem]{Lemma}

\newtheorem{proposition}[theorem]{Proposition}
\newtheorem{remark}[theorem]{Remark}

\newenvironment{proof}[1][Proof]{\textbf{#1.} }{\ \rule{0.5em}{0.5em}}

\begin{document}

\author{(Preliminary Version \#2)
\and Peter Scott\thanks{Partially supported by NSF grants DMS 034681 and 9626537}\\Mathematics Department, University of Michigan\\Ann Arbor, Michigan 48109, USA.\\email:pscott@math.lsa.umich.edu
\and Gadde A. Swarup\\Department of Mathematics and Statistics\\University of Melbourne\\Victoria 3010, Australia.\\email:gadde@ms.unimelb.edu.au }
\title{Regular Neighbourhoods and Canonical Decompositions for Groups}
\maketitle

\begin{abstract}
We find canonical decompositions for finitely presented groups which
essentially specialise to the classical JSJ-decomposition when restricted to
the fundamental groups of Haken manifolds. The decompositions that we obtain
are invariant under automorphisms of the group. A crucial new ingredient is
the concept of a regular neighbourhood of a family of almost invariant subsets
of a group. An almost invariant set is an analogue of an immersion.

\end{abstract}
\date{}

\newpage

\tableofcontents

\newpage

\textbf{{\Large Introduction}}

\medskip

This article is devoted to the study of analogues for groups of the classical
JSJ-decomposition (see Jaco and Shalen \cite{JS}, Johannson \cite{JO} and
Waldhausen \cite{Waldhausen}) for orientable Haken $3$-manifolds. The
orientability restriction is not essential but it will simplify our
discussions. This field was initiated by Kropholler \cite{K} who studied
analogous decompositions for Poincar\'{e} duality groups of any dimension
greater than $2$. But the current interest in this kind of decomposition
started with the work of Sela \cite{S1} on one-ended torsion-free hyperbolic
groups. His results were generalised by Rips and Sela \cite{RS}, Bowditch
\cite{B1}\cite{B2}, Dunwoody and Sageev \cite{D-Sageev}, and Fujiwara and
Papasoglu \cite{FP}, but none of these results yields the classical
JSJ-decomposition when restricted to the fundamental group of an orientable
Haken manifold. In this paper, we give a new approach to this subject, and we
give decompositions for finitely presented groups which essentially specialise
to the classical JSJ-decomposition when restricted to the fundamental groups
of Haken manifolds. All of the results of \cite{K}, \cite{S1}, \cite{RS},
\cite{B1}, \cite{B2}, and \cite{D-Sageev} for virtually polycyclic groups
follow from the results in this paper. We also indicate some extensions of the
results on the Algebraic Annulus and Torus Theorems in \cite{SS2}, \cite{B1}
and \cite{D-Swenson}. It should be remarked that even though we obtain
canonical decompositions for all finitely presented groups, most of the time
these decompositions are trivial. This is analogous to the fact that any
finitely generated group possesses a free product decomposition, but this
decomposition is trivial whenever the given group is freely indecomposable. We
should also remark that many of the ideas in this paper and the above
mentioned papers can be traced back to the groundbreaking work of Stallings on
groups with infinitely many ends \cite{Stallings1}\cite{Stallings2}.

We will focus on extending what we consider to be the most important aspects
of the topological JSJ-decomposition. Our choice of the crucial property of
the classical JSJ-decomposition is the Enclosing Property of the
characteristic submanifold which may be described briefly as follows. See
section \ref{charsub} for a more detailed discussion. For an orientable Haken
$3$-manifold $M$, Jaco and Shalen \cite{JS} and Johannson \cite{JO} proved
that there is a family $\mathcal{T}$ of disjoint essential annuli and tori
embedded in $M$, unique up to isotopy, and with the following properties. The
manifolds obtained by cutting $M$ along $\mathcal{T}$ are simple or are
Seifert fibre spaces or $I$-bundles over surfaces. The Seifert and $I$-bundle
pieces of $M$ are said to be \textit{characteristic}, and any essential map of
the annulus or torus into $M$ can be properly homotoped to lie in a
characteristic piece. This is called the Enclosing Property of $\mathcal{T}$.
The characteristic submanifold $V(M)$ of $M$ consists essentially of the union
of the characteristic pieces of the manifold obtained from $M$ by cutting
along $\mathcal{T}$.

The existence of the JSJ-decomposition of an orientable Haken $3$-manifold $M$
implies that the fundamental group $G$ of $M$ is the fundamental group of a
graph $\Gamma$ of groups whose edge groups are all isomorphic to $\mathbb{Z}$
or $\mathbb{Z}\times\mathbb{Z}$, and whose vertex groups are the fundamental
groups of simple manifolds or of Seifert fibre spaces or of surfaces. The
uniqueness up to isotopy of the splitting family $\mathcal{T}$ implies that
$\Gamma$ is unique. Further, the Enclosing Property implies that any subgroup
of $G$ which is represented by an essential annulus or torus in $M$ is
conjugate into a characteristic vertex group. All the previous algebraic
analogues of the topological JSJ-decomposition consist of producing a graph of
groups structure $\Gamma$ for a given group $G$ with the edge groups of
$\Gamma$ being of some specified type and with some ``characteristic''
vertices. The algebraic analogue of the topological Enclosing Property is the
property that certain ``essential'' subgroups of $G$ must be conjugate into
one of the characteristic vertex groups of $\Gamma$, but the precise meaning
of the word ``essential'' varies depending on the authors. In all of these
cases, an essential subgroup of $G$ is of the same abstract type as the edge
groups of $\Gamma$. The first such result was by Kropholler \cite{K}, who
considered the special case when $G\;$is a Poincar\'{e} duality group of
dimension $n$ and the edge groups of $\Gamma$ are virtually polycyclic (VPC)
groups of Hirsch length $n-1$. (For brevity, we will refer to the length
rather than the Hirsch length of a VPC group throughout this paper.) In his
case, a VPC subgroup $H$ of length $n-1$ is called essential if $e(G,H)\geq2$.
This corresponds to considering all $\pi_{1}$-injective maps of closed
$(n-1)$-dimensional manifolds into a $n$-manifold rather than considering just
embeddings of such manifolds. Note that a VPC group of length at most $2$ is
virtually abelian, so that when $n=3$ his result is closely related to the
topology of $3$-manifolds. In fact, Kropholler \cite{K1} used his results in
\cite{K} to give a new proof of the existence of the JSJ-decomposition for
closed $3$-manifolds.

In most of the papers which came after \cite{K}, a subgroup $H$ of $G$ is
called essential if $G$ possesses a splitting over $H$. Such subgroups
correspond to embedded codimension-$1$ manifolds in a manifold. Sela
\cite{Sela} considered the case when $G$ is a torsion-free word hyperbolic
group and the edge groups of $\Gamma$ are infinite cyclic. Rips and Sela
\cite{RS} generalised this to the case where $G$ is a torsion-free finitely
presented group. The edge groups are again infinite cyclic. Dunwoody and
Sageev \cite{D-Sageev} considered the case when $G$ is a finitely presented
group and the edge groups of $\Gamma$ are slender groups, subject to the
constraint that if $H$ is an edge group, then $G$ admits no splitting over a
subgroup of infinite index in $H$. Fujiwara and Papasoglu \cite{FP} considered
the case when $G$ is a finitely presented group and the edge groups of
$\Gamma$ are small groups, subject to the weaker constraint that if $H$ is an
edge group, then no splitting of $G$ over $H$ can cross strongly a splitting
of $G$ over a subgroup of $H$ of infinite index. (See section \ref{prelim} for
a discussion of crossing and strong crossing.) Finally, Bowditch \cite{B1}
considered the case when $G$ is a word hyperbolic group and the edge groups of
$\Gamma$ are two-ended (which is equivalent to being virtually infinite
cyclic). As in the previous cases, the essential subgroups of $G$ are
two-ended, but Bowditch includes all such subgroups $H$ of $G$ such that
$e(G,H)\geq2$. This corresponds to considering all essential annuli in a
$3$-manifold rather than just embedded ones. In this case, Bowditch proved an
existence and uniqueness result, precisely analogous to the $3$-manifold
theory in the atoroidal case.

The above results are often referred to vaguely but collectively as the
JSJ-decomposition of a finitely presented group. While these results are
commonly regarded as being an algebraic analogue of the topological JSJ
theory, none of them recovers the topological result when applied to the
fundamental group of an orientable Haken $3$-manifold. We list some reasons
for this.

\begin{itemize}
\item None of them has as strong a uniqueness property as the topological
JSJ-decomposition, apart from Bowditch's work \cite{B1} for word hyperbolic
groups. In particular, there is no invariance under automorphisms of the group
$G$ except when $G$ is word hyperbolic.

\item The topological JSJ-decomposition involves both annuli and tori. Apart
from \cite{FP}, none of the algebraic theories can simultaneously handle
splittings over free abelian groups $H\;$and $K$ of different ranks, and
\cite{FP} can only handle this in certain cases.

\item In all the previous work apart from that of Bowditch and Kropholler only
splittings are considered, whereas in the topological JSJ theory, singular
annuli and tori play a crucial role.

\item The algebraic theories only consider strong crossing of splittings,
whereas weak crossing is a key ingredient in the topology.

\item The Enclosing Property of the characteristic submanifold is stronger
than the condition stated earlier that any subgroup of $G=\pi_{1}(M)$ which is
represented by an essential annulus or torus in $M$ is conjugate into a
characteristic vertex group. For tori there is no difference, but if one
considers an embedded annulus, one has the extra information that the graph of
groups structure for $G$ determined by $\mathcal{T}$ can be refined by adding
an edge corresponding to the extra surface. Such a statement is not part of
the results in any of the previous work in this area.
\end{itemize}

Now we will discuss the above points in more detail.

A feature one would expect from any canonical decomposition of a group $G$ is
some sort of invariance under automorphisms of $G$, but in previous work this
was present only in the case when $G$ is word hyperbolic. This invariance has
been exploited by Sela in the case of word hyperbolic groups \cite{S1}, and by
Johannson in the classical case \cite{JO}, for several striking applications.

The characteristic pieces of a $3$-manifold $M$ are of three types, namely
$I$-bundles (which must meet $\partial M$), Seifert fibre spaces which meet
$\partial M$, and Seifert fibre spaces which are in the interior of $M$. We
refer to the two types which meet $\partial M$ as peripheral. For the
discussion in this paragraph, we restrict our attention to essential annuli in
$M$. Then the collection of all peripheral characteristic pieces of $M$ has
the Enclosing Property for all such annuli and no subcollection has this
property. But if one applies any of the above algebraic results to splittings
of $G=\pi_{1}(M)$ over infinite cyclic subgroups, one obtains only analogues
of the $I$-bundle pieces, and the peripheral Seifert fibre spaces are split up
in an arbitrary fashion. Note that by splitting a peripheral Seifert fibre
space, one loses the structure of this piece and hence loses the topological
Enclosing Property. A further point to note is that the above algebraic
splittings may not even yield all the $I$-bundle pieces of $M$. For suppose
that $M\;$has an $I$-bundle piece homeomorphic to $F\times I$, where $F$ is
homeomorphic to the thrice punctured sphere. Thus the only embedded
incompressible annuli in $F\times I$ are the three boundary components. Now
all the relevant algebraic results, apart from those of Bowditch in the case
when $G$ is word hyperbolic, consider only infinite cyclic subgroups of $G$
over which $G$ splits, and this corresponds to restricting attention to
embedded annuli. Thus these results will yield a graph of groups structure for
$G$ with three edges corresponding to the three boundary components of $F$,
but the vertex with associated group $\pi_{1}(F)$ will not be regarded as characteristic.

There are two important special cases of the topological JSJ-decomposition of
a $3$-manifold $M$. One occurs when $M$ is closed, in which case $\mathcal{T}$
consists of tori only and the Enclosing Property applies only to maps of tori
into $M$. The other occurs when $M$ is atoroidal, meaning that $\pi_{1}(M)$
contains no non-peripheral $(\mathbb{Z}\times\mathbb{Z})$-subgroup, in which
case $\mathcal{T}$ consists of annuli only and the Enclosing Property applies
only to maps of annuli into $M$. These two special cases seem to have guided
the development of all previous algebraic analogues of the topological
JSJ-decomposition. In particular, when trying to describe all splittings of a
group over subgroups of a given type, for example infinite cyclic, previous
authors looked for a decomposition described by splittings over subgroups of
the same type. However, if we return to the general topological situation and
consider only essential annuli, we observe that the collection of peripheral
characteristic pieces of $M$ may well have some frontier tori. This means that
even if one wishes initially to consider only splittings of a group over
infinite cyclic subgroups, one is naturally led to consider splittings over
more complicated subgroups as well. Surprisingly, we will see that, in
general, these more complicated groups need not even be finitely generated.

We believe that our ideas in this paper handle all the above problems. We
obtain decompositions of all finitely presented groups which are unique, and
hence invariant under automorphisms, and which essentially specialise to the
classical JSJ-decomposition. In particular, our ideas can handle
simultaneously splittings over free abelian groups of many different ranks. We
obtain more than one such decomposition and there seems to be a number of
further questions about finding refinements and properties of these
decompositions. Our decompositions arise in a simple and natural way, whereas
the previous constructions were all rather indirect.

Here is an introduction to our ideas. As mentioned before, our choice of the
crucial feature of the classical JSJ-decomposition is the Enclosing Property
for immersions. This property implies that the characteristic submanifold
$V(M)$ of a Haken $3$-manifold $M$ contains a representative of every homotopy
class of an essential annulus or torus in $M$. We will say that it
\textit{encloses} every essential annulus and torus in $M$. An important
observation is that, except for a few special cases, $V(M)$ is a regular
neighbourhood of some (finite) collection of essential annuli and tori in $M$.
In particular, this implies that $V(M)$ is minimal (up to isotopy) among
incompressible submanifolds of $M$ which enclose every essential annulus and
torus in $M$. Thus it seems natural to think of $V(M)$ as a regular
neighbourhood of all the essential annuli and tori in $M$. The peripheral
pieces of the characteristic submanifold can be thought of as a regular
neighbourhood of all the essential annuli only. Our main results can be
thought of as algebraic versions of these statements. In order to explain our
ideas further, we need to discuss the algebraic analogues of immersed annuli
and tori and the algebraic analogue of a regular neighbourhood.

An analogue of a $\pi_{1}$-injective immersion in codimension one has been
studied by group theorists for some time. If $H$ denotes the image in $G$ of
the fundamental group of the codimension-$1$ manifold, this analogue is a
subset of $G$ called a $H$-almost invariant set or an almost invariant set
over $H$. The appropriate notions of intersection and disjointness for such
sets were introduced by Scott in \cite{Scott:Intersectionnumbers}. These
notions were further developed in \cite{SS} and the necessary definitions and
results will be recalled in section \ref{prelim}. We should mention that
almost invariant subsets of a group $G$ can appear disconnected in the Cayley
graph of $G$ so that the analogy with immersions may seem a little forced.
However, this is an artifact of the particular choice of generators made when
constructing the Cayley graph. So long as $H$ is finitely generated, one can
always change generators to make any given $H$-almost invariant set connected.
The notion of intersection number developed in
\cite{Scott:Intersectionnumbers} and the main theorems 2.5 and 2.8 of
\cite{SS} do strongly suggest that an almost invariant set is the appropriate
analogue of an immersion.

The key new idea of this paper is an algebraic version of regular
neighbourhood theory. We describe an algebraic regular neighbourhood of a
family of almost invariant subsets of a group $G$. This is a graph of groups
structure for $G$, with the property that certain vertex groups ``enclose''
the given almost invariant sets. The idea of a regular neighbourhood of two
splittings was developed by Fujiwara and Papasoglu \cite{FP} in special cases
(it can be seen that their enclosing technique yields the same result as our
regular neighbourhood construction in these cases), but from our point of view
enclosing almost invariant sets is more basic. Any such technique will
automatically enclose splittings. In our algebraic construction of regular
neighbourhoods, as well as several other techniques, we have greatly benefited
from the two papers of Bowditch \cite{B1}\cite{B2}. Bowditch's use of pretrees
showed us how to enclose almost invariant sets under very general conditions.
See our construction in section \ref{regnbhds:construction}. In the case of
word hyperbolic groups, Bowditch \cite{B1} was effectively the first to
enclose such sets although he does not use this terminology. He also showed
that the characteristic vertex groups of his graph of groups have the
structure of finite-by-Fuchsian groups. Bowditch's techniques were further
developed by Dunwoody and Swenson in \cite{D-Swenson} and Swenson
\cite{Swenson}. Finally in \cite{B2}, Bowditch extended these techniques still
further to study simultaneously splittings of finitely presented groups over
two-ended subgroups and enclosing groups for such subgroups. This paper is
closest to our approach but there are some important differences. It turns out
that the decompositions obtained by Bowditch in \cite{B2} differ from those
which we obtain in this paper in sections \ref{JSJforsmallcommensurisers} and
\ref{JSJforlargecommensurisers}. We really want to enclose almost invariant
sets (which are the algebraic analogue of immersions of codimension-$1$
manifolds) whereas Bowditch does not have a clear analogue of an immersion.
Bowditch uses what he calls an `axis' and this can be taken as an analogue of
an immersion in some cases (see section \ref{prelim} for more details).
Secondly, the intersection number we studied in
\cite{Scott:Intersectionnumbers} and \cite{SS} seems to involve the right
notion of crossing for almost invariant sets, whereas the notion of crossing
used by Bowditch does not recognise what we call weak crossings. Taking care
of these difficulties seems to make our decompositions more canonical and also
corresponds better with the topological situation.

This paper is organised as follows. In section \ref{charsub}, we recall the
basic properties of the characteristic submanifold of a Haken $3$-manifold. In
section \ref{prelim}, we recall some of the algebraic concepts and results
that we need from our paper \cite{SS}. In sections \ref{regnbhds:construction}%
, \ref{enclosing}, \ref{regnbhds:enclosing} and \ref{regnbhds:uniqueness}, we
develop our general theory of regular neighbourhoods of almost invariant
subsets of a group. We give a precise definition of a regular neighbourhood,
and prove that it is always unique. We show that they exist for any finite
family of almost invariant subsets each of which is over a finitely generated
group, although they need not exist in general. This theory seems to be useful
for studying splittings under more general conditions than those we consider
for JSJ-decompositions. In section \ref{regnbhds:uniqueness}, we note a
strengthening of a theorem of Niblo. The results of these first five sections
are very general and apply to almost invariant subsets of any finitely
generated group.

In sections \ref{coendswhencommensuriserissmall} up to
\ref{JSJforlargecommensurisers}, we construct our first canonical
decomposition. We restrict our attention to a one-ended, finitely presented
group $G$ and almost invariant subsets over two-ended subgroups of $G$. Our
decomposition is a graph of groups structure for $G$ which is a regular
neighbourhood of all such almost invariant subsets of $G$. The restriction to
finitely presented groups is necessary because we use certain accessibility results.

In section \ref{coendswhencommensuriserissmall}, we consider two-ended
subgroups of $G$ whose commensuriser in $G$ is ``small''. We discuss the
properties of a regular neighbourhood of finitely many almost invariant
subsets over such subgroups of $G$. This means that we consider a graph of
groups structure for $G$ with certain vertices which enclose all these almost
invariant subsets. When such almost invariant subsets cross strongly, the
enclosing groups can be identified as finite-by-Fuchsian groups by the work of
Bowditch and others. For weak crossings, it follows from our general theory of
regular neighbourhoods that the enclosing groups themselves are two-ended. We
formulate a slightly nonstandard accessibility result that we need for this argument.

In section \ref{coendswhencommensuriserislarge}, we consider two-ended
subgroups of $G$ whose commensuriser in $G$ may be ``large'', and prove the
following technical result. Let $H\;$be such a subgroup of $G$, and let $B(H)$
denote the Boolean algebra of all nontrivial almost invariant subsets of $G$
which are over subgroups commensurable with $H$. We show that $B(H)$ is
finitely generated over the commensuriser of $H$ in $G$. The proof depends on
standard accessibility results and on techniques of Dunwoody and Roller
\cite{D-Roller}.

Using the results of sections \ref{coendswhencommensuriserissmall} and
\ref{coendswhencommensuriserislarge}, we obtain, in sections
\ref{JSJforsmallcommensurisers} and \ref{JSJforlargecommensurisers}, a regular
neighbourhood of all the almost invariant subsets of $G$ which are over
two-ended subgroups, i.e. a natural decomposition of any finitely presented
group $G$ which encloses all such almost invariant subsets. Our uniqueness
result for regular neighbourhoods implies that this decomposition is unique
and also is invariant under automorphisms of $G$. The most remarkable point
about this decomposition is that although the graph of groups structure
obtained is finite, in general not all the vertex and edge groups will be
finitely generated. We end section \ref{JSJforlargecommensurisers} by deducing
the existence of a regular neighbourhood of all the splittings of $G$ which
are over two-ended subgroups. This may seem a more natural object, but it only
seems possible to prove its existence by first considering the above regular
neighbourhood of all almost invariant subsets of $G$ which are over two-ended subgroups.

In section \ref{examples}, we discuss several examples of this decomposition
and compare it with the topological JSJ-decomposition for a $3$-manifold and
with Bowditch's decomposition in \cite{B2}. We show that when $G$ is the
fundamental group of a Haken manifold $M$, the enclosing groups obtained in
section \ref{JSJforlargecommensurisers} correspond to the peripheral
characteristic submanifold. In particular, the decomposition of $G$
corresponds to the full characteristic submanifold when $M$ is atoroidal.

In section \ref{JSJforVPCofgivenrank}, we generalise all the preceding results
as follows. If $G$ is a one-ended, finitely presented group which does not
split over any VPC subgroup of length less than $n$, we construct a regular
neighbourhood of all the almost invariant subsets of $G$ which are over VPC
subgroups of $G$ of length $n$, i.e. a natural decomposition of $G$ which
encloses all such almost invariant subsets. In\ the case when $n=2$, this
result is precisely analogous to the JSJ-decomposition of a closed Haken
$3$-manifold. For $n\geq2$, the special case of this result when $G$ is a
Poincar\'{e} duality group of dimension $n+1$ recovers the results of
Kropholler in \cite{K}. Our results are slightly more general because they
apply to all Poincar\'{e} duality groups, whereas Kropholler's results apply
only to Poincar\'{e} duality groups such that any VPC subgroup has finitely
generated centraliser. An example due to Mess \cite{Mess} shows that this
condition is not always satisfied.

In section \ref{JSJforVPCoftworanks}, we use all the preceding ideas to
construct an algebraic analogue of the whole characteristic submanifold of a
$3$-manifold. Given a one-ended, finitely presented group $G$, let $E_{k}$
denote the collection of all the almost invariant subsets of $G$ over VPC
subgroups of $G$ of length $k$. We start with the graph of groups structure
for $G$ which we gave in section \ref{JSJforlargecommensurisers}. This is a
regular neighbourhood of all the elements of $E_{1}$ and we will denote it by
$\Gamma_{1}$. Next we consider how to enclose elements of $E_{2}$. The closest
analogue to the topology is obtained by considering only those elements of
$E_{2}$ which do not cross any element of $E_{1}$. These are called
$1$-canonical. We show that $\Gamma_{1}$ can be refined to a graph of groups
structure $\Gamma_{1,2}$ for $G$ by adding new vertices which enclose all the
$1$-canonical elements of $E_{2}$. In the case when $G$ is the fundamental
group of a Haken $3$-manifold, we show that this corresponds to the
topological decomposition. This uses our work in \cite{SS3}. It may seem
unsatisfactory that we do not find a decomposition of $G$ with vertex groups
which enclose all elements of $E_{1}$ and $E_{2}$, but our work in \cite{SS3}
shows that if such a decomposition exists it cannot be a refinement of
$\Gamma_{1}$. We discuss this at the end of section \ref{examples}. If $G$
does not split over any VPC subgroups of length less than $n$, as in section
\ref{JSJforVPCofgivenrank}, then we start with the graph of groups structure
for $G$ which we gave in that section. This is a regular neighbourhood of all
the elements of $E_{n}$ and we will denote it by $\Gamma_{n}$. We\ say that an
element of $E_{k}$ is $n$\textit{-canonical} if it crosses no element of
$E_{i}$, for $i\leq n$. As in the case when $n=1$, we show that $\Gamma_{n}$
can be refined to a graph of groups structure $\Gamma_{n,n+1}$ for $G$ by
adding new vertices which enclose all the $n$-canonical elements of $E_{n+1}$.

In section \ref{JSJforabeliangroupsofallranks}, we discuss the natural
question of whether one can continue in the same way. We want to refine
$\Gamma_{1,2}$ to $\Gamma_{1,2,3}$ by adding new vertices which enclose all
the $2$-canonical elements of $E_{3}$. We show that this can be done if we
restrict attention to virtually abelian subgroups of $G$. Letting $A_{k}$
denote the corresponding subset of $E_{k}$, we show that this procedure can be
repeated to obtain, for every $n$, a decomposition $\Gamma_{1,2,...,n}$ of $G$
with vertices which enclose all the $(k-1)$-canonical elements of $A_{k}$,
$1\leq k\leq n$. In some cases, this sequence stabilises at some finite stage,
so that we obtain a decomposition $\Gamma_{\infty}$ of $G$ with vertices which
enclose all the $(k-1)$-canonical elements of $A_{k}$, for all $k\geq1$.
Surprisingly, this procedure does not work for VPC subgroups. We give a simple
example to show that in general there is no decomposition analogous to
$\Gamma_{1,2,3}$ for VPC subgroups of $G$.

In section \ref{non-canonicaldecompositions}, we discuss how the
JSJ-decompositions of previous authors can be related to ours, and in section
\ref{extensions}, we briefly discuss possible extensions of our results to
more general classes of groups.

Acknowledgments: The first author gratefully acknowledges the support of the
Hebrew University, the University of Melbourne, and the University of Warwick.
His visit to Warwick was supported by the EPSRC. Both authors acknowledge the
support of the Universit\'{e} Paul Sabatier in Toulouse. They were supported
by the ARC while in Melbourne, and by the CNRS while in Toulouse.

\section{The characteristic submanifold\label{charsub}}

In this section, we will give a brief summary of the theory, emphasising those
points which are closely related to the algebraic theory in this paper.

Let $M$ be an orientable Haken $3$-manifold. A map of a surface $F$ into $M$
is \textit{proper} if it sends $\partial F$ into $\partial M$. A proper map of
an orientable surface $F$ into $M$ is \textit{incompressible} if it is
$\pi_{1}$-injective. A map of the torus into $M$ is called \textit{essential}
if is incompressible and not homotopic into $\partial M$. A proper map of the
annulus $A$ into $M$ is called \textit{essential} if it is incompressible and
not properly homotopic into $\partial M$. Finally, a codimension-$0$
submanifold of $M$ is \textit{incompressible} if each frontier component is
incompressible in $M$. An incompressible codimension-$0$ submanifold $W$ of
$M$ is \textit{simple} in $M$ if any essential map of the annulus or torus
into $M$ which has image in $W$ can be properly homotoped into the frontier of
$W$.

From now on we need to assume that $M\;$has incompressible boundary. Jaco and
Shalen \cite{JS} and Johannson \cite{JO} proved that there is a family
$\mathcal{T}$ of disjoint essential annuli and tori embedded in $M$, unique up
to isotopy, and with the following properties. The manifolds obtained by
cutting $M$ along $\mathcal{T}$ are simple in $M$ or are Seifert fibre spaces
or $I$-bundles over surfaces. In fact, $\mathcal{T}$ can be characterised as
the minimal family of annuli and tori with this property. The Seifert and
$I$-bundle pieces of $M$ are said to be \textit{characteristic}, and any
essential map of the annulus or torus into $M$ can be properly homotoped to
lie in a characteristic piece. This is called the Enclosing Property of
$\mathcal{T}$. The characteristic submanifold $V(M)$ of $M$ consists
essentially of the union of the characteristic pieces of the manifold obtained
from $M$ by cutting along $\mathcal{T}$. However, if two characteristic pieces
of $M$ have a component $S$ of $\mathcal{T}$ in common, we add a second copy
of $S$ to the family $\mathcal{T}$, thus separating the two characteristic
pieces by a copy of $S\times I$, which we regard as a non-characteristic piece
of $M$. Similarly, if two non-characteristic pieces of $M$ have a component
$S$ of $\mathcal{T}$ in common, we add a second copy of $S$ to the family
$\mathcal{T}$, thus separating the two non-characteristic pieces by a copy of
$S\times I$, which we regard as a characteristic piece of $M$. This is clearly
needed if $V(M)$ is to have the Enclosing Property. Thus the frontier of
$V(M)$ is usually not equal to $\mathcal{T}$. Some annuli or tori in
$\mathcal{T}$ may appear twice in the frontier of $V(M)$. This discussion
brings out a somewhat confusing fact about the characteristic submanifold,
which is that both $V(M)\;$and its complement can have components which are
homeomorphic to $S\times I$, where $S$ is an annulus or torus. One other basic
point to note is that it is quite possible that $\mathcal{T}$ is empty. In
this case, either $V(M)$ is empty or equal to $M$. Thus $M$ is a Seifert fibre
space or admits no essential annuli and tori.

In order to complete this description of $V(M)$, we need to say a little more
about its frontier. If $W$ is a component of $V(M)$ which is a $I$-bundle over
a compact surface $F$, then the frontier of $W$ in $M$ is the restriction of
the bundle to $\partial F$, which is homeomorphic to $\partial F\times I$. If
$W$ is a Seifert fibre space component of $V(M)$, then there is a Seifert
fibration on $W$ such that the frontier of $W$ in $M$ consists of vertical
annuli and tori.

We emphasise that the Enclosing Property applies to any essential map of the
annulus or torus into $M$ and not just to essential embeddings. A related
concept which is important for our approach to JSJ-decompositions is the idea
of a canonical surface in $M$. The concept of canonical surface is not
discussed in the original memoirs \cite{JS} and \cite{JO}. It emerged from
\cite{Leeb-Scott} and \cite{NS} and is further developed in \cite{SS3}. In
\cite{SS3}, an embedded essential annulus or torus $S$ in $M$ is called
\textit{canonical} if any essential map of the annulus or torus into $M$ can
be properly homotoped to be disjoint from $S$. The Enclosing Property clearly
implies that any annulus or torus in $\mathcal{T}$ is canonical. In
\cite{SS3}, we showed that the family of isotopy classes of all canonical
annuli and tori in $M$ is equal to $\mathcal{T}$. In \cite{NS}, Neumann and
Swarup considered a slightly different version of this idea. They defined an
embedded essential annulus or torus in $M$ to be canonical if every
\emph{embedded} essential annulus or torus in $M$ can be properly homotoped to
be disjoint from it. We let $\mathcal{T}_{e}$ denote the family of isotopy
classes of essential annuli and tori in $M$ which are canonical in this sense.
The Enclosing Property of $\mathcal{T}$ implies that $\mathcal{T}$ is
contained in $\mathcal{T}_{e}$. They showed that their family $\mathcal{T}%
_{e}$ is not, in general, the same as the family $\mathcal{T}$, but they were
able to describe the differences and thereby give a new derivation of the
classical JSJ-decomposition. Thus $\mathcal{T}_{e}$ determines a canonical
decomposition of $M$ which is finer than that determined by $\mathcal{T}$.
Note that $\mathcal{T}_{e}-\mathcal{T}$ consists of annuli only, which Neumann
and Swarup call matched annuli. They list the possibilities for such annuli in
Lemma 3.4 of \cite{NS}. However, their list is not quite correct. One case
which they give occurs when $M$ is the Seifert fibre space $W_{p,q}$ which is
constructed by gluing two solid tori together along an annulus $A$ which has
degree $p$ in one solid torus and degree $q$ in the other solid torus, where
$p,q\geq2$. Thus $V(M)$ equals $M$ in this case. They assert that $A$ lies in
$\mathcal{T}_{e}$. But this is not true when $p=q=2$, although it is true for
all other values of $p$ and $q$. The easiest way to see this is to note that
$W_{2,2}$ can also be viewed as the twisted $I$-bundle over the Klein bottle
$K$. (Note that there is a unique $I$-bundle over $K$ with orientable total
space.) For recall that $K$ contains a circle $C$ which cuts $K$ into two
Moebius bands. The restriction of the $I$-bundle to each Moebius band is a
solid torus and the restriction of the $I$-bundle to $C$ is the annulus $A$.
Thus the twisted $I$-bundle over $K$ is homeomorphic to $W_{2,2}$. Now $K$
also has a non-separating two-sided simple closed curve $D$, and $C$ and $D$
cannot be homotoped apart. The restriction of the $I$-bundle to $D$ is an
annulus $B$ in $W_{2,2}$, and it follows that $A$ and $B$ cannot be homotoped
apart. Thus $A$ does not lie in $\mathcal{T}_{e}$ in the case $p=q=2$.

As discussed in the introduction, the guiding idea behind this paper is that
$V(M)$ should be thought of as a regular neighbourhood of the family of all
essential annuli and tori in $M$. By this we mean that every such map is
properly homotopic into $V(M)$ and that $V(M)$ is minimal, up to isotopy,
among all incompressible submanifolds of $M$ with this property. It will be
convenient to say that the collection of all such maps \textit{fills} $V(M)$
when $V(M)$ has this minimality property. The word ``fill'' is used in the
same way to describe certain curves on a surface. A subtle point which arises
here is that there are exceptional cases where $V(M)$, as defined by
Jaco-Shalen and Johannson, is not filled by the collection of all essential
annuli and tori in $M$. In these cases our algebraic decomposition does not
quite correspond to the topological JSJ-decomposition.

Let $V^{\prime}(M)$ denote the submanifold of $M$ which encloses every
essential annulus and torus in $M$ and is filled by them. We will see that
$V^{\prime}(M)$ is only slightly different from $V(M)$. The algebraic
decomposition which we produce in this paper corresponds to $V^{\prime}(M)$
rather than to $V(M)$. Clearly $V^{\prime}(M)$ is an incompressible
submanifold of $V(M)$, so that its frontier in $M$ must also consist of
essential annuli and tori in $M$. It follows that the frontier components of
$V^{\prime}(M)$ are precisely the canonical annuli and tori in $M$, and so are
all isotopic into $\mathcal{T}$. It follows that the difference between
$V^{\prime}(M)$ and $V(M)$ is essentially that certain components of $V(M)$
are discarded. Here is a description of the exceptional cases, most of which
are solid tori. A solid torus component $W$ of $V(M)$ will fail to be filled
by annuli essential in $M$ when its frontier consists of $3$ annuli each of
degree $1$ in $W$, or when its frontier consists of $1$ annulus of degree $2$
in $W$, or when its frontier consists of $1$ annulus of degree $3$ in $W$.
Another exceptional case occurs when a component $W$ of $V(M)$ lies in the
interior of $M$ and is homeomorphic to the manifold $W_{2,2}$ which, as
discussed above, can also be viewed as the twisted $I$-bundle over the Klein
bottle. Any incompressible torus in $W_{2,2}$ is homotopic into the boundary,
so that $W_{2,2}$ is not filled by tori which are essential in $M$. In all
these cases, one obtains $V^{\prime}(M)$ by replacing $W$ by a regular
neighbourhood of its frontier in $M$.

\section{Preliminaries\label{prelim}}

We start by introducing the idea of an almost invariant subset of a finitely
generated group $G$. Throughout this paper, we will always assume that $G$ is
finitely generated, but we will sometimes need to consider subgroups which are
not finitely generated. We emphasise here that all the results of this section
apply to the case when subgroups are not finitely generated, unless it is
specifically stated that subgroups must be finitely generated. We will need
several definitions which we take from \cite{SS}, but see
\cite{Scott:Intersectionnumbers} for a discussion.

\begin{definition}
Two sets $P$ and $Q$ are \textsl{almost equal} if their symmetric difference
$P-Q\cup Q-P$ is finite. We write $P\overset{a}{=}Q$.
\end{definition}

\begin{definition}
If a group $G$ acts on the right on a set $Z$, a subset $P$ of $Z$ is
\textsl{almost invariant} if $Pg\overset{a}{=}P$ for all $g$ in $G$. An almost
invariant subset $P$ of $Z$ is \textsl{nontrivial} if $P$ and its complement
$Z-P$ are both infinite. The complement $Z-P$ will be denoted simply by
$P^{\ast}$, when $Z$ is clear from the context.
\end{definition}

This idea is connected with the theory of ends of groups via the Cayley graph
$\Gamma$ of $G$ with respect to some finite generating set of $G$. (Note that
in this paper groups act on the left on covering spaces and, in particular,
$G$ acts on its Cayley graph on the left.) Using $\mathbb{Z}_{2}$ as
coefficients, we can identify $0$-cochains and $1$-cochains on $\Gamma$ with
sets of vertices or edges. A subset $P$ of $G$ represents a set of vertices of
$\Gamma$ which we also denote by $P$, and it is a beautiful fact, due to Cohen
\cite{Cohen}, that $P$ is an almost invariant subset of $G$ if and only if
$\delta P$ is finite, where $\delta$ is the coboundary operator in $\Gamma$.
Thus $G$ has a nontrivial almost invariant subset if and only if the number of
ends $e(G)$ of $G$ is at least $2$. If $H$ is a subgroup of $G$, we let
$H\backslash G$ denote the set of cosets $Hg$ of $H$ in $G$, i.e. the quotient
of $G$ by the left action of $H$. Of course, $G$ will no longer act on the
left on this quotient, but it will still act on the right. Thus we have the
idea of an almost invariant subset of $H\backslash G$. Further, $P$ is an
almost invariant subset of $H\backslash G$ if and only if $\delta P$ is
finite, where $\delta$ is the coboundary operator in the graph $H\backslash
\Gamma$. Thus $H\backslash G$ has a nontrivial almost invariant subset if and
only if the number of ends $e(G,H)$ of the pair $(G,H)$ is at least $2$.
Considering the pre-image $X$ in $G$ of an almost invariant subset $P$ of
$H\backslash G$ leads to the following definitions.

\begin{definition}
If $G$ is a finitely generated group and $H$ is a subgroup, then a subset $X$
of $G$ is $H$\textsl{-almost invariant} if $X$ is invariant under the left
action of $H$, and simultaneously $H\backslash X$ is an almost invariant
subset of $H\backslash G$. We may also say that $X$ is \textsl{almost
invariant over }$H$. In addition, $X$ is a \textsl{nontrivial} $H$-almost
invariant subset of $G$, if the quotient sets $H\backslash X$ and $H\backslash
X^{\ast}$ are both infinite.
\end{definition}

\begin{remark}
Note that if $X$ is a nontrivial $H$-almost invariant subset of $G$, then
$e(G,H)$ is at least $2$, as $H\backslash X$ is a nontrivial almost invariant
subset of $H\backslash G$.
\end{remark}

\begin{definition}
If $G$ is a group and $H$ is a subgroup, then a subset $W$ of $G$ is
$H$\textsl{-finite} if it is contained in the union of finitely many left
cosets $Hg$ of $H$ in $G$.
\end{definition}

\begin{definition}
If $G$ is a group and $H$ is a subgroup, then two subsets $V$ and $W$ of $G$
are $H$\textsl{-almost equal} if their symmetric difference is $H$-finite.
\end{definition}

It will also be convenient to avoid this rather clumsy terminology sometimes,
particularly when the group $H$ is not fixed, so we make the following definition.

\begin{definition}
If $X$ is a $H$-almost invariant subset of $G\;$and $Y$ is a $K$-almost
invariant subset of $G$, and if $X$ and $Y$ are $H$-almost equal, then we will
say that $X$ and $Y$ are \textsl{equivalent} and write $X\sim Y$.
\end{definition}

\begin{remark}
Note that $H$ and $K$ must be commensurable, so that $X$ and $Y$ are also
$K$-almost equal and $(H\cap K)$-almost equal.

A more elegant, but equivalent formulation is that $X$ is equivalent to $Y$ if
and only if each is contained in a bounded neighbourhood of the other.

Equivalence is important because usually one is interested in an equivalence
class of almost invariant subsets of a group rather than a specific such subset.
\end{remark}

A \textit{splitting} of a group $G$ is an expression of $G$ as a HNN extension
$A\ast_{C}$ or as an amalgamated free product $A\ast_{C}B$, where $A\neq C\neq
B$. Thus a splitting of $G$ always describes a nontrivial decomposition. If
one thinks of a splitting of a group as an algebraic analogue of the
topological notion of an embedded $\pi_{1}$-injective and codimension-$1$
submanifold, then almost invariant sets should be thought of as analogues of
immersions of such manifolds. We can describe the connection between these
ideas as follows. Let $M$ be a closed manifold with fundamental group $G$ and
consider a codimension-$1$ manifold $S$ immersed in $M$ such that the induced
map of fundamental groups is injective. Let $H$ denote the image of $\pi
_{1}(S)$ in $G$, and let $M_{H}$ denote the cover of $M$ such that $\pi
_{1}(M_{H})=H$. For simplicity, we will assume that the lift of $S$ to $M_{H}$
is an embedding, whose image we will also denote by $S$. Then $S$ must
separate $M_{H}$, and we let $A$ denote the closure of one side of $S$ in
$M_{H}$. Let $\widetilde{M}$ denote the universal cover of $M$, and let
$\widetilde{S}$ and $\widetilde{A}$ denote the pre-images in $\widetilde{M}$
of the submanifolds $S$ and $A$ of $M_{S}$. Thus $\widetilde{S}$ is a copy of
the universal cover of $S$. Next pick a generating set for $G$ and represent
it by a bouquet of circles embedded in $M$. We will assume that the wedge
point of the bouquet does not lie on the image of $S$. The pre-image of this
bouquet in $\widetilde{M}$ will be a copy of the Cayley graph $\Gamma$ of $G$
with respect to the chosen generating set. The pre-image in $M_{H}$ of the
bouquet will be a copy of the graph $H\backslash\Gamma$. Let $P$ denote the
set of all vertices of $H\backslash\Gamma$ which lie in $A$. Then $P$ has
finite coboundary, as $\delta P$ equals exactly those edges of $H\backslash
\Gamma$ which cross $S$. Hence $P$ is an almost invariant subset of
$H\backslash G$. If $X$ denotes the set of vertices in $\Gamma$ which lie in
$\widetilde{A}$, then $X$ is the pre-image of $P$ and so is a $H$-almost
invariant subset of $G$. Note that if we replace $A$ by its complement, then
$P$ is replaced by $P^{\ast}$. If we choose a different generating set for $G$
or a different embedding of the bouquet of circles in $M$, the $H$-almost
invariant subsets $P$ and $P^{\ast}$ will change, but the new sets will be
equivalent to $P$ or to $P^{\ast}$. Thus we have associated a $H$-almost
invariant set to the given immersion $S$ and this set is unique up to
equivalence and complementation.

If one has a splitting of a group $G$ over a subgroup $H$, the above
discussion leads to a natural way to associate to this splitting a standard
$H$-almost invariant subset of $G$ which is essentially unique (up to
complementation). But it is simpler and clearer to work in a more
combinatorial setting. For this we recall the basic result of Bass-Serre
theory \cite{Serre}. This tells us that an expression of a group $G$ as
$A\ast_{H}B$ or as $A\ast_{H}$ is equivalent to an action of $G$ on a tree $T$
without inversions so that $G\backslash T$ has a single edge and the edge
stabilisers are conjugates of $H$. (Throughout this paper, we will only
consider $G$-trees on which $G$ acts without inversions, i.e. if an element of
$G$ preserves an edge, it also fixes that edge pointwise. This is only a very
minor restriction. For if $G$ acts on a tree $T$ with inversions, let
$T^{\prime}$ denote the tree obtained from $T$ by dividing each edge into two
edges. There is a natural action of $G$ induced on $T^{\prime}$ which clearly
has no inversions.)

First we consider a general action without inversions of a group $G$ on a tree
$T$. Recall that there is a natural partial order on the oriented edges of a
tree $T$, given by saying that if $s$ and $t$ are oriented edges of $T$, then
$s\leq t$ if and only if there is an oriented path in $T$ which starts with
$s$ and ends with $t$. For any action without inversions of a group $G$ on a
tree $T$, and any edge $s$ of $T$, we have a natural partition of $G$ into two
sets $X_{s}=\{g:gs\geq s$ or $g\overline{s}\geq s\}$ and $X_{s}^{\ast
}=\{g:gs<s$ or $g\overline{s}<s\}$. We will show below that if $S$ denotes the
stabiliser of $s$, then $X_{s}$ and $X_{s}^{\ast}$ are both $S$-almost
invariant. Although this partition of $G$ is natural, it is not quite right
for our purposes, because it is not equivariant under the action of $G$ on
$T$. In fact, if $t$ denotes the edge $ks$ of $T$, then $X_{t}$ is equal to
$kX_{s}k^{-1}$, whereas we would like it to be equal to $kX_{s}$. We resolve
this problem in the following way. We fix some $G$-equivariant map
$\varphi:G\rightarrow V(T)$. By this we mean that $k\varphi(g)=\varphi(kg)$,
for all elements $g$ and $k$ of $G$. Of course such a map is determined once
we choose some vertex of $T$ to be $\varphi(e)$, where $e$ denotes the
identity element of $G$. Now an oriented edge $s$ of $T$ determines a natural
partition of $V(T)$ into two sets, namely the vertices of the two subtrees
obtained by removing the interior of $s$ from $T$. Let $Y_{s}$ denote the
collection of all the vertices of the subtree which contains the terminal
vertex $v$ of $s$, and let $Y_{s}^{\ast}$ denote the complementary collection
of vertices. Then $s$ determines a natural (in terms of our choice of
$\varphi$) partition of $G$ into two sets, namely $Z_{s}=\varphi^{-1}(Y_{s})$
and $Z_{s}^{\ast}=\varphi^{-1}(Y_{s}^{\ast})$. Clearly, these sets are
equivariant, i.e. if $t$ denotes the edge $ks$ of $T$, then $Z_{t}$ is equal
to $kZ_{s}$. We will show that if $S$ denotes the stabiliser of $s$, then
$Z_{s}$ and $Z_{s}^{\ast}$ are both $S$-almost invariant. Further, $Z_{s}$ is
$S$-almost equal to the set $X_{s}$ defined above. It follows that although
$Z_{s}$ depends on our choice of the map $\varphi$, the equivalence class of
$Z_{s}$ is independent of this choice.

\begin{lemma}
\label{pre-imageofhalftreeisalmostinvariant}Let $T$ be a $G$-tree, and let $s$
be an oriented edge of $T$ with stabiliser $S$.

\begin{enumerate}
\item Then the subset $X_{s}=\{g:gs\geq s$ or $g\overline{s}\geq s\}$ of $G$
is $S$-almost invariant.

\item Let $\varphi:G\rightarrow V(T)$ be a $G$-equivariant map, and define
$Z_{s}$ as above. Then $Z_{s}$ is $S$-almost invariant and is $S$-almost equal
to $X_{s}$. If $\varphi(e)$ is the terminal vertex $v$ of $s$, then
$Z_{s}=X_{s}$.
\end{enumerate}
\end{lemma}

\begin{remark}
Note that this result does not require that $G$ be finitely generated.
\end{remark}

\begin{proof}
1) We need to show that $hX_{s}=X_{s}$, for all $h$ in $S$, and that $X_{s}h$
and $X_{s}$ are $S$-almost equal for all $h$ in $G$.

If $gs\geq s$, then $hgs\geq hs$ which equals $s$, for all $h$ in $S$.
Similarly if $g\overline{s}\geq s$, then $hg\overline{s}\geq hs=s$, for all
$h$ in $S$. Thus $hX_{s}\subset X_{s}$, for all $h$ in $S$. Hence $h^{-1}%
X_{s}\subset X_{s}$, for all $h$ in $S$, so that $hX_{s}=X_{s}$, for all $h$
in $S$, as required.

Now consider an element $k$ of $X_{s}h-X_{s}$. Thus there is $g$ in $X_{s}$
such that $k=gh$ does not lie in $X_{s}$. This means that the edge $s$ lies
between the edges $gs$ and $ks=ghs$. Applying $g^{-1}$, we see that the edge
$g^{-1}s$ lies between $s$ and $hs$. Thus $g^{-1}s$ is one of the finitely
many edges of $T$ between $s$ and $hs$, so that $g^{-1}$ lies in some finite
union of cosets $g_{i}S$ of $S$ in $G$. Hence $g$ lies in the union of cosets
$Sg_{i}^{-1}$ and so $k$ also lies in some finite union of cosets $Sk_{i}$.

Next consider an element $k$ of $X_{s}-X_{s}h$. Thus $k$ lies in $X_{s}$ and
$kh^{-1}$ does not. Hence $s$ lies between $ks$ and $kh^{-1}s$ so that
$k^{-1}s$ lies between $s$ and $h^{-1}s$. As in the preceding paragraph, it
follows that $k$ lies in some finite union of cosets $Sk_{j}$ of $S$.

It follows from the previous two paragraphs that $X_{s}h$ and $X_{s}$ are
$S$-almost equal for all $h$ in $G$, as required.

2) Consider a $G$-equivariant map $\varphi:G\rightarrow V(T)$ such that
$\varphi(e)=w$, and let $Z_{s}$ denote $\varphi^{-1}(Y_{s})$. As
$\varphi(e)=w$, we have $\varphi(g)=gw$. It follows that $Z_{s}=\{g\in G:gw\in
Y_{s}\}$. Now recall that $v$ denotes the terminal vertex of $s$. It is easy
to check that $X_{s}=\{g\in G:gv\in Y_{s}\}$, so that $Z_{s}=X_{s}$ if
$\varphi(e)=v$.

Let $k$ be an element of $Z_{s}-X_{s}$. Thus $kw\in Y_{s}$ and $kv\notin
Y_{s}$. Hence $s$ lies between $kv$ and $kw$, so that $k^{-1}s$ lies between
$v$ and $w$. As in part 1), it follows that $k$ lies in some finite union of
cosets $Sk_{i}$ of $S$. Similarly, $X_{s}-Z_{s}$ is contained in a finite
union of cosets of $S$, so that $Z_{s}$ is $S$-almost equal to $X_{s}$, as required.
\end{proof}

In terms of the above, there is now an easy and natural way to associate a
$H$-almost invariant subset of $G$ to a splitting $\sigma$ of $G$ over $H$.
Given $\sigma$, let $s$ denote an oriented edge of $T$ with stabiliser $H$.
Choose the $G$-equivariant map $\varphi:G\rightarrow V(T)$ so that
$\varphi(e)$ is an endpoint of $s$, and then take the $H$-almost invariant
subset $Z_{s}$ or its complement. This description involves three choices,
namely the choices of the edge $s$, its orientation and the choice of an
endpoint of $s$. The choice of $s$ will not alter the almost invariant sets
obtained, so we end up with precisely four $H$-almost invariant subsets of $G$
which are naturally associated to the given splitting over $H$. In \cite{SS},
we gave a different, but equivalent, description of these four sets and called
them the standard $H$-almost invariant subsets of $G$ associated to $\sigma$.
Two of these sets have a particularly nice property, which will play an
important role later on. If we choose $\varphi(e)$ to be the terminal vertex
of $s$, and let $X$ denote the set $Z_{s}$, then $X=\{g\in G:gv\in
Y_{s}\}=\{g\in G:gX^{(\ast)}\subset X\}$, where we use $X^{(\ast)}$ to denote
a set which might be $X$ or $X^{\ast}$. The same equation holds if $X$ denotes
the set $Z_{s}^{\ast}$.

The next definition makes precise the notion of crossing of almost invariant
sets. This is an algebraic analogue of crossing of codimension-$1$ manifolds,
but it ignores ``inessential'' crossings.

\begin{definition}
Let $X$ be a $H$-almost invariant subset of $G$ and let $Y$ be a $K$-almost
invariant subset of $G$. We will say that $Y$ \textsl{crosses} $X$ if each of
the four sets $X\cap Y$, $X^{\ast}\cap Y$, $X\cap Y^{\ast}$ and $X^{\ast}\cap
Y^{\ast}$ is not $H$-finite. Thus each of the four sets projects to an
infinite subset of $H\backslash G$.
\end{definition}

The motivation for the above definition is that when one of the four sets is
empty, we clearly have no crossing, and if one of the four sets is ``small'',
then we have ``inessential crossing''. Note that $Y$ may be a translate of $X$
in which case such crossing corresponds to the self-intersection of a single immersion.

\begin{remark}
It is shown in \cite{Scott:Intersectionnumbers} that if $X$ and $Y$ are
nontrivial, then $X\cap Y$ is $H$-finite if and only if it is $K$-finite. It
follows that crossing of nontrivial almost invariant subsets of $G$ is
symmetric, i.e. that $X$ crosses $Y$ if and only if $Y$ crosses $X$. We will
often write $X^{(\ast)}\cap Y^{(\ast)}$ instead of listing the four sets
$X\cap Y$, $X^{\ast}\cap Y$, $X\cap Y^{\ast}$ and $X^{\ast}\cap Y^{\ast}$.
\end{remark}

\begin{definition}
Let $U$ be a $H$-almost invariant subset of $G$ and let $V$ be a $K$-almost
invariant subset of $G$. We will say that $U\cap V$ is \textsl{small} if it is
$H$-finite.
\end{definition}

\begin{remark}
This terminology will be extremely convenient, particularly when we want to
discuss translates $U$ and $V$ of $X$ and $Y$, as we do not need to mention
the stabilisers of $U$ or of $V$. However, the terminology is symmetric in $U$
and $V$ and makes no reference to $H$ or $K$, whereas the definition is not
symmetric and does refer to $H$, so some justification is required. If $U$ is
also $H^{\prime}$-almost invariant for a subgroup $H^{\prime}$ of $G$, then
$H^{\prime}$ must be commensurable with $H$. Thus $U\cap V$ is $H$-finite if
and only if it is $H^{\prime}$-finite. In addition, the preceding remark tells
us that $U\cap V$ is $H$-finite if and only if it is $K$-finite. This provides
the needed justification of our terminology.
\end{remark}

The term crossing has often been used in the literature for a somewhat
different concept which we call strong crossing. As the name suggests, strong
crossing implies crossing but the converse need not be true. We will now
define this notion. Let $G$ be a finitely generated group and let $H$ and $K$
be subgroups of $G$. Let $X$ be a nontrivial $H$-almost invariant subset of
$G$ and let $Y$ be a nontrivial $K$-almost invariant subset of $G$. It will be
convenient to think of $\delta X$ as a set of edges in $\Gamma$ or as a set of
points in $G$, where the set of points will simply be the collection of
endpoints of all the edges of $\delta X$.

\begin{definition}
We say that $Y$\textsl{\ crosses }$X$\textsl{\ strongly} if both $\delta Y\cap
X$ and $\delta Y\cap X^{\ast}$ project to infinite sets in $H\backslash G$. If
$Y$ crosses $X$ but not strongly, we say that $Y$ \textsl{crosses }$X$\textsl{
weakly}.
\end{definition}

\begin{remark}
\label{strongcrossingimpliescrossing}These definitions are independent of the
choice of generators for $G$ which is used to define $\Gamma$. Clearly, if $Y$
crosses $X$ strongly, then $Y$ crosses $X$. Note that $Y$ does not cross $X$
strongly if and only if $\delta Y$ is contained in a bounded neighbourhood of
$X$ or $X^{\ast}$.
\end{remark}

An interesting point about strong crossing of $X$ and $Y$ is that it depends
only on the subgroups $H$ and $K$. More precisely, we have the following result.

\begin{lemma}
\label{strongcrossingisdeterminedbygroups}Let $G$ be a finitely generated
group and let $H$ and $K$ be subgroups of $G$. Let $X$ and $X^{\prime}$ be
nontrivial $H$-almost invariant subsets of $G$ and let $Y$ and $Y^{\prime}$ be
nontrivial $K$-almost invariant subsets of $G$. Then $Y$ crosses $X$ strongly
if and only if $Y^{\prime}$ crosses $X^{\prime}$ strongly.
\end{lemma}

\begin{proof}
By Remark \ref{strongcrossingimpliescrossing}, $Y$ does not cross $X$ strongly
if and only if $\delta Y$ is contained in a bounded neighbourhood of $X$ or
$X^{\ast}$. As $H\backslash\delta Y$ and $H\backslash\delta Y^{\prime}$ are
both finite, it follows that $\delta Y$ is contained in a bounded
neighbourhood of $\delta Y^{\prime}$ and vice versa. Thus $Y$ crosses $X$
strongly if and only if $Y^{\prime}$ crosses $X$ strongly. By reversing the
roles of $X$ and $Y$, we immediately obtain the required result.
\end{proof}

\begin{definition}
Let $\sigma_{1}$ and $\sigma_{2}$ be splittings of $G$ over $C_{1}$ and
$C_{2}$, and let $X_{i}$ be one of the standard $C_{i}$-almost invariant
subsets of $G$ associated to the splitting $\sigma_{i}$, for $i=1$, $2$. Then
$\sigma_{1}$\textsl{ crosses }$\sigma_{2}$\textsl{ (strongly) }if $X_{1}$
crosses $X_{2}$ (strongly).
\end{definition}

\begin{remark}
As the standard $C_{i}$-almost invariant subsets of $G$ associated to the
splitting $\sigma_{i}$ are all $C_{i}$-almost equal or $C_{i}$-almost
complementary, this definition does not depend on the choice of the $X_{i}$'s.
\end{remark}

Next we give a simple example to show that strong crossing is not symmetric,
in general.

\begin{example}
Consider an essential two-sided simple closed curve $S$ on a compact surface
$F$ which intersects a simple arc $L$ transversely in a single point. Let $G$
denote $\pi_{1}(F)$, and let $H$ and $K$ respectively denote the subgroups of
$G$ carried by $S$ and $L$, so that $H$ is infinite cyclic and $K$ is trivial.
Then $S$ and $L$ each define a splitting of $G$ over $H\;$and $K$
respectively. Let $X$ and $Y$ denote associated standard $H$-almost invariant
and $K$-almost invariant subsets of $G$. These correspond to submanifolds of
the universal cover of $F$ bounded respectively by a line $\widetilde{S}$
lying above $S$ and by a compact interval $\widetilde{L}$ lying above $L$,
such that $\widetilde{S}$ meets $\widetilde{L}$ transversely in a single
point. Clearly, $X$ crosses $Y$ strongly but $Y$ does not cross $X$ strongly.
\end{example}

If $\sigma_{1}$ and $\sigma_{2}$ are splittings of $G$ over finitely generated
subgroups $C_{1}$ and $C_{2}$ respectively, Sela introduced in \cite{S1} the
following notion of crossing of $\sigma_{1}$ and $\sigma_{2}$. He says that
$\sigma_{1}$ is \textit{hyperbolic with respect to }$\sigma_{2}$ if $C_{1}$ is
not conjugate into a vertex group of the splitting $\sigma_{2}$. It is easy to
show that this idea is the same as strong crossing, and we give the proof below.

\begin{lemma}
If $\sigma_{1}$ and $\sigma_{2}$ are splittings of a finitely generated group
$G$ over finitely generated subgroups, then $\sigma_{1}$ is \textit{hyperbolic
with respect to }$\sigma_{2}$ if and only if $\sigma_{1}$ crosses $\sigma_{2}$ strongly.
\end{lemma}

\begin{proof}
Consider the $G$-tree $T_{2}$ corresponding to the splitting $\sigma_{2}$ and
let the amalgamating group of the splitting $\sigma_{i}$ be $H_{i}$. Consider
the action of $H_{1}$ on $T_{2}$. It is immediate from the definition that
$\sigma_{1}$ is hyperbolic with respect to $\sigma_{2}$ if and only if $H_{1}$
does not fix a vertex of $T_{2}$. As we are assuming that $H_{1}$ is finitely
generated, it fixes some vertex of $T_{2}$ if and only if every element of
$H_{1}$ fixes some vertex of $T_{2}$. Thus $\sigma_{1}$ is hyperbolic with
respect to $\sigma_{2}$ if and only if some element of $H_{1}$ fixes no vertex
of $T_{2}$ and thus has an axis. We claim that this implies that $\sigma_{1}$
crosses $\sigma_{2}$ strongly. Let $X_{i}$ denote a $H_{i}$-almost invariant
subset of $G$ associated to the splitting $\sigma_{i}$ as discussed just after
Lemma \ref{pre-imageofhalftreeisalmostinvariant}. As the quotient
$H_{i}\backslash\delta X_{i}$ is finite, it follows that $\delta X_{i}$ must
lie within a bounded distance of $H_{i}$. Now let $s$ denote an oriented edge
of $T_{2}$ with stabiliser $H_{2}$, let $v$ denote the terminal vertex of $s$,
and define $\varphi:G\rightarrow V(T)$ by $\varphi(g)=gv$. Then $\varphi
(\delta X_{2})=s$, and $\varphi(H_{1})$ contains points arbitrarily far from
$s$ and on each side of $s$. It now follows that $X_{1}$ crosses $X_{2}$
strongly so that $\sigma_{1}$ crosses $\sigma_{2}$ strongly as claimed. On the
other hand, if $\sigma_{1}$ is not hyperbolic with respect to $\sigma_{2}$,
then $H_{1}$ fixes a vertex of $T_{2}$ and it follows that $\sigma_{1}$ cannot
cross $\sigma_{2}$ strongly, because the same considerations show that $\delta
X_{1}$ must lie within a bounded distance of $X_{2}$ or $X_{2}^{\ast}$.
\end{proof}

\begin{remark}
A key point in the above argument is that if every element of $H_{1}$ fixes
some vertex of $T_{2}$, then $H_{1}$ itself fixes some vertex. If $H_{1}$ is
not finitely generated, this can fail, but a result of Tits, Lemma 3.4 of
\cite{Tits}, can be used instead to show that the above lemma holds for
splittings over infinitely generated subgroups.
\end{remark}

Sela showed in \cite{S1} that if $C_{1}$ and $C_{2}$ are two-ended and if $G$
does not split over a finite group, then his crossing, and hence our strong
crossing, is symmetric.

The following technical result plays an important role in the theory of almost
invariant sets.

\begin{lemma}
\label{finitenumberofdoublecosets}Let $G$ be a finitely generated group with
finitely generated subgroups $H$ and $K$, a non-trivial $H$-almost invariant
subset $X$ and a non-trivial $K$-almost invariant subset $Y$. Then $\{g\in
G:gX$ and $Y$ are not nested\} consists of a finite number of double cosets $KgH.$
\end{lemma}

\begin{proof}
Let $\Gamma$ denote the Cayley graph of $G$ with respect to some finite
generating set for $G$. Let $P$ denote the almost invariant subset
$H\backslash X$ of $H\backslash G$ and let $Q$ denote the almost invariant
subset $K\backslash Y$ of $K\backslash G$. Recall from the start of this
section, that if we identify $P$ with the $0$-cochain on $H\backslash\Gamma$
whose support is $P$, then $P$ is an almost invariant subset of $H\backslash
G$ if and only if $\delta P$ is finite. Thus $\delta P$ is a finite collection
of edges in $H\backslash\Gamma$ and similarly $\delta Q$ is a finite
collection of edges in $K\backslash\Gamma$. Now let $C$ denote a finite
connected subgraph of $H\backslash\Gamma$ such that $C$ contains $\delta P$
and the natural map $\pi_{1}(C)\rightarrow H$ is onto, and let $E$ denote a
finite connected subgraph of $K\backslash\Gamma$ such that $E$ contains
$\delta Q$ and the natural map $\pi_{1}(E)\rightarrow K$ is onto. Thus the
pre-image $D$ of $C$ in $\Gamma$ is connected and contains $\delta X$, and the
pre-image $F$ of $E$ in $\Gamma$ is connected and contains $\delta Y.$ Let
$\Delta$ denote a finite subgraph of $D$ which projects onto $C$, and let
$\Phi$ denote a finite subgraph of $F$ which projects onto $E$. If $gD$ meets
$F$, there must be elements $h$ and $k$ in $H$ and $K$ such that $gh\Delta$
meets $k\Phi$. Now $\{\gamma\in G:\gamma\Delta$ meets $\Phi\}$ is finite, as
$G$ acts freely on $\Gamma$. It follows that $\{g\in G:gD$ meets $F\}$
consists of a finite number of double cosets $KgH.$

The result would now be trivial if $X$ and $Y$ were each the vertex set of a
connected subgraph of $\Gamma$. As this need not be the case, we need to make
a careful argument as in the proof of Lemma 5.10 of
\cite{Scott-Wall:Topological}. Consider $g$ in $G$ such that $gD$ and $F$ are
disjoint. We will show that $gX$ and $Y$ are nested. As $D$ is connected, the
vertex set of $gD$ must lie entirely in $Y$ or entirely in $Y^{\ast}.$ Suppose
that the vertex set of $gD$ lies in $Y$. For a set $S$ of vertices of $\Gamma
$, let $\overline{S}$ denote the maximal subgraph of $\Gamma$ with vertex set
equal to $S$. Each component $W$ of $\overline{X}$ and $\overline{X^{\ast}}$
contains a vertex of $D$. Hence $gW$ contains a vertex of $gD$ and so must
meet $Y$. If $gW$ also meets $Y^{\ast}$, then it must meet $F$. But as $F$ is
connected and disjoint from $gD$, it lies in a single component $gW$. It
follows that there is exactly one component $gW$ of $\overline{gX}$ and
$\overline{gX^{\ast}}$ which meets $Y^{\ast}$, so that we must have $gX\subset
Y$ or $gX^{\ast}\subset Y$. Similarly, if $gD$ lies in $Y^{\ast}$, we will
find that $gX\subset Y^{\ast}$ or $gX^{\ast}\subset Y^{\ast}$. It follows that
in either case $gX$ and $Y$ are nested as required.
\end{proof}

Now we come to the definition of the intersection number of two almost
invariant sets.

\begin{definition}
\label{defnofintersectionnumber}Let $H$ and $K$ be subgroups of a finitely
generated group $G$. Let $P$ denote a nontrivial almost invariant subset of
$H\backslash G$, let $Q$ denote a nontrivial almost invariant subset of
$K\backslash G$ and let $X$ and $Y$ denote the pre-images of $P$ and $Q$
respectively in $G$. Then the \textsl{intersection number} $i(P,Q)$ of $P$ and
$Q$ equals the number of double cosets $KgH$ such that $gX$ crosses $Y$.
\end{definition}

\begin{remark}
\label{almostequalsetshavesameintersectionnumber}The following facts about
intersection numbers are proved in \cite{Scott:Intersectionnumbers}.

\begin{enumerate}
\item Intersection numbers are symmetric, i.e. $i(P,Q)=i(Q,P)$.

\item $i(P,Q)$ is finite when $G$, $H$ and $K$ are all finitely generated.
This follows immediately from Lemma \ref{finitenumberofdoublecosets}.

\item If $P^{\prime}$ is an almost invariant subset of $H\backslash G$ which
is almost equal to $P$ or to $P^{\ast}$ and if $Q^{\prime}$ is an almost
invariant subset of $K\backslash G$ which is almost equal to $Q$ or to
$Q^{\ast}$, then $i(P^{\prime},Q^{\prime})=i(P,Q)$.
\end{enumerate}
\end{remark}

One can also define the intersection number of two splittings to be the
intersection number of the almost invariant sets associated to the splittings.
Now if two curves on a surface have intersection number zero, they can be
isotoped to be disjoint. There is a natural algebraic analogue of this fact.
We define a collection of $n$ splittings of a group $G$ to be
\textit{compatible} if $G$ can be expressed as the fundamental group of a
graph of groups with $n$ edges, such that the edge splittings of the graph are
conjugate to the given splittings. The following result is a slight rewording
of Theorem 2.5 of \cite{SS}.

\begin{theorem}
\label{Theorem2.5ofSS}Let $G$ be a finitely generated group with $n$
splittings over finitely generated subgroups. Then the splittings are
compatible if and only if each pair of splittings has intersection number
zero. Further, in this situation, the graph of groups structure on $G$
obtained from these splittings is unique up to isomorphism.
\end{theorem}

In \cite{Bass}, Bass gives a discussion of isomorphisms of graphs of groups.
In particular, if two graphs of groups structures on a group $G$ are
isomorphic, they have isomorphic underlying graph, and edges and vertices
which correspond under the isomorphism carry conjugate subgroups of $G$.

This discussion of intersection numbers leads naturally to the concept which
we call a canonical splitting. Recall from section \ref{charsub} that an
embedded essential annulus or torus $F$ in a $3$-manifold $M$ is called
\textit{canonical} if any essential map of the annulus or torus into $M$ can
be properly homotoped to be disjoint from $F$. An equivalent formulation of
this definition is that an embedded essential annulus or torus $F$ in a
$3$-manifold $M$ is \textit{canonical} if any essential map of the annulus or
torus into $M$ has intersection number zero with $F$ (where the intersection
number is defined as in \cite{FHS}). In \cite{SS3}, we showed that the
canonical annuli and tori in a Haken $3$-manifold $M$ are the same (up to
isotopy) as the frontier of the characteristic submanifold, thus yielding the
classical JSJ-decomposition.

There is another natural approach to this idea. As discussed at the end of
section \ref{charsub}, Neumann and Swarup \cite{NS} defined annuli and tori
embedded in a $3$-manifold $M$ to be canonical if they have intersection
number zero with every \emph{embedded} essential annulus or torus in $M$. Each
of these approaches has natural algebraic generalisations. Generalising our
idea of canonical would involve considering splittings of a group over
subgroups isomorphic to $\mathbb{Z}$ or $\mathbb{Z}\times\mathbb{Z}$ which
have intersection number zero with many almost invariant sets. Generalising
Neumann and Swarup's idea of canonical would involve considering splittings of
a group over subgroups isomorphic to $\mathbb{Z}$ or $\mathbb{Z}%
\times\mathbb{Z}$ which have intersection number zero with many splittings.
For our purposes, the first idea turns out to be most useful, but it seems to
be important to consider a much larger class of subgroups. The following is
the algebraic definition which we derive from the above discussion.

\begin{definition}
\label{definecanonical}Let $G$ be a one-ended finitely generated group and let
$X$ be a nontrivial almost invariant subset over a subgroup $H$ of $G$.

For $n\geq1$, we will say that $X$ is $n$\textsl{-canonical} if $X$ has zero
intersection number with any nontrivial almost invariant subset of any
$K\backslash G$, for which $K$ is VPC of length at most $n$.

For $n\geq1$, we will say that $X$ is $n$\textsl{-canonical with respect to
abelian groups }if $X$ has zero intersection number with any nontrivial almost
invariant subset of any $K\backslash G$, for which $K$ is virtually free
abelian of rank at most $n$.

If $X$ is associated to a splitting $\sigma$ of $G$ and $X$ is $n$-canonical
(with respect to abelian groups), we will say that $\sigma$ is $n$%
\textsl{-canonical} (\textsl{with respect to abelian groups}).

If $H$ is virtually infinite cyclic, and $X$ is $1$-canonical, we will often
say simply that $X$ is \textsl{canonical}.
\end{definition}

\begin{remark}
If $H$ is not finitely generated, we will only use these ideas when $X$ is
associated to a splitting over $H$.

If $n\leq2$, then $X$ is $n$-canonical if and only if it is $n$-canonical with
respect to abelian groups.
\end{remark}

The definitions above make perfectly good sense when $G$ has more than one
end, and when $n=0$. However we show in Lemma \ref{never0-canonical} that if a
group $G$ has infinitely many ends, then almost invariant subsets of $G$ are
never $0$-canonical. Of course, an almost invariant subset of $G$ which is not
$0$-canonical is certainly not $n$-canonical for any $n$.

Many other related ideas can be defined by changing the class of groups in
which $K$ lies. For example, one could insist that $K$ be free abelian, or
that $K$ has length equal to some fixed number $k$. In particular, in
\cite{SS3}, we gave a similar definition, but restricted $H$ to be infinite
cyclic, $K$ to be free abelian, $X$ to be associated to a splitting, and $n$
to equal $1$ or $2$. In this situation, almost invariant subsets of
$K\backslash G$ correspond to (possibly singular) annuli or tori in a
$3$-manifold. One could also restrict attention to those almost invariant
subsets of $K\backslash G$ which are associated to splittings. This is
analogous, in the $3$-manifold situation, to considering the Enclosing
Property only for embedded annuli or tori, which is effectively what Neumann
and Swarup were doing in \cite{NS}.

In \cite{SS3}, we also considered the connection between the canonical annuli
and tori in $M$ and the $2$-canonical splittings of $G=\pi_{1}(M)$ over
subgroups isomorphic to $\mathbb{Z}$ or $\mathbb{Z}\times\mathbb{Z}$. We
showed that every such $2$-canonical splitting of $G$ arises from a canonical
annulus or torus in $M$, and that every canonical annulus in $M$ determines a
$2$-canonical splitting of $G$. Further every canonical torus in $M$
determines a $1$-canonical splitting of $G$. However, we also showed that
often $M$ will have canonical tori which determine splittings of $G$ which are
not $2$-canonical.

In order to start developing our algebraic theory of regular neighbourhoods in
the next section, we will need to consider collections of almost invariant
subsets of a given group. The following terminology will be useful.

\begin{definition}
A collection $E$ of subsets of $G$ which are closed under complementation is
called \textsl{nested} if for any pair $U$ and $V$ of sets in the collection,
one of the four sets $U^{(\ast)}\cap V^{(\ast)}$ is empty.

If each element $U$ of $E$ is a $H_{U}$-almost invariant subset of $G$ for
some subgroup $H_{U}$ of $G$, we will say that $E$ is \textsl{almost
nested}\textit{\ }if for any pair $U$ and $V$ of sets in the collection, one
of the four sets $U^{(\ast)}\cap V^{(\ast)}$ is small.
\end{definition}

If one is given a $H$-almost invariant subset $X$ of $G$, it is natural to ask
whether there is a $K$-almost invariant set $Y$ which is equivalent to $X$ and
is associated to a splitting of $G$ over $K$. This is analogous to asking
whether a given codimension-$1$ immersion in a manifold can be homotoped to
cover an embedding. We state below Theorem 1.17 of \cite{SS}, which we will
need later. To prove this result, we used the almost nested assumption to
construct a tree with $G$-action. There are more general results of this type
in section 2 of \cite{SS}, which will be extended and used later on.

\begin{theorem}
\label{Theorem1.12ofSS}

\begin{enumerate}
\item \textit{Let }$H$\textit{ be a finitely generated subgroup of a finitely
generated group }$G$\textit{. Let }$X$\textit{ be a nontrivial }%
$H$\textit{-almost invariant subset of }$G$\textit{ such that }%
$E=\{gX,gX^{\ast}:g\in G\}$\textit{ is almost nested and if two of the four
sets }$X^{(\ast)}\cap gX^{(\ast)}$\textit{ are small, then at least one of
them is empty. Then }$G$\textit{ splits over the stabiliser }$H^{\prime}%
$\textit{ of }$X$\textit{ and }$H^{\prime}$\textit{ contains }$H$\textit{ as a
subgroup of finite index. Further, one of the }$H^{\prime}$\textit{-almost
invariant sets }$Y$\textit{ determined by the splitting is }equivalent\textit{
to }$X$.

\item \textit{Let }$H_{1},\ldots,H_{k}$\textit{ be finitely generated
subgroups of a finitely generated group }$G$\textit{. Let }$X_{i}$\textit{ be
a nontrivial }$H_{i}$\textit{-almost invariant subset of }$G$\textit{ such
that }$E=\{gX_{i},gX_{i}^{\ast}:1\leq i\leq k,g\in G\}$\textit{ is almost
nested. Suppose further that, for any pair of elements }$U$\textit{ and }%
$V$\textit{ of }$E$, \textit{if two of the four sets }$U^{(\ast)}\cap
V^{(\ast)}$\textit{ are small, then at least one of them is empty. Then }%
$G$\textit{ can be expressed as the fundamental group of a graph of groups
whose }$i$-th\textit{ edge corresponds to a conjugate of a splitting of }%
$G$\textit{ over the stabiliser }$H_{i}^{\prime}$\textit{ of }$X_{i}$,
\textit{and }$H_{i}^{\prime}$\textit{ contains }$H_{i}$\textit{ as a subgroup
of finite index.}
\end{enumerate}
\end{theorem}

We end this section by briefly discussing Bowditch's notion of coends. It
turns out that this notion had been considered earlier by Ross Geoghegan
\cite{RG} under the name `filtered coends'. In \cite{B2}, Bowditch considers a
pair $(G,H)$ of finitely generated groups where $G$ is one-ended and $H$ is
two-ended. Let $G$ act properly discontinuously and cocompactly on a locally
finite one-ended graph $\Gamma$, for example the Cayley graph of $G$. Let
$S(H)$ denote the set of $H$-invariant connected subgraphs of $\Gamma$ with
finite quotient. If $A\in S(H)$, let $C(A)$ denote the collection of all
complementary components of $A$ and let $C_{\infty}(A)$ denote those
components which are not contained in any element of $S(H)$. If $A$ and $B$
are elements of $S(H)$ such that $A\subset B$, then there is a surjective map
from $C_{\infty}(B)$ to $C_{\infty}(A)$. Let $E(H)$ denote the inverse limit
in the category of topological spaces, where $A$ and $B$ range over $S(H)$.
Then $E(H)$ is a compact and totally disconnected space and an element of
$E(H)$ is called a coend of $H$. It is clear from this discussion that the
concept of coends extends to any pair of finitely generated groups $(G,H)$.
Ross Geoghegan has shown that the number of coends is the same as the number
of relative ends considered earlier by Kropholler and Roller \cite{KR} (See
also the introduction of \cite{B2}).

We note that it is clear from the definition that the number of coends of $H$
is at least as large as the number of ends of the pair $(G,H)$. But it can
easily be substantially larger. An interesting example occurs when $G$ is the
fundamental group of a closed non-orientable surface $F$, not $P^{2}$, and $H$
is the subgroup of $G$ carried by a simple but one-sided circle on $F$. For
then $e(G,H)=1$, as the cover of $F$ determined by $H$ is an open Moebius
band, but if $H_{2}$ denotes the subgroup of $H$ of index $2$, then
$e(G,H_{2})=2$, and it is clear that $H$, like $H_{2}$, has two coends.

Bowditch calls an element of $S(H)$ an \textit{axis} of $H$ if it satisfies
some further technical conditions. An axis $A$ is called \textit{proper} if
$C_{\infty}(A)$ has at least two elements. If a proper axis corresponding to a
subgroup $K$ crosses $A$, then Bowditch shows that $H$ and $K$ both have two
coends (see section 10 of \cite{B1}). If the number of coends of $H$ is $2$
and $H$ does not interchange the coends, then $H$ defines an essentially
unique (up to equivalence and complementation) $H$-almost invariant subset.
Thus in this case, crossing of axes is the same as strong crossing of the
corresponding almost invariant sets. The arguments used by Bowditch in
\cite{B2} are geared to the above case where $G$ is one-ended and $H$ is
two-ended, but his ideas work just as well in the case where $H$ is VPC and
$G$ does not split over any subgroup of length less than that of $H$.

\section{Algebraic Regular Neighbourhoods:
Construction\label{regnbhds:construction}}

In this and the following three sections, we will discuss our algebraic
analogue of the topological idea of a regular neighbourhood. It would be
possible to start by giving an abstract definition of an algebraic regular
neighbourhood, and then to prove existence and uniqueness results. But the
precise definition is somewhat technical, so we will start with our
construction of an algebraic regular neighbourhood. We will consider a
finitely generated group $G$ with finitely generated subgroups $H_{1}%
,\ldots,H_{n}$. For $i=1,\ldots,n$, let $X_{i}$ be a nontrivial $H_{i}$-almost
invariant subset of $G$. In this section, we will construct a finite graph of
groups structure $\Gamma(X_{1},\ldots,X_{n}:G)$ for $G$. In section
\ref{regnbhds:uniqueness}, we will define abstractly what constitutes an
algebraic regular neighbourhood of the $X_{i}$'s in $G$, and then will prove
that $\Gamma(X_{1},\ldots,X_{n}:G)$ satisfies the requirements of our
definition. We will also prove that any algebraic regular neighbourhood of the
$X_{i}$'s in $G$ is naturally isomorphic to $\Gamma(X_{1},\ldots,X_{n}:G)$.
These results are the algebraic analogue of the existence and uniqueness
results for regular neighbourhoods in topology.

While the restriction to finite families of almost invariant sets is very
natural, regular neighbourhoods of infinite families will play an important
role in this paper. At the end of this section, we will briefly discuss how to
modify our construction when the number of $X_{i}$'s is infinite. At the end
of section \ref{regnbhds:enclosing}, we also discuss what happens to our
construction if some of the $H_{i}$'s are not finitely generated. It is a
surprising fact that, even in this situation, our theory of algebraic regular
neighbourhoods goes through in certain cases.

In order to introduce our ideas, consider a connected manifold $M$ and let
$\mathcal{T}$ be a compact (possibly disconnected) codimension-$1$ two-sided
manifold properly embedded in $M$. (The reader will not miss anything by
thinking of $M$ as a surface, and $\mathcal{T}$ as a collection of circles and
arcs.) If each component of $\mathcal{T}$ is $\pi_{1}$-injective in $M$, then
$\mathcal{T}$ defines a graph of groups structure $\Gamma$ on $G=\pi_{1}(M)$
whose underlying graph is dual to $\mathcal{T}$. The edge groups of $\Gamma$
are the fundamental groups of the components of $\mathcal{T}$ and the vertex
groups are the fundamental groups of the components of the complement of
$\mathcal{T}$. If the components of $\mathcal{T}$ are not all $\pi_{1}%
$-injective, then $\mathcal{T}$ still determines a graph of groups structure
for $G$, with the same underlying graph, but the edge and vertex groups are
obtained from the above by replacing each group by its image in $G$.

Now suppose that we consider a compact (possibly disconnected) codimension-$0$
submanifold $N$ of $M$. We will associate to this the graph of groups
decomposition $\Gamma$ of $G$ determined as above by the frontier of $N$ in
$M$. (If $M$ is closed, this frontier is the same as the boundary of $N$.) The
vertices of $\Gamma$ correspond to components of $N\;$and of $M-N$ and each
edge of $\Gamma$ joins a vertex of one type to a vertex of the other type.
Thus $\Gamma$ is naturally a bipartite graph. Throughout this paper, we will
denote the collection of vertices of $\Gamma$ which correspond to components
of $N$ by $V_{0}(\Gamma)$, or simply $V_{0}$ if the context is clear. The
remaining vertices will be denoted by $V_{1}(\Gamma)$ or simply $V_{1}$. If we
consider the pre-image $\widetilde{N}$ of $N$ in the universal cover
$\widetilde{M}$ of $M$, the dual graph to $\partial\widetilde{N}$ is a tree
$T$ on which $G\;$acts with quotient $\Gamma$, and the vertex and edge groups
of $\Gamma$ are simply the vertex and edge stabilisers for the action of $G$
on $T$. Again $T$ is naturally bipartite with some vertices corresponding to
components of $\widetilde{N}$ and some vertices corresponding to components of
$\widetilde{M}-\widetilde{N}$.

In the previous paragraph, we discussed how any subsurface $N$ of $M$
corresponds to a bipartite graph of groups structure for $\pi_{1}(M)$. In what
follows, we will be interested in the situation where $N$ is a regular
neighbourhood of the union of a finite collection of codimension-$1$ manifolds
$C_{\lambda}$ properly immersed in a manifold $M$, and in general position.
Our aim is to understand what topological conditions on $N$ make it a regular
neighbourhood, and to translate these into algebraic conditions on a bipartite
graph of groups. The precise definition is given in Definition
\ref{defnofalgregnbhd}.

The following two conditions on $N$ clearly hold, and these are the conditions
upon which we will base our definition of an algebraic regular neighbourhood.

\begin{enumerate}
\item Each $C_{\lambda}$ is contained in $N$, and each component of $N$
contains some $C_{\lambda}$.

\item If $C$ is any codimension-$1$ manifold embedded in $M$ and disjoint from
each $C_{\lambda}$, then $C$ is homotopic into $M-N$.
\end{enumerate}

These conditions are not sufficient to characterise $N$ up to isotopy even
when $M$ is $2$-dimensional, but they do contain much of the information
needed to determine $N$.

Recall that each essential $C_{\lambda}$ determines a nontrivial $H_{\lambda}%
$-almost invariant subset of $G=\pi_{1}(M)$. Our algebraic regular
neighbourhood of the $X_{i}$'s in $G$ will be a bipartite graph of groups
structure $\Gamma(X_{1},\ldots,X_{n}:G)$ for $G$ which satisfies two
conditions analogous to 1) and 2) above. The first condition is that the
$V_{0}$-vertices of $\Gamma$ ``enclose'' the $X_{i}$'s in a certain sense. The
second condition is that splittings of $G$ which have intersection number zero
with each $X_{i}$ are enclosed by the $V_{1}$-vertices of $\Gamma$. The
definition and basic properties of enclosing are given in section
\ref{enclosing}. In addition, there are extra conditions which are analogous
to those needed to characterise $N$ completely in the case when $M$ is $2$-dimensional.

In order to understand our construction of $\Gamma(X_{1},\ldots,X_{n}:G)$,
consider further the topological situation where $N$ is a regular
neighbourhood of the $C_{\lambda}$'s. For simplicity, suppose that each
$C_{\lambda}$ lifts to an embedding in $M_{\lambda}$, the cover of $M$ whose
fundamental group equals that of $C_{\lambda}$, and let $S_{\lambda}$ denote
the pre-image in $\widetilde{M}$ of $C_{\lambda}$ in $M_{\lambda}$. Let
$\Sigma$ denote the collection of all the translates of all the $S_{\lambda}%
$'s, and let $\left|  \Sigma\right|  $ denote the union of all the elements of
$\Sigma$. Thus $\left|  \Sigma\right|  $ is the complete pre-image in
$\widetilde{M}$ of the union of the images of all the $C_{\lambda}$'s and
$\widetilde{N}$ is a regular neighbourhood of $\left|  \Sigma\right|  $. The
fact that $N\;$is a regular neighbourhood of the $C_{\lambda}$'s implies that
the inclusion of the union of the $C_{\lambda}$'s into $N$ induces a bijection
between components, and an isomorphism between the fundamental groups of
corresponding components. It follows that the inclusion of $\left|
\Sigma\right|  $ into $\widetilde{N}$ also induces a bijection between
components. Thus the $V_{0}$-vertices of $T$ correspond to the components of
$\left|  \Sigma\right|  $, and the $V_{1}$-vertices of $T$ correspond to the
components of $\widetilde{M}-\left|  \Sigma\right|  $. If two elements $S$ and
$S^{\prime}$ of $\Sigma$ belong to the same component of $\left|
\Sigma\right|  $, there must be a finite chain $S=S_{0},S_{1},...,S_{n}%
=S^{\prime}$ of elements of $\Sigma$ such that $S_{i}$ intersects $S_{i+1}$,
for each $i$. Thus the elements of $\Sigma$ which form a component of $\left|
\Sigma\right|  $ are an equivalence class of the equivalence relation on
$\Sigma$ generated by saying that two elements of $\Sigma$ are related if they
intersect. This is what we want to encode in the algebraic setting, except
that as we will be dealing with almost invariant sets, we will want to ignore
all ``inessential'' intersections.

Now return to our finitely generated group $G$ with finitely generated
subgroups $H_{1},\ldots,H_{n}$ and nontrivial $H_{i}$-almost invariant subsets
$X_{i}$ of $G$. Let $E=\{gX_{i},gX_{i}^{\ast}:g\in G,1\leq i\leq n\}$. At the
end of this section, we will also discuss situations where we are given an
infinite family of $X_{i}$'s. Previously we said that the $X_{i}$'s are the
algebraic analogue of the immersed $C_{\lambda}$'s. However, it is neater to
consider the pair $\{X,X^{\ast}\}$ as a single object, and regard this as the
algebraic analogue of an immersion. We will denote the unordered pair
$\{X,X^{\ast}\}$ by $\overline{X}$, and will say that $\overline{X}$ crosses
$\overline{Y}$ if $X$ crosses $Y$. Then our algebraic analogue of the set
$\Sigma$ is the set $\overline{E}=\{g\overline{X_{i}}:g\in G,1\leq i\leq n\}$.
Our analogue of the equivalence relation on $\Sigma$ which described the
components of $\left|  \Sigma\right|  $ is the equivalence relation on
$\overline{E}$ which is generated by saying that two elements $A$ and $B$ of
$\overline{E}$ are related if they cross. We call an equivalence class a
\textit{cross-connected component} (CCC) of $\overline{E}$, and denote the
equivalence class of $A$ by $[A]$. Note that this is a purely combinatorial
definition. The use of the word component is simply to emphasise the analogies
with the topological situation of the preceding paragraph. We will denote the
collection of all CCC's in $\overline{E}$ by $P$.

Now we want to construct in a natural way a bipartite $G$-tree $T$ with $P$ as
its set of $V_{0}$-vertices. Note that if one has a tree, there is a natural
idea of betweenness for vertices. We will reverse this process and construct
the required tree $T$ starting from an idea of betweenness on the set $P$. In
order to define this idea of betweenness, we will first introduce a partial
order on $E$ using ideas from the proof of Theorem \ref{Theorem1.12ofSS}.
There we defined a partial order on the set $E$, and constructed a $G$-tree
whose edges were the elements of $E$ by showing that the partial order
satisfied certain conditions. Unlike the situation of Theorem
\ref{Theorem1.12ofSS}, we cannot expect to construct a tree with the elements
of $E$ as its edges. For this would imply that each $X_{i}$ determined a
splitting of $G$. However, the partial order will still play a crucial role in
our situation.

If $U$ and $V$ are two elements of $E$ such that $U\subset V$, then our
partial order will have $U\leq V$. But we also want to define $U\leq V$ when
$U$ is ``nearly'' contained in $V$. Precisely, we want $U\leq V$ if $U\cap
V^{\ast}$ is small. However, an obvious difficulty arises when two of
$U^{(\ast)}\cap V^{(\ast)}$ are small, as we have no way of deciding between
two possible inequalities. It turns out that we can avoid this difficulty if
we know that whenever two of $U^{(\ast)}\cap V^{(\ast)}$ are small, then one
of them is empty. Thus we consider the following condition on $E$:

Condition (*): If $U$ and $V$ are in $E$, and two of $U^{(\ast)}\cap
V^{(\ast)}$ are small, then one of $U^{(\ast)}\cap V^{(\ast)}$ is empty.

If $E$ satisfies Condition (*), we will say that the family $X_{1}%
,\ldots,X_{n}$ is in \textit{good position}.

Assuming that this condition holds, we can define a relation $\leq$ on $E$ by
saying that $U\leq V$ if and only if $U\cap V^{\ast}$ is empty or is the only
small set among the four sets $U^{(\ast)}\cap V^{(\ast)}$. Despite the
seemingly artificial nature of this definition, one can show that $\leq$ is a
partial order on $E$. This is not entirely trivial, but the proof is in Lemma
2.4 of \cite{Scott:Annulus} and repeated more group theoretically in Lemma
1.14 of \cite{SS}. Condition (*) plays a key role in the proof. If $U\leq V$
and $V\leq U$, it is easy to see that we must have $U=V$, using the fact that
$E\;$satisfies Condition (*). Most of the proof of Lemma 1.14 of \cite{SS} is
devoted to showing that $\leq$ is transitive. We note here that the argument
that $\leq$ is transitive does not require that the $H_{i}$'s be finitely generated.

In general, the family $X_{1},\ldots,X_{n}$ need not be in good position, but
it turns out that this does not cause any problems. We will discuss this just
after Theorem \ref{Pisapretree}.

Before we go any further, we need to discuss the idea of a pretree. As already
mentioned, the vertices of a tree possess a natural idea of betweenness. The
idea of a pretree formalises this. A \textit{pretree} consists of a set $P$
together with a ternary relation on $P$ denoted $xyz$ which one should think
of as meaning that $y$ is strictly between $x$ and $z$. The relation should
satisfy the following four axioms:

\begin{itemize}
\item (T0) If $xyz$, then $x\neq z$.

\item (T1) $xyz$ implies $zyx$.

\item (T2) $xyz$ implies not $xzy$.

\item (T3) If $xyz$ and $w\neq y$, then $xyw$ or $wyz$.
\end{itemize}

A pretree is said to be \textit{discrete}, if, for any pair $x$ and $z$ of
elements of $P$, the set $\{y\in P:xyz\}$ is finite. Clearly, the vertex set
of any simplicial tree forms a discrete pretree with the induced idea of
betweenness. It is a standard result that a discrete pretree $P$ can be
embedded in a natural way into the vertex set of a simplicial tree $T$ so that
the notion of betweenness is preserved. We briefly describe the construction
of $T$ following Bowditch's papers \cite{B1} and \cite{B2}. For the proofs,
see section 2 of \cite{B1} and section 2 of \cite{B2}. These are discussed in
more detail in \cite{B4} and \cite{AN}. For any pretree $P$, we say that two
elements of $P$ are \textit{adjacent} if there are no elements of $P\;$between
them. We define a \textit{star} in $P$ to be a maximal subset of $P$ which
consists of mutually adjacent elements. (This means that any pair of elements
of a star are adjacent.) We now enlarge the set $P$ by adding in all the stars
of $P$ to obtain a new set $V$. One can define a pretree structure on $V$
which induces the original pretree structure on $P$. A star is adjacent in $V$
to each element of $P$ that it contains. Next we give $V$ the discrete
topology and add edges to $V$ to obtain a graph $T$ with $V$ as its vertex
set. For each pair of adjacent elements of $V$, we simply add an edge which
joins them. If $P$ is discrete, then it can be shown that $T$ is a tree with
vertex set $V$. It follows easily from this construction that if a group $G$
acts on the original pretree $P$, this action extends naturally to an action
of $G$ on the simplicial tree $T$. Moreover, $G$ will act without inversions
on $T$. This will then give a graph of groups decomposition for $G$, though
this decomposition would be trivial if $G$ fixed a vertex of $T$. The tree $T$
is clearly bipartite with vertex set $V(T)$ expressed as the disjoint union of
$V_{0}(T)$ and $V_{1}(T)$, where $V_{0}(T)$ equals $P$ and $V_{1}(T)$ equals
the collection of stars in $P$.

Note that if we start with a tree $T^{\prime}$, let $P$ denote its vertex set
with the induced idea of betweenness and then construct the tree $T$ as above,
then $T$ is not the same as $T^{\prime}$. In fact, $T$ is obtained from
$T^{\prime}$ by subdividing every edge into two edges.

Now we return to our discussion of the set $E$ which we still assume satisfies
Condition (*). We have the partial order $\leq$ on $E$ and we write $U<V$ if
$U\leq V$ but $U$ is not equal to $V$. If $U<Z<V$, we will say that $Z$ is
\textit{between} $U$ and $V$. We summarise below some elementary facts about
$E$.

\begin{lemma}
\label{posatisfiesDunwoody}If $E$ satisfies Condition (*), then $E$ together
with $\leq$ satisfies the following conditions.

\begin{enumerate}
\item If $U$, $V\in E$ and $U\leq V$, then $V^{\ast}\leq U^{\ast}$,

\item If $U$, $V\in E$, there are only finitely many $Z\in E$ such that $Z$ is
between $U$ and $V$,

\item If $U$, $V\in E$, one cannot have $U\leq V$ and $U\leq V^{\ast}$.
\end{enumerate}
\end{lemma}

\begin{proof}
As $U^{\ast}$ denotes the complement of $U$, the first part of this lemma is clear.

Since $E$ consists of translates of a finite number of almost invariant sets
over finitely generated subgroups, it is a standard fact that there are only
finitely many elements of $E$ between $U$ and $V$ (See Lemma 2.6 in
\cite{Scott:Annulus} or Lemma 1.15 in \cite{SS}).

Finally, if $U$ and $V$ lie in $E$, one cannot have $U\leq V$ and $U\leq
V^{\ast}$. For this would imply that $U\cap V^{\ast}$ and $U\cap V$ are each
small, so that $U$ itself would be small, contradicting the assumption that
each $X_{i}$ is a nontrivial $H_{i}$-almost invariant subset of $G$.
\end{proof}

Condition 2) of the above lemma will play an important role in our later
discussions, so we will give this condition a name.

\begin{definition}
Let $E$ be a partially ordered set. We will say that $E$ is \textsl{discrete}
if for any elements $U$ and $V$ of $E$, there are only finitely many elements
of $E$ between $U$ and $V$.
\end{definition}

Next we use these properties of our partial order on $E$ to define a notion of
betweenness on the set $\overline{E}$. Recall that $\overline{X}$ denotes the
unordered pair $\{X,X^{\ast}\}$ and that $\overline{E}$ denotes the set of all
translates of all the $\overline{X_{i}}$'s.

\begin{definition}
Let $L$, $M$, $N\in\overline{E}$. We say that $M$ is \textsl{between} $L$ and
$N$ if there exist $U\in L$, $V\in M$, and $W\in N$ such that $U<V<W$, and we
write $LMN$ or $\overline{U}\,\overline{V}\,\overline{W}$ with $U$, $V$, $W$
chosen as above.
\end{definition}

Note that it is clear that if $\overline{U}\,\overline{V}\,\overline{W}$
holds, then $\overline{W}\,\overline{V}\,\overline{U}$ also holds.

Recall that we say that $\overline{X}$ crosses $\overline{Y}$ if $X$ crosses
$Y$. This generates an equivalence relation on $\overline{E}$, whose
equivalence classes we call \textit{cross-connected components} (\textit{CCC)}%
. We denote the equivalence class of $\overline{X}$ by $[\overline{X}]$, and
denote the collection of all CCC's in $\overline{E}$ by $P$.

We extend the above idea of betweenness in $\overline{E}$ to one in $P$, as follows.

\begin{definition}
Let $A$, $B$ and $C$ be distinct cross-connected components of $\overline{E}$.
We say that $B$ is \textsl{between} $A$ and $C$ and write $ABC$ if there exist
$\overline{U}\in A$, $\overline{V}\in B$ and $\overline{W}\in C$ such that
$\overline{V}$ is between $\overline{U}$ and $\overline{W}$, i.e.
$\overline{U}\,\overline{V}\,\overline{W}$.
\end{definition}

In order for this definition to be useful, we need to know that it is
independent of the choices of $\overline{U}$ and $\overline{W}$. This is what
we prove in Corollary \ref{orderonclumpsmakessense} below. We need two small
results first.

\begin{lemma}
If $U$, $V$ and $Z$ are elements of $E$ such that $U\leq V$ and $Z$ crosses
$U$ but $Z$ does not cross $V$, then either $Z\leq V$ or $Z^{\ast}\leq V$.
\end{lemma}

\begin{proof}
As $Z$ crosses $U$, none of $Z^{(\ast)}\cap U^{(\ast)}$ is small. Since $U\leq
V$, it follows that $Z\cap V$ and $Z^{\ast}\cap V$ are not small. As $Z$ does
not cross $V$, either $Z\cap V^{\ast}$ or $Z^{\ast}\cap V^{\ast}$ is small.
Hence either $Z\leq V$ or $Z^{\ast}\leq V$ as claimed.
\end{proof}

\begin{lemma}
If $U$, $V$ and $Z$ are elements of $E$ such that $U<V$ and $Z$ crosses $U$
but $Z$ does not cross $V$, then either $Z<V$ or $Z^{\ast}<V$.
\end{lemma}

\begin{proof}
This follows from the preceding lemma, since if either $Z$ or $Z^{\ast}$ were
equal to $V$, we would have one of the inequalities $U<Z$ or $U<Z^{\ast}$,
which would contradict the assumption that $Z$ crosses $U$.
\end{proof}

\begin{corollary}
\label{orderonclumpsmakessense}Let $A$, $B$ and $C$ be distinct
cross-connected components of $\overline{E}$, and suppose that $\overline{U}$
and $\overline{U^{\prime}}$ lie in $A$, $\overline{V}$ lies in $B$, and
$\overline{W}$ and $\overline{W^{\prime}}$ lie in $C$. If $\overline
{U}\,\overline{V}\,\overline{W}$, then $\overline{U^{\prime}}\,\overline
{V}\,\overline{W^{\prime}}$.
\end{corollary}

\begin{proof}
It is easy to reduce this to the case when $\overline{U^{\prime}}$ crosses
$\overline{U}$ and $\overline{W}=\overline{W^{\prime}}$. We can also assume
that $U<V<W$. Now we write $U^{\prime}=Z$ and apply the preceding lemma. This
tells us that either $Z<V$ or $Z^{\ast}<V$. By the definition of betweenness,
this implies that $\overline{V}$ is between $\overline{U^{\prime}}$ and
$\overline{W}$, as required.
\end{proof}

Now we are ready to show that if $P$ denotes the collection of all CCC's of
$\overline{E}$ equipped with the relation of betweenness defined above, then
$P$ is a discrete pretree.

\begin{theorem}
\label{Pisapretree}Let $G$ denote a finitely generated group, and let
$H_{1},\ldots,H_{n}$ be finitely generated subgroups of $G$. For each $1\leq
i\leq n$, let $X_{i}$ be a nontrivial $H_{i}$-almost invariant subset of $G$,
and suppose that the $X_{i}$'s are in good position, so that the set $E$ of
all translates of all the $X_{i}$'s satisfies Condition (*). Form the set
$\overline{E}$ as above, and consider the collection $P$ of all
cross-connected components of $\overline{E}$ equipped with the relation of
betweenness introduced above. Then the following statements hold:

\begin{enumerate}
\item $P$ is a pretree, and $G$ acts on $P$ in a natural way.

\item The pretree $P$ is discrete, and the quotient $G\backslash P$ is finite.
Further, the stabilisers of elements of $P$ under this $G$-action are finitely generated.

\item As $P$ is discrete, it can be embedded in a natural way into the vertex
set of a $G$-tree $T$ so that the quotient $G\backslash T$ is a bipartite
graph of groups $\Gamma(X_{1},\ldots,X_{n}:G)$. This graph is finite, and the
$V_{0}$-vertex groups are finitely generated.
\end{enumerate}
\end{theorem}

\begin{remark}
\label{remarkthatedgegroupsneednotbefg}We will say that $\Gamma(X_{1}%
,\ldots,X_{n}:G)$ is a regular neighbourhood of the $X_{i}$'s in $G$. The
$V_{0}$-vertex groups of $\Gamma(X_{1},\ldots,X_{n}:G)$ are called the
\textit{enclosing groups} of the $X_{i}$'s. We will formally define a regular
neighbourhood in section \ref{regnbhds:uniqueness}, and will define enclosing
in section \ref{enclosing}. Recall that $\Gamma(X_{1},\ldots,X_{n}:G)$ is the
algebraic analogue of a regular neighbourhood of a finite family of immersed
codimension-$1$ submanifolds of a manifold with fundamental group $G$, and the
enclosing groups correspond to the fundamental groups of the components of the
regular neighbourhood.

Note that even though the enclosing groups of $\Gamma(X_{1},\ldots,X_{n}:G)$
are finitely generated, the edge groups and the $V_{1}$-vertex groups need not
be finitely generated, even when $G$ is finitely presented. We give examples
of this phenomenon in Example \ref{exampleofregnbhdwithnonfgedgegroups}.

The fact that $\Gamma(X_{1},\ldots,X_{n}:G)$ is finite can be strengthened
greatly. We will see in the next section that $T$ is a minimal $G$-tree.
\end{remark}

\begin{proof}
1) The action of $G$ on itself by left multiplication induces an action of $G$
on $E$ and hence on $\overline{E}$. As this action preserves crossing, it is
immediate that it induces an action on $P$.

Now we verify that $P$ satisfies the four axioms (T0)-(T3) for a pretree. For
the convenience of the reader, we give these axioms again.

(T0) If $xyz$, then $x\neq z$.

(T1) $xyz$ implies $zyx$.

(T2) $xyz$ implies not $xzy$.

(T3) If $xyz$ and $w\neq y$, then $xyw$ or $wyz$.

Axioms (T0) and (T1) are immediate from our definition of betweenness.

To prove (T2), suppose that $A$, $B$ and $C$ are elements of $P$ such that we
have both $ABC$ and $ACB$. As $ABC$ holds, Corollary
\ref{orderonclumpsmakessense} tells us that there is $\overline{V}\in B$ such
that we have $\overline{U}\,\overline{V}\,\overline{W}$ for any $\overline
{U}\in A$ and $\overline{W}\in C$. As $ACB$ also holds, there is $\overline
{W}\in C$ such that $\overline{U}\,\overline{W}\,\overline{V}$ holds for any
$\overline{U}\in A$. We choose some $\overline{U}\in A$. For these particular
choices of $U$, $V$ and $W$, we have both $\overline{U}\,\overline
{V}\,\overline{W}$ and $\overline{U}\,\overline{W}\,\overline{V}$. As
$\overline{U}\,\overline{V}\,\overline{W}$ holds, we can arrange that $U<V<W$,
by replacing sets by their complement if needed. As $\overline{U}%
\,\overline{W}\,\overline{V}$ holds, there exist $X\in\overline{U}$,
$Y\in\overline{V}$, $Z\in\overline{W}$ such that $X<Z<Y$. Now consider the
inequalities $U<V$ and $X<Y$, and recall that $X$ equals $U$ or $U^{\ast}$,
and $Y$ equals $V$ or $V^{\ast}$. It is easy to see that the only possibility
is that $X=U$ and $Y=V$. For example, if we had $X=U^{\ast}$ and $Y=V$, the
inequalities $U<V$ and $X<Y$ would imply that $U<V$ and $U^{\ast}<V$, which is
impossible. Similarly, the inequalities $U<W$ and $X<Z$ imply that $X=U$ and
$Z=W$. But now the inequality $V<W$ implies that $Y<Z$ which contradicts the
inequality $Z<Y$. This completes the proof of (T2).

We next verify (T3). Suppose that we have $A$, $B$, $C$, $D\in P$ with $ABC$
and $D\neq B$. We must show that $ABD$ or $DBC$. Choose $\overline{U}\in A$,
$\overline{V}\in B$, $\overline{W}\in C$ so that $U<V<W$, and choose
$\overline{X}\in D$. The result is trivial if $D$ equals $A$ or $C$, so we
will assume that $D$ is not one of $A$, $B$ or $C$. Thus $X$ does not cross
any of $U$, $V$ or $W$. Since $X$ does not cross $U$, one of $X^{(\ast)}\cap
U^{(\ast)}$ is small. Thus, we can, by interchanging $X$ and $X^{\ast}$ if
necessary, arrange that either $U<X$ or $X<U$. If $X<U$, then $X<V<W$ so that
$DBC$ holds, and we are done. If $U<X$, then we compare $X$ and $V$. Since
they do not cross, we should have one of the four inequalities $X<V$, $V<X$,
$X<V^{\ast}$ or $V^{\ast}<X$. If $X<V$, then as above we have $DBC$, and if
$V<X$, we have $U<V<X$ so that $ABD$ holds. The inequality $X<V^{\ast}$ is
impossible, as we already have $U<X$ and $U<V$. The inequality $V^{\ast}<X$
implies that $X^{\ast}<V$, which again implies $DBC$. This completes the proof
that $P$ is a pretree.

2) As the given family of $X_{i}$'s is finite, part (2) of Lemma
\ref{posatisfiesDunwoody} tells us that $P$ is discrete. There is a natural
bijection between the given family of $X_{i}$'s and the quotient
$G\backslash\overline{E}$. As the quotient $G\backslash P$ is a quotient of
$G\backslash\overline{E}$, it follows that $G\backslash P$ is finite as
required. Finally, we need to show that the stabilisers of elements of $P$ are
finitely generated. We start by noting that the stabiliser $K_{i}$ of $X_{i}$
contains $H_{i}$ as a subgroup of finite index and so must be finitely
generated. Now let $v$ denote an element of $P$, and consider those $X_{i}$'s
which have a translate in $v$. By renumbering, we can assume that $X_{i}$ has
a translate $Y_{i}$ in $v$ for $1\leq i\leq k$. We let $S_{i}$ denote the
stabiliser of $Y_{i}$. As $S_{i}$ is conjugate to $K_{i}$, it is also finitely
generated. For $1\leq i,j\leq k$, we define $\Sigma_{ij}=\{g\in G:gY_{i}$
crosses $Y_{j}\}$. Lemma 1.16 of \cite{SS} tells us that $\Sigma_{ij}$ is
contained in a finite number of double cosets $S_{i}gS_{j}$. Now it is easy to
show that $Stab(v)$, the stabiliser of $v$, is generated by the $S_{i}$'s,
$1\leq i\leq k$, together with the $\Sigma_{ij}$'s, $1\leq i,j\leq k$. It
follows that $Stab(v)$ is finitely generated.

3) Recall that a discrete pretree $P$ can be embedded in a natural way into
the vertex set of a tree $T$, and that an action of $G$ on $P$ which preserves
betweenness will automatically extend to an action without inversions on $T$.
Also $T$ is a bipartite tree with vertex set $V(T)=V_{0}(T)\cup V_{1}(T)$,
where $V_{0}(T)$ equals $P$, and $V_{1}(T)$ equals the collection of all stars
in $P$. It follows that the quotient $G\backslash T$ is naturally a bipartite
graph of groups $\Gamma$ with $V_{0}$-vertex groups conjugate to the
stabilisers of elements of $P$ and $V_{1}$-vertex groups conjugate to the
stabilisers of stars in $P$. Further, part 2) implies that $\Gamma$ has only
finitely many $V_{0}$-vertices, and that each such vertex carries a finitely
generated group. Now we can show that $\Gamma$ must be a finite graph. For as
$G$ is finitely generated, there is a finite subgraph $\Gamma_{1}$ of $\Gamma$
with fundamental group $G$. Thus $\Gamma-\Gamma_{1}$ has finitely many
components each of which must be a tree. Hence if $\Gamma$ were infinite, one
of these trees would be infinite and so $\Gamma$ would have infinitely many
$V_{0}$-vertices. This contradiction completes the proof of the lemma.
\end{proof}

Before proceeding to consider the general situation when the family
$X_{1},\ldots,X_{n}$ is not in good position, we discuss some simple examples.

\begin{example}
Let $G$ be any group which splits over a subgroup $H$. This implies that $G$
acts on a tree $T$ with quotient a graph $\Gamma$ with a single edge which
yields the given splitting of $G$. If we let $X$ denote a $H$-almost invariant
subset of $G$ associated to the given splitting and apply the preceding
construction to the set $E$ of all translates of $X$, we obtain a new graph of
groups structure $\Gamma^{\prime}$ for $G$ which is obtained from $\Gamma$ by
subdividing its edge into two edges. The reason for this is that the
translates of $X$ correspond to the oriented edges of $T$ and the partial
order $\leq$ on $E$ corresponds to the natural ordering of oriented edges of
$T$. It is easy to see that in this case each CCC of $\overline{E}$ consists
of a single translate of $X$, so that the collection $P$ of all CCC's
corresponds to the edges of $T$ and has the corresponding partial order.

In the particular case when $G$ is the fundamental group of a surface $M$ and
the splitting is given by a simple arc or closed curve on $M$, then
$\Gamma^{\prime}$ is the graph of groups structure for $G$ which corresponds
exactly to the topological regular neighbourhood of $C$.
\end{example}

\begin{example}
Let $M$ be a surface and let $C$ denote a finite family of essential simple
arcs or closed curves on $M$.

\begin{enumerate}
\item If the arcs and curves in $C$ are all disjoint, then $C$ defines a graph
of groups structure $\Gamma$ for $G=\pi_{1}(M)$ such that the underlying graph
of $\Gamma$ is dual to $C$. As in the previous example, our algebraic regular
neighbourhood construction yields a graph of groups structure $\Gamma^{\prime
}$ for $G$ which is obtained from $\Gamma$ by subdividing each edge into two
edges, which again corresponds exactly to a topological regular neighbourhood
of $C$.

\item If $C$ denotes two essential simple arcs or closed curves on $M$ which
have minimal possible intersection, then again it is true that our algebraic
regular neighbourhood construction yields a graph of groups structure
$\Gamma^{\prime}$ for $G$ which corresponds exactly to a topological regular
neighbourhood of $C$.

\item However, if $C$ denotes three essential simple arcs or closed curves on
$M$ such that each pair has minimal possible intersection, the algebraic
regular neighbourhood may be a little different from the topological one. For
example, cutting $M$ along $C$ may yield a disc component. In this case, the
dual graph for the regular neighbourhood would have a terminal vertex carrying
the trivial group, but our algebraic regular neighbourhood construction can
never yield such a vertex.
\end{enumerate}
\end{example}

In the preceding construction and the proof of Theorem \ref{Pisapretree}, we
used the assumption that the family $X_{1},\ldots,X_{n}$ was in good position.
This condition need not always be satisfied, so we need to discuss how to
modify our regular neighbourhood construction to handle the general situation.

Suppose that we consider any finite family of nontrivial $H_{i}$-almost
invariant subsets $X_{i}$ of $G$. Recall that the basic idea we used was that
of a cross-connected component (CCC) of $\overline{E}$. We can consider the
equivalence relation generated by crossing of elements of $\overline{E}$
whether or not the family $X_{1},\ldots,X_{n}$ is in good position. Thus we
can always define the family $P$ of all CCC's of $\overline{E}$ and there will
always be a natural action of $G$ on $P$. The importance of good position was
that it enabled us to define the inequality $\leq$ on $E$ and hence to define
the relation of betweenness on $P$. Suppose that we have distinct elements
$U\;$and $V$ of $E$ such that two of $U^{(\ast)}\cap V^{(\ast)}$ are small,
but neither is empty. This means that, when we attempt to define the
inequality $\leq$ on $E$, the elements $U\;$and $V$ are not comparable.
However, note that $U\;$and $V$ must be equivalent up to complementation. Thus
if there is an element $W$ of $E$ which crosses $U$, then $W$ also crosses
$V$, so that $\overline{U}$ and $\overline{V}$ will lie in the same CCC. We
will say that the family $X_{1},\ldots,X_{n}$ is \textit{in good enough
position }if whenever we find incomparable elements $U$ and $V$ of $E$ which
do not cross, there is some element $W$ of $E$ which crosses them. It is easy
to see that all the preceding discussion in this section applies essentially
unchanged if the family $X_{1},\ldots,X_{n}$ is in good enough position. The
point is that any pair of incomparable elements already lie in the same CCC,
and so we never need to be able to compare them.

Next we consider the case when the family $X_{1},\ldots,X_{n}$ is not in good
enough position. We will say that an element of $E$ which crosses no element
of $E$ is \textit{isolated in }$E$. Note that this condition depends on the
set $E$, but we will often omit the phrase ``in $E$'' when the context is
clear. As any translate of an isolated element is also isolated, such elements
can occur only if the original family $X_{1},\ldots,X_{n}$ contains elements
which are isolated in $E$. By re-labelling, we can arrange that $X_{1}%
,\ldots,X_{k}$ are the only isolated elements of the $X_{i}$'s, for some $k$
between $1$ and $n$. The first part of the following result tells us that we
can replace the isolated $X_{i}$'s by equivalent sets such that the new family
is in good enough position, and so we can define a regular neighbourhood of
the $X_{i}$'s to be a regular neighbourhood for the new family. The second
part of this result tells us that the resulting graph of groups structure for
$G$ is independent of the choices made, and so can reasonably be denoted by
$\Gamma(X_{1},\ldots,X_{n}:G)$. Thus this result completes our regular
neighbourhood construction.

\begin{lemma}
\label{candefineregnbhdwhennotingoodenoughposn}Let $G$ denote a finitely
generated group, and let $H_{1},\ldots,H_{n}$ be finitely generated subgroups
of $G$. For each $i\geq1$, let $X_{i}$ be a nontrivial $H_{i}$-almost
invariant subset of $G$, such that $X_{1},\ldots,X_{k}$ are the only isolated
elements of the $X_{i}$'s.

\begin{enumerate}
\item For each $i$, $1\leq i\leq k$, there is an almost invariant set $Z_{i}$
equivalent to $X_{i}$, such that the translates of all the $Z_{i}$'s are
nested and the family $Z_{1},\ldots,Z_{k},X_{k+1},\ldots,X_{n}$ is in good
enough position.

\item The regular neighbourhood $\Gamma(Z_{1},\ldots,Z_{k},X_{k+1}%
,\ldots,X_{n}:G)$ is independent of the choices of the $Z_{i}$'s.
\end{enumerate}
\end{lemma}

\begin{proof}
If $X_{i}$ and $X_{j}$ are two isolated elements of $E$ with stabilisers
$H_{i}$ and $H_{j}$ respectively, then the almost invariant subsets
$H_{i}\backslash X_{i}$ of $H_{i}\backslash G$ and $H_{j}\backslash X_{j}$ of
$H_{j}\backslash G$ have intersection number zero with each other, and each
has self-intersection number zero. In \cite{SS}, we discussed almost invariant
sets with intersection number zero, and the main results of that paper are
exactly what we need to understand the present situation. As $H_{i}\backslash
X_{i}$ has self-intersection number zero, Theorem 2.8 of \cite{SS} tells us
that $X_{i}$ is equivalent to $Y_{i}$ such that $Y_{i}$ is associated to a
splitting of $G$ over a subgroup $K_{i}$ commensurable with $H_{i}$. We will
replace each of $X_{1},\ldots,X_{k}$, by $Y_{1},\ldots,Y_{k}$, chosen as above
so that each $Y_{i}$ is associated to a splitting $\sigma_{i}$ of $G$ over
$K_{i}$. In topological terms, this corresponds to starting with some closed
curves $C_{i}$ on a surface such that each $C_{i}$ is homotopic to some power
of a simple closed curve $S_{i}$, and then replacing each $C_{i}$ by $S_{i}$.
Note that as $Y_{i}$ is equivalent to $X_{i}$, the splittings $\sigma
_{1},\ldots,\sigma_{k}$ have intersection number zero with each other. Now
Theorem \ref{Theorem2.5ofSS} tells us that these splittings are compatible.
This means that we can replace the $Y_{i}$'s by equivalent almost invariant
sets $Z_{i}$ over $K_{i}$, whose translates are nested. It follows that the
new family $Z_{1},\ldots,Z_{k},X_{k+1},\ldots,X_{n}$ is automatically in good
enough position. In topological terms, replacing the $Y_{i}$'s by the $Z_{i}%
$'s corresponds to starting with some simple closed curves on a surface such
that each pair has intersection number zero, and then replacing the curves by
homotopic but disjoint simple closed curves. This completes the proof of the
first part of the lemma.

In order to prove the second part of the lemma, it will be convenient to
consider first the graph of groups structure $\Gamma(Z_{1},\ldots,Z_{k}:G)$.
When we replaced $X_{i}$ by $Y_{i}$, we obtained a splitting $\sigma_{i}$ to
which $Y_{i}$ is associated. Lemma 2.3 of \cite{SS} implies that if two
splittings of $G$ have equivalent associated almost invariant sets, then the
splittings are conjugate. Thus the splitting $\sigma_{i}$ is unique up to
conjugacy. Now Theorem \ref{Theorem2.5ofSS} tells us that $\Gamma(Z_{1}%
,\ldots,Z_{k}:G)$ is determined by the conjugacy classes of the splittings
$\sigma_{i}$. It follows that $\Gamma(Z_{1},\ldots,Z_{k}:G)$ is independent of
the choices of the $Z_{i}$'s.

A more useful way of putting this is the following. Suppose that $W_{1}%
,\ldots,W_{k}$ are chosen in the same way as $Z_{1},\ldots,Z_{k}$. Then the
natural bijection between the set $F(Z)$ of all the translates of the $Z_{i}%
$'s and the set $F(W)$ of all the translates of the $W_{i}$'s is order
preserving. (If some of these splittings are conjugate, we may need to permute
some of the $W_{i}$'s to achieve this.) Now we consider the sets $E(Z)$ and
$E(W)$ obtained from $E$ by replacing each $X_{i}$ by $Z_{i}$ or by $W_{i}$,
$1\leq i\leq k$. The natural bijection between these two sets extends the
bijection between $F(Z)$ and $F(W)$. It is the identity on the translates of
$X_{k+1},\ldots,X_{n}$ and so trivially preserves the partial order on the
translates of these elements. Thus it remains to check that our bijection is
order preserving when we compare a translate of one of these elements with a
translate of some $Z_{j}$. Suppose, for example, that $gZ_{j}\leq hX_{l}$,
where $l>k$. We know that $gW_{j}$ and $hX_{l}$ must be comparable, as the
family $W_{1},\ldots,W_{k},X_{k+1},\ldots,X_{n}$ is in good enough position.
As $Z_{j}$ and $W_{j}$ are equivalent, we must have either $gW_{j}\leq hX_{l}$
or $hX_{l}\subset gW_{j}$. If the second case occurs then $gZ_{j}$, $gW_{j}$
and $hX_{l}$ must all be equivalent. But this would imply that $X_{l}$ was
also isolated, which contradicts our definition of $k$ and the fact that
$l>k$. It follows that $gW_{j}\leq hX_{l}$. This shows that the natural
bijection between the sets $E(Z)$ and $E(W)$ is order preserving, and hence
that $\Gamma(Z_{1},\ldots,Z_{k},X_{k+1},\ldots,X_{n}:G)$ is naturally
isomorphic to $\Gamma(W_{1},\ldots,W_{k},X_{k+1},\ldots,X_{n}:G)$ as required.
\end{proof}

The above result raises the more general question of how our regular
neighbourhood construction changes when one replaces $X_{i}$'s, which need not
be isolated, by equivalent sets. This is an important question because usually
one is not very interested in a particular almost invariant set but rather in
its equivalence class. Our next result shows that our regular neighbourhood
construction depends only on the equivalence classes of the $X_{i}$'s.

\begin{lemma}
\label{equivalentsetshavesameregnbhd}Let $G$ denote a finitely generated
group, and let $H_{1},\ldots,H_{n}$ be finitely generated subgroups of $G$.
For each $i\geq1$, let $X_{i}$ be a nontrivial $H_{i}$-almost invariant subset
of $G$ and let $W_{i}$ be a nontrivial $K_{i}$-almost invariant subset of $G$
which is equivalent to $X_{i}$. Then $\Gamma(X_{1},\ldots,X_{n}:G)$ and
$\Gamma(W_{1},\ldots,W_{n}:G)$ are naturally isomorphic.
\end{lemma}

\begin{proof}
Let $E$ denote the set of all translates of the $X_{i}$'s and let $F$ denote
the set of all translates of the $W_{i}$'s. As the stabilisers of $X_{i}$ and
$W_{i}$ are commensurable but need not be equal, there is no natural map
between $E$ and $F$, but we will show that there is a natural bijection
between the CCC's of $\overline{E}$ and $\overline{F}$, which preserves betweenness.

As in the preceding lemma, we will assume that $X_{1},\ldots,X_{k}$ are the
only isolated elements of the $X_{i}$'s. Thus $W_{1},\ldots,W_{k}$ are the
only isolated elements of the $W_{i}$'s. As in the proof of the preceding
lemma, we can replace $X_{1},\ldots,X_{k}$ by equivalent sets $Z_{1}%
,\ldots,Z_{k}$ such that the translates of the $Z_{i}$'s are nested. It is
then automatic that the family $Z_{1},\ldots,Z_{k},X_{k+1},\ldots,X_{n}$ is in
good enough position. We can also replace $W_{1},\ldots,W_{k}$ by the same
sets $Z_{1},\ldots,Z_{k}$. The preceding lemma allows us to identify
$\Gamma(X_{1},\ldots,X_{n}:G)$ with $\Gamma(Z_{1},\ldots,Z_{k},X_{k+1}%
,\ldots,X_{n}:G)$ and to identify $\Gamma(W_{1},\ldots,W_{n}:G)$ with
$\Gamma(Z_{1},\ldots,Z_{k},W_{k+1},\ldots,W_{n}:G)$. Thus, by changing
notation, we can suppose that in our original families $X_{1},\ldots,X_{n}$
and $W_{1},\ldots,W_{n}$, we have $X_{i}=W_{i}$, for $1\leq i\leq k$, and the
set of all translates of all the $X_{i}$'s, $1\leq i\leq k$, is nested.

Now if $U$ is an isolated element of $E$, the CCC which contains $\overline
{U}$ consists only of $\overline{U}$. Conversely, if a CCC consists of a
single element $\overline{U}$ of $\overline{E}$, then $U$ must be isolated. We
call such a CCC \textit{isolated}. We trivially have a bijection between the
isolated CCC's of $\overline{E}$ and the isolated CCC's of $\overline{F}$.
Further, this bijection preserves betweenness for isolated CCC's, possibly
after permuting equivalent $X_{i}$'s.

Now consider a non-isolated $X_{i}$. If $gX_{i}$ is equivalent to $X_{i}$,
then $g\overline{X_{i}}$ must also lie in the CCC $[\overline{X_{i}}]$ of
$\overline{E}$. As $W_{i}$ is equivalent to $X_{i}$, it also follows that
$g\overline{W_{i}}$ must lie in the CCC $[\overline{W_{i}}]$ of $\overline{F}%
$. If $kX_{j}$ crosses $X_{i}$, then $kW_{j}$ crosses $W_{i}$. It follows that
if $S_{ij}$ denotes the collection of elements $s$ of $G$ such that $sX_{j}$
lies in $[\overline{X_{i}}]$, then $S_{ij}$ must also equal the collection of
elements $s$ of $G$ such that $sW_{j}$ lies in $[\overline{W_{i}}]$. This
yields a natural bijection between the non-isolated CCC's of $\overline{E}$
and the non-isolated CCC's of $\overline{F}$. Further, it is clear that this
bijection preserves betweenness for non-isolated CCC's.

It follows that we have a natural bijection between the CCC's of $\overline
{E}$ and the CCC's of $\overline{F}$, and that this preserves betweenness
except possibly when we consider a mixture of isolated CCC's and non-isolated
CCC's. As in the proof of the previous lemma, it is easy to see that
betweenness must be preserved here also. It follows that $\Gamma(X_{1}%
,\ldots,X_{n}:G)$ and $\Gamma(W_{1},\ldots,W_{n}:G)$ are naturally isomorphic,
as required.
\end{proof}

Our regular neighbourhood construction always expresses $G$ as the fundamental
group of a graph $\Gamma(X_{1},\ldots,X_{n}:G)$ of groups, but this graph may
consist of a single point. This occurs precisely when the set $\overline{E}$
has only one cross-connected component, so that $P$ and hence $T$ consists of
a single point. We give two examples of situations where this will occur. Our
first example is from the topology of $3$-manifolds and is due to Rubinstein
and Wang \cite{R-W}.

\begin{example}
In this case $n=1$, and we denote $X_{1}$ by $X$ and the stabiliser of $X$ by
$H$. The group $G$ is the fundamental group of a closed graph manifold $M$,
and the subgroup $H$ is isomorphic to the fundamental group of a closed
surface $F$. One can choose a $\pi_{1}$-injective map $f:F\rightarrow M$ so
that $f_{\ast}\pi_{1}(F)=H$, and $F$ lifts to an embedding in the cover
$M_{F}$ of $M$ with fundamental group equal to $H$. Considering one side of
$F$ in $M_{F}$ determines the $H$-almost invariant subset $X$ of $G$.
Rubinstein and Wang show that, for many choices of the manifold $M$, the
surface $F$ cannot lift to an embedding in any finite cover of $M$. They do
this by showing that the pre-image of $F$ in the universal cover
$\widetilde{M}$ of $M$ consists of a family of embedded planes such that any
two cross in the sense of \cite{FHS}. This implies that any two distinct
translates of $X$ cross, so that the set $\overline{E}$ has only one
cross-connected component.
\end{example}

Our second example is also rather special.

\begin{example}
Let $H$ denote any finitely generated group, and let $G$ denote the group
$H\times\mathbb{Z}$. Let $X$ denote the $H$-almost invariant subset of $G$
associated to the splitting of $G$ as the HNN-extension $H\ast_{H}$. Thus the
translates of $X$ by $G$ are all equivalent to $X$. In particular, none of
these translates cross each other. Now suppose that there is a subgroup $K$ of
$G$ and a $K$-almost invariant subset $Y$ of $G$ such that $Y$ crosses $X$.
Then $Y$ crosses every translate of $X$, and hence also $X$ crosses every
translate of $Y$. It follows that if we let $E$ denote the set of all
translates of $X$ and $Y$, then the set $\overline{E}$ has only one
cross-connected component. A simple way to generate such examples of $K$ and
$Y$ is to choose $H=A\ast_{C}B$, choose $K=C\times\mathbb{Z}$ and to choose
$Y$ to be associated to the splitting $G=(A\times\mathbb{Z})\ast_{K}%
(B\times\mathbb{Z})$.
\end{example}

We have now completed discussing our regular neighbourhood construction when
one is given a finite family of nontrivial almost invariant subsets of a group
$G$. We will end this section by discussing what happens if one is given an
infinite collection of such subsets, as this will play an important role in
this paper. At first glance this may seem to be a very unreasonable thing to
consider. In topology, one never discusses regular neighbourhoods of infinite
collections of submanifolds. But if one considers a subsurface $N$ of a
surface $M$, then $N$ contains curves representing each element of $\pi
_{1}(N)$, and we want to regard $N$ as a regular neighbourhood of this
infinite family of curves. Of course, such an idea cannot make sense for
arbitrary infinite families of curves. For example, an infinite collection of
disjoint essential circles in $M$ cannot reasonably be said to have a regular neighbourhood.

Now let $G$ denote a finitely generated group with a family of finitely
generated subgroups $\{H_{\lambda}\}_{\lambda\in\Lambda}$. For each
$\lambda\in\Lambda$, let $X_{\lambda}$ denote a nontrivial $H_{\lambda}%
$-almost invariant subset of $G$. We will proceed as we did earlier in this
section, and just note the differences in the case when $\Lambda$ is infinite.

As before, we let $E$ denote the collection of all translates of the
$X_{\lambda}$'s and their complements. First we will assume that the
$X_{\lambda}$'s are in good position, i.e. that $E$ satisfies Condition (*).
This allows us to define the partial order $\leq$ on $E$ exactly as before.
The first and crucial difference between the infinite case and the finite case
occurs when we consider Lemma \ref{posatisfiesDunwoody}. While parts 1) and 3)
still hold, part 2) need not hold, i.e. $E$ need not be discrete, so that
there may be elements $U$ and $V$ of $E$ with infinitely many elements of $E$
between them. The analogous situation occurs in topology if one considers
infinitely many disjoint simple closed curves on a surface. However, we can
still define the set $\overline{E}$ of pairs $\{X,X^{\ast}\}$ for $X\in E$,
and can define $P$ to be the collection of all CCC's of $\overline{E}$.
Further the arguments following Lemma \ref{posatisfiesDunwoody} still apply,
and show that the idea of betweenness can be defined on $P$ as before. Now we
consider the proof of Theorem \ref{Pisapretree}. The argument for the first
part which asserts that $P$ is a pretree remains correct. But the second part
which asserts that $P$ is discrete depends on the discreteness of $E$, and so
$P$ may not be discrete. In fact, if there are infinitely many $X_{\lambda}$'s
which are all equivalent, then $P$ will not be discrete. However, if it
happens that $P$ is discrete, then as before $P$ can be embedded in a $G$-tree
$T$ and $G\backslash T$ is a graph of groups structure for $G$ which we call a
regular neighbourhood of the $X_{\lambda}$'s, and denote by $\Gamma
(\{X_{\lambda}\}_{\lambda\in\Lambda}:G)$. Of course, the $V_{0}$-vertex groups
of this graph need not be finitely generated, and it appears that it may have
infinitely many $V_{0}$-vertices. However, as $G$ is finitely generated, it
will follow from Proposition \ref{Tisminimal} that $\Gamma(\{X_{\lambda
}\}_{\lambda\in\Lambda}:G)$ is always a finite graph.

Having dealt with the case when the $X_{\lambda}$'s are in good position, we
say that the family $\{X_{\lambda}\}_{\lambda\in\Lambda}$ is \textit{in good
enough position }if whenever we find incomparable elements $U$ and $V$ of $E$
which do not cross, there is some element $W$ of $E$ which crosses them. As in
the finite case, the above discussion applies equally well if the $X_{\lambda
}$'s are in good enough position. If this condition does not hold, we consider
the proof of Lemma \ref{candefineregnbhdwhennotingoodenoughposn}. The argument
works exactly as before, so long as $E$ has only finitely many $G$-orbits of
isolated elements. If $E$ has infinitely many such, then $P$ cannot be
discrete. For if $P$ is discrete, then $\Gamma(\{X_{\lambda}\}_{\lambda
\in\Lambda}:G)$ can be defined and must be finite, by Proposition
\ref{Tisminimal}. But $\Gamma(\{X_{\lambda}\}_{\lambda\in\Lambda}:G)$ would
have to have infinitely many isolated $V_{0}$-vertices, one for each $G$-orbit
of isolated elements. Note that it would be simpler to replace the condition
that $E$ has only finitely many $G$-orbits of isolated elements by the
condition that only finitely many of the $X_{\lambda}$'s are isolated. But
this second condition is more restrictive than needed.

Finally the proof of Lemma \ref{equivalentsetshavesameregnbhd} still applies
to show that our construction depends only on the equivalence classes of the
$X_{\lambda}$'s.

We summarise our conclusions as follows. Suppose that we are given a finitely
generated group $G$, finitely generated subgroups $H_{\lambda}$, $\lambda
\in\Lambda$, and a nontrivial $H_{\lambda}$-almost invariant subset
$X_{\lambda}$ of $G$. Suppose that $E$ has only finitely many $G$-orbits of
isolated elements, so that after replacing the $X_{\lambda}$'s by equivalent
sets we can assume that they are in good enough position. Then one can define
the idea of betweenness on the set $P\;$of all CCC's of $\overline{E}$, and
$P$ is a pretree. If $P$ is discrete, one can construct an algebraic regular
neighbourhood $\Gamma(\{X_{\lambda}\}_{\lambda\in\Lambda}:G)$ of the
$X_{\lambda}$'s, and this depends only on the equivalence classes of the
$X_{\lambda}$'s.

\section{Enclosing\label{enclosing}}

In this section we will consider graphs of groups in general. We will discuss
the idea of a vertex of a graph of groups enclosing an almost invariant set.
In the following section, we will apply these ideas to our regular
neighbourhood construction.

Let $\Gamma$ be a graph of groups and write $\pi_{1}(\Gamma)$ for the
fundamental group of $\Gamma$. We emphasise that although we will mostly
consider situations where $\pi_{1}(\Gamma)$ is finitely generated, we will not
assume that the edge and vertex groups of $\Gamma$ are finitely generated
unless this is specifically stated. To avoid some degeneracy phenomena, we
will often assume that $\Gamma$ is minimal, meaning that for any proper
connected subgraph $K$, the natural inclusion of $\pi_{1}(K)$ into $\pi
_{1}(\Gamma)$ is not an isomorphism. Note that if $\Gamma$ is minimal and
$\pi_{1}(\Gamma)$ is finitely generated, then $\Gamma$ must be finite. If
$\Gamma^{\prime}$ is a (not necessarily connected) subgraph of $\Gamma$, we
say that a graph of groups structure $\Gamma_{1}$ is \textit{obtained from
}$\Gamma$\textit{ by collapsing }$\Gamma^{\prime}$ if the underlying graph of
$\Gamma_{1}$ is obtained from $\Gamma$ by collapsing each component of
$\Gamma^{\prime}$ to a point. In addition, if $p:\Gamma\mathcal{\rightarrow
}\Gamma_{1}$ denotes the natural projection map, we require that each vertex
$v$ of $\Gamma_{1}$ has associated group equal to $\pi_{1}(p^{-1}(v))$. These
conditions imply that $\pi_{1}(\Gamma)$ and $\pi_{1}(\Gamma_{1})$ are
naturally isomorphic. Two important special cases of this construction occur
when we consider an edge $e$ of $\Gamma$. If the subgraph $\Gamma^{\prime}$
equals $e$, we say that $\Gamma_{1}$ is \textit{obtained from }$\Gamma
$\textit{ by collapsing }$e$. If the subgraph $\Gamma^{\prime}$ equals the
complement of the interior of $e$, then $\Gamma_{1}$ has a single edge which
determines a splitting $\sigma$ of $\pi_{1}(\Gamma)$, and we call $\sigma$
\textit{the splitting of }$\pi_{1}(\Gamma)$\textit{ associated to} $e$. Note
that so long as we assume that $\Gamma$ is minimal then $\sigma$ really is a
splitting, i.e. a nontrivial decomposition of $\pi_{1}(\Gamma)$. Such
splittings of $G$ will be referred to as the edge splittings of $\Gamma$.

It will be very convenient to introduce some terminology to describe the
process which is the reverse of collapsing an edge.

\begin{definition}
\label{defnofenclosingsplitting}If a graph of groups structure $\Gamma_{1}$
for a group $G$ is obtained from a graph of groups structure $\Gamma$ by
collapsing an edge $e$, and if $e$ projects to the vertex $v_{1}$ of
$\Gamma_{1}$, we will say that $\Gamma$ is a \textsl{refinement} of
$\Gamma_{1}$ obtained by splitting at the vertex $v_{1}$.

We will also say that the vertex $v_{1}$ of $\Gamma_{1}$ \textsl{encloses}%
\textit{ the splitting }$\sigma$ associated to $e$.
\end{definition}

To understand the reasons for our terminology, the reader should consider a
subsurface $N$ of a surface $M$, and let $\Gamma_{1}$ be the graph of groups
structure for $G=\pi_{1}(M)$ determined by $\partial N$. Let $C$ be some
simple closed curve in $N$, and let $\Gamma$ be the graph of groups structure
for $G$ determined by $\partial N\cup C$. If $e$ denotes the edge of $\Gamma$
which corresponds to $C$, then $\Gamma_{1}$ is obtained from $\Gamma$ by
collapsing $e$, and the vertex $v_{1}$ of $\Gamma_{1}$ corresponds to the
component of $N$ which contains $C$. Thus saying that $v_{1}$ encloses the
splitting of $G$ associated to $C$ mirrors the fact that $N$ contains $C$.

Next we introduce a little more terminology. We will say that a vertex $v$ of
$\Gamma$ is \textit{redundant} if it has valence at most two, it is not the
vertex of a loop, and each edge group includes by an isomorphism into the
vertex group at $v$. If $\Gamma$ is minimal, then a redundant vertex must have
valence two. Clearly, these two edges determine conjugate edge splittings of
$G$. Conversely, it is easy to see that if $\Gamma$ has two edges with
conjugate edge splittings, then these edges are the end segments of a path all
of whose interior vertices are redundant. If $\Gamma$ has a redundant vertex
$v$, we can amalgamate the two edges incident to $v$ into a single edge to
obtain a new graph of groups structure for $G$. If $\Gamma$ is finite, we can
repeat this to obtain a graph of groups structure $\Gamma^{\prime}$ for $G$
with no redundant vertices. Clearly $\Gamma$ is obtained from $\Gamma^{\prime
}$ by subdividing some edges.

In Definition \ref{defnofenclosingsplitting}, we defined what it means for a
vertex $v$ of a graph of groups $\Gamma$ to enclose a splitting of $G=\pi
_{1}(\Gamma)$. Next we want to extend this notion to define what it means for
$v$ to enclose an almost invariant subset of $G$. This is meant to be an
analogue of the topological idea of a subsurface containing a possibly
singular curve. In order to avoid problems with conjugates, it is better to
consider a $G$-tree $T$ rather than the quotient graph of groups $\Gamma$. The
condition that $\Gamma=G\backslash T$ be minimal in the sense above is
equivalent to the condition that $T$ have no proper $G$-invariant subtree.
Such a $G$-tree is also called minimal. Note that any $G$-tree possesses a
minimal subtree $T_{0}$. If $G$ fixes more than one vertex of $T$, then
$T_{0}$ is not unique, but otherwise $T_{0}$ is unique and can be described
simply as the intersection of all the $G$-invariant subtrees of $T$. Note that
a minimal $G$-tree has no vertices of valence one. For if a $G$-tree has such
vertices, one can obtain a proper subtree by simply removing each such vertex
together with the interior of the incident edge.

We recall some notation from section \ref{prelim}. An oriented edge $s$ of a
tree $T$ determines a natural partition of $V(T)$ into two sets, namely the
vertices of the two subtrees obtained by removing the interior of $s$ from
$T$. Let $Y_{s}$ denote the collection of all the vertices of the subtree
which contains the terminal vertex $v$ of $s$, and let $Y_{s}^{\ast}$ denote
the complementary collection of vertices. If a group $G$ acts without
inversions on $T$, and $\varphi:G\rightarrow V(T)$ is a $G$-equivariant map,
then we have the sets $Z_{s}=\varphi^{-1}(Y_{s})$ and $Z_{s}^{\ast}%
=\varphi^{-1}(Y_{s}^{\ast})$. Lemma \ref{pre-imageofhalftreeisalmostinvariant}
shows that, if $S$ denotes the stabiliser of $s$, then $Z_{s}$ is $S$-almost
invariant, and its equivalence class is independent of the choice of the map
$\varphi$.

We now define enclosing of almost invariant sets.

\begin{definition}
Let $A$ be a nontrivial $H$-almost invariant subset of a group $G$ and let $v$
a vertex of a $G$-tree $T$. Choose a $G$-equivariant map $\varphi:G\rightarrow
V(T)$. For each edge $s$ of $T$ this determines the subsets $Z_{s}$ and
$Z_{s}^{\ast}$ of $G$ as above. We say that the vertex $v$ \textsl{encloses}
$A$, if for all edges $s$ of $T$ which are incident to $v$ and directed
towards $v$, we have $A\cap Z_{s}^{\ast}$ or $A^{\ast}\cap Z_{s}^{\ast}$ is small.
\end{definition}

\begin{remark}
This definition is independent of the choice of $\varphi$, because changing
$\varphi$ replaces each $Z_{s}$ by an equivalent almost invariant set. Also if
$B$ is an almost invariant set equivalent to $A$, then $v$ encloses $A$ if and
only if $v$ encloses $B$.

Note that it does not make much sense to consider enclosing of trivial almost
invariant subsets of $G$, because any such subset of $G$ would automatically
be enclosed by every vertex of $T$.

It would seem more natural to say that $v$ encloses $A$, if for all edges $s$
of $T$ which are incident to $v$ and directed towards $v$, we have $A^{\ast
}\geq Z_{s}^{\ast}$ or $A\geq Z_{s}^{\ast}$, but we want to ensure that any
set equivalent to $Z_{s}$ or $Z_{s}^{\ast}$ is enclosed by $v$, and not all
such sets are comparable with $Z_{s}$. See part 3) of Lemma
\ref{somefactsaboutenclosing} for precise statements.
\end{remark}

It will also be convenient to define enclosing by a vertex of a graph of groups.

\begin{definition}
Let $A$ be a nontrivial $H$-almost invariant subset of a group $G$, let $T$ be
a $G$-tree and let $\Gamma$ denote the associated graph of groups structure
for $G$ with underlying graph $G\backslash T$. We say that a vertex $w$ of
$\Gamma$ \textsl{encloses} $A$ if there is a vertex $v$ of $T$ which encloses
$A$ and projects to $w$.
\end{definition}

\begin{remark}
If $w$ encloses $A$, then it also encloses any translate of $A$ and any almost
invariant set equivalent to $A$.
\end{remark}

We now have two natural ideas of what it means for a vertex $w$ of $\Gamma$ to
enclose a splitting $\sigma$ of $G$ over a subgroup $H$. One is given in
Definition \ref{defnofenclosingsplitting}, and the other is that a $H$-almost
invariant set associated to $\sigma$ is enclosed by $w$. In Lemma
\ref{twoideasofenclosingaresameifHisfg}, we will show that these ideas are
equivalent when $H$ is finitely generated. In Lemma
\ref{twoideasofenclosingaresame} we will be able to show that this equivalence
holds even when $H$ is not finitely generated. We will need the following
basic properties of enclosing. Again we emphasise that we are not assuming
that the subgroup $H$ of $G$ is finitely generated.

\begin{lemma}
\label{somefactsaboutenclosing}Let $A$ be a nontrivial $H$-almost invariant
subset of a group $G$ and let $v$ a vertex of a $G$-tree $T$. Then the
following statements all hold:

\begin{enumerate}
\item $A$ is enclosed by $v$ if and only if $A^{\ast}$ is enclosed by $v$.

\item If $s$ is an edge of $T$ with stabiliser $S$, and if $Z_{s}$ is a
nontrivial $S$-almost invariant subset of $G$, then $Z_{s}$ is enclosed by
each of the two vertices to which $s$ is incident.

\item $A$ is enclosed by $v$ if and only if either

\begin{enumerate}
\item $A$ is equivalent to $Z_{s}$ or $Z_{s}^{\ast}$ for some edge $s$
incident to $v$, or

\item we have $A^{\ast}\geq Z_{s}^{\ast}$ or $A\geq Z_{s}^{\ast}$ for every
edge $s$ incident to $v$ and oriented towards $v$.
\end{enumerate}

\item If $A$ is enclosed by $v$, then for any edge $t$ of $T$ which is
oriented towards $v$, we have $A\cap Z_{t}^{\ast}$ or $A^{\ast}\cap
Z_{t}^{\ast}$ is small.

\item If $B$ is another nontrivial $H$-almost invariant subset of $G$, and
both $A$ and $B$ are enclosed by $v$, then $A\cup B$, $A\cap B$ and $A+B$ are
also $H$-almost invariant and, if nontrivial, each is enclosed by $v$.

\item If $s$ is an edge of $T$, and $s$ has stabiliser $S$, then $Z_{s}$ is a
nontrivial $S$-almost invariant subset of $G$ if and only if $s$ lies in the
minimal subtree of $T$.

\item If the vertex $v$ is not fixed by $G$, then $v$ encloses some nontrivial
$H$-almost invariant subset of $G$ if and only if $v$ lies in the minimal
subtree of $T$.

\item If $A$ is enclosed by two distinct vertices $u$ and $v$ of $T$, then $A$
is equivalent to $Z_{s}$ or to $Z_{s}^{\ast}$ for each edge $s$ on the path
$\lambda$ joining $u$ and $v$. Further $\lambda$ is contained in the minimal
subtree $T_{0}$ of $T$, and each interior vertex of $\lambda$ has valence $2$
in $T_{0}$.
\end{enumerate}
\end{lemma}

\begin{proof}
1) and 2) are both trivial.

3) First suppose that $A$ is enclosed by $v$, so that for all edges $s$ of $T$
which are incident to $v$ and directed towards $v$, we have $A\cap Z_{s}%
^{\ast}$ or $A^{\ast}\cap Z_{s}^{\ast}$ is small. If $A$ and $Z_{s}$ are not
comparable for some $s$, then two of the four sets $A^{(\ast)}\cap
Z_{t}^{(\ast)}$ must be small so that $A$ is equivalent to $Z_{s}$ or
$Z_{s}^{\ast}$. Otherwise, we have $A^{\ast}\geq Z_{s}^{\ast}$ or $A\geq
Z_{s}^{\ast}$ for every edge $s$ incident to $v$ and oriented towards $v$, as required.

Conversely, if $A^{\ast}\geq Z_{s}^{\ast}$ or $A\geq Z_{s}^{\ast}$ for every
edge $s$ incident to $v$ and oriented towards $v$, then $A\cap Z_{s}^{\ast}$
or $A^{\ast}\cap Z_{s}^{\ast}$ is small for each such edge $s$, so that $A$ is
enclosed by $v$. And if $A$ is equivalent to $Z_{s}$ or $Z_{s}^{\ast}$ for
some edge $s$ incident to $v$, then part 2) implies that $A$ is enclosed by
$v$. The result follows.

4) Consider the oriented path in $T$ which joins $t$ to $v$. Let $s$ be the
edge of this path incident to $v$. As $s$ is oriented towards $v$ and $A$ is
enclosed by $v$, it follows that $A\cap Z_{s}^{\ast}$ or $A^{\ast}\cap
Z_{s}^{\ast}$ is small. As $t\leq s$ in the natural ordering on oriented edges
of $T$, it follows that $Z_{s}\subset Z_{t}$, so that $A\cap Z_{t}^{\ast}$ or
$A^{\ast}\cap Z_{t}^{\ast}$ is small, as required.

5) It is clear that $A\cup B$, $A\cap B$ and $A+B$ are $H$-almost invariant.
Let $s$ be an edge of $T$ incident to $v$ and oriented towards $v$. First we
suppose that neither $A$ nor $B$ is equivalent to $Z_{s}$ or $Z_{s}^{\ast}$.
Thus part 3) tells us that there are four cases as $A^{\ast}\geq Z_{s}^{\ast}$
or $A\geq Z_{s}^{\ast}$, and $B^{\ast}\geq Z_{s}^{\ast}$ or $B\geq Z_{s}%
^{\ast}$.

If $A^{\ast}\geq Z_{s}^{\ast}$ and $B^{\ast}\geq Z_{s}^{\ast}$, then $(A\cup
B)^{\ast}=A^{\ast}\cap B^{\ast}\geq Z_{s}^{\ast}$, so that $(A\cup B)^{\ast
}\geq Z_{s}^{\ast}$. As $(A\cap B)^{\ast}\supset(A\cup B)^{\ast}$ and
$(A+B)^{\ast}\supset(A\cup B)^{\ast}$, we see that also $(A\cap B)^{\ast}\geq
Z_{s}^{\ast}$, and $(A+B)^{\ast}\geq Z_{s}^{\ast}$.

If $A^{\ast}\geq Z_{s}^{\ast}$ and $B\geq Z_{s}^{\ast}$, then $A\cup B\supset
B\geq Z_{s}^{\ast}$, $(A\cap B)^{\ast}\supset A^{\ast}\geq Z_{s}^{\ast}$, and
$A+B\supset A^{\ast}\cap B\geq Z_{s}^{\ast}$. Thus $A\cup B\geq Z_{s}^{\ast}$,
$(A\cap B)^{\ast}\geq Z_{s}^{\ast}$, and $A+B\geq Z_{s}^{\ast}$. If $A\geq
Z_{s}^{\ast}$ and $B^{\ast}\geq Z_{s}^{\ast}$, we obtain the same inequalities
by reversing the roles of $A$ and $B$.

If $A\geq Z_{s}^{\ast}$ and $B\geq Z_{s}^{\ast}$, then $A\cap B\geq
Z_{s}^{\ast}$. As $A\cup B\supset A\cap B$ and ($A+B)^{\ast}\supset A\cap B$,
we see that also $A\cup B\geq Z_{s}^{\ast}$ and ($A+B)^{\ast}\geq Z_{s}^{\ast
}$.

As the above applies to every such edge $s$, it follows that in each of the
four cases, each of $A\cup B$, $A\cap B$ and $A+B$ is enclosed by $v$.

Finally, if $A$ is equivalent to $Z_{s}$ or $Z_{s}^{\ast}$, we can replace $A$
by $Z_{s}$ or $Z_{s}^{\ast}$, as appropriate, and similarly if $B$ is
equivalent to $Z_{s}$ or $Z_{s}^{\ast}$. Then we obtain the same inequalities
as above. Thus in all cases, each of $A\cup B$, $A\cap B$ and $A+B$ is
enclosed by $v$.

6) Let $T_{0}$ denote the minimal subtree of $T$. If $G$ fixes more than one
vertex of $T$ so that $T_{0}$ is not unique, we let $T_{0}$ denote one of the
vertices fixed by $G$.

First suppose that $s$ does not lie in $T_{0}$. Note that in the special case
when $T_{0}$ is not unique, this condition is automatic as then $T_{0}$ has no
edges. We will show that $Z_{s}$ must be trivial. We can assume that $s$ is
oriented towards $T_{0}$. As we are free to choose the map $\varphi
:G\rightarrow V(T)$, we will choose it so that $\varphi(e)$ lies in $T_{0}$.
This implies that $\varphi(G)$ also lies in $T_{0}$. Hence $Z_{s}=G$, so that
$Z_{s}$ is a trivial $S$-almost invariant set as claimed.

Now suppose that $s$ lies in $T_{0}$. Note that in this case, $T_{0}$ must be
unique. Lemma 3.3 of \cite{Scott:Intersectionnumbers} tells us that $s$ lies
in $T_{0}$ if and only if there exists an element $g$ of $G$ such that $s$ and
$gs$ are distinct and coherently oriented. (Two oriented edges in $T$ are
coherently oriented if there is an oriented path in $T$ which begins with one
and ends with the other.) By repeatedly applying $g$ or $g^{-1}$, we see that
on each side of $s$ in $T$ there are translates of $s$ which are arbitrarily
far from $s$. It follows that $\varphi(G)$, which equals the orbit of
$\varphi(e)$, contains points in $Y_{s}$ and in $Y_{s}^{\ast}$ which are
arbitrarily far from $s$. This immediately implies that $Z_{s}$ must be
nontrivial, as required.

7) We choose $T_{0}$ as in part 6), and recall that we are assuming that $v$
is not fixed by $G$. It follows that if $v$ lies in $T_{0}$, there is an edge
$s$ of $T_{0}$ incident to $v$. Part 6) implies that $Z_{s}$ is a nontrivial
$S$-almost invariant set and part 2) implies that $Z_{s}$ is enclosed by $v$.

If $v$ does not lie in $T_{0}$ and $v$ encloses some nontrivial $H$-almost
invariant subset $A$ of $T$, we will obtain a contradiction. Let $s$ denote
the edge of $T$ which is incident to $v$ and on the path joining $v$ to
$T_{0}$ and we choose $s$ to be oriented towards $v$. After choosing
$\varphi:G\rightarrow V(T)$ such that $\varphi(e)=v$, we will have
$Z_{s}^{\ast}=G$. But as $A$ is enclosed by $v$, we have $A\cap Z_{s}^{\ast}$
or $A^{\ast}\cap Z_{s}^{\ast}$ is small, which implies that $A$ or $A^{\ast}$
is small. It follows that $A$ is a trivial $H$-almost invariant set, which is
the required contradiction.

8) Let $\lambda$ denote the path in $T$ which joins $u$ and $v$. Let $l$
denote the edge of $\lambda$ incident to $u$ and oriented towards $u$, and let
$m$ denote the edge of $\lambda$ incident to $v$ and oriented towards $v$. Our
choice of orientations on $l$ and $m$ implies that $Z_{l}^{\ast}\supset Z_{m}%
$. As $A\;$is enclosed by $u$, we know that $A\cap Z_{l}^{\ast}$ or $A^{\ast
}\cap Z_{l}^{\ast}$ is small. Without loss of generality, we can assume that
$A\cap Z_{l}^{\ast}$ is small. As $Z_{l}^{\ast}\supset Z_{m}$, it follows that
$A\cap Z_{m}$ is small. As $A\;$is enclosed by $v$, we know that $A\cap
Z_{m}^{\ast}$ or $A^{\ast}\cap Z_{m}^{\ast}$ is small. But if $A\cap
Z_{m}^{\ast}$ were small, then $A$ would itself be small, which contradicts
our hypothesis that $A$ is nontrivial. It follows that $A^{\ast}\cap
Z_{m}^{\ast}$ must be small. As $A\cap Z_{m}$ and $A^{\ast}\cap Z_{m}^{\ast}$
are both small, it follows that $A$ is equivalent to $Z_{m}^{\ast}$.
Similarly, $A$ must be equivalent to $Z_{l}$. If $n$ denotes any edge of
$\lambda$ oriented towards $u$, we have the inclusions $Z_{l}\subset
Z_{n}\subset Z_{m}^{\ast}$, and it immediately follows that $A$ is also
equivalent to $Z_{n}$. The fact that $\lambda$ is contained in the minimal
subtree $T_{0}$ of $T$ follows at once from part 5).

Finally, let $w$ denote an interior vertex of $\lambda$, and let $r$ and $s$
denote the edges of $\lambda$ incident to $w$ and oriented towards $w$. The
above discussion implies that $Z_{r}$ and $Z_{s}^{\ast}$ are equivalent. In
particular, their stabilisers are commensurable. Let $K$ denote $Stab(r)\cap
Stab(s)$, so that $Z_{r}$ and $Z_{s}^{\ast}$ are equivalent $K$-almost
invariant subsets of $G$. Thus $Z_{r}\cap Z_{s}$ is $K$-finite. If there is
another edge $t$ of $T$ incident to $w$, we orient $t$ towards $w$ also. As
$K$ fixes $w$, each edge $k(t)$ is also incident to $w$ and oriented towards
$w$. Thus $Z_{k(t)}^{\ast}$ is disjoint from $Z_{r}^{\ast}$ and $Z_{s}^{\ast}%
$, and so is contained in $Z_{r}\cap Z_{s}$. Hence $\cup_{k\in K}%
Z_{k(t)}^{\ast}$ is also $K$-finite. As the $Z_{k(t)}^{\ast}$ are disjoint
from each other, it follows that $Z_{t}^{\ast}$ is $Stab(t)$-finite.
Equivalently $Z_{t}^{\ast}$ is a trivial almost invariant set over $Stab(t)$.
Now part 6) implies that $t$ does not lie in the minimal subtree $T_{0}$ of
$T$. Hence every interior vertex of $\lambda$ has valence $2$ in $T_{0}$ which
completes the proof of part 8).
\end{proof}

Now we are ready to show that the two ideas of enclosing a splitting are equivalent.

\begin{lemma}
\label{twoideasofenclosingaresameifHisfg}Let $T$ be a $G$-tree, let $\Gamma$
denote the graph of groups structure for $G$ given by the quotient
$G\backslash T$, and let $w$ denote a vertex of $\Gamma$. Suppose that $A$ is
associated to a splitting $\sigma$ of $G$ over $H$. Then the following
statements hold.

\begin{enumerate}
\item If $\sigma$ is enclosed by $w$, then $A$ is enclosed by $w$.

\item If $H$ is finitely generated and $A$ is enclosed by $w$, then $\sigma$
is enclosed by $w$.
\end{enumerate}
\end{lemma}

\begin{remark}
In Lemma \ref{twoideasofenclosingaresame} we will show that part 2) holds even
when $H$ is not finitely generated.
\end{remark}

\begin{proof}
1) Suppose that $\sigma$ is enclosed by $w$. Thus there is a graph of groups
$\Gamma_{1}$ which is a refinement of $\Gamma$ obtained by splitting at $w$ so
that the extra edge has $\sigma$ as its associated edge splitting. Let
$q:\Gamma_{1}\rightarrow\Gamma$ denote the projection map, and let $e$ denote
the extra edge of $\Gamma_{1}$ so that $q(e)=w$. Let $v$ denote a vertex in
the pre-image of $w$ in the $G$-tree $T$, and let $T_{1}$ denote the universal
covering $G$-tree of $\Gamma_{1}$. Thus $T_{1}$ is obtained from $T$ by
splitting at each vertex in the orbit of $v$. Let $e$ also denote the extra
edge inserted at $v$. Then $T_{1}$ has an induced projection $p:T_{1}%
\rightarrow T$ such that $p(ge)=v$. Pick a $G$-equivariant map $\varphi
_{1}:G\rightarrow T_{1}$. Thus $\varphi=p\circ\varphi_{1}:G\rightarrow T$ is
also $G$-equivariant. Consider the set $E_{1}$ of all the sets $Z_{s}$, for
each oriented edge $s$ of $T_{1}$. If $s$ is not equal to $ge$ for any $g$,
then $p(s)$ is an oriented edge of $T$ and $Z_{s}=Z_{p(s)}$. There is a
translate $ge$ of $e$ such that $Z_{ge}$ is equivalent to $A$ or $A^{\ast}$.
We will choose $e$ so that this translate is $e$ itself. As $E_{1}$ is nested,
we know that for any edge $s$ of $T_{1}-\{e\}$, oriented towards $e$, we have
$Z_{e}\supset Z_{s}^{\ast}$ or $Z_{e}^{\ast}\supset Z_{s}^{\ast}$. Turning to
the tree $T$, it follows that either $A$ is equivalent to $Z_{s}$ or
$Z_{s}^{\ast}$ for some edge $s$ incident to $v$, or we have $A^{\ast}\geq
Z_{s}^{\ast}$ or $A\geq Z_{s}^{\ast}$ for every edge $s$ incident to $v$ and
oriented towards $v$. Now part 3) of Lemma \ref{somefactsaboutenclosing}
implies that $A$ is enclosed by $v$ and hence by $w$.

2) Suppose that $A$ is enclosed by the vertex $v$ of $T$. Recall that the
collection $E$ of all those $Z_{s}$ and $Z_{s}^{\ast}$ which are nontrivial
$S$-almost invariant subsets of $G$ is nested. We now enlarge this set to a
set $F$ by adding $A$ and $A^{\ast}$ and all their translates. As $A$ is
associated to a splitting, its translates are nested. Part 3) of Lemma
\ref{somefactsaboutenclosing} implies that either $A$ is comparable with every
$Z_{s}$ or that $A$ is equivalent to some $Z_{s}$. In\ the first case, we can
apply Theorem \ref{Theorem1.12ofSS} to obtain the required refinement of the
graph of groups $G\backslash T$. (This is where we use the assumption that $H$
is finitely generated.) In the second case, the required refinement can be
constructed by simply subdividing the edge $s$ such that $A$ is equivalent to
$Z_{s}$.
\end{proof}

Recall from part 4) of Lemma \ref{somefactsaboutenclosing} that if $A$ is a
nontrivial $H$-almost invariant subset of a group $G$ which is enclosed by a
vertex $v$ of a $G$-tree $T$, then for each edge $s$ of $T$ which is directed
towards $v$, we have $A\cap Z_{s}^{\ast}$ or $A^{\ast}\cap Z_{s}^{\ast}$ is
small. Thus we have a naturally defined $H$-invariant partition of the edges
of $T$, where one set consists of those $s$ with $A\cap Z_{s}^{\ast}$ small
and the other set consists of those $s$ with $A^{\ast}\cap Z_{s}^{\ast}$
small. This induces a $H$-invariant partition of all the vertices of $T-\{v\}$
as in the following definition.

\begin{definition}
\label{defnofA-sideofv}Let $A$ be a nontrivial $H$-almost invariant subset of
a group $G$ and let $v$ a vertex of a $G$-tree $T$. Suppose that $A$ is
enclosed by $v$, and that $s$ is an edge of $T$ which is directed towards $v$.
We will say that $s$ \textsl{lies on the }$A$\textsl{-side of }$v$ if
$A^{\ast}\cap Z_{s}^{\ast}$ is small.

We will say that a vertex $w$ of $T-\{v\}$ \textsl{lies on the }%
$A$\textsl{-side of }$v$ if the path from $w$ to $v$ ends in an edge $s$ which
lies on the $A$-side of $v$. The collection of all vertices of $T-\{v\}$ which
lie on the $A$-side of $v$ will be denoted by $\Sigma_{v}(A)$ or by
$\Sigma(A)$ if the context is clear.
\end{definition}

It is easy to see that if an edge $s$ of $T$ lies on the $A$-side of $v$, then
the same holds for every edge (and hence every vertex) in the path joining $s$
to $v$. Also if a vertex $w$ of $T-\{v\}$ lies on the $A$-side of $v$, then
the path from $w$ to $v$ consists entirely of edges which lie on the $A$-side
of $v$. Thus the ideas of an edge and a vertex being on the $A$-side of $v$
are compatible. Clearly, every vertex of $T-\{v\}$ lies in $\Sigma_{v}(A)$ or
$\Sigma_{v}(A^{\ast})$, so that these two sets partition the vertices of
$T-\{v\}$. Also $\Sigma_{v}(A)$ and $\Sigma_{v}(A^{\ast})$ are each clearly
$H$-invariant, under the left action of $H$.

To understand the reason for our terminology, think of $A$ as determined by a
closed curve on a surface $M$ which lies inside a subsurface $N$ of $M$, think
of $G$ as being $\pi_{1}(M)$, and think of $T$ as being the $G$-tree
determined by $N$, so that the picture in $T$ corresponds to the picture in
the universal cover of $M$.

It is natural to ask if we can replace $A$ by an equivalent almost invariant
subset $B$ of $G$ such that for all edges $s$ of $T$ which are directed
towards $v$, we have $B\cap Z_{s}$ or $B^{\ast}\cap Z_{s}$ is empty.
Equivalently, can we replace $A$ by $B$ which is nested with respect to every
$Z_{s}$? We will show that this is indeed the case.

Suppose that we are given a nontrivial $H$-almost invariant subset $A$ of $G$
which is enclosed by a vertex $v$ of $T$. The following result shows how to
replace $A$ by a subset $B(A)$ which is nested with respect to every $Z_{s}$.

We define%
\begin{align*}
B(A)  &  =\varphi^{-1}(\Sigma_{v}(A))\cup(A\cap\varphi^{-1}(v)),\\
C(A)  &  =\varphi^{-1}(\Sigma_{v}(A^{\ast}))\cup(A^{\ast}\cap\varphi^{-1}(v)).
\end{align*}

Note that these definitions are clearly equivariant, i.e. $B(kA)=kB(A)$, for
all $k$ in $G$.

\begin{lemma}
\label{enclosedsetscanbechosennested}Let $G$ be a finitely generated group
with a finitely generated subgroup $H$, and a nontrivial $H$-almost invariant
subset $A$. Suppose that $T$ is a $G$-tree, and that $\varphi$ is a
$G$-equivariant map from $G$ to $V(T)$ such that the vertex $\varphi(e)$
encloses some nontrivial $K$-almost invariant subset $U$ of $G$. If $A$ is
enclosed by a vertex $v$ of $T$, then $C(A)=B(A)^{\ast}$, and $B(A)$ is
$H$-almost invariant and is equivalent to $A$. Further, $B(A)$ is nested with
respect to $Z_{s}$, for every oriented edge $s$ of $T$.
\end{lemma}

\begin{remark}
The technical hypothesis on $\varphi(e)$ is automatically satisfied if it lies
in the minimal subtree of $T$, by part 7) of Lemma
\ref{somefactsaboutenclosing}. If $T$ is obtained by our regular neighbourhood
construction, Lemma \ref{XiisenclosedbyitsV0-vertex} shows that this
hypothesis will be satisfied by simply choosing $\varphi(e)$ to be any $V_{0}%
$-vertex of $T$.
\end{remark}

\begin{proof}
To simplify notation, we will write $\Sigma(A)$ in place of $\Sigma_{v}(A)$
throughout this proof.

It is clear from their definitions that $B$ and $C$ are disjoint and that
$B\cup C=G$. Thus $C=B^{\ast}$. It is also clear that $HB=B$, because
$H(\Sigma(A))=\Sigma(A)$. Finally, it is also clear that if $s$ is any edge of
$T$ which is directed towards $v$, and lies on the $A$-side of $v$, then
$B(A)\supset\varphi^{-1}(\Sigma(A))\supset Z_{s}^{\ast}$. If $s$ lies on the
$A^{\ast}$-side of $v$, then $C(A)\supset Z_{s}^{\ast}$. Hence $B(A)$ is
nested with respect to $Z_{s}$, for every oriented edge $s$ of $T$.

It remains to show that $B(A)$ is $H$-almost invariant and is equivalent to
$A$.

Let $w$ denote $\varphi(e)$. As $\varphi(g)=gw$, we see that $\varphi
^{-1}(\Sigma(A))=\{g\in G:gw\in\Sigma(A)\}$. Our hypothesis on $\varphi(e)$
implies that $w$ encloses some nontrivial $K$-almost invariant subset $U$ of
$G$. Hence $\varphi^{-1}(\Sigma(A))\subset\{g\in G:gU^{(\ast)}<A\}$. Now Lemma
\ref{containedinboundednbhdofA} implies that $\varphi^{-1}(\Sigma(A))$ lies in
a bounded neighbourhood of $A$. It follows that $B$ itself is contained in a
bounded neighbourhood of $A$. Similarly $B^{\ast}$ is contained in a bounded
neighbourhood of $A^{\ast}$. It follows that $\delta B$ lies in a bounded
neighbourhood of $\delta A$. As $A$ is $H$-almost invariant, we know that
$\delta A$ projects to a finite subset of the quotient graph $H\backslash
\Gamma$. It follows that $\delta B$ also projects to a finite subset of
$H\backslash\Gamma$, so that $B$ projects to an almost invariant subset of
$H\backslash G$ and hence is $H$-almost invariant. As $B$ is contained in a
bounded neighbourhood of $A$, and $B^{\ast}$ is contained in a bounded
neighbourhood of $A^{\ast}$, it now follows that $B$ is equivalent to $A$,
which completes the proof of the lemma.
\end{proof}

Before our final result of this section, we will need the following simple proposition.

\begin{proposition}
\label{trivialitylemma}Let $G$ be a finitely generated group, and let $X$ be a
$H$-almost invariant subset of $G$ which is contained in a proper subgroup $K$
of $G$. Then $X$ is trivial.
\end{proposition}

\begin{proof}
As $K$ is a proper subgroup of $G$, there is an element $g\in G-K$. As $X$ is
$H$-almost invariant, we must have $X$ and $Xg$ being $H$-almost equal. The
assumption that $X$ is contained in $K$ implies that $X$ and $Xg$ are
disjoint. It follows that $X$ is $H$-finite and hence trivial.
\end{proof}

Now we can prove the following useful result.

\begin{corollary}
\label{nontrivialpartition}Let $G$ be a finitely generated group with finitely
generated subgroups $H$ and $K$, and let $T$ be a minimal $G$-tree which is
not a single point. Let $U$ be a nontrivial $H$-almost invariant subset of $G$
and let $V$ be a nontrivial $K$-almost invariant subset of $G$.

\begin{enumerate}
\item If $U$ is enclosed by a vertex $v$ of $T$, then both $\Sigma_{v}(U)$,
the $U$-side of $v$, and $\Sigma_{v}(U^{\ast})$, the $U^{\ast}$-side of $v$,
are nonempty, so that $U$ determines a nontrivial partition of the vertices of
$T-\{v\}$.

\item If $U$ is enclosed by a vertex $v$ of $T$, then $H$ is contained in
$Stab(\Sigma_{v}(U))$ with finite index.

\item If $U$ and $V$ are enclosed by a vertex $v$ of $T$, and if they
determine the same partition of the vertices of $T-\{v\}$, then $U$ and $V$
are equivalent.
\end{enumerate}
\end{corollary}

\begin{proof}
First fix $v$ and choose the $G$-equivariant map $\varphi:G\rightarrow V(T)$
so that $\varphi(e)=v$.

1) By applying Lemma \ref{enclosedsetscanbechosennested}, we can assume that
$U=B(U)=\varphi^{-1}(\Sigma_{v}(U))\cup(U\cap\varphi^{-1}(v))$. If $\Sigma
_{v}(U)$ is empty, this implies that $U\subset\varphi^{-1}(v)$. Now
$\varphi^{-1}(v)=Stab(v)$, and the assumption that $T$ is minimal and not a
single point implies that $Stab(v)$ must be a proper subgroup of $G$. Thus
Lemma \ref{trivialitylemma} implies that $U$ must be trivial. We conclude that
$\Sigma_{v}(U)$ cannot be empty, and similarly that $\Sigma_{v}(U^{\ast})$
cannot be empty. Thus $U$ determines a nontrivial partition of the vertices of
$T-\{v\}$, as claimed.

2) Clearly $H\subset Stab(\Sigma_{v}(U))$. Further, it is clear that
$\Sigma_{gv}(gU)=g\Sigma_{v}(U)$, for any element $g$ of $G$. Hence if $g\in
Stab(\Sigma_{v}(U))$, then $\Sigma_{gv}(gU)=\Sigma_{v}(U)$. Let $U_{1}$ be a
translate of $U$ enclosed by a vertex of $T$ on the $U$-side of $v$, and let
$U_{2}$ be a translate of $U$ enclosed by a vertex of $T$ on the $U^{\ast}%
$-side of $v$. Then, for any $g\in Stab(\Sigma_{v}(U))$, the set
$g\overline{U}$ lies between $\overline{U_{1}}$ and $\overline{U_{2}}$. As the
number of translates of $\overline{U}$ between $\overline{U_{1}}$ and
$\overline{U_{2}}$ is finite, it follows that the orbit of $U$ under the
action of $Stab(\Sigma_{v}(U))$ is finite, so that $H$ has finite index in
$Stab(\Sigma_{v}(U))$.

3) Apply Lemma \ref{enclosedsetscanbechosennested} so that we can assume that
$U=B(U)=\varphi^{-1}(\Sigma_{v}(U))\cup(U\cap\varphi^{-1}(v))$, and
$V=B(V)=\varphi^{-1}(\Sigma_{v}(V))\cup(V\cap\varphi^{-1}(v))$. Suppose that
$\Sigma_{v}(U)=\Sigma_{v}(V)$. It follows that the symmetric difference $U+V$
of $U\;$and $V$ is contained in $\varphi^{-1}(v)=Stab(v)$. As part 2) implies
that the stabilisers $H_{U}$ and $H_{V}$ of $U\;$and $V$ are commensurable,
both $U$ and $V$ are $H$-almost invariant, where $H=H_{U}\cap H_{V}$. It
follows that $U+V$ is also $H$-almost invariant. As in part 1), $Stab(v)$ must
be a proper subgroup of $G$, and now Proposition \ref{trivialitylemma} shows
that $U+V$ must be trivial, so that $U\;$and $V$ are equivalent as claimed.
\end{proof}

\section{Algebraic Regular Neighbourhoods: Enclosing\label{regnbhds:enclosing}}

In this section, we will apply the results of the previous section to graphs
of groups which are obtained by the regular neighbourhood construction. Let
$G$ denote a finitely generated group with a family of finitely generated
subgroups $\{H_{\lambda}\}_{\lambda\in\Lambda}$. For each $\lambda\in\Lambda$,
let $X_{\lambda}$ denote a nontrivial $H_{\lambda}$-almost invariant subset of
$G$. In section \ref{regnbhds:construction}, we discussed how to construct the
regular neighbourhood $\Gamma(\{X_{\lambda}\}_{\lambda\in\Lambda}:G)$ by
producing a bipartite $G$-tree $T$ whose $V_{0}$-vertices are the CCC's of
$\overline{E}$. Recall that the construction always works when $\Lambda$ is
finite. In order to show that this construction and our ideas about enclosing
all fit together, we need to prove that the $V_{0}$-vertices of $\Gamma
(\{X_{\lambda}\}_{\lambda\in\Lambda}:G)$ enclose the given $X_{\lambda}$'s.

First we will need the following general result. A special case of this result
is proved in \cite{SS}, but is not formulated as a separate statement, so we
give here a brief description of the proof of the result. It is well known to experts.

\begin{lemma}
\label{containedinboundednbhdofA}Let $G$ be a finitely generated group with
finitely generated subgroups $H$ and $K$, a nontrivial $H$-almost invariant
subset $A$ and a nontrivial $K$-almost invariant subset $U$. Then $\{g\in
G:gU^{(\ast)}\leq A\}$ is contained in a bounded neighbourhood of $A$ in the
Cayley graph of $G$.
\end{lemma}

\begin{proof}
The starting point is Lemma \ref{finitenumberofdoublecosets} which tells us
that $\{g\in G:gU$ and $A$ are not nested\} consists of a finite number of
double cosets $HgK$. (This uses the finite generation of $G$, $H$ and $K$.)
Now suppose that $gU\leq A$ but $gU$ is not contained in $A$. As $gU\cap
A^{\ast}$ is small, i.e. it projects to a finite subset in $H\backslash G$, it
follows that $gU$ is contained in a bounded neighbourhood of $A$. Now the fact
that $g$ must lie in a finite number of double cosets $HgK$ implies that there
is a uniform bound on the size of this neighbourhood of $A$. This means that
there is a number $d$ such that for all $g$ such that $gU^{(\ast)}\leq A$, we
have $gU^{(\ast)}$ is contained in a $d$-neighbourhood of $A$. Now let $W$
denote $\{g\in G:gU^{(\ast)}\leq A\}$. The preceding discussion shows that if
$g\in W$, then $g\delta A$ lies in the $(d+1)$-neighbourhood of $A$. If we let
$c$ denote the distance of the identity of $G$ from $\delta U$, then $g$ must
lie in the $(c+d+1)$-neighbourhood of $A$, so that $W$ lies in the
$(c+d+1)$-neighbourhood of $A$, as required.
\end{proof}

Now we can bring together our ideas about enclosing.

\begin{lemma}
\label{XiisenclosedbyitsV0-vertex}Let $G$ denote a finitely generated group
with a family of finitely generated subgroups $\{H_{\lambda}\}_{\lambda
\in\Lambda}$. For each $\lambda\in\Lambda$, let $X_{\lambda}$ denote a
nontrivial $H_{\lambda}$-almost invariant subset of $G$, and suppose that the
regular neighbourhood $\Gamma(\{X_{\lambda}\}_{\lambda\in\Lambda}:G)$ can be
constructed. Let $T$ denote the bipartite tree constructed in section
\ref{regnbhds:construction} in order to define this regular neighbourhood.

If $v$ is a $V_{0}$-vertex of $T$, and the corresponding CCC of $\overline{E}$
contains an element $\overline{U}$ of $\overline{E}$, then $v$ encloses $U$.
\end{lemma}

\begin{proof}
Recall that by replacing any isolated $X_{\lambda}$'s by suitably chosen
equivalent sets, we can suppose that the $X_{\lambda}$'s are in good enough
position. Now we start by choosing an equivariant map $\varphi:G\rightarrow
V(T)$ such that $\varphi(e)=v$. Thus $\varphi(G)$ is contained in $V_{0}(T)$,
which we can identify with the pretree $P\ $consisting of all the CCC's of
$\overline{E}$. Let $s$ be an edge of $T$ which is oriented towards $v$. We
will show that $U\cap Z_{s}^{\ast}$ or $U^{\ast}\cap Z_{s}^{\ast}$ is small.
Recall that $Z_{s}=\varphi^{-1}(Y_{s})$, where $Y_{s}$ denotes the collection
of all the vertices of $T$ which lie on the terminal vertex side of $s$. Thus
$Y_{s}$ includes $v$. As $\varphi(e)=v$, it follows that $\varphi
^{-1}(w)=\{g\in G:gv=w\}$ for any vertex $w$ of $T$. Thus $Z_{s}^{\ast
}=\varphi^{-1}(Y_{s}^{\ast})=\{g\in G:gv\in Y_{s}^{\ast}\}$. Hence if $g$ and
$h$ lie in $Z_{s}^{\ast}$, then $v$ does not lie between $gv$ and $hv$. Recall
that the idea of betweenness which we defined on the set $P$ of all CCC's of
$\overline{E}$ is the same as the idea of betweenness for the $V_{0}$-vertices
of $T$. Thus the CCC $[\overline{U}]$ does not lie between $g[\overline{U}]$
and $h[\overline{U}]$. Our definition of betweenness for $P$ implies that
$\overline{U}$ does not lie between $g\overline{U}$ and $h\overline{U}$. As
$Y_{s}^{\ast}$ does not include $v$, we know that $g$ cannot fix $v$, so that
$gU$ and $U$ are in distinct CCC's. Thus they are comparable using our partial
order $\leq$ on the elements of $E$. Similarly $hU\;$and $U$ are comparable.
By replacing $U$ by $U^{\ast}$ if necessary, we can arrange that $gU^{(\ast
)}\leq U$, and the fact that $\overline{U}$ does not lie between
$g\overline{U}$ and $h\overline{U}$ implies that we must also have
$hU^{(\ast)}\leq U$. We conclude that, by replacing $U$ by $U^{\ast}$ if
necessary, we can arrange that $gU^{(\ast)}\leq U$, for every $g\in
Z_{s}^{\ast}$. Now Lemma \ref{containedinboundednbhdofA} implies that
$Z_{s}^{\ast}$ lies in a bounded neighbourhood of $U$, so that $U^{\ast}\cap
Z_{s}^{\ast}$ is small. It follows that for any edge $s$ of $T$ which is
oriented towards $v$, we have $U\cap Z_{s}^{\ast}$ or $U^{\ast}\cap
Z_{s}^{\ast}$ is small, so that $U$ is enclosed by $v$ as required. This
completes the proof of the lemma.
\end{proof}

Next we will apply Lemma \ref{enclosedsetscanbechosennested} to the regular
neighbourhood $\Gamma(\{X_{\lambda}\}_{\lambda\in\Lambda}:G)$. Let $T$ denote
the universal covering $G$-tree of $\Gamma(\{X_{\lambda}\}_{\lambda\in\Lambda
}:G)$. As each $X_{\lambda}$ is enclosed by a vertex of $T$, this lemma tells
us how to replace each $X_{\lambda}$ by an equivalent set $B(X_{\lambda})$,
which is nested with respect to $Z_{s}$, for every edge $s$ of $T$. Note that
before we can define $Z_{s}$, we should first choose a $G$-equivariant map
$\varphi$ from $G$ to $V(T)$ such that the vertex $\varphi(e)$ encloses some
nontrivial $K$-almost invariant subset $U$ of $G$. Lemma
\ref{XiisenclosedbyitsV0-vertex} implies that we can choose $\varphi(e)$ to be
any $V_{0}$-vertex of $\Gamma(\{X_{\lambda}\}_{\lambda\in\Lambda}:G)$. Recall
from Lemma \ref{equivalentsetshavesameregnbhd} that replacing each
$X_{\lambda}$ by an equivalent set does not alter the regular neighbourhood.
Thus we obtain the following interesting result which can be thought of as
asserting that we can replace the $X_{\lambda}$'s by equivalent sets which are
in ``very good position''.

\begin{corollary}
Let $G$ be a finitely generated group with a family of finitely generated
subgroups $\{H_{\lambda}\}_{\lambda\in\Lambda}$. For each $\lambda\in\Lambda$,
let $X_{\lambda}$ denote a nontrivial $H_{\lambda}$-almost invariant subset of
$G$, and suppose that the regular neighbourhood $\Gamma(\{X_{\lambda
}\}_{\lambda\in\Lambda}:G)$ can be constructed. Then we can replace each
$X_{\lambda}$ by an equivalent almost invariant set $B(X_{\lambda})$ so that
if $U$ and $V$ are any two elements of the set $E=\{gX_{\lambda},gX_{\lambda
}^{\ast}:g\in G\}$, then either $B(U)\;$and $B(V)$ are nested or $U\;$and $V$
lie in the same CCC of $\overline{E}$.
\end{corollary}

\begin{proof}
If $U$ and $V$ lie in distinct CCC's of $\overline{E}$, they are enclosed by
distinct vertices of $T$. Let $s$ denote an edge of $T$ on the path joining
these two vertices. As $B(U)\;$and $B(V)$ are nested with respect to every
$Z_{s}$, it follows that $\overline{Z_{s}}$ lies between $\overline{B(U)}$ and
$\overline{B(V)}$, so that the two sets $B(U)\;$and $B(V)$ are nested, as claimed.
\end{proof}

The examples at the end of section \ref{regnbhds:construction} show that our
regular neighbourhood construction can yield a graph of groups $\Gamma$
consisting of a single point, but they also suggest that this is quite
unusual. A more delicate question is whether the graph of groups decomposition
of $G$ which we obtain can decompose $G$ trivially when $\Gamma$ is not a
point. The answer is that this cannot happen. In fact, we can show the far
stronger result that $\Gamma(\{X_{\lambda}\}_{\lambda\in\Lambda}:G)$ is always minimal.

\begin{proposition}
\label{Tisminimal}Let $G$ be a finitely generated group with a family of
finitely generated subgroups $\{H_{\lambda}\}_{\lambda\in\Lambda}$. For each
$\lambda\in\Lambda$, let $X_{\lambda}$ denote a nontrivial $H_{\lambda}%
$-almost invariant subset of $G$, and suppose that the regular neighbourhood
$\Gamma(\{X_{\lambda}\}_{\lambda\in\Lambda}:G)$ can be constructed. Let $E$,
$P$ and $T$ be as in our construction of the regular neighbourhood
$\Gamma(\{X_{\lambda}\}_{\lambda\in\Lambda}:G)$. Then $T$ is a minimal
$G$-tree, so that $\Gamma(\{X_{\lambda}\}_{\lambda\in\Lambda}:G)$ is also minimal.
\end{proposition}

\begin{proof}
Recall from Lemma \ref{enclosedsetscanbechosennested}, that by replacing the
$X_{\lambda}$'s by equivalent sets, we can arrange that each element of $E$ is
nested with respect to every $Z_{s}$. We will assume that this has been done.

Let $T_{0}$ denote the minimal subtree of $T$. If $G$ fixes more than one
vertex of $T$ so that $T_{0}$ is not unique, we let $T_{0}$ denote one of the
vertices fixed by $G$. We will show that $T_{0}=T$. As in the proof of part 7)
of Lemma \ref{somefactsaboutenclosing}, any vertex of $T$ which encloses a
nontrivial almost invariant subset of $G$ must lie in $T_{0}$. Hence $T_{0}$
contains every $V_{0}$-vertex of $T$. Now recall that each $V_{1}$-vertex of
$T$ corresponds to a star in the pretree $P$, and any such star must contain
at least two points of $P$. Hence each $V_{1}$-vertex is joined by edges of
$T$ to at least two $V_{0}$-vertices. (It is possible that $T$ has a single
$V_{0}$-vertex, but in this case $T$ consists only of this vertex.) As $T$ is
a tree, it follows that $T_{0}=T$ as required.
\end{proof}

We next examine the special features of isolated elements of $E$ in the case
when the elements of $E$ are in good enough position. Recall that if $A$ is an
isolated element of $E$, then the translates of $A$ are nested and since $gA$
is also isolated in $E$, each $gA$ determines exactly one cross-connected
component of $\overline{E}$, which we call isolated. The following result
characterises isolated CCC's of $\overline{E}$ in terms of the $G$-tree $T$.

\begin{proposition}
\label{isolatedvertices}A $V_{0}$-vertex $v$ of $T$ has valence two if and
only if $v$ corresponds to an isolated CCC of $\overline{E}$.
\end{proposition}

\begin{proof}
Suppose that $v$ is a $V_{0}$-vertex of valence two and let $s$ and $t$ denote
the edges which are incident to $v$. Let $U$ and $V$ denote elements of $E$
enclosed by $v$. Part 1) of Corollary \ref{nontrivialpartition} implies that
$s$ and $t$ must lie one on the $U$-side of $v$ and the other on the $U^{\ast
}$-side of $v$, and the same holds with $U$ replaced by $V$. Thus $U$ and $V$
determine the same partition of the vertices of $T-\{v\}$. Now part 3) of
Corollary \ref{nontrivialpartition} shows that $U$ and $V$ must be equivalent,
up to complementation. As required, it now follows that the CCC which
corresponds to $v$ contains exactly one element of $\overline{E}$.

To prove the converse, suppose that $A$ is an isolated element of
$\overline{E}$. Consider the $V_{0}$-vertex $v$ which is the CCC of $A$ and,
in the pretree $P$, let $S$ be a star which contains $v$. Thus $S$ can also be
identified with a $V_{1}$-vertex $v_{1}$ in $T$ which is adjacent to $v$. Let
$v_{2}\neq v$ be a $V_{0}$-vertex corresponding to another CCC in $S$. This
means that both $v_{2}$ and $v$ are adjacent to $v_{1}$ in the tree $T$. If
$U$ is an element of $E$ which is enclosed by $v_{2}$, then we must have
$U^{(\ast)}>A$ or $U^{(\ast)}<A$. If $U^{(\ast)}<A$, then we have the same
inequality for any other CCC in $S$ other than $v$. For otherwise, $v$ would
lie between two of these vertices. Conversely, if $[\overline{U}]=u$ is any
vertex of a star which contains $v$ and if $U^{(\ast)}<A$, then $u$ must lie
in $S$ since $v$ cannot lie between $u$ and $v_{2}$. Thus, there are only two
stars which contain $v$. Equivalently, there are only two edges of $T$
incident to $v$.
\end{proof}

Before proceeding, we will need the following simple fact about splittings of groups.

\begin{lemma}
\label{ascendingHNN}Let $G$ be a group with a splitting $\sigma$ over a
subgroup $H$. Let $X$ be one of the $H$-almost invariant subsets of $G$
associated to $\sigma$. Then either there is an element $g$ of $G$ such that
$gX\subset X^{\ast}$ or $\sigma$ is a HNN extension $G=H\ast_{H}$ in which at
least one of the inclusions of the edge group in the vertex group is an isomorphism.
\end{lemma}

\begin{remark}
Such an extension is often called an ascending HNN extension.
\end{remark}

\begin{proof}
The given splitting $\sigma$ determines a $G$-tree $T$ such that the quotient
$G\backslash T$ has a single edge. We pick an orientation of this edge,
thereby fixing a $G$-invariant orientation on every edge of $T$. In terms of
our previous notation, there is an edge $s$ of $T$ such that $X$ or $X^{\ast}$
is equivalent to $Z_{s}$. We will assume that $X$ is equivalent to $Z_{s}$.
Let $v$ denote the vertex of $T$ at the initial end of $s$, and let $gs$
denote another edge incident to $v$. Then $gX\subset X^{\ast}$ or $gX^{\ast
}\subset X^{\ast}$, depending on whether $gs$ points away from $v$ or towards
$v$. Suppose that there is no element $g$ of $G$ such that $gX\subset X^{\ast
}$. Then $gX^{\ast}\subset X^{\ast}$, and $gs$ must point towards $v$. It
follows that $s$ is the only edge of $T$ which is incident to $v$ and points
away from $v$. Hence the stabilisers of $s$ and of $v$ are equal, which
implies the result of the lemma.
\end{proof}

Our next result is crucial for understanding algebraic regular neighbourhoods.
It is the algebraic analogue of the topological fact that if $N$ is a regular
neighbourhood of a finite collection $C_{\lambda}$ of closed curves on a
surface $M$, and if $C$ is a closed curve which is disjoint from each
$C_{\lambda}$, then we can homotop $C$ into $M-N$. In our result the curve $C$
is replaced by a $H$-almost invariant subset $X$ of $G$, and the conclusion is
that if $X$ crosses no element of $E$, then $X$ is enclosed by a $V_{1}%
$-vertex of $\Gamma(\{X_{\lambda}\}_{\lambda\in\Lambda}:G)$. As usual, we
argue with the $G$-tree $T$ rather than with the graph of groups $\Gamma$ itself.

\begin{proposition}
\label{XdoesnotcrossanyXiimpliesXisenclosedbyaV1vertex}Let $G$ be a finitely
generated group with a family of finitely generated subgroups $\{H_{\lambda
}\}_{\lambda\in\Lambda}$. For each $\lambda\in\Lambda$, let $X_{\lambda}$
denote a nontrivial $H_{\lambda}$-almost invariant subset of $G$, and suppose
that the regular neighbourhood $\Gamma(\{X_{\lambda}\}_{\lambda\in\Lambda}:G)$
can be constructed. Further suppose that the $X_{\lambda}$'s are in good
enough position, and that there is more than one CCC, so that the pretree $P$
is not a single point. Let $T$ denote the $G$-tree with $V_{0}(T)=P$. Let $X$
be a nontrivial $H$-almost invariant subset of $G$ which does not cross any
element of $E$. Then the following statements hold.

\begin{enumerate}
\item If $H$ is finitely generated, then $X$ is enclosed by a $V_{1}$-vertex
of $T$.

\item If $X$ is a standard $H$-almost invariant set associated to a splitting
of $G$ over $H$, which need not be finitely generated, then $X$ is enclosed by
a $V_{1}$-vertex of $T$.
\end{enumerate}
\end{proposition}

\begin{proof}
We will give the argument assuming that the $X_{\lambda}$'s are in good
position, so that $E$ satisfies Condition (*) and the relation $\leq$ can be
defined on $E$. However, the reader can easily verify that, as in section
\ref{regnbhds:construction}, the argument goes through if the $X_{\lambda}$'s
are in good enough position, because we will never need to compare elements of
$E$ which belong to the same CCC.

If $X$ is equivalent to $Z_{s}$ for any edge $s$ of $T$, the result follows
from part 2) of Lemma \ref{somefactsaboutenclosing}. If $X$ is equivalent to
an element $U$ of $E$, our hypothesis on $X$ implies that $U$ must be isolated
in $E$. This implies that $X$ is equivalent to $Z_{s}$ for each of the two
edges of $T$ incident to the $V_{0}$-vertex which corresponds to the CCC
$[\overline{U}]$. In particular, the result follows in this case also. So we
will assume that $X$ is not equivalent to any $Z_{s}$ nor to any element of
$E$. This implies that the relation $\leq$ on $E$ can be extended to the set
obtained by adding in all translates of $X$ and $X^{\ast}$. Note that $\leq$
is a partial order on this larger set even if the stabiliser of $X$ is not
finitely generated, because the proof of Lemma 1.14 of \cite{SS} still applies.

Our first step is to show that $X$ is sandwiched between two elements of $E$,
i.e. that there are elements $U_{1}$ and $U_{2}$ of $E$ such that
$U_{1}<X<U_{2}$. Let $Y$ denote an element of $E$. As $X$ crosses no element
of $E$, we know that, for each element $g$ of $G$, one of the four
inequalities $gY\leq X$, $gY^{\ast}\leq X$, $gY\leq X^{\ast}$, $gY^{\ast}\leq
X^{\ast}$ must hold. If the stabiliser $H$ of $X$ is finitely generated, then
Lemma \ref{containedinboundednbhdofA} tells us that $\{g\in G:gY^{(\ast)}\leq
X\}$ is contained in a bounded neighbourhood of $X$, and that $\{g\in
G:gY^{(\ast)}\leq X^{\ast}\}$ is contained in a bounded neighbourhood of
$X^{\ast}$. As $G$ is the union of these two sets it follows that neither is
empty so that there are elements $U_{1}$ and $U_{2}$ of $E$ such that
$U_{1}<X<U_{2}$, as required.

If $H$ is not finitely generated, we use the hypothesis that $X$ is associated
to a splitting $\sigma$ of $G$. Suppose that $X$ is not sandwiched between two
elements of $E$. Then, by replacing $X$ by $X^{\ast}$, if necessary, we have
that for every element $U$ of $E$, either $U<X$ or $U^{\ast}<X$. As $E$ is
$G$-invariant, it follows that for every element $U$ of $E$, and for every
element $g$ of $G$, either $U<gX$ or $U^{\ast}<gX$. Now suppose that there is
an element $g$ of $G$ such that $gX\subset X^{\ast}$. This implies that for
every element $U$ of $E$, either $U<X^{\ast}$ or $U^{\ast}<X^{\ast}$, which
contradicts the fact that either $U<X$ or $U^{\ast}<X$. It remains to handle
the situation where there is no element $g$ of $G$ such that $gX\subset
X^{\ast}$. Then Lemma \ref{ascendingHNN} shows that $\sigma$ is an ascending
HNN extension $G=H\ast_{H}$. Let $f:G\rightarrow\mathbb{Z}$ denote the natural
homomorphism associated to this HNN extension, and let $U$ be an element of
$E$ such that $U<X$. As $U$ does not cross $X$, it follows, in particular,
that the image in $\mathbb{Z}$ of the coboundary $\delta U$ of $U$ is bounded
above or below. Hence the stabiliser of $U$ must be contained in $\ker(f)$. It
follows that the image of $\delta U$ in $\mathbb{Z}$ must be finite. As the
image in $\mathbb{Z}$ of $\delta gX$ is finite for all $g$ in $G$, there must
be an element $g$ of $G$ such that $gX\subset U$. But we know that $U<gX$ or
$U^{\ast}<gX$. The first implies that $U$ is equivalent to $gX$, which
contradicts our assumption that the stabiliser $H$ of $X$ is not finitely
generated, and the second implies that $U^{\ast}<U$ which is also a
contradiction. This contradiction completes the proof that, in all cases, $X$
must be sandwiched between two elements $U_{1}$ and $U_{2}$ of $E$.

By considering the path in $T$ which joins the $V_{0}$-vertices $[\overline
{U_{1}}]$ and $[\overline{U_{2}}]$, it is easy to see that there is a $V_{1}%
$-vertex $v$ with two $V_{0}$-vertices $v_{1}$ and $v_{2}$ adjacent to $v$,
and an element $V_{i}$ of $E$ enclosed by $v_{i}$, such that $V_{1}<X<V_{2}$.
We will show that $X$ is enclosed by $v$.

If $U$ is any element of $E$, then either $U^{(\ast)}<X$ or $U^{(\ast
)}<X^{\ast}$ but not both. Further if $V$ lies in the same CCC as $U$, then
the same inequality must hold as for $U$. Thus we obtain a partition of the
$V_{0}$-vertices of $T-\{v\}$ into two subsets $\Phi$ and $\Phi^{\ast}$, where
the vertices of $\Phi$ enclose those elements $U$ of $E$ such that $U^{(\ast
)}<X$ and the vertices of $\Phi^{\ast}$ enclose those elements $U$ of $E$ such
that $U^{(\ast)}<X^{\ast}$. If $w$ lies in $\Phi$, then every $V_{0}$-vertex
on the path from $w$ to $v$ also lies in $\Phi$. This enables us to define a
partition of all the vertices of $T-\{v\}$ into two subsets $\Psi$ and
$\Psi^{\ast}$. We will say that a vertex $w$ of $T-\{v\}$ lies in $\Psi$ if
the last $V_{0}$-vertex on the path from $w$ to $v$ lies in $\Phi$. Thus
$\Phi$ is contained in $\Psi$. Note that if we already knew that $X$ was
enclosed by $v$, then $\Psi$ and $\Psi^{\ast}$ would constitute the $X$-side
of $v$ and the $X^{\ast}$-side of $v$. If $U<X$ and $h\in H$, then $hU<X$.
Thus $\Psi$ and $\Psi^{\ast}$ are $H$-invariant. In particular, $H$ must be a
subgroup of $Stab(v)$.

Now we choose $\varphi(e)=v$ and define $B(X)=\varphi^{-1}(\Psi)\cup
(X\cap\varphi^{-1}(v))$, and $C(X)=\varphi^{-1}(\Psi^{\ast})\cup(X^{\ast}%
\cap\varphi^{-1}(v))$, as in the proof of Lemma
\ref{enclosedsetscanbechosennested}. Clearly $B(X)$ and $C(X)$ partition $G$,
so that $C(X)=B(X)^{\ast}$. If $s$ is an edge of $T$ which is directed towards
$v$, then $Z_{s}^{\ast}\subset B(X)$ or $Z_{s}^{\ast}\subset C(X)$. It is also
clear that $HB=B$, because $\Psi$, $X$ and $v$ are all $H$-invariant.

We will show that $B(X)$ is $H$-almost invariant and is equivalent to $X$. It
will then follow that $Z_{s}^{\ast}\leq X$ or $Z_{s}^{\ast}\leq X^{\ast}$, for
every edge $s$ of $T$ which is directed towards $v$, so that $X$ is enclosed
by $v$.

If $H$ is finitely generated, we will proceed very much as in the proof of
Lemma \ref{enclosedsetscanbechosennested}. As $\varphi(e)=v$, we see that
$\varphi(g)=gv$, so that $\varphi^{-1}(\Psi)=\{g\in G:gv\in\Psi\}$. As
$V_{1}<X<V_{2}$, we have $gV_{1}<gX<gV_{2}$, for all $g\in G$. If $gv$ lies in
$\Psi$, then $gv_{1}$ and $gv_{2}$ must also lie in $\Psi$. Thus we have the
inequalities $gV_{i}^{(\ast)}<X$, for $i=1$, $2$. It follows that one of the
inequalities $gX^{(\ast)}<X$ also holds. Thus $\varphi^{-1}(\Psi)=\{g\in
G:gv\in\Psi\}$ is contained in $\{g\in G:gX^{(\ast)}<X\}$, which is contained
in a bounded neighbourhood of $X$, by Lemma \ref{containedinboundednbhdofA}.
It follows that $B(X)$ itself is contained in a bounded neighbourhood of $X$.
Similarly $C(X)=B(X)^{\ast}$ is contained in a bounded neighbourhood of
$X^{\ast}$. Exactly as in the proof of Lemma
\ref{enclosedsetscanbechosennested} it follows that $B(X)$ is $H$-almost
invariant and is equivalent to $X$.

If $H$ is not finitely generated, we will choose $X$ carefully and then show
that $X=B(X)$. Recall our assumption that $X$ is associated to a splitting of
$G$ over $H$. As discussed after Lemma
\ref{pre-imageofhalftreeisalmostinvariant}, we can choose $X$ to satisfy
$X=\{g\in G:gX^{(\ast)}\subset X\}$. As the translates of $X$ are nested, we
have $\{g\in G:gX^{(\ast)}\subset X\}=\{g\in G:gX^{(\ast)}\leq X\}$ which
allows us to express $X$ as the disjoint union $\{g\in G:g\notin
Stab(v),gX^{(\ast)}<X\}\cup(X\cap Stab(v))$.

As $\varphi(e)=v$, we have $\varphi^{-1}(v)=Stab(v)$, so that the second terms
in our expressions for $X$ and $B(X)$ are equal. Now we consider the first
terms. As above, we have $\varphi^{-1}(\Psi)=\{g\in G:gv\in\Psi\}\subset\{g\in
G:g\notin Stab(v),gX^{(\ast)}<X\}$.

Now suppose that $g\notin Stab(v)$ and $gX<X$. As $g\notin Stab(v)$, the
vertices $gv$, $gv_{1}$ and $gv_{2}$ of $T$ must all lie in $\Psi$ or all in
$\Psi^{\ast}$. As $gV_{1}<gX<X$, we see that $gv_{1}$ must lie in $\Psi$. Thus
$gv\in\Psi$, so that $g\in\varphi^{-1}(\Psi)$. A similar argument applies if
$g\notin Stab(v)$ and $gX^{\ast}<X$. Thus $\{g\in G:g\notin Stab(v),gX^{(\ast
)}<X\}\subset\varphi^{-1}(\Psi)$. We conclude that $\varphi^{-1}(\Psi)=\{g\in
G:g\notin Stab(v),gX^{(\ast)}<X\}$, so that $X=B(X)$ as claimed.
\end{proof}

Recall from the definition of betweenness on the pretree of CCC's of
$\overline{E}$ that if three $V_{0}$-vertices $v_{1}$, $v_{2}$ and $v_{3}$ lie
on a path in $T$ with $v_{2}$ between $v_{1}$ and $v_{3}$, then there is an
element $X$ of $E$ which is enclosed by $v_{2}$ and such that for any elements
$Y$ and $Z$ of $E$ with $Y$ enclosed by $v_{1}$ and $Z$ enclosed by $v_{3}$,
we have $\overline{Y}\,\overline{X}\,\overline{Z}$. Hence if a $V_{0}$-vertex
$v_{2}$ lies between edges $s$ and $t$ which point towards $v_{2}$, there is
an element $X$ of $E$ enclosed by $v_{2}$ such that $Z_{s}^{\ast}\leq X\leq
Z_{t}$. The following result is an immediate consequence and is the analogous
result for $V_{1}$-vertices.

\begin{proposition}
Let $Y$ be a nontrivial $H$-almost invariant subset of $G$, and let $Z$ be a
nontrivial $K$-almost invariant subset of $G$. If $H$ is not finitely
generated, suppose in addition that $Y$ is a standard $H$-almost invariant set
associated to a splitting of $G$ over $H$, and similarly for $Z$. If $Y$ and
$Z$ are enclosed by distinct $V_{1}$-vertices $v_{1}$ and $v_{3}$, then there
is a $V_{0}$-vertex $v_{2}$ and an element $X$ of $E$ which is enclosed by
$v_{2}$ such that $Y^{(\ast)}\leq X\leq Z^{(\ast)}$.
\end{proposition}

\begin{proof}
Let $v_{2}$ be any $V_{0}$-vertex on the path joining $v_{1}$ and $v_{3}$. Let
$s$ and $t$ be the edges of this path which are incident to $v_{2}$ and point
towards $v_{2}$, labelled so that $s$ is nearer to $v_{1}$ than is $t$. As we
pointed out above, there is an element $X$ of $E$ enclosed by $v_{2}$ such
that $Z_{s}^{\ast}\leq X\leq Z_{t}$. As $s$ points away from $v_{1}$, we have
$Y^{(\ast)}\leq Z_{s}^{\ast}$, and as $t$ points away from $v_{3}$, we have
$Z^{(\ast)}\geq Z_{t}$. It follows that $Y^{(\ast)}\leq Z_{s}^{\ast}\leq X\leq
Z_{t}\leq Z^{(\ast)}$, so that $Y^{(\ast)}\leq X\leq Z^{(\ast)}$, as required.
\end{proof}

The following result will be useful when we consider taking regular
neighbourhoods of increasing finite families of almost invariant subsets of a
fixed group.

\begin{lemma}
\label{edgesplittingswhenyouenlargeE}Let $G$ be a finitely generated group
with a finite family of finitely generated subgroups $H_{1},\ldots,H_{n}$. For
$1\leq i\leq n$, let $X_{i}$ denote a nontrivial $H_{i}$-almost invariant
subset of $G$. Let $E_{n}$ denote the set of all translates of the $X_{i}$'s
and their complements, let $P_{n}$ denote the pretree of all CCC's of
$\overline{E_{n}}$, let $T_{n}$ denote the associated $G$-tree and let
$\Gamma_{n}$ denote the corresponding graph of groups structure for $G$, so
that $\Gamma_{n}=G\backslash T_{n}$.

Let $m<n$, and let $f$ denote the natural map from $P_{m}$ to $P_{n}$. If $A$
is a $H$-almost invariant subset of $G$ enclosed by a $V_{0}$-vertex $v$ of
$T_{m}$, then $A$ is enclosed by the $V_{0}$-vertex $f(v)$ of $T_{n}$.
\end{lemma}

\begin{proof}
Lemma \ref{enclosedsetscanbechosennested} implies that each edge splitting of
$\Gamma_{n}$ has zero intersection number with each $X_{i}$, $1\leq i\leq n$.
Proposition \ref{XdoesnotcrossanyXiimpliesXisenclosedbyaV1vertex} now implies
that each edge splitting of $\Gamma_{n}$ is enclosed by some $V_{1}$-vertex of
$\Gamma_{m}$. Hence there is a common refinement $\Gamma_{m,n}$ of $\Gamma
_{n}$ and $\Gamma_{m}$ obtained by splitting the $V_{1}$-vertices of
$\Gamma_{m}$ using the edge splittings of $\Gamma_{n}$. (Thus the number of
edges of $\Gamma_{m,n}$ is the sum of the number of edges of $\Gamma_{n}$ and
of $\Gamma_{m}$.) Associated to the construction of $\Gamma_{m,n}$ there is a
natural quotient map $p_{m}:T_{m,n}\rightarrow T_{m}$. There is also a natural
quotient map $p_{n}:T_{m,n}\rightarrow T_{n}$ obtained by collapsing those
edges of $T_{m,n}$ which were originally in $T_{m}$. We choose a
$G$-equivariant map $\varphi:G\rightarrow T_{m,n}$. As both $p_{m}$ and
$p_{n}$ are $G$-equivariant, this determines $G$-equivariant maps $\varphi
_{m}$ and $\varphi_{n}$ to $T_{m}$ and $T_{n}$ respectively. Now we consider
the $V_{0}$-vertex $v$ of $T_{m}$. The construction of $T_{m,n}$ means that
the pre-image of $v$ is a single vertex $w$. Further the pre-image of a small
neighbourhood of $v$ maps homeomorphically into $T_{m}$.

Let $A$ be a $H$-almost invariant subset of $G$ enclosed by the vertex $v$ of
$T_{m}$. Thus for any edge $s$ of $T_{m,n}$ incident to $w$ and oriented
towards $w$, if $t$ denotes $p_{m}(s)$, then we have $A\cap Z_{t}^{\ast}$ or
$A^{\ast}\cap Z_{t}^{\ast}$ is small. As the pre-image in $T_{m,n}$ of a small
neighbourhood of $v$ maps homeomorphically into $T_{m}$, our definition of
$\varphi_{m}:G\rightarrow T_{m}$ implies that $Z_{s}=Z_{t}$, so that we have
$A\cap Z_{s}^{\ast}$ or $A^{\ast}\cap Z_{s}^{\ast}$ is small, i.e. $A$ is
enclosed by the vertex $w$ of $T_{m,n}$. Hence part 4) of Lemma
\ref{somefactsaboutenclosing} shows that $A\cap Z_{s}^{\ast}$ or $A^{\ast}\cap
Z_{s}^{\ast}$ is small for every edge $s$ of $T_{m,n}$ which is oriented
towards $w$, whether or not $s$ is incident to $w$. Let $r$ denote an edge of
$T_{n}$ which is incident to $p_{n}(w)$ and oriented towards $p_{n}(w)$. The
construction of $p_{n}$ implies that there is a unique edge $s$ of $T_{m,n}$,
which need not be incident to $w$, such that $p_{n}(s)=r$. The edge $s$ is
automatically oriented towards $w$. Now the definition of $\varphi_{n}$
implies that $Z_{r}=Z_{s}$, so that we have $A\cap Z_{r}^{\ast}$ or $A^{\ast
}\cap Z_{r}^{\ast}$ is small, as required. Hence $A$ is enclosed by the vertex
$p_{n}(w)$ of $T_{n}$. If $p_{n}(w)=f(v)$, we have completed the proof of the lemma.

Now suppose that $p_{n}(w)\neq f(v)$, and let $U$ denote an element of $E_{m}$
which belongs to the CCC of $\overline{E_{m}}$ corresponding to the $V_{0}%
$-vertex $v$ of $T_{m}$. Lemma \ref{twoideasofenclosingaresameifHisfg} tells
us that $U\;$is enclosed by $v$, so that the above argument shows that $U$
must also be enclosed by the vertex $p_{n}(w)$ of $T_{n}$. But we know that
$U$ is enclosed by the $V_{0}$-vertex $f(v)$ of $T_{n}$, as it belongs to the
CCC of $\overline{E_{n}}$ associated to $f(v)$. As these vertices of $T_{n}$
are distinct, part 8) of Lemma \ref{somefactsaboutenclosing} tells us that $U$
is equivalent to $Z_{s}$ or to $Z_{s}^{\ast}$ for each edge $s$ on the path
joining them. As each element of $E_{n}$ is enclosed by some vertex of $T_{n}%
$, it follows that no element of $E_{n}$ crosses any $Z_{s}$, so that $U$ must
be an isolated element of $E_{n}$. Hence $U$ is an isolated element of $E_{m}%
$, and so Proposition \ref{isolatedvertices} tells us that $v$ has valence $2$
in $T_{m}$. It follows from Corollary \ref{nontrivialpartition} that $A$ is
equivalent to $U$ or $U^{\ast}$, so that again $A$ is enclosed by the vertex
$f(v)$ of $T_{n}$, as required. The lemma follows.
\end{proof}

Next we apply the preceding results on regular neighbourhoods to obtain new
results about more general graphs of groups.

First we can use the argument of Proposition
\ref{XdoesnotcrossanyXiimpliesXisenclosedbyaV1vertex} to show that the two
ideas of enclosing a splitting of $G$ over a subgroup $H$ are equivalent even
when $H$ is not finitely generated.

\begin{lemma}
\label{twoideasofenclosingaresame}Suppose that $A$ is associated to a
splitting $\sigma$ of $G$ over $H$, and $T$ is a $G$-tree. Let $\Gamma$ denote
the graph of groups structure for $G$ given by the quotient $G\backslash T$,
and let $w$ denote the image of $v$ in $\Gamma$. Then $A$ is enclosed by $w$
if and only if $\sigma$ is enclosed by $w$.
\end{lemma}

\begin{proof}
In Lemma \ref{twoideasofenclosingaresameifHisfg}, we showed that if $\sigma$
is enclosed by $w$ then $A$ is enclosed by $w$. We also proved the converse in
the case when $H$ is finitely generated. It remains to prove the converse in
the case when $H\;$is not finitely generated.

Suppose that $A$ is enclosed by the vertex $w$ of $\Gamma$. This means that
$A$ is enclosed by some translate in $T$ of $v$, and we will assume that $A$
is enclosed by $v$ itself. As always there is a special case if $v$ is fixed
by $G$, but the result is trivial in this case as $v$ encloses all splittings
and all almost invariant subsets of $G$. So we will now assume that $v$ is not
fixed by $G$. Part 7) of Lemma \ref{somefactsaboutenclosing} shows that $v$
lies in the minimal subtree of $T$. By replacing $T$ by its minimal subtree,
we can assume that $T$ itself is minimal.

Now the argument at the end of the proof of Proposition
\ref{XdoesnotcrossanyXiimpliesXisenclosedbyaV1vertex} applies to show that
$A=B(A)$. As $B(A)$ is clearly nested with respect to every $Z_{s}$, the set
$E$ of all $Z_{s}$, for all oriented edges of $T$, together with all the
translates of $A$ and $A^{\ast}$, forms a nested set $F$ of subsets of $G$.
Thus $F$ is a $G$-set partially ordered by inclusion. We want to apply
Dunwoody's construction in \cite{D-D} to $F$ to obtain a $G$-tree $T^{\prime}$
whose oriented edges naturally correspond to the elements of $F$. This will
yield a graph of groups structure $\Gamma^{\prime}=G\backslash T^{\prime}$ for
$G$ which is a refinement of $\Gamma$, with one extra edge whose edge
splitting is $\sigma$. This will show that $\sigma$ is enclosed by $w$ as
required. In order for his construction to be applicable, we need to know that
$F$ is discrete. Recall that the elements of $E$ correspond to the edges of
$T$ and have the corresponding partial order. It follows at once that $E$ is
discrete. Also the fact that $A$ is associated to a splitting implies that the
set of all translates of $A$ is discrete. (In fact, the translates of $A$
correspond to the edges of the $G$-tree determined by the splitting.) As $T$
is minimal, each vertex has valence at least two, so that $A$, and hence any
translate of $A$, is sandwiched between two elements of $E$. Also any element
of $E$ lies between two translates of $A$ or $A^{\ast}$. Combining these facts
shows that $F$ is discrete, as required.
\end{proof}

The preceding results allow us to give a surprising generalisation of Theorem
\ref{Theorem2.5ofSS}, which we proved in \cite{SS}. In that theorem, we showed
that if $G$ is a finitely generated group with $n$ splittings over finitely
generated subgroups, then the splittings are compatible if and only if each
pair of splittings has intersection number zero. Now we will show that the
hypothesis that the splittings be over finitely generated subgroups of $G$ can
be removed. The precise result we obtain is the following.

\begin{theorem}
\label{splittingswithintersectionnumberzeroarecompatible}Let $G$ be a finitely
generated group with $n$ splittings over possibly infinitely generated
subgroups. Then the splittings are compatible if and only if each pair of
splittings has intersection number zero. Further, in this situation, the graph
of groups structure on $G$ obtained from these splittings is unique up to isomorphism.
\end{theorem}

\begin{proof}
For $1\leq i\leq n$, let $\sigma_{i}$ be a splitting of $G$ over a subgroup
$H_{i}$, and let $X_{i}$ be an associated $H_{i}$-almost invariant subset of
$G$.

We start by discussing the existence proof. In order to use the previous
arguments directly, we will proceed by induction on $n$. Thus we need to
consider the situation where we have $k$ compatible splittings $\sigma
_{1},\ldots,\sigma_{k}$ and then another splitting $\sigma_{k+1}$ which has
intersection number zero with each of $\sigma_{1},\ldots,\sigma_{k}$. Thus $G$
has a graph of groups decomposition with $k$ edges and the edge splittings are
conjugate to the $\sigma_{i}$'s, $1\leq i\leq k$. By subdividing each edge
into two, we obtain the regular neighbourhood $\Gamma(X_{1},\ldots,X_{k}:G)$.
Now part 2) of Proposition
\ref{XdoesnotcrossanyXiimpliesXisenclosedbyaV1vertex} shows that $X_{k}$ is
enclosed by a $V_{1}$-vertex of $\Gamma(X_{1},\ldots,X_{k}:G)$, and then Lemma
\ref{twoideasofenclosingaresame} implies that $\sigma_{1},\ldots,\sigma_{k+1}$
are compatible as required.

In order to prove the uniqueness, we want to use the proof of the second part
of Theorem 1.12 of \cite{SS}. This argument never directly uses the hypothesis
that the $H_{i}$'s are finitely generated. It does quote one result, Lemma 2.3
of \cite{SS}, and this in turn uses Lemma 2.2 of \cite{SS}. Both lemmas are
about splittings of a finitely generated group $G$, and they each contain the
hypothesis that the splittings be over finitely generated subgroups. However,
this hypothesis is not needed and the proofs work perfectly well even when the
splittings are over infinitely generated subgroups of $G$. Thus the proof of
the second part of Theorem 1.12 of \cite{SS} yields the required result.
\end{proof}

The following result is an immediate consequence of the uniqueness part of
Theorem \ref{splittingswithintersectionnumberzeroarecompatible}.

\begin{theorem}
\label{graphswithsameedgesplittingsareiso}Let $\Gamma_{1}$ and $\Gamma_{2}$ be
minimal graphs of groups structures for a finitely generated group $G$. If
each $\Gamma_{i}$ has no redundant vertices, and if they have the same
conjugacy classes of edge splittings, then $\Gamma_{1}$ and $\Gamma_{2}$ are isomorphic.
\end{theorem}

\begin{remark}
As $G$ is finitely generated, and $\Gamma_{1}$ and $\Gamma_{2}$ are minimal,
it follows that $\Gamma_{1}$ and $\Gamma_{2}$ are finite.

Note that this result needs no assumptions on the edge groups involved.
\end{remark}

We end this section by discussing how to generalise our theory of regular
neighbourhoods of almost invariant subsets to the case of almost invariant
subsets over infinitely generated subgroups. We will use the preceding results
to show that this can be done provided that such sets are associated to
splittings. While this may seem a little exotic, it is a natural extension
once one realises that the edge groups in a regular neighbourhood may be
infinitely generated even if the given almost invariant sets are all over
finitely presented groups. (See Example
\ref{exampleofregnbhdwithnonfgedgegroups}.) Also no more work is needed. We
simply need to use the results we already have. Here is an outline of the theory.

As always we start with a finitely generated group $G$, a family
$\{H_{\lambda}\}_{\lambda\in\Lambda}$ of subgroups of $G$, and for each
$\lambda\in\Lambda$, a nontrivial $H_{\lambda}$-almost invariant subset
$X_{\lambda}$ of $G$. We no longer assume that every $H_{\lambda}$ is finitely
generated. Instead we assume that if $H_{\lambda}$ is not finitely generated,
then $X_{\lambda}$ is associated to a splitting of $G$ over $H_{\lambda}$. Now
consider the construction in section \ref{regnbhds:construction}. We will
proceed exactly as we did there, and note the differences.

First we let $E$ denote the collection of all translates of the $X_{\lambda}%
$'s and their complements, and assume that the $X_{\lambda}$'s are in good
position, i.e. that $E$ satisfies Condition (*). This allows us to define the
partial order $\leq$ on $E$ exactly as before. For the proof that $\leq$ is a
partial order given in Lemma 1.14 of \cite{SS} does not use the finite
generation of the groups involved. As in our discussion at the end of section
\ref{regnbhds:construction}, $E$ may not be discrete. However, we can still
define the set $\overline{E}$ of pairs $\{X,X^{\ast}\}$ for $X\in E$, and can
define $P$ to be the collection of all CCC's of $\overline{E}$. Further the
arguments of section \ref{regnbhds:construction} still apply, and show that
the idea of betweenness can be defined on $P$ as before and this makes $P$
into a pretree. Of course, $P$ need not be discrete, but if it is, then as
before $P$ can be embedded in a $G$-tree $T$ and $G\backslash T$ is a graph of
groups structure for $G$ which we denote by $\Gamma(\{X_{\lambda}%
\}_{\lambda\in\Lambda}:G)$. This is a regular neighbourhood of the
$X_{\lambda}$'s. Of course, the $V_{0}$-vertex groups of this graph need not
be finitely generated, but the proof of Proposition \ref{Tisminimal} can still
be applied by using later results from this section to show that
$\Gamma(\{X_{\lambda}\}_{\lambda\in\Lambda}:G)$ is minimal and so must be a
finite graph.

Having dealt with the case when the $X_{\lambda}$'s are in good position, we
say that the family $\{X_{\lambda}\}_{\lambda\in\Lambda}$ is \textit{in good
enough position }if whenever we find incomparable elements $U$ and $V$ of $E$
which do not cross, there is some element $W$ of $E$ which crosses them. As
before, the above discussion applies equally well if the $X_{\lambda}$'s are
in good enough position. If this condition does not hold, we need to insist
that $E$ has only finitely many $G$-orbits which consist of isolated elements.
Now we recall the proof of Lemma \ref{candefineregnbhdwhennotingoodenoughposn}%
. This consisted of first replacing each isolated $X_{\lambda}$ by an
equivalent almost invariant set $Y_{\lambda}$ associated to a splitting, and
then replacing the $Y_{\lambda}$'s by equivalent almost invariant sets
$Z_{\lambda}$ such the collection of translates of the $Z_{\lambda}$'s and
their complements is nested. In the present situation some isolated
$X_{\lambda}$ may be over an infinitely generated group $H_{\lambda}$, but in
this case our hypothesis is that $X_{\lambda}$ is already associated to a
splitting, so we can simply take $Y_{\lambda}$ to equal $X_{\lambda}$. For the
second stage, Theorem \ref{splittingswithintersectionnumberzeroarecompatible}
tells us that we can still find the required sets $Z_{\lambda}$. Finally the
proof of Lemma \ref{equivalentsetshavesameregnbhd} still applies. We conclude
that if $E$ has only finitely many $G$-orbits which consist of isolated
elements, we can arrange that the $X_{\lambda}$'s are in good enough position
and then can define the pretree $P$. Further, if $P$ is discrete, one can
construct an algebraic regular neighbourhood $\Gamma(\{X_{\lambda}%
\}_{\lambda\in\Lambda}:G)$ of the $X_{\lambda}$'s which depends only on the
equivalence classes of the $X_{\lambda}$'s.

Finally, if the given family is finite, say $X_{1},\ldots,X_{n}$ we claim that
$P$ is discrete, so that $\Gamma(X_{1},\ldots,X_{n}:G)$ always exists in this
case. By renumbering we can suppose that $H_{1},\ldots,H_{k}$ are the only
non-finitely generated $H_{i}$'s. Let $E_{k}$ denote the set of all translates
of $X_{1},\ldots,X_{k}$ and their complements, and let $F$ denote the set of
all translates of $X_{k+1},\ldots,X_{n}$ and their complements. By replacing
$X_{1},\ldots,X_{k}$ by equivalent almost invariant sets, as discussed above,
we can arrange that $E_{k}$ is nested. Then we already know that $E_{k}$ and
$F$ are each discrete. As at the end of the proof of Lemma
\ref{twoideasofenclosingaresame}, any element of $E_{k}$ lies between two
elements of $F$ and vice versa. It follows that $E$ itself is discrete, so
that $P$ is discrete as required.

\section{Algebraic Regular Neighbourhoods: Uniqueness and
Applications\label{regnbhds:uniqueness}}

In this section, we finally give our precise definition of a regular
neighbourhood, and prove existence and uniqueness results which correspond
closely to the situation in topology. Then we discuss further generalisations
and some applications. The rest of this paper uses heavily the existence and
uniqueness of regular neighbourhoods.

As usual, let $G$ be a finitely generated group with a family of subgroups
$\{H_{\lambda}\}_{\lambda\in\Lambda}$. For each $\lambda\in\Lambda$, let
$X_{\lambda}$ denote a nontrivial $H_{\lambda}$-almost invariant subset of
$G$. In the preceding sections, we discussed how to construct a bipartite
graph of groups structure $\Gamma(\{X_{\lambda}\}_{\lambda\in\Lambda}:G)$. In
section \ref{regnbhds:construction}, we showed that this construction always
works if $\Lambda$ is finite and each $H_{\lambda}$ is finitely generated.
However Example \ref{exampleofregnbhdwithnonfgedgegroups} will show that the
edge groups of $\Gamma(X_{1},\ldots,X_{n}:G)$ need not be finitely generated
even if $G$ and each $H_{i}$ is finitely presented. Thus we need to consider
splittings over non-finitely generated subgroups of $G$. This is why we
formulate our definitions without assuming that the $H_{\lambda}$'s are
finitely generated, although in most of this paper, we will restrict to the
case when the $H_{\lambda}$'s are finitely generated.

We start by defining an algebraic regular neighbourhood of a family of
$X_{\lambda}$'s in $G$.

\begin{definition}
\label{defnofalgregnbhd} Let $G$ be a finitely generated group with a family
of subgroups $\{H_{\lambda}\}_{\lambda\in\Lambda}$. For each $\lambda
\in\Lambda$, let $X_{\lambda}$ denote a nontrivial $H_{\lambda}$-almost
invariant subset of $G$. Let $E$ denote the set of all translates of the
$X_{\lambda}$'s and their complements. Then an \textsl{algebraic regular
neighbourhood of the }$X_{\lambda}$\textsl{ 's in }$G$ is a bipartite graph of
groups structure $\Gamma$ for $G$ such that the following conditions hold:

\begin{enumerate}
\item Each $X_{\lambda}$ is enclosed by some $V_{0}$-vertex of $\Gamma$, and
each $V_{0}$-vertex of $\Gamma$ encloses some $X_{\lambda}$.

\item If $\sigma$ is a splitting of $G$ over a subgroup $H$ (which need not be
finitely generated) such that $\sigma$ does not cross any element of $E$, then
$\sigma$ is enclosed by some $V_{1}$-vertex of $\Gamma$.

\item $\Gamma$ is minimal.

\item Let $T$ denote the universal covering $G$-tree of $\Gamma$. We will say
that a $V_{0}$-vertex $v$ of $T$ is \textsl{isolated} if there is an element
$X$ of $E$ such that $v$ encloses $X$ and encloses only those elements of $E$
which are equivalent to $X$.

If $X$ is an isolated element of $E$ which is enclosed by a $V_{0}$-vertex $v$
of $T$, then either $v$ is isolated or $v$ encloses some non-isolated element
of $E$. Further there is a bijection $f$ from the isolated elements of $E$ to
the isolated vertices of $T$, such that $f(X_{\lambda})$ encloses $X_{\lambda
}$.
\end{enumerate}
\end{definition}

\begin{remark}
As $G$ is finitely generated, Condition 3) implies that $\Gamma$ is finite.
Clearly this condition is necessary if we want to prove any uniqueness result.

If Condition 4) holds, then the number of $G$-orbits of isolated elements of
$E$ must be finite.

Condition 4) is subtle, and only comes into play when there are isolated
$X_{\lambda}$'s. To understand the problems here, we again consider a finite
family of immersed curves $C_{\lambda}$ in a surface $M$, and consider the
problem of characterising their regular neighbourhood $N$. Condition 1) above
is analogous to asserting that each $C_{\lambda}$ is homotopic into $N$, and
that for each component of $N$, some $C_{\lambda}$ is homotopic into it.
Condition 2) is analogous to asserting that any simple curve on $M$ which has
intersection number zero with each $C_{\lambda}$ is homotopic into $M-N$. If
we assume that the $C_{\lambda}$'s are $\pi_{1}$-injective in $M$, and also
assume that each component of $\partial N$ is $\pi_{1}$-injective in $M$, then
these conditions characterise $N$, unless some $C_{\lambda}$ is homotopic to
an embedding disjoint from all the other $C_{\lambda}$'s. In this case, $N$
clearly has some annulus components, and the two conditions above do not
completely determine $N$, because they do not control the number of such
components. First, one can always add annulus components to $N$ parallel to
other such components without affecting the above two conditions. Second if
there are two $C_{\lambda}$'s which are simple and disjoint from each other
and from all other $C_{\lambda}$'s, so that $N$ has two parallel annulus
components, the subsurface of $M$ obtained from $N$ by simply deleting one of
these annulus components will still satisfy the above two conditions. If some
$C_{\lambda}$'s are arcs, there is another possible problem. Suppose that we
have three disjoint simple arcs $C_{1}$, $C_{2}$ and $C_{3}$ in $M$, such that
a component of $M$ cut along the $C_{i}$'s is a disc $D$ with copies of
$C_{1}$, $C_{2}$ and $C_{3}$ in its boundary. Let $N$ denote a regular
neighbourhood of the $C_{i}$'s, so that $N$ consists of three discs. If we
enlarge $N$ by adding $D$, the new submanifold still satisfies the two
conditions above which correspond to conditions 1) and 2) of our definition of
an algebraic regular neighbourhood. This example has nothing to do with the
triviality of the groups involved. Taking the product of this example with the
circle $S^{1}$ yields three annuli $C_{i}\times S^{1}$ in a $3$-manifold
$M\times S^{1}$, and we can enlarge a regular neighbourhood of these three
annuli by adding the solid torus $D\times S^{1}$. This example is related to
some subtleties in the topological JSJ-decomposition of a Haken $3$-manifold.
\end{remark}

Before we prove our existence and uniqueness results, we note the following
result which is an immediate consequence of our definition of a regular
neighbourhood and the fact that if a vertex of a $G$-tree encloses an almost
invariant subset $X$ of $G$, it also encloses any almost invariant subset $Y$
of $G$ which is equivalent to $X$.

\begin{lemma}
\label{regnbhdofXisisregnbhdofWis}Let $G$ be a finitely generated group with a
family of subgroups $\{H_{\lambda}\}_{\lambda\in\Lambda}$. For each
$\lambda\in\Lambda$, let $X_{\lambda}$ denote a nontrivial $H_{\lambda}%
$-almost invariant subset of $G$, and let $W_{\lambda}$ be a nontrivial
$K_{\lambda}$-almost invariant subset of $G$ which is equivalent to
$X_{\lambda}$. Assume further, that if $X_{\lambda}$ and $X_{\mu}$ are
isolated and $X_{\mu}=gX_{\lambda}$ , then $W_{\mu}=gW_{\lambda}$. Then a
bipartite graph of groups structure $\Gamma$ for $G$ is a regular
neighbourhood of the $X_{\lambda}$'s if and only if it is a regular
neighbourhood of the $W_{\lambda}$'s.
\end{lemma}

\begin{proof}
Checking the first three conditions of Definition \ref{defnofalgregnbhd} is
trivial. The technical hypothesis of the lemma about isolated elements is
required in order to ensure that the number of $G$-orbits of isolated elements
does not change when we replace the $X_{\lambda}$'s by the $W_{\lambda}$'s.
This allows us to check the fourth condition.
\end{proof}

One more definition will be very useful later in this paper. The above lemma
shows that a regular neighbourhood is really determined by a collection of
equivalence classes of almost invariant subsets of $G$. Thus it will be
convenient to define a regular neighbourhood of a family of such equivalence classes.

\begin{definition}
\label{defnofregnbhdofequivalenceclasses}Let $G$ be a finitely generated group
with a family $\mathcal{F}$ of equivalence classes of almost invariant
subsets. Then \textsl{an algebraic regular neighbourhood of} $\mathcal{F}$
\textsl{in} $G$ is an algebraic regular neighbourhood of a family of almost
invariant subsets of $G$ obtained by picking a representative of each
equivalence class in $\mathcal{F}$, subject to the condition that if $A$ and
$B$ are elements of $\mathcal{F}$ such that $B=gA$, for some $g$ in $G$, then
the representatives $X$ and $Y$ chosen for $A$ and $B$ must satisfy $Y=gX$.
\end{definition}

\begin{remark}
The reason for requiring equivariance in the choice of representatives is to
ensure that each equivalence class in $\mathcal{F}$ has a unique
representative. This condition is not needed unless $\mathcal{F}$ has isolated elements.

Let $F$ denote the collection of all the almost invariant subsets of $G$ which
represent elements of $\mathcal{F}$, and let $\Gamma$ denote an algebraic
regular neighbourhood of $\mathcal{F}$ in $G$. If no element of $F$ is
isolated, then $\Gamma$ is also an algebraic regular neighbourhood of the
collection $F$. However, if some element $X$ of $F$ is isolated, then the
collection $F$ does not have a regular neighbourhood, because it will contain
infinitely many distinct elements equivalent to $X$.
\end{remark}

Now we are ready to prove existence and uniqueness results for algebraic
regular neighbourhoods. For our existence result, we need to restrict to
finite families of almost invariant subsets of $G$, but our uniqueness result
does not need this restriction.

\begin{theorem}
(Existence of algebraic regular neighbourhoods) Let $G$ be a finitely
generated group, and for each $1\leq i\leq n$, let $H_{i}$ be a subgroup of
$G$, and let $X_{i}$ be a nontrivial $H_{i}$-almost invariant subset of $G$.
If $H_{i}$ is not finitely generated, then we assume that $X_{i}$ is
associated to a splitting of $G$ over $H_{i}$.

Then $\Gamma(X_{1},\ldots,X_{n}:G)$ is an algebraic regular neighbourhood of
the $X_{i}$'s in $G$.
\end{theorem}

\begin{remark}
It is natural to ask what happens if one considers an empty family of almost
invariant subsets of $G$. In this case, it will be convenient to say that the
graph of groups structure $\Gamma$ for $G$ which consists of a single $V_{0}%
$-vertex labelled $G$ is the regular neighbourhood.
\end{remark}

\begin{proof}
Recall that we constructed $\Gamma(X_{1},\ldots,X_{n}:G)$ in section
\ref{regnbhds:construction} in the case when every $H_{i}$ is finitely
generated, and constructed it at the end of the previous section in the
general case. We will write $\Gamma$ for $\Gamma(X_{1},\ldots,X_{n}:G)$.
Clearly $\Gamma$ satisfies Condition 4) of Definition \ref{defnofalgregnbhd},
as the $V_{0}$-vertices of $T$ are precisely the CCC's of $\overline{E}$.

Suppose that each $H_{i}$ is finitely generated. We showed in Lemma
\ref{XiisenclosedbyitsV0-vertex} that $\Gamma$ satisfies Condition 1). We
showed in Proposition \ref{XdoesnotcrossanyXiimpliesXisenclosedbyaV1vertex}
that $\Gamma$ satisfies Condition 2), and we showed in Lemma \ref{Tisminimal}
that $T$ is a minimal $G$-tree so that $\Gamma$ satisfies Condition 3). Thus
$\Gamma(X_{1},\ldots,X_{n}:G)$ is an algebraic regular neighbourhood of the
$X_{i}$'s in $G$. Note that in Proposition
\ref{XdoesnotcrossanyXiimpliesXisenclosedbyaV1vertex} we proved a result which
is stronger than Condition 2) in the case when $H$ is finitely generated, as
it applies to $H$-almost invariant subsets of $G$ which need not be associated
to splittings.

Now suppose that some $H_{i}$ is not finitely generated. The proof of
Proposition \ref{XdoesnotcrossanyXiimpliesXisenclosedbyaV1vertex} never uses
the hypothesis that the $H_{i}$'s are finitely generated. Thus $\Gamma$
satisfies Condition 2). The proof in Lemma \ref{Tisminimal} that $T$ is a
minimal $G$-tree needed the finite generation of the $H_{i}$'s only when Lemma
\ref{enclosedsetscanbechosennested} was used. But the proof of the last part
of Proposition \ref{XdoesnotcrossanyXiimpliesXisenclosedbyaV1vertex} shows
that Lemma \ref{enclosedsetscanbechosennested} remains true in our situation
because of our assumption that if $H_{i}$ is not finitely generated then
$X_{i}$ is associated to a splitting. Thus again $T$ is a minimal $G$-tree, so
that $\Gamma$ satisfies Condition 3).

It remains to show that $\Gamma$ satisfies Condition 1) of Definition
\ref{defnofalgregnbhd}. Let $U$ be an element of $E$ such that $Stab(U)$ is
not finitely generated, and let $v$ denote the $V_{0}$-vertex of $T$
determined by the CCC of $\overline{E}$ which contains $\overline{U}$. We must
show that $U$ is enclosed by $v$. As $U$ is associated to a splitting of $G$,
Let $s$ be an edge of $T$ which is incident to $v$ and oriented towards $v$,
and consider the proof of Lemma \ref{XiisenclosedbyitsV0-vertex}. This proof
needed the finite generation of the $H_{i}$'s only at the end when it used
Lemma \ref{containedinboundednbhdofA}. Thus the proof of Lemma
\ref{XiisenclosedbyitsV0-vertex} shows that, by replacing $U$ by $U^{\ast}$ if
needed, we can arrange that $gU^{(\ast)}\leq U$, for every $g\in Z_{s}^{\ast}%
$. Our assumption that $U$ is associated to a splitting of $G$ now implies
that $gU^{(\ast)}\subset U$, for every $g\in Z_{s}^{\ast}$. Further our
discussion after Lemma \ref{pre-imageofhalftreeisalmostinvariant} shows that
we can choose $U$ so that $U=\{g\in G:gU^{(\ast)}\subset U\}$. It follows that
$Z_{s}^{\ast}\subset U$. As this holds for every such edge $s$, it follows
that $U$ is enclosed by $v$ as required. This completes the proof that
$\Gamma$ satisfies Condition 1), and hence is an algebraic regular
neighbourhood of the $X_{i}$'s in $G$.
\end{proof}

\begin{theorem}
\label{uniquenessofregbhds}(Uniqueness of algebraic regular neighbourhoods)
Let $G$ be a finitely generated group with a family of subgroups
$\{H_{\lambda}\}_{\lambda\in\Lambda}$. For each $\lambda\in\Lambda$, let
$X_{\lambda}$ denote a nontrivial $H_{\lambda}$-almost invariant subset of
$G$. If $\Gamma_{1}$ and $\Gamma_{2}$ are algebraic regular neighbourhoods of
the $X_{\lambda}$ 's in $G$, then they are naturally isomorphic.
\end{theorem}

\begin{remark}
If $\Lambda$ is finite, then the construction of section
\ref{regnbhds:construction} yields a regular neighbourhood $\Gamma
(\{X_{\lambda}\}_{\lambda\in\Lambda}:G)$. It follows that all regular
neighbourhoods of the $X_{\lambda}$ 's in $G$ are isomorphic to $\Gamma
(\{X_{\lambda}\}_{\lambda\in\Lambda}:G)$. However it seems conceivable that,
when $\Lambda$ is infinite, the $X_{\lambda}$'s could possess a regular
neighbourhood but that the construction of section \ref{regnbhds:construction}
does not yield a regular neighbourhood because the pretree $P$ is not discrete.
\end{remark}

\begin{proof}
Let $\Gamma$ denote any algebraic regular neighbourhood of the $X_{\lambda}$'s
in $G$, let $T$ denote the universal covering $G$-tree of $\Gamma$, and let
$E$ denote the set of all translates of the $X_{\lambda}$'s and their
complements. Note that the existence of $\Gamma$ implies that $E$ contains
only finitely many $G$-orbits of isolated elements. Now Lemma
\ref{regnbhdofXisisregnbhdofWis} implies that we will lose nothing by assuming
that the $X_{\lambda}$'s are in good enough position. We need to prove some
general facts.

Suppose that an element $U$ of $E$ is enclosed by distinct vertices $v_{1}$
and $v_{2}$ of $T$. Part 8) of Lemma \ref{somefactsaboutenclosing} tells us
that $U$ is equivalent to $Z_{s}$ or to $Z_{s}^{\ast}$ for each edge $s$ on
the path joining $v_{1}$ and $v_{2}$. As each element of $E$ is enclosed by
some vertex of $T$, no element of $E$ crosses any $Z_{s}$. Hence $U$ must be
an isolated element of $E$. It follows that if $U$ is a non-isolated element
of $E$, then it is enclosed by a unique vertex of $T$. Next suppose that
$U_{1}$ and $U_{2}$ are elements of $E$ which are enclosed by distinct
vertices $v_{1}$ and $v_{2}$ of $T$. Then $U_{1}$ and $U_{2}$ do not cross, as
$\overline{Z_{s}}$ lies between $\overline{U_{1}}$ and $\overline{U_{2}}$, for
any edge $s$ of $T$ between $v_{1}$ and $v_{2}$. It follows that if $U$ is a
non-isolated element of $E$, then $U\;$is enclosed by a unique vertex $v$ of
$T$, and any element of $\overline{E}$ which lies in the CCC containing
$\overline{U}$ must also be enclosed by $v$.

Now let $\Gamma_{1}$ and $\Gamma_{2}$ denote two algebraic regular
neighbourhoods of the $X_{\lambda}$'s in $G$. We will first suppose that
$\Gamma_{1}$ and $\Gamma_{2}$ do not have any redundant vertices. In
particular, it follows that $\Gamma_{1}$ and $\Gamma_{2}$ do not have isolated
vertices. Hence Condition 4) implies that no $X_{\lambda}$ is isolated. Now it
also follows that each $X_{\lambda}$ is enclosed by a unique vertex in each
graph of groups. We will show that $\Gamma_{1}$ and $\Gamma_{2}$ have the same
conjugacy classes of edge splittings, which will then imply that they are
isomorphic by Theorem \ref{graphswithsameedgesplittingsareiso} as required.
Note that $V_{0}$-vertices of $\Gamma_{1}$ must correspond to $V_{0}$-vertices
of $\Gamma_{2}$ under this isomorphism because the $V_{0}$-vertices enclose
the $X_{\lambda}$'s and the $V_{1}$-vertices do not. Thus the isomorphism of
$\Gamma_{1}$ and $\Gamma_{2}$ automatically preserves their bipartite structure.

Let $\sigma$ be an edge splitting of $\Gamma_{2}$ over a subgroup $H$ of $G$.
(Recall that $\sigma$ is defined by collapsing $\Gamma_{2}$ with the interior
of an edge removed.) If $X$ is a $H$-almost invariant subset of $G$ associated
to $\sigma$, then $X$ does not cross any translate of any $X_{\lambda}$. Thus
Condition 2) for $\Gamma_{1}$ implies that $X$ is enclosed by some $V_{1}%
$-vertex of $\Gamma_{1}$. Now Lemma \ref{twoideasofenclosingaresame} shows
that we can refine the graph of groups $\Gamma_{1}$ by splitting at this
$V_{1}$-vertex using the edge splitting $\sigma$. If the new graph of groups
has a redundant vertex, this can only be because $\Gamma_{1}$ already had an
edge splitting conjugate to $\sigma$. Now let $\Gamma_{12}$ denote the graph
of groups structure for $G$ obtained from $\Gamma_{1}$ by splitting at $V_{1}%
$-vertices using those edge splittings of $\Gamma_{2}$ which are not conjugate
to any edge splitting of $\Gamma_{1}$. Thus $\Gamma_{12}$ has no redundant
vertices. Similarly let $\Gamma_{21}$ be obtained from $\Gamma_{2}$ by
splitting at $V_{1}$-vertices using those edge splittings of $\Gamma_{1}$
which are not conjugate to any edge splitting of $\Gamma_{2}$, so that
$\Gamma_{21}$ also has no redundant vertices. As $\Gamma_{12}$ and
$\Gamma_{21}$ have exactly the same conjugacy classes of edge splittings,
Theorem \ref{splittingswithintersectionnumberzeroarecompatible} tells us that
they are isomorphic. Now consider the universal covering $G$-trees $T_{12}$
and $T_{21}$. Although $T_{12}$ and $T_{21}$ are not bipartite, they still
have the property that each element of $E$ is enclosed by some vertex.
Further, as each element of $E$ is non-isolated, there is a unique vertex of
$T_{12}$ which encloses every element $U$ of $E$ which lies in a given CCC of
$\overline{E}$, and there is a unique vertex of $T_{21}$ with the same
property. Thus the $G$-isomorphism between $T_{12}$ and $T_{21}$ preserves
these vertices. For each such vertex, consider the edge splittings associated
to the incident edges. In one case they are all edge splittings of $\Gamma
_{1}$, and in the other case they are all edge splittings of $\Gamma_{2}$. The
$G$-isomorphism between the trees implies that each such edge splitting of
$\Gamma_{1}$ is conjugate to some edge splitting of $\Gamma_{2}$, and
conversely. Thus $\Gamma_{1}$ and $\Gamma_{2}$ have the same conjugacy classes
of edge splittings, as required.

Now consider the case when $\Gamma_{1}$ and $\Gamma_{2}$ may have redundant
vertices. Note that a $V_{0}$-vertex will be redundant if and only if it is
isolated. But even if there are no isolated $X_{\lambda}$'s, it is possible
for a $V_{1}$-vertex to be redundant. Now Condition 4) implies that a
non-isolated $V_{0}$-vertex must enclose a non-isolated element of $E$, and
hence is the unique vertex which encloses this element of $E$. We want to
apply the construction of the previous paragraph but first we need to remove
all the redundant vertices of $\Gamma_{1}$ and $\Gamma_{2}$, by amalgamating
suitable segments to a single edge. The resulting graphs of groups $\Gamma
_{1}^{\prime}$ and $\Gamma_{2}^{\prime}$ are no longer bipartite. But the
non-isolated $V_{0}$-vertices do not get removed. The argument of the previous
paragraph applies to show that $\Gamma_{1}^{\prime}$ and $\Gamma_{2}^{\prime}$
are isomorphic. Now $\Gamma_{1}$ and $\Gamma_{2}$ are obtained from this
common graph of groups structure by subdividing some edges, and Condition 4)
implies that the same edges get subdivided the same number of times, so that
$\Gamma_{1}$ and $\Gamma_{2}$ must be isomorphic. As before, this isomorphism
must preserve the non-isolated $V_{0}$-vertices of $\Gamma_{1}$ and
$\Gamma_{2}$, so it follows that the isomorphism must preserve the bipartite
structure, except possibly when every $V_{0}$-vertex is isolated. The same
argument applies to $V_{1}$-vertices also. We conclude that the isomorphism
must preserve the bipartite structure, except possibly when every vertex is
isolated. In this case, $\Gamma_{1}$ and $\Gamma_{2}$ will each be a circle,
and it is trivial to change the isomorphism to one which does preserve the
bipartite structure.
\end{proof}

Now we can use our results about regular neighbourhoods to give the result
promised in section \ref{prelim}, that almost invariant subsets of a group
with infinitely many ends are never $0$-canonical. Note that regular
neighbourhoods are not essential for this argument. They are simply convenient.

\begin{lemma}
\label{never0-canonical}Let $G$ be a finitely generated group with infinitely
many ends. If $X$ is any nontrivial $H$-almost invariant subset of $G$, where
$H$ is finitely generated, or $X$ is associated to a splitting of $G$ over
$H$, then $X$ is not $0$-canonical.
\end{lemma}

\begin{proof}
As $G$ has infinitely many ends, it admits a splitting $\sigma$ over a finite
subgroup $K$. Let $Y$ denote a standard $K$-almost invariant set associated to
$\sigma$. Suppose that $X$ is a nontrivial $H$-almost invariant subset of $G$
which is $0$-canonical. In particular, $X$ has zero intersection number with
$Y$. Now we consider the regular neighbourhood $\Gamma(X,Y:G)$, and its
universal covering $G$-tree $T$. As $\sigma$ is a splitting and has
intersection number zero with $X$, it follows that $Y$ is isolated and so
there is a corresponding isolated $V_{0}$-vertex of $T$. Let $v$ be a $V_{0}%
$-vertex of $T$ which encloses $X$, and pick edges $s$ and $t$ of $T$ with
finite stabiliser such that one lies in $\Sigma(X)$ and the other lies in
$\Sigma(X^{\ast})$. Choose $s$ and $t$ to be oriented away from $v$. Then $X$
crosses $Z_{s}\cup Z_{t}$, showing that $X$ is not $0$-canonical, as claimed.
\end{proof}

We can also give our example, promised in Remark
\ref{remarkthatedgegroupsneednotbefg}, of a finitely presented group $G$ which
splits over finitely presented subgroups $H$ and $K$, such that the algebraic
regular neighbourhood $\Gamma$ of these splittings has some edge and vertex
groups which are not finitely generated. We start with a general construction.

\begin{example}
\label{generalconstruction}This construction will give many examples of a
group $G$ which has two splittings with intersection number $1$, and also
yields the regular neighbourhood of these two splittings. Our construction is
based on the following topological picture. Consider two arcs $l$ and $m$
embedded properly in a surface $M$ so that each arc separates $M$ and the two
arcs meet transversely in a single point $w$. Thus a regular neighbourhood $N$
of the union of the two arcs has four boundary arcs and $M-N$ has four
components. Let $\Gamma$ denote the graph of groups structure for $G$
determined by $\partial N$. Thus $\Gamma$ is a tree which has a single vertex
$v_{0}$ corresponding to $N$, has four edges corresponding to the components
of $\partial N$ and four other vertices corresponding to the components of
$M-N$. The vertex $v_{0}$ and the four edges all carry the trivial group. We
will use this simple picture to guide us in constructing a group $G$ which
corresponds to $\pi_{1}(M)$ and possesses splittings over subgroups $H$ and
$K$ which correspond to the arcs $l$ and $m$. The group $H\cap K$ will
correspond to the point $w$. Finally the regular neighbourhood of the two
splittings will be a graph of the same combinatorial type as $\Gamma$, with a
single $V_{0}$-vertex corresponding to $v_{0}$.

To understand our idea, consider constructing $M$ by starting with a point
$w$, adding the four halves of $l$ and $m$ to obtain $l\cup m$, then
constructing $N$, and finally adding in the four remaining pieces of $M$. We
will follow a similar procedure, but we will use spaces with nontrivial
fundamental groups.

Start with a group $C$ and groups $A$, $B$, $D$ and $E$ which each contain $C$
as a proper subgroup. We will assume that the intersection of any two of these
groups is $C$. Let $H=A\ast_{C}B$ and let $K=D\ast_{C}E$. Think of $C$ as
corresponding to the point $w$, the groups $A$ and $B$ as corresponding to the
two halves of $l$, and the groups $D$ and $E$ as corresponding to the two
halves of $m$. Thus $H$ corresponds to $l$ and $K$ corresponds to $m$. Let
$L_{0}$ denote $H\ast_{C}K$, which corresponds to $\pi_{1}(N)$. Now the four
components of $\partial N$ naturally have corresponding groups which are
$A\ast_{C}D$, $D\ast_{C}B$, $B\ast_{C}E$ and $E\ast_{C}A$, each of which is a
subgroup of $L_{0}$. Denote these by $L_{1}$, $L_{2}$, $L_{3}$ and $L_{4}$
respectively. For $i=1,2,3,4$, pick a group $G_{i}$ which properly contains
$L_{i}$, and think of the $G_{i}$'s as corresponding to the components of
$M-N$. Now we define the group $G$ to be the fundamental group of the graph of
groups $\Gamma$ which is a tree with a vertex $v_{0}$ with associated group
$L_{0}$, with four edges attached to $v_{0}$ which carry the $L_{i}$'s, and
with the four other vertices carrying the $G_{i}$'s. Let $v_{i}$ denote the
vertex which carries $G_{i}$.

To understand this construction topologically, one needs to build a space with
fundamental group $G$ which mimics the structure of our initial example $M$.
We pick spaces $M_{A}$, $M_{B}$, $M_{C}$, $M_{D}$, and $M_{E}$ with
fundamental groups $A$, $B$, $C$, $D$ and $E$ respectively, such that each
space contains $M_{C}$ and the intersection of any two equals $M_{C}$. Let $Z$
denote the union of these spaces, so that $\pi_{1}(Z)=L_{0}$. Then for each
$i\geq1$, we choose a space $M_{i}$ with fundamental group $L_{i}$, take its
product with the unit interval and glue one end to $Z$ using the inclusion of
$L_{i}$ into $L_{0}$. Finally, for each $G_{i}$, we choose a space with
fundamental group $G_{i}$ and glue the other end of $M_{i}\times I$ to it
using the inclusion of $L_{i}$ in $G_{i}$. The resulting space $M$ has
fundamental group $G$. Further, its structure clearly yields splittings
$\sigma$ and $\tau$ over $H$ and $K$ respectively. For $\sigma$ is the
splitting of $G$ obtained by ``cutting along'' $M_{A}\cup M_{B}$, and $\tau$
is the splitting of $G$ obtained by ``cutting along'' $M_{D}\cup M_{E}$.
Consider the pre-image $\widetilde{Z}$ of $Z$ in the universal cover
$\widetilde{M}$ of $M$. If $G$ is finitely generated, then it is easy to see
that the pretree constructed combinatorially in section
\ref{regnbhds:construction} is the same as the pretree of components of
$\widetilde{Z}$. (The main point to notice is that the stabiliser of a
component of $\widetilde{Z}$ equals $L_{0}$, which also equals the stabiliser
of the corresponding CCC. Thus $Z$ is in `good position'.) It follows that
$\sigma$ and $\tau$ have intersection number $1$, and that $\Gamma$ is their
regular neighbourhood in $G$.
\end{example}

By making interesting choices of the groups involved in the above
construction, one can give many interesting examples. Here is an important
example which explains why we have spent so much time considering splittings
over non-finitely generated subgroups.

\begin{example}
\label{exampleofregnbhdwithnonfgedgegroups}We give here an example of a
finitely presented group $G$ which splits over finitely presented subgroups
$H$ and $K$, such that the algebraic regular neighbourhood $\Gamma$ of these
splittings has an edge group and a vertex group which is not finitely
generated. We do this by making choices of the groups involved in the
construction of Example \ref{generalconstruction}. The edge group $L_{1}$ and
the vertex group $G_{1}$ of the regular neighbourhood are not finitely generated.

In Example \ref{generalconstruction}, choose $C=F_{\infty}$, the free group of
countably infinite rank, choose $B$ and $E$ to be $F_{2}$, the free group of
rank $2$, and choose $A$ and $D$ to be $F_{2}\ast C$. The inclusions of $C$ in
$A$ and $D$ are the obvious ones. The inclusions of $C$ in $B$ and $E$ can be
chosen in any reasonable way. A good example would be to map the $i$-th basis
element of $C$ to $u^{-i}vu^{i}$, where $u$ and $v$ are the basis elements of
$F_{2}$. Thus $H=A\ast_{C}B$ and $K=D\ast_{C}E$ are each isomorphic to
$(F_{2}\ast C)\ast_{C}F_{2}$ which is simply $F_{4}$, the free group of rank
$4$. In order to complete the construction of the group $G$, we need to choose
the groups $G_{i}$. Recall that the groups $L_{1}$, $L_{2}$, $L_{3}$ and
$L_{4}$ are respectively isomorphic to $A\ast_{C}D$, $D\ast_{C}B$, $B\ast
_{C}E$ and $E\ast_{C}A$. This means that $L_{1}$ is not finitely generated,
though the remaining $L_{i}$'s are finitely generated. Further, $L_{2}$ and
$L_{4}$ are each isomorphic to $F_{4}$. The group $L_{3}=B\ast_{C}E=F_{2}%
\ast_{C}F_{2}$ is finitely generated but is not finitely presented. It is
trivial to choose finitely presented groups $G_{2}$ and $G_{4}$ which properly
contain $L_{2}$ and $L_{4}$ respectively. We can use Higman's Embedding
Theorem \cite{Higman} to find a finitely presented group $G_{3}$ which
contains $L_{3}$. Finally, we choose $G_{1}$ to be any group of the form
$P\ast_{Q}L_{1}$, where $Q$ is finitely generated and $P$ is finitely
presented. Now we claim that $G$ is finitely presented. To see this, we build
up $G$ in stages. Recall that, by construction, each of $G_{2}$, $G_{3}$ and
$G_{4}$ is finitely presented. Thus the group $G_{2}\ast_{B}G_{3}$ is finitely
presented. Hence the group $\left(  G_{2}\ast_{B}G_{3}\right)  \ast_{E}G_{4}$
is finitely presented. Denote this last group by $G_{5}$. Then $G=G_{1}%
\ast_{L_{1}}G_{5}=\left(  P\ast_{Q}L_{1}\right)  \ast_{L_{1}}G_{5}=P\ast
_{Q}G_{5}$, which is finitely presented because $P$ and $G_{5}$ are finitely
presented and $Q$ is finitely generated.
\end{example}

\begin{remark}
In the above example, the group $G$ has subgroups which are finitely generated
but not finitely presented. This raises the question of whether examples such
as the above can exist for a group $G$ which is coherent, i.e. finitely
generated subgroups are finitely presented. We have no ideas about how to
answer this question.
\end{remark}

We next consider an application which strengthens a result of Niblo in
\cite{Niblo} on the existence of splittings of a given group. Given a
nontrivial $H$-almost invariant subset $X$ of a group $G$, we consider the
regular neighbourhood $\Gamma(X:G)$. Recall that Proposition \ref{Tisminimal}
implies that $\Gamma(X:G)$ is minimal. It follows that unless it consists of a
single vertex, then any edge will yield a splitting of $G$. We define the
subgroup $S(X)$ of $G$, to be the stabiliser of the CCC of $\overline{E}$
which contains $\overline{X}$ in the construction of $\Gamma(X:G)$, so that
$S(X)$ is the vertex group for the corresponding $V_{0}$-vertex of
$\Gamma(X:G)$. Thus we immediately deduce the following result.

\begin{corollary}
Let $G$ be a finitely generated group with finitely generated subgroup $H$,
and let $X$ be a nontrivial $H$-almost invariant subset of $G$. If $S(X)$ is
not equal to $G$, then $G$ splits over a subgroup of $S(X)$.
\end{corollary}

\begin{remark}
As $S(X)$ was defined in terms of the regular neighbourhood $\Gamma(X:G)$,
Lemma \ref{regnbhdofXisisregnbhdofWis} implies that if $X$ and $Y$ are
equivalent, then $S(X)=S(Y)$.
\end{remark}

In order to understand the implication of this, we need to consider the group
$S(X)$ more carefully. Let $v$ denote the $V_{0}$-vertex of $T$ which
corresponds to the CCC of $\overline{E}$ which contains $\overline{X}$. If $X$
is not isolated in the set $E$ of all translates of $X$ and $X^{\ast}$, then
the argument at the end of the proof of Theorem \ref{Pisapretree} shows that
$S(X)$ is generated by $H$ and $\{g\in G:gX$ crosses $X\}$. In \cite{Niblo},
Niblo defined a group $T(X)$ which is the subgroup of $G$ generated by $H$ and
$\{g\in G:gX$ and $X$ are not nested\}. He proved that if $T(X)\neq G$, then
$G$ splits over a subgroup of $T(X)$. Clearly $S(X)$ is a subgroup of $T(X)$
in this case, so our result implies his when $X$ is not isolated. If $X$ is
isolated, then Theorem 2.8 of \cite{SS} implies that $G$ splits over a
subgroup commensurable with $H$.

An interesting related result due to Dunwoody and Roller \cite{D-Roller} is
that if $S(X)$ is contained in $Comm_{G}(H)$, then $G$ splits over a subgroup
commensurable with $H$ even if $G=Comm_{G}(H)$. See \cite{SS} and
{\cite{Niblo} for alternative proofs of this fact}.

The above Corollary was proved by considering the regular neighbourhood
$\Gamma(X:G)$ and using the fact that $S(X)$ is a vertex group. If $G$ has a
$H$-almost invariant subset $X$ and a $K$-almost invariant subset $Y$ such
that $H\backslash X$ and $K\backslash Y$ have intersection number zero, then
considering the regular neighbourhood $\Gamma(X,Y:G)$ yields the following
result, which strengthens another result of Niblo in \cite{Niblo}.

\begin{corollary}
Let $G$ be a finitely generated group with finitely generated subgroups $H$
and $K$. If $G\;$has a $H$-almost invariant subset $X$ and a $K$-almost
invariant subset $Y$ such that $H\backslash X$ and $K\backslash Y$ have
intersection number zero, then $G$ has a minimal graph of groups decomposition
with at least two edges, in which two of the edge groups are a subgroup of
$S(X)$ and a subgroup of $S(Y)$.
\end{corollary}

Now we can use this corollary to give another generalisation of Theorem
\ref{Theorem2.5ofSS}. In that result, we showed that if $G$ is a finitely
generated group with splittings $\sigma$ and $\tau$ over finitely generated
subgroups $H$ and $K$ of $G$ such that $\sigma$ and $\tau$ have zero
intersection number, then $\sigma$ and $\tau$ are compatible. In Theorem
\ref{splittingswithintersectionnumberzeroarecompatible} of this paper, we
generalised Theorem \ref{Theorem2.5ofSS} to the case of splittings over
subgroups $H$ and $K$ which need not be finitely generated. A natural question
is whether an analogous result holds for almost invariant sets which are not
associated to splittings.

An equivalent formulation of Theorem \ref{Theorem2.5ofSS} is that if $X$ is a
$H$-almost invariant subset of $G$ and $Y$ is a $K$-almost invariant subset of
$G$ such that each is associated to a splitting of $G$ and $H\backslash X$ and
$K\backslash Y$ have intersection number zero, then $X$ and $Y\;$are
equivalent to subsets $X^{\prime}$ and $Y^{\prime}$ of $G$ such that the set
$E$ of all translates of $X^{\prime}$ and $Y^{\prime}$ is nested. We will
prove the following result which is the natural analogue when $X$ and $Y$ are
not associated to splittings.

\begin{lemma}
Let $G$ be a finitely generated group with finitely generated subgroups $H$
and $K$. Let $X$ be a $H$-almost invariant subset of $G$ and let $Y$ be a
$K$-almost invariant subset of $G$ such that $H\backslash X$ and $K\backslash
Y$ have intersection number zero. Then $X$ and $Y\;$are equivalent to subsets
$X^{\prime}$ and $Y^{\prime}$ of $G$ such that any translate of $X^{\prime}$
and any translate of $Y^{\prime}$ are nested.
\end{lemma}

\begin{remark}
We now have two generalisations of Theorem \ref{Theorem2.5ofSS}, one of which
allows splittings over non-finitely generated subgroups and the other replaces
splittings over finitely generated subgroups $H$ and $K$ by almost invariant
subsets over $H$ and $K$. It seems natural to ask if there is a common
generalisation for arbitrary almost invariant sets over arbitrary subgroups of
$G$. We have no methods for answering such a question but we think that such a
generalisation is unlikely.
\end{remark}

\begin{proof}
Let $\Gamma$ denote $\Gamma(X,Y:G)$, as constructed in section
\ref{regnbhds:construction}. Then $\Gamma$ has two $V_{0}$-vertices, one of
which encloses $X$ and the other encloses $Y$. Now, Lemma
\ref{enclosedsetscanbechosennested} shows that $X$ and $Y\;$are equivalent to
subsets $X^{\prime}$ and $Y^{\prime}$ of $G$ such that $X^{\prime}$ is nested
with respect to all the $Z_{s}$ and $Z_{s}^{\ast}$ and so is $Y^{\prime}$. In
particular, the claim follows.
\end{proof}

\section{Coends when the commensuriser is small
\label{coendswhencommensuriserissmall}}

The rest of this paper consists of understanding regular neighbourhoods of
certain families of almost invariant subsets of a group and showing that in
several interesting cases, certain infinite families of such subsets possess a
regular neighbourhood. For the results about infinite families of almost
invariant subsets, we will need to assume that $G$ is finitely presented, but
most of our results about finite families work without this additional
restriction. We will be interested in a one-ended, finitely generated group
$G$ and almost invariant subsets of $G$ which are over virtually polycyclic
(VPC) subgroups. In this and the following three sections, we will mainly be
interested in the case of VPC subgroups of length $1$. However many of our
results are valid for VPC subgroups of length $n\geq1$ assuming that $G$ has
no almost invariant subsets over VPC subgroups of length $<n$, and in this
section we will prove several technical results in that generality for later
use. Any VPC group of length $1$ is virtually infinite cyclic, or equivalently
two-ended, and for brevity we will use the phrase two-ended in what follows.
In section \ref{JSJforlargecommensurisers}, we will show that if $G$ is
finitely presented, there is a regular neighbourhood of all the equivalence
classes of nontrivial almost invariant subsets of $G$ over all two-ended
subgroups. It turns out that if $H$ is a two-ended subgroup of $G$ such that
there is a nontrivial almost invariant subset over $H$, then the commensuriser
$Comm_{G}(H)$ of $H$ in $G$ plays an important role. We will analyse the role
of the commensuriser in this and the next section.

We will say that a subgroup $H$ of $G$ \textit{has small commensuriser in} $G$
if $Comm_{G}(H)$ contains $H$ with finite index. In this section we will
consider a one-ended finitely generated group $G$ and almost invariant subsets
over two-ended subgroups. To study these, we need Bowditch's results from
\cite{B2} as well as a non-standard accessibility result. We start by quoting
a result in \cite{B2}, but reformulated in the language of almost invariant
sets (see \cite{SS2} for a similar result for hyperbolic groups).

\begin{proposition}
\label{twocoends} Let $G$ be a one-ended finitely generated group, and let $X$
and $Y$ be nontrivial almost invariant subsets over two-ended subgroups $H$
and $K$. If $X$ crosses $Y$ strongly, then $Y$ crosses $X$ strongly and the
number of coends of both $H$ and $K$ is two.
\end{proposition}

We will prove a more general result in section \ref{JSJforVPCoftworanks}. The
following result tells us what happens when $X$ and $Y$ cross weakly.

\begin{proposition}
\label{commensurable} Let $G$ be a one-ended finitely generated group, and let
$X$ and $Y$ be nontrivial almost invariant subsets over two-ended subgroups
$H$ and $K$. If $X$ crosses $Y$ weakly, then $H$ and $K$ are commensurable.
\end{proposition}

\begin{proof}
As $H$ and $K$ are virtually infinite cyclic, either they are commensurable or
$H\cap K$ is finite. We will suppose that $H\cap K$ is finite and derive a contradiction.

By Proposition \ref{twocoends}, as $X$ crosses $Y$ weakly, we know that $Y$
crosses $X$ weakly. Thus one of $\delta Y\cap X^{(\ast)}$ is $H$-finite and
one of $\delta X\cap Y^{(\ast)}$ is $K$-finite. By changing notation if
necessary, we can arrange that $\delta Y\cap X$ is $H$-finite and $\delta
X\cap Y$ is $K$-finite. As $\delta X$ is $H$-finite and $\delta Y$ is
$K$-finite, it follows that each of $\delta Y\cap X$ and $\delta X\cap Y$ is
both $H$-finite and $K$-finite. Thus they are both $(H\cap K)$-finite. Now
consider the coboundary $\delta(X\cap Y)$. Every edge in this coboundary meets
$\delta Y\cap X$ or $\delta X\cap Y$. Hence $\delta(X\cap Y)$ is also $(H\cap
K)$-finite. As $X\cap Y$ is clearly invariant under the left action of $H\cap
K$, it is $(H\cap K)$-almost invariant. As we are assuming that $H\cap K$ is
finite, this means that $G$ has more than one end, which is the required contradiction.
\end{proof}

We now recall from \cite{SS} that a pair of finitely generated groups $(G,H)$
is of \textit{surface type} if $e(G,H^{\prime})=2$ for every subgroup
$H^{\prime}$ of finite index in $H$ and $e(G,H^{\prime})=1$ for every subgroup
$H^{\prime}$ of infinite index in $H$. It follows that $(G,H)$ has two coends.
Conversely, suppose that $(G,H)$ has two coends. Then $H$ has a subgroup
$H_{1}$ of index at most $2$ such that $e(G,H_{1})=2$ and hence the pair
$(G,H_{1})$ is of surface type. The following result will be useful.

\begin{proposition}
\label{twocoendscrossing} Let $G$ be a finitely generated group with finitely
generated subgroups $H$ and $K$, a nontrivial $H$-almost invariant subset $X$,
and a nontrivial $K$-almost invariant subset $Y$. Suppose also that the number
of coends of $H$ in $G$ is $2$. Then $Y$ crosses $X$ if and only if $Y$
crosses $X$ strongly.
\end{proposition}

\begin{proof}
By replacing $H$ by the above subgroup $H_{1}$ of index at most $2$, we can
assume that the pair $(G,H)$ is of surface type. Now Proposition 3.7 of
\cite{SS} tells us that if $(G,H)$ is a pair of finitely generated groups of
surface type, $X$ is a nontrivial $H$-almost invariant subset of $G$ and $Y$
is a nontrivial $K$-almost invariant subset of $G$, then $Y$ crosses $X$ if
and only if $Y$ crosses $X$ strongly. The result follows.
\end{proof}

The following result summarises the above in the form which we will need.

\begin{proposition}
\label{crossingsareallstrongorallweak}Let $G$ be a one-ended finitely
generated group and let $\{X_{\lambda}\}_{\lambda\in\Lambda}$ be a family of
nontrivial almost invariant subsets over two-ended subgroups of $G$. As usual,
let $E$ denote the set of all translates of the $X_{\lambda}$'s and their
complements. Form the pretree $P$ of cross-connected components (CCC's) of
$\overline{E}$ as in the construction of regular neighbourhoods in section
\ref{regnbhds:construction}. Then the following statements hold:

\begin{enumerate}
\item The crossings in a CCC of $\overline{E}$ are either all strong or are
all weak.

\item In a CCC with all crossings weak, the stabilisers of the corresponding
elements of $E$ are all commensurable.
\end{enumerate}
\end{proposition}

\begin{proof}
1) If $X$ and $Y$ are elements of $E$ which cross strongly, then Propositions
\ref{twocoendscrossing} and \ref{twocoends} imply not only that $Y$ must cross
$X$ strongly, but the same applies to any other element of $E$ which crosses
$X$. Hence all crossings in the CCC determined by $X$ and $Y$ are strong. It
follows that the crossings in a CCC of $\overline{E}$ are either all strong or
are all weak, as required.

2) If a CCC has weak crossing, then Proposition \ref{commensurable} implies
that any two elements of this CCC have commensurable stabilisers, which proves
the required result.
\end{proof}

In a minimal graph of groups decomposition $\Gamma$ of a group $G$, we will
say that a vertex $v$ is \textit{of finite-by-Fuchsian type}, or that the
associated vertex group $G(v)$ is of finite-by-Fuchsian type, if $G(v)$ is a
finite-by-Fuchsian group, where the Fuchsian group is not finite nor
two-ended, and there is exactly one edge of $\Gamma$ which is incident to $v$
for each peripheral subgroup $K$ of $G(v)$ and this edge carries $K$. If
$G=G(v)$, then the Fuchsian quotient group corresponds to a closed orbifold.
We should note that usually a Fuchsian group means a discrete group of
isometries of the hyperbolic plane, but in this paper, it will be convenient
to include also discrete groups of isometries of the Euclidean plane. As we
are excluding finite and two-ended Fuchsian groups, the extra groups this
includes are all virtually $\mathbb{Z}\times\mathbb{Z}$.

Now suppose that the family of $X_{\lambda}$'s in the above proposition is
finite. Thus their regular neighbourhood $\Gamma(\{X_{\lambda}\}_{\lambda
\in\Lambda}:G)$ exists by Theorem \ref{Pisapretree}. Consider a $V_{0}$-vertex
$v$ of $\Gamma(\{X_{\lambda}\}_{\lambda\in\Lambda}:G)$ which comes from a CCC
of $\overline{E}$ in which all crossing is strong. If $X_{\lambda}$ lies in
the CCC corresponding to $v$ and is not isolated, Lemma \ref{twocoends} shows
that $H_{\lambda}$ has two coends in $G$. Bowditch \cite{B2} and
Dunwoody-Swenson \cite{D-Swenson} have shown that the enclosing group $G(v)$
is of finite-by-Fuchsian type. Bowditch deals only with the case when
$H_{\lambda}$ has two coends in $G$ for each $\lambda\in\Lambda$. In this
case, our construction of regular neighbourhoods coincides with Bowditch's
construction of enclosing groups as both use the same pretree. In fact our
construction of regular neighbourhoods is suggested by Bowditch's use of
pretrees in \cite{B1} and \cite{B2}. We state his result in the form that we
will use later. We will also need similar results for VPC subgroups of length
greater than $1$, but we will simply observe that the
Bowditch-Dunwoody-Swenson arguments go through in general. The following
result is contained in Propositions 7.1 and 7.2 of Bowditch's paper \cite{B2},
and in the JSJ-decomposition theorem of \cite{D-Swenson} but we have
reformulated it using our regular neighbourhood terminology.

\begin{theorem}
\label{strongcrossingimpliesFuchsiantype}Let $G$ be a one-ended finitely
generated group with a finite family of two-ended subgroups $\{H_{\lambda
}\}_{\lambda\in\Lambda}$. For each $\lambda\in\Lambda$, let $X_{\lambda}$
denote a nontrivial $H_{\lambda}$-almost invariant subset of $G$, let $E$
denote the set of all translates of the $X_{\lambda}$'s and their complements,
and let $\Gamma$ denote the regular neighbourhood of the $X_{\lambda}$'s. Let
$X$ denote an element of $E$, let $H$ denote its stabiliser, and let $v$
denote a vertex of $\Gamma$ which encloses $X$.

Suppose $H$ has two coends and that there exists an element of $E$ which
crosses $X$. Then the vertex group $G(v)$ is of finite-by-Fuchsian type, and
$H$ is not commensurable with a peripheral subgroup of $G(v)$.
\end{theorem}

\begin{remark}
If $\Gamma$ consists of a single vertex, so that $G=G(v)$, then $G$ must
itself be of finite-by-Fuchsian type.
\end{remark}

In a Fuchsian group $\Sigma$, any two-ended subgroup has small commensuriser
unless $\Sigma$ is virtually $\mathbb{Z}\times\mathbb{Z}$. When combined with
the above theorem, this yields the following result.

\begin{corollary}
\label{strongcrossingimpliessmallcomensuriser}Let $G$ be a one-ended, finitely
generated group and suppose that $X$ is a nontrivial almost invariant subset
over a two-ended subgroup $H$. Suppose some almost invariant set over a
two-ended subgroup $K$ crosses $X$ strongly. Then either $G$ is of
finite-by-Fuchsian type or $H$ has small commensuriser in $G$. Further, if $G$
has a finite normal subgroup $K$ with Fuchsian quotient $\Sigma$, then either
$H$ has small commensuriser in $G$, or $\Sigma$ is virtually $\mathbb{Z}%
\times\mathbb{Z}$.
\end{corollary}

Later we will want to consider infinite families of such almost invariant
subsets of $G$. We want to do this by taking increasing finite families of
such sets and showing that the graphs of groups structures for $G$ obtained in
this way must stabilise. If all crossings of such subsets of $G$ are strong,
one can use Theorem \ref{strongcrossingimpliesFuchsiantype} to show that this
happens for homological reasons. However if weak crossings occur this argument
does not work. We will need to assume that $G$ is finitely presented and to
use variants of previous accessibility results.

Recall that, in a graph of groups decomposition, we call a vertex redundant if
it has valence at most two, it is not the vertex of a loop, and each edge
group includes by an isomorphism into the vertex group. Recall that a vertex
is reducible if it has two incident edges, it is not the vertex of a loop, and
one of the incident edge groups includes by an isomorphism into that vertex
group. For a finitely presented group $G$, the main result of \cite{B-F} gives
a bound on the complexity of reduced, nontrivial graphs of groups
decompositions of $G$ with all edge groups being small groups. We will use
their result to prove the following.

\begin{theorem}
\label{graphsstabilise}Let $G$ be a finitely presented group and suppose that
$G$ does not split over VPC groups of length less than $n$. For each positive
integer $k$, let $\Gamma_{k}$ be a graph of groups decomposition of $G$
without redundant vertices and with all edge groups VPC groups of length $n$,
and suppose that for each $k$, $\Gamma_{k+1}$ is a refinement of $\Gamma_{k}$.
Then the sequence $\Gamma_{k}$ stabilises.
\end{theorem}

\begin{proof}
Since $G$ is finitely presented, Bestvina-Feighn's accessibility result
\cite{B-F} implies the theorem provided the $\Gamma_{k}$'s are reduced. Thus
we only have to bound the length of chains of splittings of $G$ over
descending subgroups (unfoldings in the language of \cite{RS}). This was done
in \cite{RS} when $G$ is finitely generated and torsion-free and $n=1$. In
\cite{B3}, Bowditch gave a much simpler argument using tracks when $G$ is
finitely presented and $n=1$. We give our argument which is similar to
Bowditch's but which was arrived at independently.

Suppose that we have an infinite sequence of splittings of $G$ over descending
VPC subgroups $H_{i}$ of length $n$. Thus $\cap_{i\geq1}H_{i}$ is VPC of
length $n-1$. We fix a finite $2$-complex $K$ with fundamental group $G$ and
universal cover $\widetilde{K}$. For each $m\geq1$, there is a $G$-tree
$T_{m}$ which corresponds to the first $m$ splittings, and $T_{m+1}$ is a
refinement of $T_{m}$, i.e. there is a natural collapsing map $q_{m+1}%
:T_{m+1}\rightarrow T_{m}$. We now pick $G$-equivariant linear maps
$p_{m}:\widetilde{K}\rightarrow T_{m}$ such that $p_{m}=q_{m+1}p_{m+1}$, let
$W_{m}$ denote the midpoints of the edges of $T_{m}$ and consider $p_{m}%
^{-1}(W_{m})$. This is a $G$-invariant pattern in $\widetilde{K}$, which
projects to a finite pattern $L_{m}$ in $K$. By construction of the maps
$p_{m}$, we have $L_{m}\subset L_{m+1}$. Each component of $L_{m}$ carries a
subgroup of $H_{m}$. Since $G$ does not split over a VPC subgroup of length
less than $n$, $L_{m+1}-L_{m}$ must have at least one component $C_{m}$ with
stabiliser which is VPC\ of length $n$. Now Dunwoody (see \cite{D-D}) showed
that there is an upper bound on the number of non-parallel disjoint tracks one
can have in $K$. In particular, it follows that the $C_{m}$'s carry only
finitely many distinct subgroups of $G$, and hence that the descending
sequence of $H_{i}$'s must stabilise, which is a contradiction.
\end{proof}

Let $\Gamma$ be a graph of groups decomposition of a group $G$ without
redundant vertices. We will say that $\Gamma$ is \textit{maximal with respect
to two-ended subgroups}, if whenever a vertex encloses a splitting over a
two-ended subgroup, this splitting is already an edge splitting of $\Gamma$.
This means that if we form a refinement of $\Gamma$ by splitting at some
vertex so that the extra edge splitting is over a two-ended group, this
refinement must have a redundant vertex. The above result in particular
implies the following.

\begin{corollary}
\label{fpGhasmaximaldecompovertwo-endedgroups}A one-ended, finitely presented
group has maximal decompositions with respect to two-ended subgroups.
\end{corollary}

Similarly we call a decomposition of $G$ maximal with respect to VPC groups of
length $\leq n$, if it cannot be refined without introducing redundant
vertices by splitting at a vertex along a VPC group of length $\leq n$.

The proof of Theorem \ref{graphsstabilise} applies essentially unchanged to
yield the following result.

\begin{theorem}
\label{graphsstabilise2} Let $G$ be a finitely presented group and let
$\Gamma_{0}$ be a graph of groups decomposition of $G$ which is maximal with
respect to splittings over VPC groups of length $\leq n$. For each $k\geq1$,
let $\Gamma_{k}$ be a graph of groups decomposition of $G$ without redundant
vertices, and suppose that for each $k\geq0$, $\Gamma_{k+1}$ is a refinement
of $\Gamma_{k}$. Suppose further that all the edge splittings of $\Gamma_{k}$
which are not edge splittings of $\Gamma_{0}$ are over VPC subgroups of length
$(n+1)$. Then the sequence $\Gamma_{k}$ stabilises.
\end{theorem}

Now consider a two-ended subgroup $H$ of $G$. Let $Q(H)$ denote the collection
of all almost invariant subsets of $G$ which are over subgroups of $G$
commensurable with $H$, and let $F(H)$ denote the subset of $Q(H)$ which
consists of all the trivial elements. (Recall that a $H$-almost invariant set
is trivial if it is $H$-finite.) Clearly if $X$ and $Y$ lie in $Q(H)$, then
$X\cap Y$, $X+Y$ and $X\cup Y$ also lie in $Q(H)$. Thus $Q(H)$ is a subalgebra
of the Boolean algebra of all subsets of $G$. Also $F(H)$ is an ideal in
$Q(H)$. We let $B(H)$ denote the quotient Boolean algebra $Q(H)/F(H)$. Thus
$B(H)$ is precisely the collection of equivalence classes of elements of
$Q(H)$. Note that if $H$ has small commensuriser so that $Comm_{G}(H)$
contains $H$ with finite index, then $Comm_{G}(H)$ is also two-ended

\begin{theorem}
\label{smallcommensuriserimpliesB(H)isfinite}Let $G$ be a one-ended, finitely
presented group with a two-ended subgroup $H$ with small commensuriser.
Suppose that there is an element $X$ of $Q(H)$ such that no almost invariant
set over a two-ended subgroup of $G$ crosses $X$ strongly. Then $B(H)$ is finite.
\end{theorem}

\begin{remark}
Proposition \ref{twocoends} shows that the conclusion of this result remains
true if $X$ does cross strongly some almost invariant set over a two-ended subgroup.
\end{remark}

\begin{proof}
Lemma \ref{strongcrossingisdeterminedbygroups} implies that the assumption on
$X$ holds for every element of $Q(H)$. If there is an almost invariant subset
$Y$ of $G$ which is over some two-ended subgroup and crosses $X$, then it must
do so weakly and Lemma \ref{commensurable} implies that the stabilisers of $X$
and $Y$ are commensurable. This implies that $Y$ also lies in $Q(H)$.

Given a finite subset $\{U_{j}\}_{j\in J}$ of $Q(H)$, we can form $E$, $T$ and
the regular neighbourhood $\Gamma(\{U_{j}\}_{j\in J}:G)$, which we will denote
by $\Gamma$. Consider a $V_{0}$-vertex $v$ of $T$ such that the corresponding
CCC of $\overline{E}$ contains $\overline{U}$, for some element $U$ of the
$U_{j}$'s. The other elements of $E$ which lie in this CCC must all lie in
$Q(H)$, by the preceding paragraph. Now we consider $Stab(v)$, the stabiliser
of $v$. Recall that there are a finite number of elements $X_{i}$ of $E$
enclosed by $v$ such that $Stab(v)$ is generated by the stabilisers $C_{i}$ of
$X_{i}$ and by finitely many elements $g_{ij}$ such that $g_{ij}X_{i}$ crosses
$X_{j}$. We know that the $C_{i}$ are commensurable with $H$, and that
$C_{i}^{g_{ij}}$ is commensurable with $C_{j}$ since $g_{ij}X_{i}$ crosses
$X_{j}$. Thus $Stab(v)$ commensurises $H$ and so is a subgroup of
$Comm_{G}(H)$, which we denote by $C$ in this proof for brevity. As $H$ has
small commensuriser, $C=Comm_{G}(H)$ has $2$ ends. Hence $Stab(v)$ has $2$
ends. Let $e$ denote an edge of $T$ incident to $v$, so that $Stab(e)$ is a
subgroup of $Stab(v)$ and hence of $C$. If $Stab(e)$ were finite, the fact
that $T$ is minimal would imply that $G$ splits over this finite subgroup and
hence has more than one end. As $G$ is one-ended, it follows that $Stab(e)$ is
infinite and so has finite index in $Stab(v)$. In particular $Stab(e)$ is
two-ended. As $\Gamma$ is finite, it follows that there are only finitely many
edges of $T$ incident to $v$. Recall that part 3) of Corollary
\ref{nontrivialpartition} tells us that each almost invariant subset $X$
enclosed by $v$ is determined up to equivalence by the induced partition of
the edges incident to $v$. It follows that the elements of $E$ enclosed by $v$
belong to finitely many equivalence classes.

The fact that $\Gamma$ has finitely many $V_{0}$-vertices implies that the
$V_{0}$-vertices of $T$ lie in finitely many $G$-orbits. If a $V_{0}$-vertex
$w$ of $T$ has stabiliser which is a subgroup of $C$, then the stabiliser of
$gw$ satisfies the same condition if and only if $g$ commensurises $H$ and so
lies in $C$. As the $C$-orbit of $w$ is finite, it follows that there are only
finitely many $V_{0}$-vertices of $T$ whose stabiliser is contained in $C$. As
for $v$, each such $V_{0}$-vertex can enclose only finitely many equivalence
classes of elements of $E$.

We conclude from the above discussion that if $\Gamma$ is the regular
neighbourhood of a finite set of elements of $Q(H)$, then all the edge
splittings of $\Gamma$ are over two-ended subgroups commensurable with $H$,
and the $V_{0}$-vertices of $\Gamma$ can only enclose finitely many
equivalence classes of elements of $Q(H)$.

Now suppose that $B(H)$ is infinite. We will describe how to pick a sequence
$\{U_{i}\}_{i\geq1}$ of elements of $Q(H)$ which represent distinct elements
of $B(H)$. Having chosen $U_{1},\ldots,U_{k}$, we will form their regular
neighbourhood $\Gamma_{k}$. As the $V_{0}$-vertices of $\Gamma_{k}$ can only
enclose finitely many equivalence classes of elements of $Q(H)$, there is
$U_{k+1}$ in $Q(H)$ which is not enclosed by any $V_{0}$-vertex of $\Gamma
_{k}$. This implies that when we form $\Gamma_{k+1}$, it is distinct from
$\Gamma_{k}$. Lemma \ref{edgesplittingswhenyouenlargeE} implies that each edge
splitting of $\Gamma_{k}$ is enclosed by some $V_{0}$-vertex of $\Gamma_{k+l}%
$, for any $l\geq1$. In particular, the edge splittings of $\Gamma_{k}$ are
compatible with those of $\Gamma_{k+l}$, for any $l\geq1$. Consider all the
edge splittings of $\Gamma_{1},\ldots,\Gamma_{k}$ and choose one from each
conjugacy class. Let $\Delta_{k}$ denote the graph of groups structure for $G$
whose edge splittings are the chosen ones. Such a graph of groups exists by
Theorem \ref{Theorem2.5ofSS}. It is trivial that $\Delta_{k+1}$ is a
refinement of $\Delta_{k}$. Further, our construction implies that $\Delta
_{k}$ has no redundant vertices. As the edge groups of $\Delta_{k}$ are
two-ended, the accessibility result of Theorem \ref{graphsstabilise} applies
and tells us that the sequence $\{\Delta_{k}\}$ must eventually stabilise,
i.e. there is $N$ such that $\Delta_{N}=\Delta_{n}$, for all $n\geq N$. It
follows that the $V_{0}$-vertices of $\Gamma_{N}$ enclose every $U_{i}$. But
the $U_{i}$'s represent infinitely many distinct elements of $B(H)$ which
contradicts the preceding paragraph. This contradiction shows that $B(H)$ must
be finite as required.
\end{proof}

Now we know that $B(H)$ is finite, we can prove the following result.

\begin{proposition}
\label{regnbhdofeveryelementofB(H)whenHhassmallcommensuriser}Let $G$ be a
one-ended, finitely presented group, and let $H$ be a two-ended subgroup with
small commensuriser. Let $\Gamma$ denote the regular neighbourhood of the
collection $B(H)$. Thus $\Gamma$ is the graph of groups $\Gamma_{N}$ above.
Then one of the following cases holds:

\begin{enumerate}
\item All $V_{0}$-vertices of $\Gamma$ are isolated. In this case, $\Gamma$
has at most three $V_{0}$-vertices, and each has associated group which is
commensurable with $H$.

\item There is exactly one non-isolated $V_{0}$-vertex of $\Gamma$ with
associated group $Comm_{G}(H)$, and a non-zero number of isolated $V_{0}%
$-vertices. The non-isolated $V_{0}$-vertex encloses every element of $Q(H)$.
Further each isolated $V_{0}$-vertex is joined to the non-isolated $V_{0}%
$-vertex by a path of length $2$, such that the single $V_{1}$-vertex on this
path has valence $2$.
\end{enumerate}
\end{proposition}

\begin{proof}
Choose representatives $U_{1},\ldots,U_{k}$ of the elements of $B(H)$, so that
$U_{i}$ has stabiliser $H_{i}$ which is commensurable with $H$. We will use
our usual notation from section \ref{regnbhds:construction}. Thus $E$ denotes
the set of all translates of $U_{1},\ldots,U_{k}$ and their complements, and
$T$ denotes the universal covering $G$-tree of $\Gamma$.

Let $K$ denote the intersection of the $H_{i}$'s, so that each $U_{i}$ is
$K$-almost invariant, and consider the almost invariant subsets $K\backslash
U_{i}$ of $K\backslash G$. Note that the coboundary $\delta(K\backslash
U_{i})$ is finite for each $i$. Let $\Lambda$ denote the Cayley graph of $G$
with respect to some finite generating set, and let $A$ denote a finite
connected subcomplex of the graph $K\backslash\Lambda$ which contains every
edge of each of $\delta(K\backslash U_{i})$ and which carries $K$. Consider
the inverse image $Z$ of $A$ in $\Lambda$. Then $Z$ is connected and
$K$-finite. Let $X_{1},\ldots,X_{m}$ denote the $K$-infinite components of the
complement of $Z$. Since $\delta(K\backslash U_{i})$ is contained in the
interior of $A$, we see that $\delta U_{i}$ is contained in the interior of
$Z$ and thus each $U_{i}$ is a union of some of the $X_{j}$. Hence any
$K$-almost invariant subset of $G$ is equivalent to a union of some of the
$X_{j}$'s. In particular, if a nontrivial almost invariant subset over a
subgroup commensurable with $H$ is contained in one of the $X_{j}$'s then it
must be $H$-almost equal to $X_{j}$. It follows that $X_{j}$ cannot cross any
element of $B(H)$. Thus any element of $B(H)$ is represented by some union of
the $X_{j}$'s, and those elements of $E$ which are equivalent to some $X_{j}$
must be isolated in $E$.

If $m>3$, we will show that we have case 2) of the proposition. Recall that
the CCC's of $\overline{E}$ correspond to the $V_{0}$-vertices of $T$.
Consider all unions of $k$ of the $X_{j}$'s, for all $k$ such that $2\leq
k\leq(m-2)$. It is easy to check that all elements of $E$ equivalent to such
unions lie in the same CCC. For example, $X_{1}\cup X_{2}$ and $X_{1}\cup
X_{3}$ cross each other. Thus, if we consider only those CCC's which enclose
representatives of elements of $B(H)$, there is exactly one non-isolated such
CCC, say $v$, and there are exactly $m$ isolated such CCC's, say $v_{1}%
,\ldots,v_{m}$, where $v_{j}$ encloses $X_{j}$. Clearly, $v$ is invariant
under $Comm_{G}(H)$, and as its stabiliser must commensurise $H$ by the proof
of Theorem \ref{smallcommensuriserimpliesB(H)isfinite}, it follows that
$Stab(v)=Comm_{G}(H)$. Now part 5) of Lemma \ref{somefactsaboutenclosing}
tells us that if a vertex of a $G$-tree encloses two almost invariant subsets
of $G$, then it also encloses their intersection. Thus the fact that
$X_{1}=(X_{1}\cup X_{2})\cap(X_{1}\cup X_{3})$ implies that $v$ encloses
$X_{1}$. As Proposition \ref{Tisminimal} tells us that $T$ is a minimal
$G$-tree, we can apply part 8) of Lemma \ref{somefactsaboutenclosing} which
tells us that all the vertices on the interior of the path $\lambda_{1}$ in
$T$ joining $v_{1}$ to $v$ have valence $2$, and all these vertices enclose
$X_{1}$. Lemma \ref{isolatedvertices} now implies that each $V_{0}$-vertex on
the interior of $\lambda_{1}$ is isolated. Recall that no two distinct
elements of $E$ are equivalent. This implies that $X_{1}$ can be enclosed by
only one isolated $V_{0}$-vertex of $T$. It follows that $\lambda_{1}$ has no
interior $V_{0}$-vertices, so that $\lambda_{1}$ has length $2$. This
completes the proof that $\Gamma$ has all the properties in case 2) of the proposition.

If $m\leq3$, every element of $Q(H)$ is equivalent to some $X_{j}$ or its
complement. Thus we have $m$ isolated CCC's, say $v_{1},\ldots,v_{m}$, where
$v_{j}$ encloses $X_{j}$, and so $\Gamma$ has $m$ isolated $V_{0}$-vertices as
in case 1) of the proposition.
\end{proof}

The first part of the proof of the above proposition shows that each $X_{j}$
contains only one coend of $H$ in $G$. Thus we have the following proposition
which answers a question of Bowditch \cite{B2} in the finitely presented case:

\begin{proposition}
\label{finitelymanycoends} Suppose that $G$ is one-ended and finitely
presented and that $H$ is a two-ended subgroup of $G$ with small
commensuriser. Then the number of coends of $H$ in $G$ is finite.
\end{proposition}

Note that in either case in Proposition
\ref{regnbhdofeveryelementofB(H)whenHhassmallcommensuriser}, there are
isolated elements of $E$. Such elements determine splittings of $G$ which have
intersection number zero with every element of $Q(H)$. We will call such a
splitting $H$\textit{-canonical}. Thus we have the following result.

\begin{corollary}
\label{smallcommimpliesGhasH-canonicalsplitting}Suppose that $G$ is a
one-ended, finitely presented group and assume that $G$ is not of
finite-by-Fuchsian type. Let $H$ be a two-ended subgroup of $G$ and suppose
that $H$ has small commensuriser in $G$. If $G\;$possesses a nontrivial almost
invariant subset over a subgroup commensurable with $H$, then $G$ possesses a
$H$-canonical splitting over a subgroup commensurable with $H$.
\end{corollary}

\section{Coends when the commensuriser is
large\label{coendswhencommensuriserislarge}}

Recall that $Q(H)$ denotes the collection of all almost invariant subsets of
$G$ which are over subgroups of $G$ commensurable with $H$, and $B(H)$ denotes
the collection of equivalence classes of elements of $Q(H)$. In Proposition
\ref{regnbhdofeveryelementofB(H)whenHhassmallcommensuriser}, we described the
structure of the regular neighbourhood $\Gamma$ of all elements of $B(H)$, in
the case when $H\;$has small commensuriser in $G$. In this section, we
consider a one-ended finitely generated group $G$, and a two-ended subgroup
$H$ such that $H$ has large commensuriser in $G$, i.e. $H$ has infinite index
in $Comm_{G}(H)$. Again we want to study the structure of the regular
neighbourhood $\Gamma$ of all elements of $B(H)$. However $B(H)$ may be
infinite, so that it is not clear that such a regular neighbourhood exists. We
will show that it does exist in section \ref{JSJforlargecommensurisers}. A key
point in the argument is that $B(H)$ possesses certain algebraic finiteness
properties. Recall from the previous section that $Q(H)$ and $B(H)$ are
Boolean algebras. Also $Q(H)$ is invariant under the action by left
multiplication of $Comm_{G}(H)$, and this action induces an action on $B(H)$.
Thus $B(H)$ is a Boolean algebra with a natural action of $Comm_{G}(H)$.

We will use some arguments from \cite{D-Roller}. As before, the results extend
to the case where $H$ is virtually polycyclic (VPC) of any length but we will
first discuss the case where $H$ is VPC\ of length $1$, i.e. $H$ is two-ended.
We start by recalling a result of Kropholler and Roller \cite{KR} (see also
\cite{RG}), and we sketch an argument for this which follows the proof of
Proposition 2.8 of \cite{SS}.

\begin{proposition}
\label{Comm_G(H)largeimpliescoendsare12orinfty}Let $G$ be a finitely generated
group and $H$ a two-ended subgroup with large commensuriser. Then the number
of coends of $H$ in $G$ is $1$, $2$ or infinity. The number of coends is $2$
if and only if $G$ is virtually a torus group.
\end{proposition}

\begin{proof}
If the number of coends of $H$ in $G$ is greater than $1$, there is a
nontrivial $H$-almost invariant subset $X$ of $G$. Consider the translates of
$X$ by elements of $Comm_{G}(H)$. If infinitely many of these are equivalent
to $X$, then the proof of Lemma 2.13 in \cite{SS} shows that there is a
subgroup $K_{1}$ of finite index in $G$ which contains a subgroup $H_{1}$
commensurable with $H$ such that $H_{1}$ is normal in $K_{1}$ and
$H_{1}\backslash K_{1}$ has two ends. Thus $H$ has two coends in $G$, and as
$H$ is virtually infinite cyclic, $G$ is virtually a torus group. If only
finitely many translates of $X$ by $Comm_{G}(H)$ are equivalent to $X$, then
$X$ has infinitely many distinct such translates and it follows that $H$ has
infinitely many coends in $G$.
\end{proof}

Now we can prove the following finiteness result for the Boolean algebra
$B(H)$.

\begin{theorem}
\label{Booleanalgebrafinitelygenerated} If $G$ is a one-ended, finitely
presented group and $H$ is a two-ended subgroup of $G$, then $B(H)$ is
finitely generated over $Comm_{G}(H)$.
\end{theorem}

\begin{proof}
Theorem \ref{smallcommensuriserimpliesB(H)isfinite} tells us that if $H$ has
small commensuriser, then $B(H)$ is finite. Thus the result is trivial in this
case. So we will assume that $H$ has infinite index in $Comm_{G}(H)$, and that
$B(H)$ is infinite. Note that all the crossings between elements of $B(H)$
must be weak, by Corollary \ref{strongcrossingimpliessmallcomensuriser}. The
accessibility result of Theorem \ref{graphsstabilise} tells us that there is a
finite graph of groups decomposition $\mathcal{G}$ of $G$ with all edge groups
commensurable with $H$, such that $\mathcal{G}$ cannot be properly refined
using such splittings. An alternative way of expressing this condition is to
say that if $G$ possesses a splitting over a two-ended subgroup commensurable
with $H$ which has intersection number zero with the edge splittings of
$\mathcal{G}$, then this splitting is conjugate to one of these edge
splittings. Let $X_{1},\ldots,X_{n}$ denote almost invariant subsets of $G$
associated to the edge splittings of $\mathcal{G}$, and let $H_{i}$ denote the
stabiliser of $X_{i}$. Let $A(H)$ denote the subalgebra of $B(H)$ generated
over $Comm_{G}(H)$ by the equivalence classes of the $X_{i}$'s. We will show
that $B(H)=A(H)$.

Let $E$ denote the collection of translates of the $X_{i}$'s by elements of
$Comm_{G}(H)$. Let $Y$ be an element of $Q(H)$, so that $Y$ is an almost
invariant subset of $G$ over a subgroup $K$ commensurable with $H$. We will
show that $Y$ crosses only finitely many elements of $E$. The intersection
number $i(H_{i}\backslash X{_{i}},K\backslash Y)$ is a finite number which is
the number of double cosets $KgH_{i}$ such that $gX_{i}$ crosses $Y$. If
$gX_{i}$ crosses $Y$, it must do so weakly by Corollary
\ref{strongcrossingimpliessmallcomensuriser}. Now Proposition
\ref{commensurable} tells us that the stabilisers of $gX_{i}$ and $Y$ are
commensurable, so that $H_{i}^{g}$ and $K$ are commensurable. As $H_{i}$ and
$K$ are commensurable, it follows that $g$ commensurises $H_{i}$, and hence
also commensurises $K$. If we let $L_{i}$ denote $K\cap H_{i}$, then
$g^{-1}Kg$ can be expressed as the union of cosets $g_{j}(g^{-1}Kg\cap L_{i}%
)$, for $1\leq j\leq n.$ Hence
\[
KgH_{i}=g(g^{-1}Kg)H_{i}=g\left(  \cup_{j=1}^{n}g_{j}(g^{-1}Kg\cap
L_{i})\right)  H_{i}=g\left(  \cup_{j=1}^{n}g_{j}H\right)  =\cup_{j=1}%
^{n}gg_{j}H,
\]
so that $KgH_{i}$ is the union of finitely many cosets $gH_{i}$ of $H_{i}$. It
follows that there are only finitely many translates of the $X_{i}$'s which
cross $Y$. We conclude that $Y$ crosses only finitely many of the translates
of the $X_{i}$'s by $G$ and that these must all lie in $E$.

If $Y$ crosses no elements of $E$, then $Y$ is enclosed by some vertex $v$ of
the universal covering $G$-tree $T$ of $\mathcal{G}$. In particular, the
stabiliser $K$ of $Y$ is a subgroup of $Stab(v)$. Suppose that $Y\cap Stab(v)$
is a nontrivial $K$-almost invariant subset of $Stab(v)$. Then the proof of
the main theorem of \cite{D-Roller} produces another nontrivial almost
invariant subset $U$ of $G$ such that $U$ is over a subgroup of $Stab(v)$
which is commensurable with $K$, and $U\cap Stab(v)$ is also a nontrivial
almost invariant subset of $Stab(v)$, and $U$ has self-intersection number
zero. As $U$ has self-intersection number zero, Theorem 2.8 of \cite{SS}
implies that $U$ is equivalent to an almost invariant set $W$ which is
associated to a splitting $\sigma$ of $G$. The set $U$ was constructed by
taking successively intersections of $Y^{(\ast)}$ and $cY^{(\ast)}$ where
$c\in Comm_{G}(H)\cap Stab(v)$, so that part 5) of Lemma
\ref{somefactsaboutenclosing} shows that it is enclosed by $v$. Hence $W$ is
also enclosed by $v$. Thus $\sigma$ has intersection number zero with the edge
splittings of $\mathcal{G}$. Our choice of $\mathcal{G}$ now implies that
$\sigma$ must be conjugate to one of the edge splittings of $\mathcal{G}$. In
particular, $\sigma$ does not split $Stab(v)$, so that $U\cap Stab(v)$ must be
a trivial almost invariant subset of $Stab(v)$, which is a contradiction. This
contradiction shows that if $Y$ crosses no elements of $E$, then $Y\cap
Stab(v)$ must be a trivial almost invariant subset of $Stab(v)$.

Now recall from Lemma \ref{enclosedsetscanbechosennested} that $Y$ is
equivalent to the set $B(Y)=\varphi^{-1}(\Sigma_{v}(Y))\cup(Y\cap\varphi
^{-1}(v))$, where $\varphi:G\rightarrow V(T)$ is some $G$-equivariant map. We
will choose $\varphi(e)=v$, so that $\varphi^{-1}(v)=Stab(v)$. Thus
$Y\cap\varphi^{-1}(v)$ is a trivial almost invariant subset of $Stab(v)$. It
follows that $Y$ is equivalent to $\varphi^{-1}(\Sigma_{v}(Y))$. This set is
the union of some of the $Z_{s}^{\ast}$, for edges $s$ incident to $v$ and
oriented towards $v$. As the $Z_{s}^{\ast}$ and their coboundaries are
disjoint, the coboundary of this union equals the union of the $\delta
Z_{s}^{\ast}$. As $\delta Y$ is $H$-finite, it follows that $Y$ is equivalent
to the union of a finite number of $Z_{s}^{\ast}$. Thus $Y$ lies in the
subalgebra $A(H)$ of $B(H)$.

Next suppose that $Y$ crosses some elements of $E$ and let $v$ be a vertex of
$T$ some of whose incident edges have associated edge splittings crossed by
$Y$. If $s_{1},...s_{k}$ are these edges, we denote the almost invariant
subset of $G$ associated to $s_{i}$ by $Z_{i}$ for obvious typographical
reasons. Then, we can express $Y$ as the union of $Y\cap Z_{i}^{\ast}$, $1\leq
i\leq k$, and of $W=Y\cap\left(  \bigcap_{i=1}^{k}Z_{i}\right)  $. Now part 5)
of Lemma \ref{somefactsaboutenclosing} shows that $W$ is enclosed by $v$, so
that $W$ lies in the subalgebra $A(H)$ of $B(H)$, by the preceding paragraph.
Moreover, each of $Y\cap Z_{i}^{\ast}$ crosses a smaller number of the $X_{i}%
$'s than does $Y$ and so by induction we see that each $Y\cap Z_{i}^{\ast}$
also lies in the subalgebra $A(H)$ of $B(H)$. This now implies that $Y$ itself
lies in $A(H)$ and completes the proof that $B(H)\;$is finitely generated over
$Comm_{G}(H)$.
\end{proof}

In the preceding proof, we referred several times to the arguments of Dunwoody
and Roller in \cite{D-Roller}. It will also be convenient to state one of the
main results of that paper.

\begin{theorem}
\label{D-Rollertheorem}Let $G$ be a one-ended, finitely generated group which
does not split over VPC subgroups of length $<n$, and let $H$ be a VPC
subgroup of length $n$ with large commensuriser, such that $e(G,H)\geq2$. Then
$G$ splits over some subgroup commensurable with $H$.
\end{theorem}

Note that the proof of Theorem \ref{Booleanalgebrafinitelygenerated} shows
that $B(H)$ is generated by almost invariant sets associated to splittings of
$G$ over subgroups commensurable with $H$, so that Theorem
\ref{D-Rollertheorem} follows. However, this is not a new proof of Theorem
\ref{D-Rollertheorem} as we used Theorem \ref{D-Rollertheorem} in our proof of
Theorem \ref{Booleanalgebrafinitelygenerated}.

We call a nontrivial element $X$ of $Q(H)$ \textit{special} if $X\cap
Comm_{G}(H)$ is $H$-finite. A splitting of $G$ over a subgroup commensurable
with $H$ will be called \textit{special} if one of the associated almost
invariant subsets of $G$ is a special element of $Q(H)$. Recall that a
splitting of $G$ over a subgroup commensurable with $H$ is $H$-canonical if it
has zero intersection number with every element of $Q(H)$.

\begin{proposition}
Let $G$ be a one-ended, finitely presented group which is not virtually
$\mathbb{Z}\times\mathbb{Z}$, with a two-ended subgroup $H$ which has large
commensuriser. If $Q(H)$ has a nontrivial special element, then $G$ has an
$H$-canonical special splitting over a subgroup commensurable with $H$.
\end{proposition}

\begin{remark}
If $H$ has small commensuriser, then every element of $Q(H)$ is special, so
the result follows immediately from Corollary
\ref{smallcommimpliesGhasH-canonicalsplitting}.
\end{remark}

\begin{proof}
The proof uses details from the arguments of \cite{D-Roller} with slight
improvements from \cite{D-Swenson}. Note that as $H$ has large commensuriser
and $G$ is not virtually $\mathbb{Z}\times\mathbb{Z}$, Corollary
\ref{strongcrossingimpliessmallcomensuriser} tells us that an element of
$Q(H)$ cannot cross strongly any almost invariant set over a two-ended subgroup.

The proof in \cite{D-Roller} shows that if $X$ is a special element of $Q(H)$,
then there is a splitting of $G$ over a subgroup commensurable with $H$ such
that one of the almost invariant sets associated with the splitting is
contained in $X$ (and thus special). We recall some of the argument. By
modifying $X$ by a $H$-finite set, we can arrange that $X$ does not intersect
$Comm_{G}(H)$. In this case, Lemma 5 of \cite{D-Roller} tells us that if $c\in
Comm_{G}(H)$ and $X$ and $cX$ are nested, then either $cX=X$ or $X\cap
cX=\emptyset$. Now \cite{D-Roller} produces a new element $Y$ of $Q(H)$ which
is a finite intersection of some $c_{i}X$, $c_{i}\in Comm_{G}(H)$, such that
the translates of $Y$ by all elements of $G$ are almost nested. This is enough
to produce a splitting of $G$ with associated almost invariant set equivalent
to $Y$. There is an improvement in \cite{D-Swenson} (see section 3, last
paragraph on page 622), which shows there is a subset $Z$ of $Y$ which is
equivalent to $Y$ such that the translates of $Z$ are actually nested. This
implies that $Z$ is associated to a splitting of $G$ over a subgroup
commensurable with $H$. Thus, if we start with an element of $Q(H)$ which is
disjoint from $Comm_{G}(H)$, there is a subset which is also an element of
$Q(H)$ and is associated to a splitting of $G$. Of course, this new element of
$Q(H)$ is also disjoint from $Comm_{G}(H)$.

Let $X_{1}$ denote an element of $Q(H)$ which is disjoint from $Comm_{G}(H)$,
and is associated to a splitting of $G$. If $X_{1}$ contains an element
$Y_{1}$ of $Q(H)$ which is not equivalent to $X_{1}$, the preceding paragraph
yields a subset $X_{2}$ of $Y_{1}$, which is also an element of $Q(H)$ and is
associated to a splitting of $G$. We repeat this process to obtain a
descending sequence of inequivalent elements $X_{i}$ of $Q(H)$, each of which
is disjoint from $Comm_{G}(H)$, and is associated to a splitting $\sigma_{i}$
of $G$.

Suppose that this sequence stops after finitely many steps. Then we will have
found an element $X$ of $Q(H)$ which is disjoint from $Comm_{G}(H)$ and is
associated to a splitting of $G$, such that any subset of $X$ which lies in
$Q(H)$ is equivalent to $X$. Thus $X$ is a minimal element in $B(H)$. Any such
minimal element cannot be crossed by any element in $Q(H)$, so it determines a
$H$-canonical splitting of $G$, as required.

If the sequence does not stop after finitely many steps, we will obtain an
infinite sequence $X_{i}$. We will show that this leads to a contradiction.
Let $H_{i}$ denote the stabiliser of $X_{i}$. We claim $H_{2}\subset H_{1}$.
Let $c$ be an element of $H_{2}$. As $cX_{2}=X_{2}$, we see that $cX_{1}\cap
X_{1}$ is non-empty. As $X_{1}$ is associated to a splitting, all its
translates are nested. In particular, as $c$ lies in $Comm_{G}(H_{1})$, we
have $cX_{1}=X_{1}$ or $cX_{1}\cap X_{1}=\emptyset$. As $cX_{1}\cap X_{1}$ is
non-empty, we must have $cX_{1}=X_{1}$, so that $c$ lies in $H_{1}$. It
follows that $H_{2}\subset H_{1}$. Next we claim that $X_{1}$ does not cross
any translate of $X_{2}$, so that the splittings $\sigma_{1}$ and $\sigma_{2}$
have intersection number zero. For suppose that $X_{1}$ and $gX_{2}$ cross
each other. We pointed out at the start of this proof that they must cross
weakly because $X_{1}$ lies in $Q(H)$, so that Proposition \ref{commensurable}
tells us that their stabilisers must be commensurable. As $X_{1}$ and $X_{2}$
have commensurable stabiliser, it follows that $g$ lies in $Comm_{G}(H)$. As
before this means that $gX_{1}=X_{1}$ or $gX_{1}\cap X_{1}=\emptyset$. But the
fact that $X_{1}$ and $gX_{2}$ cross implies that $X_{1}\cap gX_{2}$ is not
empty, so that $X_{1}\cap gX_{1}$ is not empty. It follows that $gX_{1}=X_{1}%
$, so that $gX_{2}\subset X_{1}$, which contradicts our assumption that
$X_{1}$ and $gX_{2}$ cross. Similarly, $H_{i+1}\subset H_{i}$, for all
$i\geq1$, and $X_{i}$ does not cross any translate of $X_{j}$, for any $i$ and
$j$. Thus the splittings $\sigma_{i}$ are all compatible. Hence we obtain an
infinite sequence of graphs of groups decompositions of $G$ each refining the
previous one. As the $\sigma_{i}$'s are distinct, this contradicts the
accessibility result in Theorem \ref{graphsstabilise}. This contradiction
shows that the sequence of $X_{i}$'s cannot be infinite, so that we obtain a
minimal $X$ as in the preceding paragraph, and hence obtain a $H$-canonical
splitting of $G$.
\end{proof}

The argument above shows that even though the number of coends of $H$ in $G$
is infinite, if $Q(H)$ has a special element, then there are elements of
$Q(H)$ which contain a finite number of coends of $H$ with the complementary
set containing all the other coends. Of course, there may not be any special
elements in $Q(H)$.

Now we will consider the following construction. Let $G$ be a one-ended,
finitely presented group, pick one representative for each element of $B(H)$
and let $E$ denote the set of all translates of this collection and their
complements. We will consider the CCC's of $\overline{E}$. As usual, these
form a pretree, but we will not be able to prove that this is discrete until
section \ref{JSJforlargecommensurisers}.

\begin{proposition}
\label{largecommensuriserimpliesoneinfiniteCCC}Let $G$ be a one-ended,
finitely presented group and suppose that $G$ is not virtually $\mathbb{Z}%
\times\mathbb{Z}$. Let $H$ be a two-ended subgroup of $G$ such that $H$ has
large commensuriser and there are nontrivial $H$-almost invariant subsets of
$G$. Consider the CCC's of $\overline{E}$ which consist of elements of $Q(H)$.
There is exactly one such CCC which is infinite, and there are finitely many
(possibly zero) other CCC's all of which are isolated. The infinite CCC has
stabiliser equal to $Comm_{G}(H)$, and encloses every element of $Q(H)$.
\end{proposition}

\begin{remark}
In the case when $G$ is virtually $\mathbb{Z}\times\mathbb{Z}$, the CCC's of
$\overline{E}$ which consist of elements of $Q(H)$ form a single CCC. This is
discussed at the start of the proof of Theorem
\ref{JSJforlargecommensurisersexist}.
\end{remark}

\begin{proof}
Theorem \ref{Booleanalgebrafinitelygenerated} shows that there are a finite
number of compatible splittings $\sigma_{1},\ldots\sigma_{n}$ of $G$ each over
a subgroup of $G$ commensurable with $H$, such that the equivalence classes of
the associated almost invariant subsets $\{X_{i},X_{i}^{\ast}:1\leq i\leq n\}$
generate $B(H)$ over $Comm_{G}(H)$. In fact the proof shows more. It shows
that any element of $B(H)$ is represented by some finite union of translates
by $Comm_{G}(H)$ of the $X_{i}$'s. As $G$ is not virtually $\mathbb{Z}%
\times\mathbb{Z}$, Proposition \ref{Comm_G(H)largeimpliescoendsare12orinfty}
tells us that $H$ has infinitely many coends in $G$. As the translates of
$X_{i}$ and $X_{i}^{\ast}$ are nested, the proof of part 2) of Proposition
\ref{regnbhdofeveryelementofB(H)whenHhassmallcommensuriser} shows that the
collection of all finite unions of more than one of these sets forms a single
CCC. It further shows that this CCC also encloses the $X_{i}$'s themselves. It
follows that this CCC has stabiliser equal to $Comm_{G}(H)$. Finally, there
are isolated CCC's corresponding to those $X_{i}$'s which are isolated (if any).
\end{proof}

\section{Canonical decompositions over two-ended groups when commensurisers
are small\label{JSJforsmallcommensurisers}}

In this section and the next, we will find canonical decompositions of a
one-ended, finitely presented group $G$ which are analogous to the topological
JSJ-decomposition of an atoroidal $3$-manifold. Our approach is similar to an
unpublished approach of Scott \cite{Unpublished} for proving the classical
results on the JSJ-decomposition. The idea is to enclose all almost invariant
subsets of $G$ which are over two-ended subgroups. However, instead of
enclosing all at once, we form regular neighbourhoods of larger and larger
finite families and use accessibility to show that the sequence obtained in
this way must stabilise. This is where we use finite presentability. The
result is a regular neighbourhood of all equivalence classes of almost
invariant subsets of $G$ over two-ended subgroups. (See Definition
\ref{defnofregnbhdofequivalenceclasses}.) Since we enclose all almost
invariant sets over two-ended groups, the decompositions that we obtain are
unique and are invariant under automorphisms of the group.

The preceding two sections contain the crucial pattern for obtaining our
canonical decompositions. We want to form a regular neighbourhood of an
infinite family of almost invariant subsets of $G$. The first step is to show
that the cross-connected components are of two types, those which contain only
strong crossings and those which contain only weak crossings. The structure of
the strong crossing components is handled by the techniques of Bowditch
\cite{B2} and of Dunwoody and Swenson \cite{D-Swenson}. If commensurisers are
small, the structure of weak crossing components is easy to describe using
regular neighbourhoods as in section \ref{coendswhencommensuriserissmall}. If
the commensuriser of a subgroup $H$ of $G$ is large, we use the fact that the
Boolean algebra $B(H)$ is finitely generated over $Comm_{G}(H)$, as in section
\ref{coendswhencommensuriserislarge}. This used the arguments of Dunwoody and
Roller \cite{D-Roller} for a special case of the annulus theorem. To obtain
our canonical decompositions, the only remaining difficulty is to show that
the pretree which we construct from the cross-connected components is
discrete. This is clear in the case when all commensurisers are small, and is
proved in general using again the fact that $B(H)$ is finitely generated over
$Comm_{G}(H)$. This is what we will do in this and the next section. We will
use the same strategy in the more general cases which we consider later.

Let $\Gamma$ be a minimal graph of groups decomposition of $G$. Recall that a
vertex $v$ of $\Gamma$ is of finite-by-Fuchsian type if $G(v)$ is a
finite-by-Fuchsian group, where the Fuchsian group is not finite nor
two-ended, and there is exactly one edge of $\Gamma$ which is incident to $v$
for each peripheral subgroup $K$ of $G(v)$ and this edge carries $K$. We will
also need the following definition.

\begin{definition}
\label{defnofsimplefortwo-endedsubgroups}Let $\Gamma$ be a minimal graph of
groups decomposition of a group $G$. A vertex $v$ of $\Gamma$ is
\textsl{simple}, if whenever $X$ is a nontrivial almost invariant subset of
$G$ over a two-ended subgroup such that $X$ is enclosed by $v$, then $X$ is
associated to an edge splitting of $\Gamma$.
\end{definition}

Note that if $v$ is simple, then part 8) of Lemma
\ref{somefactsaboutenclosing} implies that $X$ is associated to an edge
splitting for an edge of $\Gamma$ which is incident to $v$. Also note that it
is not possible for $v$ to be both simple and of finite-by-Fuchsian type. For
if a Fuchsian group is not finite nor two-ended, the corresponding
$2$-orbifold has non-peripheral loops representing elements of infinite order,
and any such yields a nontrivial almost invariant subset of $G$ over a
two-ended subgroup which is not associated to an edge splitting of $\Gamma$.

The above definition is the analogue of the topological fact that if $V(M)$
denotes the characteristic submanifold of a $3$-manifold $M$, then any
essential annulus in $M$ which lies in $M-V(M)$ is homotopic to a covering of
an annulus in the characteristic family $\mathcal{T}$. We could have defined
simple differently by insisting that the condition applies only to almost
invariant subsets which are associated to a splitting. Call this condition
``simple for splittings''. This would correspond to a topological definition
of simple which applies only to embedded annuli. Unfortunately, the word
simple is used in both these senses in the literature of $3$-manifold theory.
Note that a vertex $v$ of a graph of groups $\Gamma$ can be of
finite-by-Fuchsian type and be simple for splittings. For example, this occurs
when $G(v)$ is the fundamental group of a thrice punctured sphere.

If one applies Definition \ref{defnofsimplefortwo-endedsubgroups} to
splittings of $G$, one sees that if $v$ is simple, then $\Gamma$ cannot be
properly refined by splitting at $v$ using a splitting of $G$ over a two-ended subgroup.

Finally, note that $v$ need not be simple even if $G(v)$ is two-ended, despite
the fact that such a group does not admit any splitting over a two-ended
subgroup. Again this is clear from consideration of topological examples.
Consider a solid torus component $\Sigma$ of $V(M)$ such that $\Sigma$ has at
least four frontier annuli. Then $\Sigma$ will contain essential embedded
annuli which are not homotopic to covers of the frontier annuli. This example
shows that $v$ need not even be simple for splittings.

In this section, we will assume that whenever $H$ is a two-ended subgroup of
$G$ and $e(G,H)\geq2$, then $H$ has small commensuriser. The class of groups
for which this holds includes all word hyperbolic groups. Now we can state the
main result of this section which is our version of the JSJ-decomposition for
this class of groups.

\begin{theorem}
\label{JSJforsmallcommensurisersexists}Let $G$ be a one-ended, finitely
presented group such that whenever $H$ is a two-ended subgroup and
$e(G,H)\geq2$, then $H$ has small commensuriser. Let $\mathcal{F}$ denote the
collection of equivalence classes of all nontrivial almost invariant subsets
of $G$ which are over a two-ended subgroup.

Then the regular neighbourhood construction of section
\ref{regnbhds:construction} works and yields a regular neighbourhood
$\Gamma(\mathcal{F}:G)$. Further each $V_{0}$-vertex $v$ of $\Gamma
(\mathcal{F}:G)$ satisfies one of the following conditions:

\begin{enumerate}
\item $v$ is isolated.

\item $G(v)$ is the full commensuriser $Comm_{G}(H)$ for some two-ended
subgroup $H$, such that $e(G,H)\geq2$.

\item $v$ is of finite-by-Fuchsian type.
\end{enumerate}

$\Gamma(\mathcal{F}:G)$ consists of a single vertex if and only if either
$\mathcal{F}$ is empty or $G$ is of finite-by-Fuchsian type.
\end{theorem}

\begin{remark}
In cases 1) and 2), $G(v)$ is two-ended. If $v$ is of type 2) but not
isolated, we will say that $v$ is of small commensuriser type.
\end{remark}

\begin{proof}
We start by picking a representative for each element of $\mathcal{F}$,
subject to the condition that if $A$ and $B$ are elements of $\mathcal{F}$
such that $B=gA$, for some $g$ in $G$, then the representatives $X$ and $Y$
chosen for $A$ and $B$ must satisfy $Y=gX$. As $\mathcal{F}$ is countable, it
can be expressed as the union of an ascending sequence of finite subsets
$\mathcal{F}_{i}$, for $i\geq1$. As each $\mathcal{F}_{i}$ is finite, it has a
regular neighbourhood $\Gamma_{i}=\Gamma(\mathcal{F}_{i}:G)$. We will choose
the sequence $\mathcal{F}_{i}$ carefully, and then show that the sequence
$\Gamma_{i}$ must stabilise eventually. The resulting graph of groups
structure $\Gamma$ for $G$ will be the required regular neighbourhood
$\Gamma=\Gamma(\mathcal{F}:G)$.

We will start with some choice of the $\mathcal{F}_{i}$'s and then will modify
this sequence inductively. We will continue to denote the modified sets by
$\mathcal{F}_{i}$. We modify $\mathcal{F}_{1}$ as follows. Let $E_{1}$ denote
the collection of all the translates of the chosen representatives of
$\mathcal{F}_{1}$. Consider a $H$-almost invariant subset $X$ of $G$ which
lies in $E_{1}$. If $X$ crosses weakly some element of $\mathcal{F}$, Lemma
\ref{commensurable} tells us that this element lies in $B(H)$. Theorem
\ref{smallcommensuriserimpliesB(H)isfinite} tells us that $B(H)$ is finite.
Thus, by enlarging $\mathcal{F}_{1}$, we can arrange that whenever $E_{1}$
contains a $H$-almost invariant subset of $G$ which crosses weakly some
element of $\mathcal{F}$, then $\mathcal{F}_{1}$ contains each element of
$B(H)$. This is the final version of $\mathcal{F}_{1}$.

The sequence $\mathcal{F}_{i}$ is constructed inductively starting from
$\mathcal{F}_{1}$. Having constructed $\mathcal{F}_{i}$, we first check
whether $\mathcal{F}_{i}=\mathcal{F}$. If it does, then $\Gamma(\mathcal{F}%
_{i}:G)$ is the required regular neighbourhood $\Gamma(\mathcal{F}:G)$.
Otherwise, we need to construct $\mathcal{F}_{i+1}$. As the new $\mathcal{F}%
_{i}$ is still finite, there is some index $j>i$ such that $\mathcal{F}_{j}$
properly contains the new $\mathcal{F}_{i}$. By replacing our sequence of
subsets of $\mathcal{F}$ by a subsequence, we can suppose that $j=i+1$. Now we
enlarge $\mathcal{F}_{i+1}$ in the same way in which we enlarged
$\mathcal{F}_{1}$.

In order to show that the sequence $\Gamma_{i}$ stabilises, we first need to
describe the $V_{0}$-vertices of $\Gamma_{i}$. Recall that $\Gamma_{i}$ can be
constructed using the collection $E_{i}$ of all the translates of the chosen
representatives of $\mathcal{F}_{i}$. Proposition
\ref{crossingsareallstrongorallweak} shows that the non-isolated CCC's of
$\overline{E_{i}}$ are of two types. In a given CCC, the crossings are either
all strong or all weak. In the strong crossing case, Theorem
\ref{strongcrossingimpliesFuchsiantype} shows that the corresponding $V_{0}%
$-vertex of $\Gamma_{i}$ is of finite-by-Fuchsian type. Suppose that $E_{i}$
has an element $X$ which is $H$-almost invariant and crosses some other
element of $E_{i}$ weakly. Then Proposition
\ref{crossingsareallstrongorallweak} shows that the CCC of $\overline{E_{i}}$
which contains $\overline{X}$ contains only elements whose stabiliser is
commensurable with $H$. Now the proof of Proposition
\ref{regnbhdofeveryelementofB(H)whenHhassmallcommensuriser} tells us that the
CCC of $\overline{E_{i}}$ which contains $\overline{X}$ has stabiliser equal
to $Comm_{G}(H)$, and encloses every element of $E_{i}$ whose stabiliser is
commensurable with $H$. As $Comm_{G}(H)$ contains $H$ with finite index, it is
two-ended. It follows that every $V_{0}$-vertex group of $\Gamma_{i}$ is
either two-ended or of finite-by-Fuchsian type, and hence that every edge
group of $\Gamma_{i}$ is two-ended. Lemma \ref{edgesplittingswhenyouenlargeE}
shows that the edge splittings of $\Gamma_{i}$ are compatible with those of
$\Gamma_{j}$ for every $j>i$. Now we proceed as at the end of the proof of
Theorem \ref{smallcommensuriserimpliesB(H)isfinite}. Consider all the edge
splittings of $\Gamma_{1},\ldots,\Gamma_{k}$ and choose one from each
conjugacy class. Let $\Delta_{k}$ denote the graph of groups structure for $G$
whose edge splittings are the chosen ones. Such a graph of groups exists by
Theorem \ref{Theorem2.5ofSS}. It is trivial that $\Delta_{k+1}$ is a
refinement of $\Delta_{k}$. Further, our construction implies that $\Delta
_{k}$ has no redundant vertices. As the edge groups of $\Delta_{k}$ are
two-ended, the accessibility result of Theorem \ref{graphsstabilise} applies
and tells us that the sequence $\Delta_{k}$ must eventually stabilise, i.e.
there is $N$ such that $\Delta_{N}=\Delta_{n}$, for all $n\geq N$. It follows
that the sequence $\Gamma_{i}$ must stabilise, as required. We call this final
graph of groups $\Gamma$.

By construction, $\Gamma$ is the required regular neighbourhood $\Gamma
(\mathcal{F}:G)$. Each $V_{0}$-vertex group of $\Gamma$ is isolated, of small
commensuriser type or of finite-by-Fuchsian type because this holds for each
$\Gamma_{i}$. Finally, $\Gamma$ will consist of a single vertex if and only if
either $\mathcal{F}$ is empty or the representatives of $\mathcal{F}$ lie in a
single CCC. In the second case, the vertex group of $\Gamma$ must be two-ended
or of finite-by-Fuchsian type. As the vertex group is $G$ and we assumed that
$G$ is one-ended, it follows that $\Gamma(\mathcal{F}:G)$ consists of a single
vertex if and only if either $\mathcal{F}$ is empty or $G$ is of
finite-by-Fuchsian type.
\end{proof}

In order to understand $\Gamma(\mathcal{F}:G)$ in more detail, our next result
lists some properties which follow almost immediately from the above theorem
and the properties of an algebraic regular neighbourhood.

\begin{theorem}
\label{propertiesofJSJforsmallcommensurisers}Let $G$ be a one-ended, finitely
presented group such that whenever $H$ is a two-ended subgroup and
$e(G,H)\geq2$, then $H$ has small commensuriser. Let $\mathcal{F}$ denote the
collection of equivalence classes of all nontrivial almost invariant subsets
of $G$ which are over a two-ended subgroup.

Then the regular neighbourhood $\Gamma(\mathcal{F}:G)$ is a minimal bipartite
graph of groups decomposition of $G$ with the following properties:

\begin{enumerate}
\item Each $V_{0}$-vertex $v$ of $\Gamma(\mathcal{F}:G)$ satisfies one of the
following conditions:

\begin{enumerate}
\item $v$ is isolated.

\item $G(v)$ is the full commensuriser $Comm_{G}(H)$ for some two-ended
subgroup $H$, such that $e(G,H)\geq2$.

\item $v$ is of finite-by-Fuchsian type.
\end{enumerate}

\item the edge groups of $\Gamma(\mathcal{F}:G)$ are two-ended.

\item any element of $\mathcal{F}$ is enclosed by some $V_{0}$-vertex of
$\Gamma(\mathcal{F}:G)$, and each $V_{0}$-vertex of $\Gamma(\mathcal{F}:G)$
encloses such a subset of $G$. In particular, any splitting of $G$ over a
two-ended subgroup is enclosed by some $V_{0}$-vertex of $\Gamma
(\mathcal{F}:G)$.

\item if $X$ is an almost invariant subset of $G$ over a finitely generated
subgroup $H$, and if $X$ does not cross any element of $\mathcal{F}$, then $X$
is enclosed by a $V_{1}$-vertex of $\Gamma(\mathcal{F}:G)$.

\item if $X$ is a $H$-almost invariant subset of $G$ associated to a splitting
of $G$ over $H$, and if $X$ does not cross any element of $\mathcal{F}$, then
$X$ is enclosed by a $V_{1}$-vertex of $\Gamma(\mathcal{F}:G)$.

\item the $V_{1}$-vertex groups of $\Gamma(\mathcal{F}:G)$ are simple.

\item If $\Gamma_{1}$ and $\Gamma_{2}$ are minimal bipartite graphs of groups
structures for $G$ which satisfy conditions 3 and 5 above, they are isomorphic
provided there is a one-to-one correspondence between their isolated $V_{0}$-vertices.

\item The graph of groups $\Gamma(\mathcal{F}:G)$ is invariant under the
automorphisms of $G$.

\item The edge splittings of $\Gamma(\mathcal{F}:G)$ are precisely the
canonical splittings of $G$ over two-ended subgroups. Hence if $G$ is not of
finite-by-Fuchsian type and $\mathcal{F}$ is non-empty, then $G$ has a
canonical splitting over a two-ended subgroup.
\end{enumerate}
\end{theorem}

\begin{proof}
Part 1) holds by Theorem \ref{JSJforsmallcommensurisersexists}. Part 2) holds
because every edge group of $\Gamma(\mathcal{F}:G)$ has one end at a $V_{0}%
$-vertex, and each such vertex group is two-ended or of finite-by-Fuchsian type.

By construction, every element of $\mathcal{F}$ is enclosed by some $V_{0}%
$-vertex of $\Gamma(\mathcal{F}:G)$, and every almost invariant subset of $G$
over a two-ended subgroup is equivalent to some element of $\mathcal{F}$. Thus
part 3) follows at once. Parts 4) and 5) follow from Proposition
\ref{XdoesnotcrossanyXiimpliesXisenclosedbyaV1vertex}.

To see that part 6) holds, suppose that $X$ is an almost invariant subset of
$G$ over a two-ended subgroup which is enclosed by some $V_{1}$-vertex $v$ of
$\Gamma(\mathcal{F}:G)$. Part 3) tells us that $X$ is enclosed by some $V_{0}%
$-vertex of $\Gamma(\mathcal{F}:G)$. Now it follows from part 8) of Lemma
\ref{somefactsaboutenclosing} that $X$ is associated to an edge splitting of
$\Gamma(\mathcal{F}:G)$. This implies that $v$ is simple, so that every
$V_{1}$-vertex of $\Gamma(\mathcal{F}:G)$ is simple as required.

Part 7) is exactly the uniqueness result for regular neighbourhoods stated in
Theorem \ref{uniquenessofregbhds}.

For part 8), consider any automorphism $\alpha$ of $G$. Then $\alpha$ induces
a natural $G$-invariant action on $\mathcal{F}$. If we denote by $P$ the
pretree of CCC's of this collection, where we choose one representative for
each element of $\mathcal{F}$ as in Definition
\ref{defnofregnbhdofequivalenceclasses}, then $\alpha$ induces a $G$-invariant
automorphism of $P$ and thus defines a simplicial automorphism of the tree
$T$. Thus there is an induced automorphism of $\Gamma(\mathcal{F}:G)$. (See
\cite{Bass} for a discussion of automorphisms of graphs of groups.)

To prove part 9), we start by observing that it is clear that every edge
splitting of $\Gamma(\mathcal{F}:G)$ is canonical, i.e. it has zero
intersection number with any element of $\mathcal{F}$, because any element of
$\mathcal{F}$ is enclosed by some vertex of $\Gamma(\mathcal{F}:G)$. It
remains to show that these are the only canonical splittings of $G$ over
two-ended subgroups. Let $\sigma$ denote a canonical splitting of $G$ over a
two-ended subgroup $H$, and let $X$ denote the associated $H$-almost invariant
subset of $G$. As $X$ has intersection number zero with every element of
$\mathcal{F}$, Proposition
\ref{XdoesnotcrossanyXiimpliesXisenclosedbyaV1vertex} implies that it is
enclosed by some $V_{1}$-vertex of $T$, the universal covering $G$-tree of
$\Gamma(\mathcal{F}:G)$. But $X$ is also enclosed by a $V_{0}$-vertex of $T$.
Now part 8) of Lemma \ref{somefactsaboutenclosing} shows that $H$ must
stabilise some edge of $T$, which implies that $\sigma$ is conjugate to an
edge splitting of $\Gamma(\mathcal{F}:G)$ as required. Finally, if $G$ is not
of finite-by-Fuchsian type and $\mathcal{F}$ is non-empty, Theorem
\ref{JSJforsmallcommensurisersexists} implies that $\Gamma(\mathcal{F}:G)$
does not consist of a single vertex. Thus $\Gamma(\mathcal{F}:G)$ has at least
one edge and so $G$ has a canonical splitting over a two-ended subgroup.
\end{proof}

At this point we need to discuss further a point about the topological
JSJ-decomposition which we mentioned in section \ref{charsub}. Recall that the
frontier of the characteristic submanifold $V(M)$ of a $3$-manifold $M$ is not
quite the same as the canonical family $\mathcal{T}$ of annuli and tori in
$M$. Some of the components of $\mathcal{T}$ may appear twice in the frontier
of $V(M)$. This means that we can obtain two slightly different graphs of
groups structures for $\pi_{1}(M)$, one graph being dual to $\mathcal{T}$ and
the other dual to the frontier of $V(M)$. The algebraic decomposition
$\Gamma(\mathcal{F}:G)$ which we obtained in Theorem
\ref{JSJforsmallcommensurisersexists} is closer to the second case. However,
there is a natural algebraic object which corresponds to the first case also,
and this can be defined without any regular neighbourhood theory. Namely
consider the family of all conjugacy classes of canonical splittings of $G$
over two-ended subgroups. Any finite subset of this family will be compatible
and so will determine a graph of groups structure for $G$, by Theorem 2.5 of
\cite{SS}. The accessibility result of Theorem \ref{graphsstabilise} implies
that if we take an ascending sequence of such finite families of splittings,
the resulting sequence of graphs of groups structures will stabilise. Thus $G$
has only finitely many conjugacy classes of canonical splittings and they
determine a natural graph of groups structure $\Gamma^{\prime}$ for $G$. Note
that $\Gamma^{\prime}$ need not be bipartite. The following result gives the
connection between $\Gamma(\mathcal{F}:G)$ and $\Gamma^{\prime}$. Recall that
if a graph of groups structure $\Gamma$ for $G$ has a redundant vertex, we can
remove it by replacing the two incident edges by a single edge. If $\Gamma$ is
finite, repeating this will yield a graph of groups structure with no
redundant vertices.

\begin{theorem}
\label{decompositionusingcanonicalsplittings}Let $G$ be a one-ended, finitely
presented group such that whenever $H$ is a two-ended subgroup and
$e(G,H)\geq2$, then $H$ has small commensuriser. Suppose $G$ is not of
finite-by-Fuchsian type and possesses a nontrivial almost invariant subset
over some two-ended subgroup $H$ of $G$, so that $G$ has a canonical
splitting. Let $\Gamma^{\prime}$ denote the graph of groups structure for $G$
determined by a maximal family of non-conjugate canonical splittings of $G$. Then

\begin{enumerate}
\item The graph of groups structure $\Gamma^{\prime}$ for $G$ is obtained from
$\Gamma(\mathcal{F}:G)$ by removing all redundant vertices as above.

\item The vertex groups of $\Gamma^{\prime}$ are each simple, two-ended or of
finite-by-Fuchsian type.

\item Any element of $\mathcal{F}$ is enclosed by some vertex of
$\Gamma^{\prime}$ which is either two-ended or of finite-by-Fuchsian type. In
particular, any splitting of $G$ over a two-ended subgroup is enclosed by such
a vertex of $\Gamma^{\prime}$.
\end{enumerate}
\end{theorem}

\begin{proof}
Let $\Gamma_{1}$ denote the graph of groups structure obtained from
$\Gamma(\mathcal{F}:G)$ by removing all the redundant vertices. Thus the edge
splittings of $\Gamma_{1}$ are exactly those of $\Gamma(\mathcal{F}:G)$, but
now distinct edges of $\Gamma_{1}$ have non-conjugate splittings of $G$
associated. It follows that there is a bijection between the edge splittings
of $\Gamma^{\prime}$ and of $\Gamma_{1}$, so that these graphs of groups
structures for $G$ are isomorphic as required. This proves part 1). Now parts
2) and 3) follow immediately from Theorems
\ref{JSJforsmallcommensurisersexists} and
\ref{propertiesofJSJforsmallcommensurisers}.
\end{proof}

\begin{remark}
In the case of word hyperbolic groups, the graph of groups $\Gamma
(\mathcal{F}:G)$ obtained in Theorem \ref{JSJforsmallcommensurisersexists} is
similar to that obtained by Bowditch in \cite{B1}, but it may differ from that
in \cite{B2} when $\Gamma(\mathcal{F}:G)$ has isolated $V_{0}$-vertices
corresponding to splittings over two-ended subgroups $H$ which have at most
$3$ coends. Moreover our decomposition has redundant vertices corresponding to
canonical splittings which may not appear in Bowditch's decomposition.

The decomposition $\Gamma^{\prime}$ described in Theorem
\ref{decompositionusingcanonicalsplittings} also differs from that obtained by
Sela in \cite{S1} for similar reasons. He further decomposes some of the
vertex groups not of finite-by-Fuchsian type, and seems to take unfoldings of
some of the edge splittings considered here. Moreover some vertex groups which
are of finite-by-Fuchsian type in our terminology (for example, pairs of
pants) may be counted as simple in his decomposition. Conditions 3), 4) and 5)
of Theorem \ref{propertiesofJSJforsmallcommensurisers} are not in either of
their results.
\end{remark}

There is yet another graph of groups structure for $G$ which is also natural
and is similar to $\Gamma^{\prime}$. We define a splitting $\sigma$ of $G$ to
be \textit{splitting-canonical} if it has zero intersection number with every
splitting of $G$ over a two-ended subgroup. As for $\Gamma^{\prime}$, there
can only be only finitely many conjugacy classes of splitting-canonical
splittings of $G$ which are over two-ended subgroups, and these yield a
natural graph of groups structure $\Gamma^{\prime\prime}$ for $G$, whose edge
splittings are these splittings of $G$. The concepts of canonical and
splitting-canonical are in general different. The following example
demonstrates this, and simultaneously gives some insight into the properties
of our regular neighbourhood $\Gamma(\mathcal{F}:G)$.

\begin{example}
Let $A$ be a finitely presented group which is one-ended, admits no splitting
over a two-ended subgroup, and is not finite-by-Fuchsian. For example, $A$
could be the fundamental group of a closed hyperbolic $3$-manifold. Let $B$ be
an infinite cyclic group, let $d\geq4$ be an integer, and let $C$ be the
subgroup $dB$ of index $d$ in $B$. Finally let $G=A\ast_{C}B$. Then $G$ is
finitely presented and the regular neighbourhood $\Gamma(\mathcal{F}:G)$ is
the graph of groups with a single edge given by this splitting. The vertex
which carries $B$ is a $V_{0}$-vertex of small commensuriser type. Note that
if $d$ is $2$ or $3$, this would not be the correct description of
$\Gamma(\mathcal{F}:G)$, as the given splitting would then be isolated, so
that $\Gamma(\mathcal{F}:G)$ would have graph corresponding to $G=A\ast
_{C}C\ast_{C}B$, and the vertex carrying $C$ would be the only $V_{0}$-vertex.
This corresponds to case 1) of Proposition
\ref{regnbhdofeveryelementofB(H)whenHhassmallcommensuriser}. Note also that if
$A$ is word hyperbolic, then so is $G$.

Now let $T$ be the $G$-tree determined by the given splitting of $G$, and let
$v$ be the vertex with stabiliser $B$. Then $v$ has valence $d$ in $T$. Let
$s_{1},\ldots,s_{d}$ denote the edges of $T$ which are incident to $v$. As
usual, we let $\varphi:G\rightarrow V(T)$ be a $G$-equivariant map such that
$\varphi(e)=v$, and orient the edges $s_{1},\ldots,s_{d}$ towards $v$. Let
$Z_{i}$ denote the almost invariant subset of $G$ associated to $s_{i}$. As
$v$ has finite valence, Corollary \ref{nontrivialpartition} implies that the
almost invariant subsets of $G$ which are enclosed by $v$ are all equivalent
to some union of the $Z_{i}$'s and their complements. As each $Z_{i}$ is the
union of all the remaining $Z_{j}^{\ast}$'s, it follows that any almost
invariant subset of $G$ which is enclosed by $v$ is equivalent to some union
of the $Z_{i}^{\ast}$'s. Note that the $Z_{i}^{\ast}$'s are disjoint. Now it
is clear that two such unions fail to cross precisely when they are disjoint
or coincide. Suppose that $d$ is a composite number, say $d=ab$. Then the
above splitting of $G$ can be refined by splitting at the vertex carrying $B$
to obtain $G=A\ast_{C}aB\ast_{aB}B$. This induces a refinement of $T$, in
which the vertex $v$ of valence $d=ab$ is replaced by a tree which is the cone
on points $v_{1},\ldots,v_{a}$, and the $d$ edges which were attached to $v$
are now attached to the vertices $v_{1},\ldots,v_{a}$ with $b$ of these edges
being attached to each $v_{j}$. This shows that if we pick some $Z_{i}$, then
the union of all its translates by the subgroup of index $a$ in $B$ is
associated to a splitting of $G$. Now pick a prime $p\geq3$ and let $a=b=p$,
so that $d=p^{2}$. The above discussion shows that if an almost invariant
subset $X$ is associated to a splitting $\sigma$ of $G$ and is enclosed by
$v$, then, up to equivalence, either $X$ is a single $Z_{i}$ or $X$ is the
orbit of a single $Z_{i}$ under the action of the subgroup $pB$ of $B$. Thus
$\sigma$ must be conjugate to the original splitting over $C$, or to the
splitting over $pB$, described above. In particular, it follows that the
splitting of $G$ over $pB$ is splitting-canonical. But there are many almost
invariant subsets of $G$ enclosed by $v$ which it crosses, so it is not canonical.
\end{example}

\section{Canonical decompositions over two-ended groups when commensurisers
are large\label{JSJforlargecommensurisers}}

Again we consider a one-ended, finitely presented group $G$ and almost
invariant subsets over two-ended subgroups. This time, we do not assume that
our two-ended subgroups have small commensurisers. This leads to two
additional difficulties. When we form the regular neighbourhood of a finite
number of almost invariant subsets of $G$ which are over two-ended subgroups,
the edge groups of the regular neighbourhood may no longer be two-ended. As
pointed out in the introduction, this happens in the topological situation.
For if one wants to enclose essential annuli in a $3$-manifold $M$, one may
obtain Seifert fibre space components of $V(M)$ whose frontier has toral
components. In the case of general finitely presented groups, the edge groups
may be even more complicated. Thus if we proceed, as in the previous section,
to take regular neighbourhoods of larger and larger finite collections of
almost invariant subsets over two-ended subgroups, we will not be able to show
that our construction stabilises. Hence we are forced to consider directly
regular neighbourhoods of infinite families of almost invariant sets. This
leads to the additional problem of showing that the pretrees which appear in
the regular neighbourhood construction are discrete. For this, Theorem
\ref{Booleanalgebrafinitelygenerated} on the finite generation of certain
Boolean algebras plays a key role.

The main result of this section is our version of the JSJ-decomposition for
arbitrary finitely presented groups with one end.

\begin{theorem}
\label{JSJforlargecommensurisersexist}Let $G$ be a one-ended, finitely
presented group, and let $\mathcal{F}$ denote the collection of equivalence
classes of all nontrivial almost invariant subsets of $G$ which are over a
two-ended subgroup.

Then the regular neighbourhood construction of section
\ref{regnbhds:construction} works and yields a regular neighbourhood
$\Gamma(\mathcal{F}:G)$. Each $V_{0}$-vertex $v$ of $\Gamma(\mathcal{F}:G)$
satisfies one of the following conditions:

\begin{enumerate}
\item $v$ is isolated.

\item $v$ is of finite-by-Fuchsian type.

\item $G(v)$ is the full commensuriser $Comm_{G}(H)$ for some two-ended
subgroup $H$, such that $e(G,H)\geq2$.
\end{enumerate}

$\Gamma(\mathcal{F}:G)$ consists of a single vertex if and only if
$\mathcal{F}$ is empty, or $G$ itself satisfies one of conditions 2) or 3) above.
\end{theorem}

\begin{remark}
We will say that a $V_{0}$-vertex in case 3) above is \textsl{of commensuriser
type}, if $v$ is not isolated nor of finite-by-Fuchsian type, and is
\textsl{of large commensuriser type}, if in addition $H$ has large commensuriser.

Note that even if $G$ is finitely presented, Example \ref{example
ofCommG(H)isnotfg} shows that a vertex group of commensuriser type need not be
finitely generated.
\end{remark}

\begin{proof}
We start by picking a representative for each element of $\mathcal{F}$,
subject to the condition that if $A$ and $B$ are elements of $\mathcal{F}$
such that $B=gA$, for some $g$ in $G$, then the representatives $X$ and $Y$
chosen for $A$ and $B$ must satisfy $Y=gX$.

Before proceeding further, we consider the very special case when $G$ is
virtually $\mathbb{Z}\times\mathbb{Z}$. The equivalence classes of almost
invariant subsets of $\mathbb{Z}\times\mathbb{Z}$ which are over two-ended
subgroups correspond to all the simple closed curves on the torus. It follows
that the collection of all the chosen representatives of elements of
$\mathcal{F}$ is cross connected, so that the required regular neighbourhood
exists and consists of a single $V_{0}$-vertex with associated group $G$. As
$G$ is finite-by-Fuchsian, this proves the theorem in this case. In the
following we will assume that $G$ is not virtually $\mathbb{Z}\times
\mathbb{Z}$.

Next we let $\mathcal{F}_{0}$ denote the subset of $\mathcal{F}$ whose
elements are represented by $H$-almost invariant subsets of $G$, such that $H$
has small commensuriser. The proof of Theorem
\ref{JSJforsmallcommensurisersexists} applies to show that the regular
neighbourhood construction of section \ref{regnbhds:construction} works to
yield $\Gamma_{0}=\Gamma(\mathcal{F}_{0}:G)$. In what follows, we will express
$\mathcal{F}$ as an ascending sequence of subsets $\mathcal{F}_{i}$ of
$\mathcal{F}$, for $i\geq0$, and show that each $\mathcal{F}_{i}$ has a
regular neighbourhood $\Gamma_{i}=\Gamma(\mathcal{F}_{i}:G)$. Finally, we will
show that the sequence $\Gamma_{i}$ must stabilise eventually. The resulting
graph of groups structure $\Gamma$ for $G$ will be the required regular
neighbourhood $\Gamma=\Gamma(\mathcal{F}:G)$.

By our definition of $\mathcal{F}_{0}$, any element of $\mathcal{F-}$
$\mathcal{F}_{0}$ will be represented by a $H$-almost invariant subset $X$ of
$G$, such that $H$ has large commensuriser. From the collection of all such
subgroups $H$ of $G$, we choose one group from each conjugacy class and denote
the chosen groups by $H_{j}$, $j\geq1$. We choose $\mathcal{F}_{i+1}$ to be
the union of $\mathcal{F}_{i}$ and all translates by $G$ of elements of
$B(H_{i+1})$. Clearly the union of the $\mathcal{F}_{i}$'s equals
$\mathcal{F}$.

Let $E_{0}$ denote the collection of the chosen representatives of all the
elements of $\mathcal{F}_{0}$. Denote $H_{1}$ by $H$, so that $\mathcal{F}%
_{1}$ is the union of $\mathcal{F}_{0}$ and all translates by $G$ of elements
of $B(H)$. We let $E_{1}$ denote the collection of the chosen representatives
of all the elements of $\mathcal{F}_{1}$. Theorem
\ref{Booleanalgebrafinitelygenerated} tells us that $B(H)$ has a finite system
of generators when we regard $B(H)\;$as a\ Boolean algebra over $Comm_{G}(H)$.
We let $X_{1},\ldots,X_{n}$ be the chosen representatives of this system of
generators. Thus any element of $B(H)$ can be represented by taking finite
sums of finite intersections of the $X_{i}$'s and their complements. The proof
of Proposition \ref{largecommensuriserimpliesoneinfiniteCCC} shows that any
isolated element of $B(H)$ is a translate of some $X_{i}$ by an element of
$Comm_{G}(H)$.

In order to show that $\mathcal{F}_{1}$ has a regular neighbourhood
$\Gamma_{1}=\Gamma(\mathcal{F}_{1}:G)$, we will need to recall some more facts
from section \ref{coendswhencommensuriserislarge}. Let $G$ be a finitely
presented group and let $H$ be a two-ended subgroup of $G$ such that $H$ has
large commensuriser and $G$ has a nontrivial $H$-almost invariant subset. As
we are assuming that $G$ is not virtually $\mathbb{Z}\times\mathbb{Z}$,
Proposition \ref{Comm_G(H)largeimpliescoendsare12orinfty} tells us that the
number of coends of $H$ in $G$ is infinite. Proposition
\ref{largecommensuriserimpliesoneinfiniteCCC} tells us that the corresponding
cross-connected components (CCC's) consist of translates of a finite number of
isolated almost invariant sets and a single CCC $H_{\infty}$ which consists of
an infinite number of almost invariant sets (we will call this `the infinite
CCC' corresponding to $H$). We saw that any $K$-almost invariant set with $K$
commensurable with $H$ is either isolated or in $H_{\infty}$. Note also that
the stabiliser of $H_{\infty}$ is $Comm_{G}(H)$.

Let $P_{0}$ be the pretree of CCC's of $\overline{E_{0}}$ and $P_{1}$ that of
$\overline{E_{1}}$. We know that $P_{0}$ is discrete and want to show that
$P_{1}$ is discrete. We note that the natural map $P_{0}\rightarrow P_{1}$ is injective.

Recall that the proof that $P_{0}$ is discrete, depended crucially on the
discreteness of $E_{0}$, which holds because $E_{0}$ consists of translates of
a finite family of almost invariant sets. Recall that discreteness of a
partially ordered set $F$ means that, for any $U$, $V\in F$, there are only
finitely many $Z\in F$ such that $U\leq Z\leq V$. (See Lemma
\ref{posatisfiesDunwoody}).

To show that $P_{1}$ is discrete, we need to show that if $A$ and $C$ are
distinct CCC's of $\overline{E_{1}}$, then there are only finitely many CCC's
$B$ of $\overline{E_{1}}$ such that $ABC$. We know that there are only
finitely many finite CCC's between $A$ and $C$, because these come from
elements of $E_{0}$, which is discrete. So it remains to show that there are
only finitely many infinite CCC's $B$ between $A$ and $C$. By construction,
the only infinite CCC's in $E_{1}$ are translates of the infinite CCC
$H_{\infty}$ corresponding to the commensurability class of $H$. We choose
elements $U$ and $W$ of $E_{1}$, such that $\overline{U}\in A$ and
$\overline{W}\in C$. Suppose that $B=H_{\infty}$. The definition of
betweenness for CCC's implies that there exists an element $\overline{V}$ of
$H_{\infty}$ such that $U<V<W$. Lemma \ref{B(H)fgimpliesPisdiscrete} below
shows that we can choose $V$ in a special way. It shows that there is an
almost invariant subset $X$ of $G$ which is a translate of one of the $X_{i}%
$'s by an element of $Comm_{G}(H)$, such that $U\leq X\leq W$. Note that $X$
will represent an element of $B(H)$. If $B$ is a translate $gH_{\infty}$ of
$H_{\infty}$ such that $ABC$, then $H_{\infty}$ lies between $g^{-1}A$ and
$g^{-1}C$, and applying Lemma \ref{B(H)fgimpliesPisdiscrete} to this situation
yields an almost invariant subset $X$ of $G$ as above except that $U\leq
gX\leq W$. Thus $gX$ is a translate by an element of $gComm_{G}(H)$ of one of
the $X_{i}$'s. As only finitely many translates of the $X_{i}$'s can lie
between $U$ and $W$, it follows that there are only finitely many translates
of $H_{\infty}$ between $A$ and $C$. Thus $P_{1}$ is a discrete pretree.

Similar arguments yield an inductive proof that the pretree $P_{i}$ is
discrete, for all $i\geq0$, so that each $\mathcal{F}_{i}$ has a regular
neighbourhood $\Gamma_{i}=\Gamma(\mathcal{F}_{i}:G)$, as required.

In order to show that the sequence $\Gamma_{i}$ stabilises, we first need to
describe the $V_{0}$-vertices of $\Gamma_{i}$. The results of sections
\ref{coendswhencommensuriserissmall} and \ref{coendswhencommensuriserislarge}
show that these are of three types, isolated, finite-by-Fuchsian type, and
commensuriser type, which is where a $V_{0}$-vertex carries the group
$Comm_{G}(H)$ for some two-ended subgroup $H$ of $G$ such that $e(G,H)\geq2$.
Further, our construction implies that each $\Gamma_{i}$ has $i$ vertices of
large commensuriser type. Now we claim that each $V_{0}$-vertex of large
commensuriser type encloses a splitting of $G$ over a two-ended subgroup of
$G$. For consider the $V_{0}$-vertex $v$ determined by the infinite CCC
$H_{\infty}$. Theorem \ref{D-Rollertheorem} of Dunwoody and Roller tells us
that $G$ splits over some subgroup commensurable with $H$. As the almost
invariant subset associated to this splitting must lie in $Q(H)$, and the
proof of Lemma \ref{largecommensuriserimpliesoneinfiniteCCC} shows that every
element of $Q(H)$ is enclosed by $H_{\infty}$, it follows that $v$ encloses a
splitting of $G$ over a two-ended subgroup, as claimed. Thus we can refine
$\Gamma_{i}$ by splitting at each $V_{0}$-vertex of large commensuriser type
using such a splitting, to obtain a new graph of groups structure $\Gamma
_{i}^{\prime}$ for $G$. This construction means that if we let $f(i)$ denote
the number of those edge splittings of $\Gamma_{i}^{\prime}$ which are over a
two-ended group, then $f(i)$ is strictly increasing. Now Theorem
\ref{graphsstabilise} implies that the sequence $\Gamma_{i}^{\prime}$ must
stabilise and hence that the sequence $\Gamma_{i}$ must stabilise, as
required. We call this final graph of groups $\Gamma$. By construction,
$\Gamma$ is the required regular neighbourhood $\Gamma(\mathcal{F}:G)$. Each
$V_{0}$-vertex group of $\Gamma$ satisfies one of the three conditions in the
statement of the theorem because this holds for each $\Gamma_{i}$.

Finally, $\Gamma$ will consist of a single vertex if and only if either
$\mathcal{F}$ is empty or the representatives of $\mathcal{F}$ lie in a single
CCC. In the second case, the vertex group of $\Gamma$ must satisfy one of the
three conditions in the statement of the theorem. As the vertex group is $G$
and we assumed that $G$ is one-ended, it follows that $\Gamma(\mathcal{F}:G)$
consists of a single vertex if and only if $\mathcal{F}$ is empty, or $G$ is
of finite-by-Fuchsian type, or $G$ equals $Comm_{G}(H)$ for some two-ended
subgroup $H$, such that $e(G,H)\geq2$.
\end{proof}

Now we prove the following result which was used in the above proof. Recall
that $X_{1},\ldots,X_{n}$ are the chosen representatives of elements of $B(H)$
which generate $B(H)$ over $Comm_{G}(H)$.

\begin{lemma}
\label{B(H)fgimpliesPisdiscrete}Let $A$, $B$ and $C$ be distinct CCC's of
$\overline{E_{1}}$, such that $B$ equals the infinite CCC $H_{\infty}$, and
$B$ lies between $A$ and $C$. Let $U$ and $W$ be almost invariant subsets of
$G$ such that $\overline{U}\in A$ and $\overline{W}\in C$. Then there is an
almost invariant subset $X$ of $G$ which is a translate of one of the $X_{i}%
$'s by an element of $Comm_{G}(H)$, such that $U\leq X\leq W$.
\end{lemma}

\begin{proof}
As $B$ lies between $A$ and $C$, there is $\overline{V}\in B$ with $U<V<W$. In
particular, $V$ represents an element of $B(H)$. As $B(H)$ is generated over
$Comm_{G}(H)$ by the $X_{i}$'s, any element of $B(H)$ can be represented as a
finite sum of finite intersections of translates of the $X_{i}$'s and their
complements. We will need to prove the following two claims.

\begin{claim}
If $Y$ and $Z$ represent elements of $B(H)$ and $U<(Y+Z)<W$, then either $U$
or $W$ represents an element of $B(H)$, or one of the four sets $Y$, $Y^{\ast
}$, $Z$, $Z^{\ast}$ lies between $U$ and $W$.
\end{claim}

\begin{claim}
If $Y$ and $Z$ represent elements of $B(H)$ and $U<(Y\cap Z)<W$, then either
$U$ or $W$ represents an element of $B(H)$, or one of the four sets $Y$,
$Y^{\ast}$, $Z$, $Z^{\ast}$ lies between $U$ and $W$.
\end{claim}

Starting from an expression of $V$ as a finite sum of finite intersections of
the $X_{i}$'s and their complements, we apply one of these two claims. If one
of the four sets $Y$, $Y^{\ast}$, $Z$, $Z^{\ast}$ lies between $U$ and $W$, we
again apply one of these two claims, and repeat this process as long as
possible. This process will eventually stop, at which point either we will see
that $U$ or $W$ represents an element of $B(H)$, or we will find an almost
invariant subset $X$ of $G$ which is a translate of one of the $X_{i}$'s by an
element of $Comm_{G}(H)$, such that $U<X<W$. If $U$ represents an element of
$B(H)$, then $U$ is enclosed by the CCC $B=H_{\infty}$. As $U$ is also
enclosed by the CCC $A$ which is distinct from $B$, part 8) of Lemma
\ref{somefactsaboutenclosing} shows that $U$ is equivalent to $Z_{s}$ for some
edge $s$ of $T$ which is incident to $H_{\infty}$. In particular, $U$
represents an isolated element of $B(H)$, which implies that $U$ is equivalent
to a translate of some $X_{i}$ by an element of $Comm_{G}(H)$. In this case,
we find the required set $X$ by simply choosing $X=U$. Similarly if $W$
represents an element of $B(H)$, we can choose $X=W$. Thus in all cases, we
have found the required set $X$.

To prove Claim 1, suppose that $U<(Y+Z)<W$. If $Y$ is isolated in $E_{1}$,
then it automatically cannot cross $U$ or $W$. If $Y$ is not isolated, then it
must lie in $H_{\infty}$. In particular, it lies in a different CCC from $U$
and $W$, and so again cannot cross $U$ or $W$. Thus in either case, one of the
four sets $U^{(\ast)}\cap Y^{(\ast)}$ is small, and one of the four sets
$W^{(\ast)}\cap Y^{(\ast)}$ is small. Similarly one of the four sets
$U^{(\ast)}\cap Z^{(\ast)}$ is small, and one of the four sets $W^{(\ast)}\cap
Z^{(\ast)}$ is small. In order to understand our argument, it will be helpful
to first suppose that the inequality $\leq$ is actual inclusion, and that $Y$
and $Z$ are each nested with respect to $U$ and to $W$. As $Y+Z=(Y\cap
Z^{\ast})\cup(Y^{\ast}\cap Z)$, and $U$ and $W$ do not cross $Y$ or $Z$, it is
easy to see that we must have $U\subset Y\cap Z^{\ast}$ or $U\subset Y^{\ast
}\cap Z$, and that we must have $W^{\ast}\subset Y\cap Z$ or $W^{\ast}\subset
Y^{\ast}\cap Z^{\ast}$. For each of the four possibilities, we will then
obtain one of the four inclusions $U\subset Y^{(\ast)}\subset W$ or $U\subset
Z^{(\ast)}\subset W$ as required. Here is the formal argument.

We consider the inequality $U<(Y+Z)$. This is equivalent to the statement that
$U\cap(Y+Z)^{\ast}$ is small. As $(Y+Z)^{\ast}=(Y\cap Z)\cup(Y^{\ast}\cap
Z^{\ast})$, it follows that $U\cap(Y\cap Z)$ and $U\cap(Y^{\ast}\cap Z^{\ast
})$ are each small. Hence $U^{\ast}\cap(Y\cap Z)$ and $U^{\ast}\cap(Y^{\ast
}\cap Z^{\ast})$ are each not small. This implies that each of $U^{\ast}\cap
Y$, $U^{\ast}\cap Z$, $U^{\ast}\cap Y^{\ast}$ and $U^{\ast}\cap Z^{\ast}$ is
not small. As $U$ does not cross $Y$ or $Z$, we know that one of the four sets
$U^{(\ast)}\cap Y^{(\ast)}$ is small and that one of the four sets $U^{(\ast
)}\cap Z^{(\ast)}$ is small. It follows that one of the two sets $U\cap
Y^{(\ast)}$ is small and that one of the two sets $U\cap Z^{(\ast)}$ is small.
If $U\cap Y$ is small, then $U\cap(Y\cap Z^{\ast})$ is small, so that we have
$U\leq(Y^{\ast}\cap Z)$. Similar arguments apply in the other three cases. We
conclude that $U\leq(Y^{\ast}\cap Z)$ or $U\leq(Y\cap Z^{\ast})$.

Similar arguments show that if $(Y+Z)<W$, then $W^{\ast}\leq(Y\cap Z)$ or
$W^{\ast}\leq(Y^{\ast}\cap Z^{\ast})$. In each of these four cases, one sees
that one of the four sets $Y$, $Y^{\ast}$, $Z$, $Z^{\ast}$ is either equal to
$U$ or $W$ or lies between $U$ and $W$. Thus either $U$ or $W$ represents an
element of $B(H)$, or one of the four sets $Y$, $Y^{\ast}$, $Z$, $Z^{\ast}$
lies between $U$ and $W$, as required.

To prove Claim 2, suppose that $U<(Y\cap Z)<W$ . Much as in the proof of Claim
1, one can show that as $W$ does not cross $Y$ or $Z$, we have $W^{\ast}%
\leq(Y\cap Z^{\ast})$, $W^{\ast}\leq(Y^{\ast}\cap Z)$ or $W^{\ast}\leq
(Y^{\ast}\cap Z^{\ast})$. It follows that we have one of $Y\leq W$ or $Z\leq
W$. Thus we have $U\leq Y\leq W$ or $U\leq Z\leq W$. As above, it follows that
either $U$ or $W$ represents an element of $B(H)$, or one of the two sets $Y$
and $Z$ lies between $U$ and $W$, as required.
\end{proof}

In order to understand $\Gamma(\mathcal{F}:G)$ in more detail, our next result
lists some properties which follow almost immediately from the above theorem
and the properties of an algebraic regular neighbourhood.

\begin{theorem}
\label{propertiesofJSJforlargecommensurisers}Let $G$ be a one-ended, finitely
presented group, and let $\mathcal{F}$ denote the collection of equivalence
classes of all nontrivial almost invariant subsets of $G$ which are over a
two-ended subgroup.

Then the regular neighbourhood $\Gamma(\mathcal{F}:G)$ is a minimal bipartite
graph of groups decomposition of $G$ with the following properties:

\begin{enumerate}
\item each $V_{0}$-vertex $v$ of $\Gamma(\mathcal{F}:G)$ satisfies one of the
following conditions:

\begin{enumerate}
\item $v$ is isolated.

\item $v$ is of finite-by-Fuchsian type.

\item $G(v)$ is the full commensuriser $Comm_{G}(H)$ for some two-ended
subgroup $H$, such that $e(G,H)\geq2$.
\end{enumerate}

Further, if $H$ is a two-ended subgroup of $G$ such that $e(G,H)\geq2$, and if
$H$ has large commensuriser, then $\Gamma(\mathcal{F}:G)$ will have a $V_{0}%
$-vertex $v$ such that $G(v)=Comm_{G}(H)$.

\item If an edge of $\Gamma(\mathcal{F}:G)$ is incident to a $V_{0}$-vertex of
type a) or b) above, then it carries a two-ended group.

\item any representative of an element of $\mathcal{F}$ is enclosed by some
$V_{0}$-vertex of $\Gamma(\mathcal{F}:G)$, and each $V_{0}$-vertex of
$\Gamma(\mathcal{F}:G)$ encloses such a subset of $G$. In particular, any
splitting of $G$ over a two-ended subgroup is enclosed by some $V_{0}$-vertex
of $\Gamma(\mathcal{F}:G)$.

\item if $X$ is an almost invariant subset of $G$ over a finitely generated
subgroup $H$, and if $X$ does not cross any element of $\mathcal{F}$, then $X$
is enclosed by a $V_{1}$-vertex of $\Gamma(\mathcal{F}:G)$.

\item if $X$ is a $H$-almost invariant subset of $G$ associated to a splitting
of $G$ over $H$, and if $X$ does not cross any element of $\mathcal{F}$, then
$X$ is enclosed by a $V_{1}$-vertex of $\Gamma(\mathcal{F}:G)$.

\item the $V_{1}$-vertex groups of $\Gamma(\mathcal{F}:G)$ are simple. (See
Definition \ref{defnofsimplefortwo-endedsubgroups}.)

\item If $\Gamma_{1}$ and $\Gamma_{2}$ are minimal bipartite graphs of groups
structures for $G$ which satisfy conditions 3) and 5) above, they are
isomorphic provided there is a one-to-one correspondence between their
isolated $V_{0}$-vertices.

\item The graph of groups $\Gamma(\mathcal{F}:G)$ is invariant under the
automorphisms of $G$.

\item The canonical splittings of $G$ over two-ended subgroups are precisely
those edge splittings of $\Gamma(\mathcal{F}:G)$ which are over two-ended
subgroups. This includes, but need not be limited to, all those edges of
$\Gamma(\mathcal{F}:G)$ which are incident to $V_{0}$-vertices whose
associated groups are of type a) or b) above.
\end{enumerate}
\end{theorem}

\begin{proof}
The description of the possible types of $V_{0}$-vertex given in part 1)
follows from Theorem \ref{JSJforlargecommensurisersexist}. Further, if $H$ is
a two-ended subgroup of $G$ such that $e(G,H)\geq2$, and if $H$ has large
commensuriser, then either $G$ is virtually $\mathbb{Z}\times\mathbb{Z}$, or
$e(G,H)$ is infinite. In the first case, $\Gamma(\mathcal{F}:G)$ consists of a
single vertex labeled $G$, and in the second case, Proposition
\ref{largecommensuriserimpliesoneinfiniteCCC} shows that $\Gamma
(\mathcal{F}:G)$ will have a $V_{0}$-vertex $v$ such that $G(v)=Comm_{G}(H)$.
Thus in either case, $\Gamma(\mathcal{F}:G)$ will have a $V_{0}$-vertex $v$
such that $G(v)=Comm_{G}(H)$.

The proofs for parts 2)-9) are the same as for parts 2)-9) of Theorem
\ref{propertiesofJSJforsmallcommensurisers}.
\end{proof}

At this point, we note the connection between the above results and the
Algebraic Annulus Theorem \cite{D-Swenson} for finitely generated groups. (See
also \cite{SS2} and \cite{B1} for the case of word hyperbolic groups.) The
proof we give below is for finitely presented groups only and is not
essentially different from that given by Dunwoody and Swenson in
\cite{D-Swenson}. We include the argument here for completeness only. In the
topological context, one can deduce the Annulus Theorem from the
JSJ-decomposition in much the same way. Clearly regular neighbourhood theory
is not essential for the proof of the Algebraic Annulus Theorem. Nor can it
yield a proof for groups which are not finitely presented.

\begin{theorem}
(Algebraic Annulus Theorem) Let $G$ be a one-ended, finitely presented group.
If $G$ has a two-ended subgroup $H$ such that $e(G,H)\geq2$, then either $G$
splits over some two-ended subgroup or $G$ is of finite-by-Fuchsian type.
\end{theorem}

\begin{proof}
The assumption implies that the set $\mathcal{F}$ in Theorem
\ref{JSJforlargecommensurisersexist} and
\ref{propertiesofJSJforlargecommensurisers} is non-empty. Applying Theorem
\ref{JSJforlargecommensurisersexist} yields the regular neighbourhood
$\Gamma=\Gamma(\mathcal{F}:G)$. If $\Gamma$ consists of a single vertex, then
$G$ is of finite-by-Fuchsian type. Otherwise, each $V_{0}$-vertex of $\Gamma$
has at least one incident edge. Any edge incident to a $V_{0}$-vertex of type
a) or b) carries a two-ended group and so yields a splitting of $G$ over such
a group. If $\Gamma$ has no such $V_{0}$-vertices, then it must have a $V_{0}%
$-vertex $v$ of type c), and Theorem \ref{D-Rollertheorem} shows that $G$
splits over some two-ended subgroup. The result follows.
\end{proof}

A key point about the preceding arguments was that we considered all almost
invariant subsets of $G$ over two-ended subgroups and did not restrict to
those which are associated to splittings. However, now we have Theorem
\ref{JSJforlargecommensurisersexist}, it is quite easy to deduce the existence
of a regular neighbourhood of this smaller collection of almost invariant
subsets. The result we obtain is the following.

\begin{theorem}
Let $G$ be a one-ended, finitely presented group, and let $\mathcal{S}$ denote
the collection of equivalence classes of all almost invariant subsets of $G$
which are associated to a splitting of $G$ over a two-ended subgroup.

Then the regular neighbourhood construction of section
\ref{regnbhds:construction} works and yields a regular neighbourhood
$\Gamma(\mathcal{S}:G)$. Each $V_{0}$-vertex $v$ of $\Gamma(\mathcal{S}:G)$
satisfies one of the following conditions:

\begin{enumerate}
\item $v$ is isolated.

\item $v$ is of finite-by-Fuchsian type.

\item $G(v)$ contains a two-ended subgroup $H$ which it commensurises, such
that $e(G,H)\geq2$.
\end{enumerate}

If $\Gamma(\mathcal{S}:G)$ consists of a single vertex, then either
$\mathcal{S}$ is empty, or $G$ itself satisfies one of conditions 2) or 3) above.
\end{theorem}

\begin{remark}
\label{remarkonregnbhdofsplittings}Note that even if $G$ is finitely
presented, Example \ref{example ofCommG(H)isnotfg} shows that a vertex group
of type 3) need not be finitely generated. Note also that if $G$ commensurises
a two-ended subgroup $H$ such that $e(G,H)\geq2$, then $\Gamma(\mathcal{S}:G)$
need not consist of a single vertex. This is in contrast with the situation of
Theorem \ref{propertiesofJSJforlargecommensurisers}. We give some simple
examples after the proof of the theorem.
\end{remark}

\begin{proof}
We start from the regular neighbourhood $\Gamma(\mathcal{F}:G)$ obtained in
Theorem \ref{JSJforlargecommensurisersexist}. Now we know that the
construction of section \ref{regnbhds:construction} works, we can consider
this construction directly. As in the proof of Theorem
\ref{JSJforlargecommensurisersexist}, we start by picking a representative for
each element of $\mathcal{F}$, subject to the condition that if $A$ and $B$
are elements of $\mathcal{F}$ such that $B=gA$, for some $g$ in $G$, then the
representatives $X$ and $Y$ chosen for $A$ and $B$ must satisfy $Y=gX$. This
determines the set $E$ of all translates of these subsets of $G$, and we let
$S$ denote the subset of $E$ consisting of almost invariant subsets of $G$
which are associated to a splitting. When we replace $E$ by $S$, we want to
describe how the CCC's and their stabilisers alter. We claim that each
$G$-orbit of CCC's of $\overline{E}$ is the union of a finite number of
$G$-orbits of CCC's of $\overline{S}$. Given this, it is not difficult to
verify that the pretree determined by $S$ must be discrete, so that the
regular neighbourhood $\Gamma(\mathcal{S}:G)$ exists. The reason for our claim
is simply that otherwise, some $V_{0}$-vertex of $\Gamma(\mathcal{F}:G)$ would
enclose an infinite number of non-conjugate compatible splittings of $G$ over
two-ended subgroups obtained by picking a splitting from each CCC of
$\overline{S}$, and this would contradict the accessibility result of Theorem
\ref{graphsstabilise}.
\end{proof}

Now we can give the examples referred to in Remark
\ref{remarkonregnbhdofsplittings}. These examples are to demonstrate that if
$G$ commensurises a two-ended subgroup $H$ such that $e(G,H)\geq2$, then
$\Gamma(\mathcal{S}:G)$ need not consist of a single vertex.

\begin{example}
Let $G_{p,q}=A\ast_{C}B$, where $A$ and $B$ are both infinite cyclic and $C$
has index $p$ in $A$ and index $q$ in $B$. Thus $G_{p,q}$ centralises, and
hence commensurises, the two-ended subgroup $C$. If $p,q\geq2$, then $G$
splits over $C$, so that $e(G,C)\geq2$. If, in addition, we exclude the case
$p=q=2$, then $\Gamma(\mathcal{S}:G_{p,q})$ does not consist of a single
vertex. For it is easy to show that, up to conjugacy, this is the only
splitting of $G_{p,q}$ over a two-ended subgroup, which implies that
$\Gamma(\mathcal{S}:G_{p,q})$ is the graph of groups associated to
$G=A\ast_{C}C\ast_{C}B$, where the vertex carrying $C$ is the only $V_{0}%
$-vertex. However, $\Gamma(\mathcal{S}:G_{2,2})$ does consist of a single vertex.

These examples are closely related to some topological examples discussed in
section \ref{charsub}. Let $W_{p,q}$ denote the $3$-manifold which is obtained
by gluing two solid tori along an annulus $A$ which has degree $p$ in one
solid torus and degree $q$ in the other. Thus $G_{p,q}=\pi_{1}(W_{p,q})$. If
$p,q\geq2$ and we exclude the case $p=q=2$, the fact that there is only one
splitting of $G_{p,q}$ over a two-ended subgroup, up to conjugacy, corresponds
to the fact that the annulus $A$ is the only embedded essential annulus in
$W_{p,q}$, up to isotopy. However, the reader should be warned that this fact
about $G_{p,q}$ does not follow from that about $W_{p,q}$, but needs its own
proof. This is because a splitting of a $3$-manifold group over an infinite
cyclic subgroup need not, in general, be induced by an embedded annulus. On
the other hand, the fact that $\Gamma(\mathcal{S}:G_{2,2})$ consists of a
single vertex follows from the fact that $W_{2,2}$ is filled by essential
embedded annuli.
\end{example}

\section{Some examples\label{examples}}

We start this section with some specific examples, and then give some more
general ones.

Our first example is of a one-ended, finitely presented group $G$ such that
the regular neighbourhood $\Gamma(G)$ of equivalence classes of all almost
invariant subsets of $G$ which are over two-ended subgroups has a $V_{0}%
$-vertex of commensuriser type with a non-finitely generated vertex group.

\begin{example}
\label{example ofCommG(H)isnotfg}We start by showing that there exists a
one-ended, finitely presented group $A$ which has an infinite cyclic subgroup
$H$ such that $Comm_{A}(H)$ is not finitely generated. To construct such a
group, we take a free group $F$ of countably infinite rank, and an infinite
cyclic group $H$. As $F$ embeds in $F_{2}$, the free group of rank $2$, we can
embed $F\times H$ in $E=(F\times H)\ast_{F}F_{2}$. Clearly $E$ is finitely
generated and recursively presented, and $Comm_{E}(H)=F\times H$. Now we embed
$E$ in a finitely presented group $L$ using Higman's Embedding Theorem
\cite{Higman}. He first constructs a certain finitely presented group $K$ and
then constructs $L$ as a HNN extension of the form $(K\times E)\ast_{H}$. It
is clear from the construction that $Comm_{L}(H)=K\times Comm_{E}(H)=K\times
F\times H$, which is not finitely generated. If $L$ is one-ended, we take
$A=L$. Otherwise, the accessibility result of Dunwoody in \cite{Dunwoody}
implies that $L$ can be expressed as the fundamental group of a graph of
groups with all edge groups finite, and all vertex groups having zero or one
end. The vertex groups must then be finitely presented. Now one of these
vertex groups must contain $K\times F\times H$, and this is the required
one-ended, finitely presented group $A$.

Now let $C$ denote $K\times F$ so that $Comm_{A}(H)=C\times H$, let $D$ denote
any nontrivial finitely presented group, and let $B$ denote $C\ast D$. We
define $G=A\ast_{C\times H}(B\times H)$. As $B=C\ast D$, we can also write
$G=A\ast_{H}(D\ast H)$, so that $G$ is finitely presented and splits over $H$.
Now $Comm_{G}(H)=B\times H$ which is not finitely generated. As $G$ splits
over $H$ and $Comm_{G}(H)$ is not finitely generated, it follows that
$\Gamma(G)$ has a $V_{0}$-vertex with associated group $Comm_{G}(H)$.

If one can choose $A$ so as not to split over any two-ended subgroup, then the
graph of groups $\Gamma(G)$ consists of a single edge which induces the
decomposition $G=A\ast_{C\times H}(B\times H)$.
\end{example}

It is natural to ask what can be said about the edge groups of $\Gamma(G)$
which are incident to a $V_{0}$-vertex of commensuriser type. We have already
pointed out that if the commensuriser vertex group is not finitely generated,
then some incident edge group must also be not finitely generated. In the case
of $3$-manifolds, the components of $V(M)$ which correspond to a $V_{0}%
$-vertex of commensuriser type are Seifert fibre spaces, and each frontier
component is a vertical annulus or torus. In particular, the incident edge
groups all contain the normal subgroup $H$ carried by a regular fibre of the
Seifert fibre space. However the following example shows that this does not
hold in general.

\begin{example}
Let $K$ and $L$ be free groups of rank at least $2$, let $H$ be an infinite
cyclic group and let $A$ and $B$ be groups which properly contain $K$ and $L$
respectively. Then define $G=A\ast_{K}(K\times H)\ast_{H}(H\times L)\ast_{L}%
B$. Then $Comm_{G}(H)=(K\ast L)\times H$. As $G$ splits over $H$, it follows
that $\Gamma(G)$ has a $V_{0}$-vertex of commensuriser type with associated
group $Comm_{G}(H)$. If we choose $A$ and $B$ to be one-ended, it is easy to
see that $G$ is also one-ended. For if $G$ splits over a finite subgroup, the
one-endedness of $A$ implies that $A$ must be conjugate into a vertex group
$G_{1}$ of $G$. In particular $K$ is conjugate into $G_{1}$ which implies that
$H$, and then $L$ must also be conjugate into $G_{1}$. As $B$ contains $L$ and
is one-ended, it follows that $B$ is also conjugate into $G_{1}$, which
contradicts the assumption that $G$ has a splitting.

If we assume that $A$ has no two-ended subgroups $D$ with $e(A,D)\geq2$ and
similarly for $B$, then we claim that $\Gamma(G)$ is the graph of groups given
by $G=A\ast_{K}\left[  (K\ast L)\times H\right]  \ast_{L}B$. Assuming this,
then it is clear that neither of the edge groups of the two edges incident to
the commensuriser vertex of $\Gamma(G)$ contains $H$.

To justify the above claim about $\Gamma(G)$, we need to show that if $C$ is a
two-ended subgroup of $G$, then any nontrivial $C$-almost invariant subset of
$G$ is enclosed by the commensuriser vertex of the above graph of groups. This
is easy to show topologically. Pick compact spaces with fundamental groups
$A$, $B$, $K$ and $L$, and use them to form a compact space with fundamental
group $G$. Then consider a covering space with fundamental group $C$.
\end{example}

Next we give a specific example, which puzzled us for many years before we
understood the theory of regular neighbourhoods. This is related to the
problem of unfolding of splittings over two-ended subgroups, which appeared to
make it very difficult to produce a truly canonical algebraic
JSJ-decomposition. Our work in this paper solves this problem by showing how
to enclose all splittings over two-ended subgroups simultaneously.

\begin{example}
Let $A$ and $B$ be finitely generated groups, and let $C$ and $D$ be infinite
cyclic subgroups of $A$ and $B$ respectively. Let $nD$ denote the subgroup of
$D$ of index $n$. Let $G$ denote the group $A\ast_{C=6D}B$, and let
$\sigma_{6}$ denote this splitting of $G$ over $6D$. For $k=1$, $2$ or $3$,
let $A_{k}$ denote $A\ast_{C=6D}kD$. Then we can also express $G$ as
$A_{k}\ast_{kD}B$, for $k=1$, $2$ or $3$. Let $\sigma_{1}$, $\sigma_{2}$ and
$\sigma_{3}$ denote these three splittings of $G$ over $D$, $2D$ and $3D$
respectively. The two splittings $\sigma_{2}$ and $\sigma_{3}$ of $G$ must
have non-zero intersection number, because otherwise they would be compatible
by Theorem \ref{Theorem2.5ofSS}, and it is easy to see that this is
impossible. We claim that the regular neighbourhood of $\sigma_{2}$ and
$\sigma_{3}$ in $G$ is the graph of groups $\Gamma$ with two edges given by
$A\ast_{C=6D}D\ast_{D}B$. This graph has a single $V_{0}$-vertex carrying $D$
and two other vertices which are $V_{1}$-vertices. One way to prove this would
be to check that the conditions in Definition \ref{defnofalgregnbhd} hold.
Certainly the vertex of $\Gamma$ which carries $D$ does enclose each of the
splittings $\sigma_{2}$ and $\sigma_{3}$. To see this for $\sigma_{2}$,
observe that the graph of groups $\Gamma_{2}$ given by $G=A\ast_{C=6D}%
2D\ast_{2D}D\ast_{D}B$ is the required refinement of $\Gamma$, and similarly
for $\sigma_{3}$. Also $\Gamma$ is minimal and the condition on isolated
$V_{0}$-vertices is vacuous, because neither $\sigma_{2}$ nor $\sigma_{3}$ is
isolated. But it is not easy to verify directly that $\Gamma$ satisfies
Condition 2) of the definition. Instead, we will directly consider the
construction of the regular neighbourhood of $\sigma_{2}$ and $\sigma_{3}$ in
$G$, which we gave in section \ref{regnbhds:construction}. This is possible in
this case, because we can directly understand the connections between the
almost invariant sets associated to the four different splittings described above.

For $k=1$, $2$, or $3$, we let $Z_{k}$ denote one of the standard almost
invariant subsets of $G$ associated to the splitting $\sigma_{k}$. We let $Z$
denote one of the standard almost invariant subsets of $G$ associated to the
original splitting $\sigma_{6}$. Finally, let $d$ denote a generator of $D$.
Consider the $G$-tree $T_{2}$ determined by the graph of groups $\Gamma_{2}$.
Let $v$ be a vertex with stabiliser $2D$. There is one edge $s$ incident to
$v$ with stabiliser $2D$ and three other edges incident to $v$ each with
stabiliser $6D$. If $t$ denotes one of these three edges, then the other two
equal $d^{2}t$ and $d^{4}t$. We orient $s$ towards $v$ and $t$ away from $v$,
and pick any $G$-equivariant function $\varphi:G\rightarrow V(T_{2})$. We
choose $Z_{s}=Z_{2}$ and $Z_{t}=Z$. It is now immediate that $Z_{2}=Z\cup
d^{2}Z\cup d^{4}Z$. Similarly, considering $\Gamma_{3}$, shows that
$Z_{3}=Z\cup d^{3}Z$. We claim that the CCC $v_{0}$ of $\overline{E}$ which
contains $\overline{Z_{2}}$ and $\overline{Z_{3}}$ consists precisely of
$\overline{Z_{2}}$, $d^{3}\overline{Z_{2}}$, $\overline{Z_{3}}$,
$d^{2}\overline{Z_{3}}$ and $d^{4}\overline{Z_{3}}$. Clearly these must all
lie in $v_{0}$, and we are claiming that no other translates of $Z_{2}$ cross
$Z_{3}$, and that no other translates of $Z_{3}$ cross $Z_{2}$. Assuming this,
it follows that the stabiliser of $v_{0}$ is simply the group generated by the
stabilisers of $Z_{2}$ and $Z_{3}$, namely $2D$ and $3D$, which is exactly
$D$. Further, we claim that for any element $a$ of $A-D$, the CCC's $v_{0}$
and $av_{0}$ are adjacent in the pretree, and that the analogous statement
holds for any element of $B-D$. It follows that the regular neighbourhood of
$\sigma_{2}$ and $\sigma_{3}$ in $G$ has a single $V_{0}$-vertex with
associated group $D$ and has a $V_{1}$-vertex with associated group $A$ and a
$V_{1}$-vertex with associated group $B$. As these groups generate $G$ and the
regular neighbourhood is a minimal graph of groups, it follows that the
regular neighbourhood must be the graph $\Gamma$ described above.

Now we prove the first claim about the composition of the CCC $v_{0}$ of
$\overline{E}$ which contains $\overline{Z_{2}}$ and $\overline{Z_{3}}$. This
required showing that no translate of $Z_{2}$ by an element of $G-D$ can cross
$Z_{3}$, and that no translate of $Z_{3}$ by an element of $G-D$ can cross
$Z_{2}$. Let $g$ be an element of $G-D$. We need to show that $gZ_{2}$ and
$Z_{3}$ are nested. Considering $T_{2}$ shows that we have either
$gZ_{2}^{(\ast)}\subset Z^{\ast}$ or $gZ_{2}^{(\ast)}\subset kZ^{(\ast)}$,
where $kZ^{(\ast)}$ is not $Z$ or $Z^{\ast}$, so that $k\notin D$. Considering
$T_{3}$ shows that $Z^{\ast}\subset Z_{3}^{\ast}$ and $kZ^{(\ast)}\subset
Z_{3}$, when $k\notin D$, so that $gZ_{2}$ and $Z_{3}$ are nested, as
required. Similarly, $gZ_{3}$ and $Z_{2}$ are nested, for any $g\in G-D$.
Hence the CCC $v_{0}$ contains the five elements claimed. Note that as $v_{0}$
contains $\overline{Z_{2}}$ and $\overline{Z_{3}}$, every CCC is a translate
of $v_{0}$ by some element of $G$.

Next we prove that $av_{0}$ is adjacent to $v_{0}$, for any element $a$ of
$A-D$. If some CCC $gv_{0}$ lies between $v_{0}$ and $av_{0}$, then there is
an element $X$ of $v_{0}$ such that $gX$ lies between $\overline{Z_{2}}$ and
$a\overline{Z_{2}}$. It is easy to see that this is impossible, using the fact
that $a$ stabilises a vertex of $T_{2}$ adjacent to $v$. We prove similarly
that $bv_{0}$ is adjacent to $v_{0}$, for any element $b$ of $B-D$.
\end{example}

Now we come to our general examples. Our first such example is when $G$ is the
fundamental group of an orientable Haken manifold. Recall from the discussion
in section \ref{charsub} that for our purposes we will consider the
submanifold $V^{\prime}(M)$ of $M$ rather than the characteristic submanifold
$V(M)$. The $V_{0}$-vertices of the decomposition $\Gamma$ of the previous
section applied to $G$ essentially correspond to the peripheral components of
$V^{\prime}(M)$. However we get extra $V_{0}$-vertices corresponding to most
of the annulus components in the frontier of the peripheral components of
$V^{\prime}(M)$. In fact, if $S$ is an annulus component of the frontier of a
peripheral component $W$ of $V^{\prime}(M)$, we get an extra $V_{0}$-vertex
corresponding to $S$ except in the case when $W$ is homeomorphic to $S\times
I$. To see this, observe that the peripheral components of $V^{\prime}(M)$
have enough immersions of the annulus (mostly embeddings) in them to make them
cross-connected. Moreover, we showed in \cite{SS3}, that the frontier
components of $V(M)$ induce splittings of $G$ which are $1$-canonical. The
non-peripheral components of $V^{\prime}(M)$ are the same as those of $V(M)$
and do not enclose $\mathbb{Z}$-almost invariant sets, since this would give a
splitting of a Seifert fibre space over $\mathbb{Z}$ relative to its boundary.
Thus the $V_{0}$-vertex groups correspond to the peripheral components of
$V^{\prime}(M)$ together with the extra annuli mentioned above. The whole
submanifold $V^{\prime}(M)$ can be obtained by the methods of section
\ref{JSJforVPCoftworanks}.

We next compare the decompositions of arbitrary finitely presented groups
obtained by Bowditch in section 15 of \cite{B2} with ours. The construction of
Bowditch in \cite{B2} is in terms of axes and does not seem to give the
enclosing properties that we described in the previous sections. Bowditch's
decomposition uses axes and it is not clear whether it is independent of the
axis chosen. It seems that different choices of axes give the same $V_{0}%
$-vertex groups, but the edge groups can be different. We give an example
where the decomposition given in the previous section may differ from that in
\cite{B2}. Start with a one-ended hyperbolic group $K$ which does not have any
splittings over two-ended subgroups and let $L$ be the HNN-extension obtained
by identifying two non-conjugate infinite cyclic subgroups $H_{1}$, $H_{2}$ of
$K$. Let $M$ be the product of a one-ended group with an infinite cyclic group
$H$. Let $G$ be the group obtained by amalgamating $L$ and $M$ along $H_{1}$
and $H$. If we take the axis corresponding to the final decomposition,
Bowditch's construction again yields the decomposition $L\ast_{H_{1}=H}M$.
However in our decomposition, there are two more $V_{0}$-vertices
corresponding to the two splittings over $H_{i}$, $i=1$, $2$. This is more
canonical, since it takes care of all splittings of $G$ over virtually cyclic
groups. Finally, we recall that there are examples when $Comm_{G}(H)$ is not
finitely generated. See Example \ref{example ofCommG(H)isnotfg}. In this case,
it is not at all clear how the edge groups in Bowditch's decomposition and
ours correspond. It is possible that the decompositions obtained by Bowditch
for big enough axes are the same as ours. Even if this is possible, it is not
clear how to prove the enclosing property for almost invariant sets without
going through some work similar to that in this paper.

Finally, we give an example from $3$-manifold topology which motivates some of
our later work.

\begin{example}
\label{canonicaltorusmaynotbealgebraicallycanonical}Let $F_{1}$ and $F_{2}$
denote two compact surfaces each with at least two boundary components. Let
$\Sigma_{i}$ denote $F_{i}\times S^{1}$, let $T_{i}$ denote a boundary
component of $\Sigma_{i}$, and construct a $3$-manifold $M$ from $\Sigma_{1}$
and $\Sigma_{2}$ by gluing $T_{1}$ to $T_{2}$ so that the given fibrations by
circles do not match. Let $T$ denote the torus $\Sigma_{1}\cap\Sigma_{2}$.
Then $T$ is a canonical torus in $M$, and the characteristic submanifold
$V(M)$ of $M$ has two components which are a copy of $\Sigma_{1}$ and a copy
of $\Sigma_{2}$. Let $H$ denote the subgroup of $G=\pi_{1}(M)$ carried by $T$.
Then there is a splitting $\sigma$ of $G$ over $H$ which has non-zero
intersection number with the splitting $\tau$ determined by $T$.

The splitting $\sigma$ is constructed as follows. Let $G_{i}$ denote $\pi
_{1}(\Sigma_{i})$, and let $C_{i}$ denote the subgroup of $G_{i}$ carried by
$T_{i}$. The starting point of our construction is that if $F$ is a compact
surface with at least two boundary components, and if $S$ denotes a boundary
circle of $F$, then $S$ carries an infinite cyclic subgroup of $\pi_{1}(F)$
which is a free factor of $\pi_{1}(F)$. Now it is easy to give a splitting of
$\pi_{1}(F)$ over $\pi_{1}(S)$, and hence a splitting of $G_{i}$ over $C_{i}$.
If each $\pi_{1}(F_{i})$ is free of rank at least $3$, then we can write
$G_{i}=A_{i}\ast_{C_{i}}B_{i}$. If we let $A$ denote the subgroup of $G$
generated by $A_{1}$ and $A_{2}$, i.e. $A=A_{1}\ast_{H}A_{2}$, and define $B$
similarly, then we can express $G$ as $A\ast_{H}B$, and it is easy to see that
this splitting $\sigma$ of $G$ has non-zero intersection number with the
splitting $\tau$ of $G$ determined by $T$. If $\pi_{1}(F_{i})$ has rank $2$,
then we can write $G_{i}=A_{i}\ast_{C_{i}}$ and a similar construction can be made.
\end{example}

The point of this example is the following. The fact that $T$ is topologically
canonical means that any essential annulus or torus in $M$ has zero
intersection number with $T$. But this example shows that the splitting $\tau$
determined by $T$ is not algebraically $2$-canonical. Now the natural next
step after the results of the previous sections would be to attempt to define
an algebraic analogue of the characteristic submanifold of a $3$-manifold to
be the regular neighbourhood of all almost invariant subsets of $G$ which are
over subgroups isomorphic to $\mathbb{Z}$ or to $\mathbb{Z}\times\mathbb{Z}$.
Suppose that this can be done and consider the case when $G$ is the
fundamental group of the manifold $M$ in the above example. Let $\Gamma$
denote the regular neighbourhood. The fact that $\tau$ crosses another
splitting over $\mathbb{Z}\times\mathbb{Z}$ implies that it cannot be an edge
splitting of $\Gamma$, so that clearly $\Gamma$ would not be the same as the
topological JSJ-decomposition of $M$. In fact, we do not know whether there is
such a regular neighbourhood. However, our results in \cite{SS3} imply that
the topological JSJ-decomposition is an algebraic regular neighbourhood of all
almost invariant subsets over $\mathbb{Z}$ and of all $1$-canonical almost
invariant subsets over $\mathbb{Z}\times\mathbb{Z}$, and this is what we will
generalise in later sections.

\section{Canonical decompositions over VPC groups of a given
rank\label{JSJforVPCofgivenrank}}

As stated at the beginning of section \ref{coendswhencommensuriserissmall},
the analogues of the results of sections \ref{coendswhencommensuriserissmall}
to \ref{JSJforlargecommensurisers} go through for almost invariant sets over
VPC groups of length $n$ assuming that $G$ does not have nontrivial almost
invariant sets over VPC groups of length $<n$ (the analogue of Proposition
\ref{twocoends} is Proposition \ref{twocoends2}, which we prove in the next
section). Note that Theorem \ref{D-Rollertheorem}, which is due to Dunwoody
and Roller, implies that the condition that $G$ does not have nontrivial
almost invariant sets over VPC groups of length $<n$ is equivalent to the
condition that $G$ does not split over such a subgroup. We will use these two
conditions interchangeably.

We will need the following definitions.

\begin{definition}
Let $\Gamma$ be a minimal graph of groups decomposition of a group $G$. A
vertex $v$ of $\Gamma$ is of \textsl{VPC-by-Fuchsian type} if $G(v)$ is a
VPC-by-Fuchsian group, where the Fuchsian group is not finite nor two-ended,
and there is exactly one edge of $\Gamma$ which is incident to $v$ for each
peripheral subgroup $K$ of $G(v)$ and this edge carries $K$. If the length of
the normal VPC subgroup of $G(v)$ is $n$, we will say that $G(v)$ is of
\textsl{length} $n$.
\end{definition}

Note that if $G=G(v)$, then the Fuchsian quotient group corresponds to a
closed orbifold.

\begin{definition}
\label{defnofn-simple}Let $\Gamma$ be a minimal graph of groups decomposition
of a group $G$. A vertex $v$ of $\Gamma$ is $n$\textsl{-simple}, if whenever
$X$ is a nontrivial almost invariant subset of $G$ over a VPC subgroup of
length at most $n$ such that $X$ is enclosed by $v$, then $X$ is associated to
an edge splitting of $\Gamma$.
\end{definition}

The following are the results which we obtain. To prove these results, we will
need generalisations of the results in section
\ref{coendswhencommensuriserislarge} on the Boolean algebra $B(H)$ to the case
when $H$ is VPC\ of rank $n$. Such results can be proved by the same methods.
This will then allow us to handle the $V_{0}$-vertices of large commensuriser
type in our regular neighbourhood. As in sections
\ref{JSJforsmallcommensurisers} and \ref{JSJforlargecommensurisers}, we first
state the basic existence result for the appropriate regular neighbourhood and
then list its properties. Recall that we use the term length instead of Hirsch
length, for brevity.

\begin{theorem}
\label{JSJforVPCofgivenrankexists}Let $G$ be a one-ended, finitely presented
group which does not split over VPC groups of length $<n$, and let
$\mathcal{F}$ denote the collection of equivalence classes of all nontrivial
almost invariant subsets of $G$ which are over VPC groups of length $n$.

Then the regular neighbourhood construction of section
\ref{regnbhds:construction} works and yields a regular neighbourhood
$\Gamma(\mathcal{F}:G)$. Each $V_{0}$-vertex $v$ of $\Gamma(\mathcal{F}:G)$
satisfies one of the following conditions:

\begin{enumerate}
\item $v$ is isolated, so that $G(v)$ is VPC of length $n$.

\item $G(v)$ is of VPC-by-Fuchsian type of length $(n-1)$.

\item $G(v)$ is the full commensuriser $Comm_{G}(H)$ for some VPC subgroup $H$
of length $n$, such that $e(G,H)\geq2$.
\end{enumerate}

$\Gamma(\mathcal{F}:G)$ consists of a single vertex if and only if
$\mathcal{F}$ is empty, or $G$ itself satisfies one of the above three conditions.
\end{theorem}

Now we list the properties of $\Gamma(\mathcal{F}:G)$.

\begin{theorem}
\label{propertiesofJSJforVPCofgivenrank}Let $G$ be a one-ended, finitely
presented group which does not split over VPC groups of length $<n$, and let
$\mathcal{F}$ denote the collection of equivalence classes of all nontrivial
almost invariant subsets of $G$ which are over VPC subgroups of length $n$.

Then the regular neighbourhood $\Gamma(\mathcal{F}:G)$ is a minimal bipartite
graph of groups decomposition of $G$ with the following properties:

\begin{enumerate}
\item each $V_{0}$-vertex $v$ of $\Gamma(\mathcal{F}:G)$ satisfies one of the
following conditions:

\begin{enumerate}
\item $v$ is isolated, so that $G(v)$ is VPC of length $n$.

\item $G(v)$ is of VPC-by-Fuchsian type of length $(n-1)$.

\item $G(v)$ is the full commensuriser $Comm_{G}(H)$ for some VPC subgroup $H$
of length $n$, such that $e(G,H)\geq2$.
\end{enumerate}

Further, if $H$ is a VPC subgroup of length $n$ such that $e(G,H)\geq2$, and
if $H$ has large commensuriser, then $\Gamma(\mathcal{F}:G)$ will have a
$V_{0}$-vertex $v$ such that $G(v)=Comm_{G}(H)$.

\item If an edge of $\Gamma(\mathcal{F}:G)$ is incident to a $V_{0}$-vertex of
type a) or b) above, then it carries a VPC group of length $n$.

\item any representative of an element of $\mathcal{F}$ is enclosed by some
$V_{0}$-vertex of $\Gamma(\mathcal{F}:G)$, and each $V_{0}$-vertex of
$\Gamma(\mathcal{F}:G)$ encloses such a subset of $G$. In particular, any
splitting of $G$ over a VPC subgroup of length $n$ is enclosed by some $V_{0}%
$-vertex of $\Gamma(\mathcal{F}:G)$.

\item if $X$ is an almost invariant subset of $G$ over a finitely generated
subgroup $H$, and if $X$ does not cross any element of $\mathcal{F}$, then $X$
is enclosed by a $V_{1}$-vertex of $\Gamma(\mathcal{F}:G)$.

\item if $X$ is a $H$-almost invariant subset of $G$ associated to a splitting
of $G$ over $H$, and if $X$ does not cross any element of $\mathcal{F}$, then
$X$ is enclosed by a $V_{1}$-vertex of $\Gamma(\mathcal{F}:G)$.

\item the $V_{1}$-vertex groups of $\Gamma(\mathcal{F}:G)$ are $n$-simple.
(See Definition \ref{defnofn-simple}.)

\item If $\Gamma_{1}$ and $\Gamma_{2}$ are minimal bipartite graphs of groups
structures for $G$ which satisfy conditions 3 and 5 above, then they are
isomorphic provided there is a one-to-one correspondence between their
isolated $V_{0}$-vertices.

\item The graph of groups $\Gamma(\mathcal{F}:G)$ is invariant under the
automorphisms of $G$.

\item The $n$-canonical splittings of $G$ over a VPC subgroup of length $n$
are precisely those edge splittings of $\Gamma(\mathcal{F}:G)$ which are over
such a subgroup. This includes, but need not be limited to, all those edges of
$\Gamma(\mathcal{F}:G)$ which are incident to $V_{0}$-vertices whose
associated groups are of type a) or b) above.
\end{enumerate}
\end{theorem}

As in section \ref{JSJforlargecommensurisersexist}, it follows that one can
also form a regular neighbourhood of only those almost invariant subsets which
are associated to splittings. This is the result we obtain.

\begin{theorem}
Let $G$ be a one-ended, finitely presented group which does not split over VPC
groups of length $<n$, and let $\mathcal{S}$ denote the collection of
equivalence classes of all almost invariant subsets which are associated to
splittings of $G$ over VPC groups of length $n$.

Then the regular neighbourhood construction of section
\ref{regnbhds:construction} works and yields a regular neighbourhood
$\Gamma(\mathcal{S}:G)$. Each $V_{0}$-vertex $v$ of $\Gamma(\mathcal{S}:G)$
satisfies one of the following conditions:

\begin{enumerate}
\item $v$ is isolated, so that $G(v)$ is VPC of length $n$.

\item $G(v)$ is of VPC-by-Fuchsian type of length $(n-1)$.

\item $G(v)$ contains a VPC subgroup $H$ of length $n$, which it
commensurises, such that $e(G,H)\geq2$.
\end{enumerate}

If $\Gamma(\mathcal{S}:G)$ consists of a single vertex, then either
$\mathcal{S}$ is empty, or $G$ itself satisfies one of conditions 2) or 3) above.
\end{theorem}

\section{Canonical splittings over VPC groups of two successive
ranks\label{JSJforVPCoftworanks}}

As usual, let $G$ denote a one-ended, finitely presented group. In this
section, we examine the problem of enclosing almost invariant subsets of $G$
over VPC subgroups of two successive lengths $n$ and $n+1$. Note that VPC
groups of length at most $2$ are virtually abelian. As discussed at the end of
section \ref{examples} in the case when $n=1$, we should not expect to be able
to enclose all almost invariant subsets of $G$ over VPC groups of lengths $1$
and $2$. The characteristic submanifold $V(M)$ of a $3$-manifold $M$ can be
regarded as a regular neighbourhood of all the essential annuli in $M$ and of
all the essential tori in $M$ which have intersection number zero with any
essential annulus. Thus, in the case $n=1$, we will show that one can enclose
all almost invariant subsets of $G$ over two-ended subgroups together with all
$1$-canonical almost invariant subsets over VPC subgroups of length $2$. This
corresponds to the classical JSJ-decomposition. For general values of $n$, we
will need to assume that $G$ does not have any nontrivial almost invariant
subsets over VPC groups of length $<n$, and we will then show that one can
enclose all almost invariant subsets of $G$ over VPC subgroups of length $n$
together with all $n$-canonical almost invariant subsets over VPC subgroups of
length $n+1$.

We will frequently use the following known facts about VPC groups (see, for
example, \cite{KR} and \cite{D-Swenson}).

\begin{lemma}
\label{VPCofcolengthonehaslargenormaliser}Let $K$ be a VPC group of length
$n+1$, with a subgroup $L$ of length $n$. Then the number of coends of $L$ in
$K$ is $2$ and moreover there is a subgroup $L^{\prime}$ of finite index in
$L$ such that $L^{\prime}$ has infinite index in its normaliser.
\end{lemma}

\begin{proof}
The number of coends of $L$ in $K$ is $2$ since both are (virtual)
Poincar\'{e} duality groups. Now, without loss of generality, we can assume
that they are both Poincar\'{e} duality groups. In order to find the required
subgroup $L^{\prime}$, we use the fact that $L$ can be separated from elements
of $K-L$. If $K$ and $L$ are both orientable, then the number of ends of the
pair $(K,L)$ is $2$. It follows from Scott's Theorem 4.1 in \cite{Scott:Ends}
on ends of pairs of groups that a subgroup $K^{\prime}$ of finite index in $K$
splits over $L$. Now the fact that the number of ends of the pair $(K^{\prime
},L)$ is $2$ implies that $L$ is normal in $K^{\prime}$ and $L\backslash
K^{\prime}$ is either infinite cyclic or infinite dihedral. In particular, $L$
has infinite index in its normaliser, so that $L^{\prime}$ equals $L$ in this
case. If $L$ is not orientable, we let $L^{\prime}$ denote the orientable
subgroup of index two, and the preceding argument shows that $L^{\prime}$ has
infinite index in its normaliser, as required.
\end{proof}

The same argument shows the following result.

\begin{lemma}
\label{VPCofcolengthoneintwogroupshaslargenormaliser}Suppose that $L$ is a VPC
group of length $n$ and that $L$ is a subgroup of two VPC groups $K_{1}$,
$K_{2}$ of length $n+1$. Then there is a subgroup $L^{\prime}$ of finite index
in $L$ and subgroups $L_{1}$ and $L_{2}$ of finite index in $K_{1}$ and
$K_{2}$ respectively so that $L^{\prime}\backslash L_{1}$ and $L^{\prime
}\backslash L_{2}$ are both infinite cyclic.
\end{lemma}

Now consider a one-ended, finitely presented group $G$ which does not have any
nontrivial almost invariant subsets over VPC groups of length $<n$, and apply
Theorem \ref{JSJforVPCofgivenrankexists} to obtain a regular neighbourhood
$\Gamma_{n}$ of $\mathcal{F}_{n}$, the equivalence classes of all almost
invariant subsets of $G$ over VPC subgroups of length $n$. Let $P_{n}$ denote
the corresponding pretree , and $T_{n}$ the corresponding tree. A
$n$-canonical almost invariant subset of $G$ over a VPC subgroup of length
$(n+1)$ does not cross any element of $\mathcal{F}_{n}$, and so must be
enclosed by some $V_{1}$-vertex of $T_{n}$, by part 1 of Proposition
\ref{XdoesnotcrossanyXiimpliesXisenclosedbyaV1vertex}. Further if two such
subsets cross, they must be enclosed by the same $V_{1}$-vertex of $T_{n}$. In
order to obtain the regular neighbourhood for which we are looking, we will
refine $\Gamma_{n}$ by splitting at some of the $V_{1}$-vertices. We will need
to consider crossing of $n$-canonical subsets of $G$ over VPC subgroups of
length $(n+1)$. First we consider strong crossing of such subsets. The
following proposition and proof are suggested by the symmetry of crossings
proved in \cite{SS2}.

\begin{proposition}
\label{twocoends2} Let $G$ be a one-ended, finitely generated group, and let
$X$ and $Y$ be $n$-canonical nontrivial almost invariant subsets of $G$ over
VPC subgroups $H$ and $K$ of length $(n+1)$. If $X$ crosses $Y$ strongly, then
$Y$ crosses $X$ strongly and the number of coends in $G$ of both $H$ and $K$
is $2$.
\end{proposition}

\begin{proof}
Let $\Delta$ be the Cayley graph of $G$ with respect to some finite system of
generators. As $H\backslash\delta X$ is finite, there is a finite connected
subcomplex of $H\backslash\Delta$ which contains $H\backslash\delta X$ and
carries the group $H$. The pre-image of this subcomplex in $\Delta$ is a
connected subcomplex $C$ which contains $\delta X$ and is $H$-finite.
Similarly, there is a connected subcomplex $D$ of $\Delta$ which contains
$\delta Y$ and is $K$-finite. Since $X$ crosses $Y$ strongly, there are points
of $\delta X$, and hence of $C$, in $Y$ and in $Y^{\ast}$ which are outside
any $d$-neighbourhood of $\delta Y$. Thus the projection of $C$ into
$K\backslash\Delta$ has at least two ends, so that $e(H,H\cap K)\geq2$. It
follows that $H\cap K$ is a VPC group of length $n$.

Lemma \ref{VPCofcolengthoneintwogroupshaslargenormaliser} tells us that we can
find a subgroup $L$ of finite index in $H\cap K$ so that $L$ has infinite
index in its normalisers in both $H$ and $K$. We claim that $L\backslash
\Delta$ has one end. For suppose that there is a nontrivial $L$-almost
invariant subset $Z$ of $G$. Since $X$ and $Y$ are $n$-canonical, we have one
of the four inequalities $Z^{(\ast)}\leq X^{(\ast)}$ and one of the four
inequalities $Z^{(\ast)}\leq Y^{(\ast)}$. By appropriately replacing some of
$X$, $Y$ and $Z$ by their complements, we may arrange that $Z\leq X$ and
$Z\leq Y$, which implies that $Z\leq X\cap Y$. As $L$ has infinite index in
its normaliser $N_{H}(L)$ in $H$, there is an infinite cyclic subgroup $J$ of
$N_{H}(L)/L$ which acts freely on $L\backslash\Delta$. As $X$ crosses $Y$
strongly, the orbit of any point of $L\backslash\Delta$ under the action of
$J$ contains points on each side of $L\backslash\delta Y$ which are
arbitrarily far from $L\backslash\delta Y$. As $L\backslash\delta Z$ is
finite, there is an element $j$ of $J$ such that $j(L\backslash\delta Z)$ is
contained in $L\backslash Y^{\ast}$. Thus there is an element $h$ of
$N_{H}(L)$ such that $h\delta Z$ is contained in $Y^{\ast}$. Thus we have one
of the inclusions $hZ^{(\ast)}\subset Y^{\ast}$. Suppose that $hZ\subset
Y^{\ast}$. As $Z\leq X$ and $hX=X$, we have $hZ\leq X\cap Y^{\ast}$. As $h$
normalises $L$, it follows that $Z\cup hZ$ is a nontrivial $L$-almost
invariant subset of $G$. As $Z\leq Y$ and $hZ\leq Y^{\ast}$, this set crosses
$Y$ which contradicts our assumption that $Y$ is $n$-canonical. If $hZ^{\ast
}\subset X\cap Y^{\ast}$, we consider $Z\cup hZ^{\ast}$ instead, to complete
the proof of the claim that $L\backslash\Delta$ has one end.

Further, both the covers $L\backslash\Delta\rightarrow H\backslash\Delta$ and
$L\backslash\Delta\rightarrow K\backslash\Delta$ have two-ended covering
groups. Let $C_{L}$ and $D_{L}$ be the images of $C$ and $D$ respectively in
$L\backslash\Delta$, so that both $C_{L}$ and $D_{L}$ are two-ended. Then also
both $(L\backslash\Delta)-C_{L}$ and $(L\backslash\Delta)-D_{L}$ have at least
two components and $C_{L}$ has two infinite pieces in different components of
$(L\backslash\Delta)-D_{L}$.

Our hypothesis implies that $H\backslash\Delta$ has more than one end. Now we
claim that $H\backslash\Delta$ has only two ends. For suppose that
$H\backslash\Delta$ has more than two ends. Since $e(H,H\cap K)=2$, the image
of $D$ in $H\backslash\Delta$ has two ends. Thus an end of $H\backslash\Delta$
is free of the image of $D$. Choose a compact set separating this end from the
image of $D$ and let $M$ be an infinite component of the complement of this
compact set which does not intersect the image of $D$. Let $N$ be a component
of the pre-image of $M$ in $L\backslash\Delta$, so that $N$ is disjoint from
$D_{L}$. Since $L\backslash\Delta$ has only one end, the coboundary of $\delta
N$ is not finite. Thus the stabiliser of $N$ projects to an infinite cyclic
group in the covering transformation group of the cover $L\backslash
\Delta\rightarrow H\backslash\Delta$. Thus $\delta N$ and $C_{L}$ are in a
bounded neighbourhood of each other. Recall that as $X$ crosses $Y$ strongly,
$C_{L}$ has points on both sides of $D_{L}$ which are arbitrarily far from
$D_{L}$. Thus there are points of $\delta N$, and hence of $N$, on both sides
of $D_{L}$ which are arbitrarily far from $D_{L}$. This is a contradiction
since $D_{L}$ and $N$ are disjoint and are both connected. It follows that the
number of ends of $H\backslash\Delta$ is $2$, as claimed.

Now $Y$ crosses $X$ since crossing is symmetric, and Proposition 3.7 of
\cite{SS} shows that $Y$ crosses $X$ strongly. We can repeat the above
argument for subgroups of finite index in $H$ and $K$ to conclude that both
$H$ and $K$ have two coends in $G$, as required.
\end{proof}

We now consider $n$-canonical subsets $X$ and $Y$ of $G$ which cross weakly.
Thus $X$ and $Y$ are enclosed by a $V_{1}$-vertex $v$ of $T_{n}$. Let $H$ and
$K$ denote the stabilisers of $X$ and $Y$ respectively. If $H$ and $K$ are
commensurable, their commensuriser will play an important role, so we will
need to compare their commensurisers in $Stab(v)$ and in $G$. We will show in
Proposition \ref{commensuriser2} that these two groups coincide in the cases
in which we are interested.

\begin{proposition}
\label{commensurable2}Let $G$ be a one-ended, finitely generated group without
nontrivial almost invariant subsets over VPC groups of length $<n$. Let $H$
and $K$ be VPC\ subgroups of $G$ of length $(n+1)$, and let $X$ and $Y$ be
nontrivial $n$-canonical subsets of $G$ over $H$ and $K$ respectively. Suppose
that $X$ crosses $Y$ weakly. Then $H$ and $K$ are commensurable.
\end{proposition}

\begin{proof}
By Proposition \ref{twocoends2}, we know that $X$ and $Y$ cross each other
weakly. Thus one of $\delta Y\cap X^{(\ast)}$ is $H$-finite and one of $\delta
X\cap Y^{(\ast)}$ is $K$-finite. By replacing $X$ and $Y$ by their complements
as needed, we can arrange that $\delta Y\cap X$ is $H$-finite and $\delta
X\cap Y$ is $K$-finite. As $\delta X$ is $H$-finite and $\delta Y$ is
$K$-finite, it follows that $\delta Y\cap X$ and $\delta X\cap Y$ are both
$H$-finite and $K$-finite. Thus they are $(H\cap K)$-finite. Now consider
$X\cap Y$. Every edge in $\delta(X\cap Y)$ lies in $\delta Y\cap X$ or $\delta
X\cap Y$. It follows that $\delta(X\cap Y)$ is also $(H\cap K)$-finite, so
that $X\cap Y$ is $(H\cap K)$-almost invariant.

If $H$ and $K$ are not commensurable, then $H\cap K$ has infinite index in
both $H$ and $K$. In particular, $H\cap K$ has length $\leq n$. Suppose first
that $H\cap K$ has length $<n$. As $G$ has no almost invariant subsets over
VPC subgroups of length $<n$, it follows that $X\cap Y$ is $(H\cap K)$-finite,
and hence $H$-finite, which contradicts our hypothesis that $X$ and $Y$ cross.
Now suppose that $H\cap K$ has length $n$, so that $e(H,H\cap K)=2$. As in the
proof of Proposition \ref{twocoends2}, we can translate $X\cap Y$ by an
element $h$ of $H$ such that $h$ commensurises $H\cap K$ and $h(X\cap Y)$ is
contained in $X\cap Y^{\ast}$. Thus $(X\cap Y)\cup h(X\cap Y)$ is an almost
invariant set over a subgroup commensurable with $H\cap K$ which crosses $Y$.
This contradicts the assumption that $Y$ is $n$-canonical. Hence $H$ and $K$
are commensurable as required.
\end{proof}

Note that the proof of Proposition \ref{twocoends2} did not use the hypothesis
that $G$ has no nontrivial almost invariant subsets over VPC groups of length
$<n$. But the proof of Proposition \ref{commensurable2} used this hypothesis
in an essential way. It was used to exclude the case when $H\cap K$ has length
$<n$. If this occurs, we might not be able to construct an almost invariant
set which crosses $Y$, because there need not be any suitable elements of $G$
which commensurise $H\cap K$.

\begin{remark}
We note that if $X$ and $Y$ are $n$-canonical and their stabilisers are
commensurable, then $X\cap Y$, $X+Y$, $X\cup Y$ are again $n$-canonical.
\end{remark}

These results show that as in Proposition \ref{crossingsareallstrongorallweak}%
, we have:

\begin{proposition}
Let $G$ be a one-ended finitely generated group and let $\{X_{\lambda
}\}_{\lambda\in\Lambda}$ be a family of nontrivial almost invariant subsets
over VPC subgroups of length $n$ and of nontrivial $n$-canonical almost
invariant subsets over VPC subgroups of length $n+1$. As usual, let $E$ denote
the set of all translates of the $X_{\lambda}$'s and their complements. Form
the pretree $P$ of cross-connected components (CCC's) of $\overline{E}$ as in
the construction of regular neighbourhoods in section
\ref{regnbhds:construction}. Then the following statements hold:

\begin{enumerate}
\item The crossings in a CCC of $\overline{E}$ are either all strong or are
all weak.

\item In a CCC with all crossings weak, the stabilisers of the corresponding
elements of $E$ are all commensurable. In a CCC with all crossings strong, the
stabilisers of the corresponding elements of $E$ have $2$ coends in $G$.
\end{enumerate}
\end{proposition}

Recall that $T_{n}$ is the universal covering $G$-tree of the regular
neighbourhood $\Gamma_{n}$ of $\mathcal{F}_{n}$, the equivalence classes of
all almost invariant subsets of $G$ over VPC subgroups of length $n$. The last
result we need to produce the required regular neighbourhood is that if $X$ is
a $n$-canonical $H$-almost invariant subset of $G$ which is enclosed by a
$V_{1}$-vertex $v$ of $T_{n}$, then the commensurisers of $H$ in $Stab(v)$ and
in $G$ are equal. This is the last part of the following proposition.

\begin{proposition}
\label{commensuriser2}Let $G$ be a one-ended, finitely presented group without
nontrivial almost invariant subsets over VPC groups of length $<n$. Let $X$
and $Y$ be $n$-canonical subsets of $G$ over VPC groups of length $n+1$, and
let $H$ and $K$ denote the stabilisers of $X$ and $Y$ respectively.

\begin{enumerate}
\item Then $X$ is enclosed by a unique $V_{1}$-vertex $v_{X}$ of $T_{n}$.

\item If $H$ and $K$ are commensurable, then $v_{X}$ and $v_{Y}$ are equal.

\item Let $v$ denote the $V_{1}$-vertex of $T_{n}$ which encloses $X$. Then
$Comm_{Stab(v)}(H)=Comm_{G}(H)$.
\end{enumerate}
\end{proposition}

\begin{proof}
Recall that a $n$-canonical almost invariant subset of $G$ over a VPC subgroup
of length $(n+1)$ must be enclosed by some $V_{1}$-vertex of $T_{n}$, by part
1 of Proposition \ref{XdoesnotcrossanyXiimpliesXisenclosedbyaV1vertex}. Now
part 1) is the special case of part 2) obtained when $Y$ equals $X$, so we
will prove part 2). Let $v_{1}$ and $v_{2}$ be $V_{1}$-vertices of $T_{n}$
which enclose $X\;$and $Y$ respectively, and suppose that $v_{1}$ and $v_{2}$
are distinct. Then there is a $V_{0}$-vertex $v$ separating $v_{1}$ from
$v_{2}$. This implies that there is an almost invariant subset $Z$ of $G$ over
a VPC group $L$ of length $n$ such that $Z$ is enclosed by $v$ and
$\overline{Z}$ lies between $\overline{X}$ and $\overline{Y}$. Note that the
stabilisers of $X$ and $Y$ are VPC groups of length $(n+1)$ and the stabiliser
of $Z$ is a length $n$ group. Since $H$ and $K$ are commensurable, $\delta X$
and $\delta Y$ lie in a bounded neighbourhood of each other. In $L\backslash
\Delta$, where $\Delta$ is the Cayley graph of $G$, the images of $\delta X$
and of $\delta Y$ are non-compact and lie essentially on different sides of
the image $L\backslash\delta Z$ which is compact. This contradicts the fact
that $\delta X$ and $\delta Y$ lie in a bounded neighbourhood of each other.
This contradiction shows that $v_{1}$ must equal $v_{2}$, as required.

For the third part of the proposition, let $g$ be an element of $Comm_{G}(H)$.
Then $gX$ is an almost invariant subset of $G$ over $H^{g}$, which is
commensurable with $H$, and $gX$ is enclosed by the $V_{1}$-vertex $gv$. Part
2) of the proposition tells us that $gv=v$, showing that we must have
$Comm_{Stab(v)}(H)=Comm_{G}(H)$, as required.
\end{proof}

With propositions \ref{twocoends2} and \ref{commensuriser2} available, we can
obtain our analogue of the JSJ-decomposition by considering almost invariant
sets over VPC groups of two successive lengths. As usual we state first the
existence result.

\begin{theorem}
Let $G$ be a one-ended, finitely presented group which does not split over VPC
subgroups of length $<n$, and let $\mathcal{F}_{n,n+1}$ denote the collection
of equivalence classes of all nontrivial almost invariant subsets of $G$ which
are over a VPC subgroup of length $n$, together with the equivalence classes
of all nontrivial $n$-canonical almost invariant subsets of $G$ which are over
a VPC subgroup of length $n+1$.

Then the regular neighbourhood construction of section
\ref{regnbhds:construction} works and yields a regular neighbourhood
$\Gamma_{n,n+1}=\Gamma(\mathcal{F}_{n,n+1}:G)$.

Each $V_{0}$-vertex $v$ of $\Gamma_{n,n+1}$ satisfies one of the following conditions:

\begin{enumerate}
\item $v$ is isolated, so that $G(v)$ is VPC of length $n$ or $n+1$.

\item $v$ is of VPC-by-Fuchsian type of length $n-1$ or $n$.

\item $G(v)$ is the full commensuriser $Comm_{G}(H)$ for some VPC subgroup $H$
of length $n$ or $n+1$, such that $e(G,H)\geq2$.
\end{enumerate}

$\Gamma_{n,n+1}$ consists of a single vertex if and only if $\mathcal{F}%
_{n,n+1}$ is empty, or $G$ itself satisfies one of the above three conditions.
\end{theorem}

Now we list the properties of the decomposition $\Gamma_{n,n+1}$.

\begin{theorem}
Let $G$ be a one-ended, finitely presented group which does not split over VPC
subgroups of length $<n$, and let $\mathcal{F}_{n,n+1}$ denote the collection
of equivalence classes of all nontrivial almost invariant subsets of $G$ which
are over a VPC subgroup of length $n$, together with the equivalence classes
of all nontrivial $n$-canonical almost invariant subsets of $G$ which are over
a VPC subgroup of length $n+1$.

Then the regular neighbourhood $\Gamma_{n,n+1}=\Gamma(\mathcal{F}_{n,n+1}:G)$
is a minimal bipartite graph of groups decomposition of $G$ with the following properties:

\begin{enumerate}
\item each $V_{0}$-vertex $v$ of $\Gamma_{n,n+1}$ satisfies one of the
following conditions:

\begin{enumerate}
\item $v$ is isolated, so that $G(v)$ is VPC of length $n$ or $n+1$.

\item $v$ is of VPC-by-Fuchsian type of length $n-1$ or $n$.

\item $G(v)$ is the full commensuriser $Comm_{G}(H)$ for some VPC subgroup $H$
of length $n$ or $n+1$, such that $e(G,H)\geq2$.
\end{enumerate}

Further if $H$ is a VPC subgroup of length $n$ such that $e(G,H)\geq2$, and if
$H$ has large commensuriser, then $\Gamma_{n,n+1}$ will have a $V_{0}$-vertex
$v$ such that $G(v)=Comm_{G}(H)$. The same holds if $H$ is VPC of length
$n+1$, so long as there exists a $n$-canonical almost invariant subset of $G$
over $H$.

\item If an edge of $\Gamma$ is incident to a $V_{0}$-vertex of type a) or b)
above, then it carries a VPC group of length $n$ or $n+1$, as appropriate.

\item any representative of an element of $\mathcal{F}_{n,n+1}$ is enclosed by
some $V_{0}$-vertex of $\Gamma_{n,n+1}$, and each $V_{0}$-vertex of
$\Gamma_{n,n+1}$ encloses such a subset of $G$. In particular, any splitting
of $G$ over a VPC subgroup of length $n$, and any $n$-canonical splitting of
$G$ over a VPC subgroup of length $n+1$ is enclosed by some $V_{0}$-vertex of
$\Gamma_{n,n+1}$.

\item if $X$ is an almost invariant subset of $G$ over a finitely generated
subgroup $H$, and if $X$ does not cross any element of $\mathcal{F}_{n,n+1}$,
then $X$ is enclosed by a $V_{1}$-vertex of $\Gamma_{n,n+1}$.

\item if $X$ is a $H$-almost invariant subset of $G$ associated to a splitting
of $G$ over $H$, and if $X$ does not cross any element of $\mathcal{F}%
_{n,n+1}$, then $X$ is enclosed by a $V_{1}$-vertex of $\Gamma_{n,n+1}$.

\item the $V_{1}$-vertex groups of $\Gamma_{n,n+1}$ are $(n+1)$-simple. In
particular, $\Gamma_{n,n+1}$ cannot be further refined by splitting at a
$V_{1}$-vertex along a VPC group of length $\leq(n+1)$.

\item If $\Gamma_{1}$ and $\Gamma_{2}$ are minimal bipartite graphs of groups
structures for $G$ which satisfy conditions 3 and 5 above, they are isomorphic
provided there is a one-to-one correspondence between their isolated $V_{0}$-vertices.

\item The graph of groups $\Gamma_{n,n+1}$ is invariant under the
automorphisms of $G$.

\item For $k=n$ or $n+1$, the $k$-canonical splittings of $G$ over a VPC
subgroup of length $k$ are precisely those edge splittings of $\Gamma_{n,n+1}$
which are over such a subgroup. This includes, but need not be limited to, all
those edges of $\Gamma_{n,n+1}$ which are incident to $V_{0}$-vertices whose
associated groups are of types a) or b) above.
\end{enumerate}
\end{theorem}

\begin{proof}
The only new point which arises is in the proof of part 6). Suppose that $X$
is a nontrivial almost invariant subset of $G$ over a VPC subgroup $H$ of
length at most $(n+1)$ such that $X$ is enclosed by a $V_{1}$-vertex of
$\Gamma_{n,n+1}$. This implies that $X$ crosses no element of $\mathcal{F}%
_{n,n+1}$, so that $X$ crosses no almost invariant subset of $G$ over a VPC
subgroup of length at most $n$. It follows that $X$ represents an element of
$\mathcal{F}_{n,n+1}$, and so is enclosed by some $V_{0}$-vertex of
$\Gamma_{n,n+1}$. It follows from part 8) of Lemma
\ref{somefactsaboutenclosing} that $X$ is associated to an edge splitting of
$\Gamma_{n,n+1}$, as required.
\end{proof}

If we specialise to the case $n=1$, and apply this result to the fundamental
group of a Haken $3$-manifold $M$, then the $V_{0}$-vertices of $\Gamma_{1,2}$
essentially correspond to the components of the submanifold $V^{\prime}(M)$,
which we discussed in section \ref{charsub}. The only difference is that
$\Gamma_{1,2}$ has extra $V_{0}$-vertices corresponding to most of the
components of the frontier of $V^{\prime}(M)$. In fact, if $S$ is a component
of the frontier of a component $W$ of $V^{\prime}(M)$, we get an extra $V_{0}%
$-vertex corresponding to $S$ except in the case when $W$ is homeomorphic to
$S\times I$. To see this, observe that the peripheral components of
$V^{\prime}(M)$ have enough immersions of the annulus to make them
cross-connected, and the interior components of $V^{\prime}(M)$ have enough
immersions of the torus to make them cross-connected. Moreover, we showed in
\cite{SS3}, that the frontier components of $V(M)$ induce splittings of $G$
which are all$\mathbb{\ }1$-canonical. This is similar to the discussion in
section \ref{examples} for the case of the canonical decomposition obtained in
section \ref{JSJforlargecommensurisers}.

We already saw from Example \ref{canonicaltorusmaynotbealgebraicallycanonical}
that if there is a regular neighbourhood of \emph{all} the almost invariant
subsets of $G$ which are over VPC subgroups of length $1$ or $2$, then it
cannot be a refinement of $\Gamma_{1}$. Thus there may be almost invariant
subsets of $G$ which are over VPC subgroups of length $2$ but are not enclosed
by any $V_{0}$-vertex of $\Gamma_{1,2}$. However, we prove below that the
stabiliser of such an almost invariant subset of $G$ must be `almost'
conjugate into a $V_{0}$-vertex of $\Gamma_{1,2}$.

In order to discuss $V_{0}$-vertices, it will be helpful to introduce some new
language. For the graph of groups $\Gamma_{n,n+1}$, there is a natural idea of
the level of a $V_{0}$-vertex $v$. If $v$ appears only after refining
$\Gamma_{n}$, then $v$ has level $n+1$. Otherwise $v$ has level $n$. In this
second case, it is natural to think that $v$ belongs to $\Gamma_{n}$, in some
sense, but as there is no map from $\Gamma_{n}$ to $\Gamma_{n,n+1}$, this is
not very precise. A more accurate way to describe the level of a $V_{0}%
$-vertex is to use the projection map from $\Gamma_{n,n+1}$ to $\Gamma_{n}$,
which is part of the definition of a refinement. This map sends each edge of
$\Gamma_{n,n+1}$ either to an edge or to a $V_{1}$-vertex of $\Gamma_{n}$.
Then a $V_{0}$-vertex of $\Gamma_{n,n+1}$ has level $n+1$ if it is sent to a
$V_{1}$-vertex, and has level $n$ otherwise.

\begin{proposition}
\label{groupsarealmostenclosed}Let $G$ be a one-ended, finitely presented
group which does not split over VPC groups of length $<n$, and let
$\Gamma_{n,n+1}$ denote the regular neighbourhood of the previous theorem. Let
$X$ be any almost invariant subset of $G$ over a VPC subgroup $H$ of length
$(n+1)$. Then either $X$ represents an element of $\mathcal{F}_{n,n+1}$, and
so is enclosed by some $V_{0}$-vertex of $\Gamma_{n,n+1}$, or some subgroup of
finite index in $H$ is conjugate into the vertex group of a $V_{0}$-vertex of
$\Gamma_{n,n+1}$ which is of large commensuriser type and level $n$. In the
second case, there is a VPC subgroup $A$ of $G$ of length $n$, such that
$e(G,A)\geq2$, and a $V_{0}$-vertex $v$ of $\Gamma_{n,n+1}$ which is of large
commensuriser type such that $G(v)=Comm_{G}(A)$, and some subgroup of finite
index in $H$ is conjugate to a subgroup of $G(v)$ which contains $A$.
\end{proposition}

\begin{proof}
If $X$ does not cross any almost invariant subset $Y$ over a VPC subgroup $K$
of length $n$, then $X$ is $n$-canonical, and so represents an element of
$\mathcal{F}_{n,n+1}$. Otherwise $X$ crosses such a set $Y$. If $X$ crosses
$Y$ strongly, then the first paragraph of the proof of Proposition
\ref{twocoends2} shows that $H\cap K$ must have length $n$ and hence be of
finite index in $K$. Now Lemma \ref{VPCofcolengthonehaslargenormaliser} tells
us that a subgroup of finite index in $H$ commensurises $K$. It follows that
$K\;$has large commensuriser so that $\Gamma_{n,n+1}$ has a $V_{0}$-vertex
group which equals $Comm_{G}(K)$, and so contains a subgroup of finite index
in $H$ as required. We can take the group $A$ to be $H\cap K$.

Now suppose that $X$ crosses $Y$ weakly. If $Y$ also crosses $X$ weakly, then
the first paragraph of the proof of Proposition \ref{commensurable2} shows
that one of $X^{(\ast)}\cap Y^{(\ast)}$, say $W$, is almost invariant over
$H\cap K$. As $X$ and $Y$ cross, $W$ will be a nontrivial almost invariant set
over $H\cap K$. It follows that the length of $H\cap K$ cannot be less than
$n$, since $G$ does not have any nontrivial almost invariant subsets over VPC
subgroups of length less than $n$. Thus $H\cap K$ has length $n$, and we can
apply the arguments in the preceding paragraph. Note that $X$ need not cross
$W$. We simply need the fact that $H$ contains $H\cap K$.

The only remaining case is when $X$ crosses $Y$ weakly and $Y$ crosses $X$
strongly. In this case, let $L=H\cap K$. By replacing $K$ by a subgroup of
finite index, we may assume that $L$ is normal in $K$ and that $L\backslash K$
is infinite cyclic. Since $X$ crosses $Y$ weakly one of $\delta X\cap
Y^{(\ast)}$ is $K$-finite. We will assume that $\delta X\cap Y$ is $K$-finite.
Now by again replacing $K$ by a subgroup of finite index, we may assume that
for a generator $k$ of $L\backslash K$, the translates of $\delta X\cap Y$ by
the powers of $k$ do not intersect. We choose $k$ so that $k(\delta X\cap
Y)\subset X\cap Y$, and consider the set $Z=X\cap Y\cap kX^{\ast}$. This set
is almost invariant over $L$. As $L$ has length $<n$, any $L$-almost invariant
subset of $G$ is trivial. Thus $Z$ is $L$-finite. In particular, $Z$ lies
within a finite distance of $\delta Y$. As $\cup_{i\geq1}k^{i}Z=X\cap Y$ and
$k$ preserves $\delta Y$, it follows that $X\cap Y$ also lies within a finite
distance of $\delta Y$, contradicting the hypothesis that $X$ and $Y$ cross.
This contradiction completes the proof.
\end{proof}

Finally, as in sections \ref{JSJforlargecommensurisers} and
\ref{JSJforVPCofgivenrank}, it follows that one can also form a regular
neighbourhood of only those almost invariant subsets which are associated to
splittings. This is the result we obtain.

\begin{theorem}
Let $G$ be a one-ended, finitely presented group which does not split over VPC
subgroups of length $<n$, and let $\mathcal{S}_{n,n+1}$ denote the collection
of equivalence classes of all almost invariant subsets which are associated to
a splitting of $G$ over a VPC subgroup of length $n$, together with the
equivalence classes of all $n$-canonical almost invariant subsets which are
associated to a splitting of $G$ over a VPC subgroup of length $n+1$.

Then the regular neighbourhood construction of section
\ref{regnbhds:construction} works and yields a regular neighbourhood
$\Gamma(\mathcal{S}_{n,n+1}:G)$.

Each $V_{0}$-vertex $v$ of $\Gamma(\mathcal{S}_{n,n+1}:G)$ satisfies one of
the following conditions:

\begin{enumerate}
\item $v$ is isolated, so that $G(v)$ is VPC of length $n$ or $n+1$.

\item $v$ is of VPC-by-Fuchsian type of length $n-1$ or $n$.

\item $G(v)$ contains a VPC subgroup $H$ of length $n$ or $n+1$, which it
commensurises, such that $e(G,H)\geq2$.
\end{enumerate}

If $\Gamma(\mathcal{S}_{n,n+1}:G)$ consists of a single vertex, then either
$\mathcal{S}_{n,n+1}$ is empty, or $G$ itself satisfies one of the above three conditions.
\end{theorem}

\section{Canonical splittings over virtually abelian
groups\label{JSJforabeliangroupsofallranks}}

Here is a summary of what we achieved in the last section. Consider a
one-ended, finitely presented group $G$ which does not have any nontrivial
almost invariant subsets over VPC groups of length $<n$. Let $\Gamma_{n}$
denote the regular neighbourhood of $\mathcal{F}_{n}$, the equivalence classes
of all nontrivial almost invariant subsets of $G$ over VPC subgroups of length
$n$, and let $\mathcal{F}_{n,n+1}$ denote $\mathcal{F}_{n}$ together with the
equivalence classes of all nontrivial $n$-canonical almost invariant subsets
of $G$ which are over a VPC subgroup of length $n+1$. We showed that there is
a regular neighbourhood $\Gamma_{n,n+1}$ of $\mathcal{F}_{n,n+1}$ which is a
refinement of $\Gamma_{n}$ obtained by splitting $\Gamma_{n}$ at some of its
$V_{1}$-vertices.

The natural next step would be to let $\mathcal{F}_{n,n+1,n+2}$ denote
$\mathcal{F}_{n,n+1}$ together with the equivalence classes of all nontrivial
$(n+1)$-canonical almost invariant subsets of $G$ which are over a VPC
subgroup of length $n+2$, and show that $\mathcal{F}_{n,n+1,n+2}$ has a
regular neighbourhood $\Gamma_{n,n+1,n+2}$ which is a refinement of
$\Gamma_{n,n+1}$ obtained by splitting $\Gamma_{n,n+1}$ at some of its $V_{1}%
$-vertices. However, the following example for the case $n=1$ shows that this
cannot be done following the pattern of the previous results. On the other
hand, we will show in this section, that such refinements always exist if we
restrict our attention to almost invariant sets over virtually abelian groups,
and that the process can be repeated up to any given rank. This seems to
indicate that there may be geometric differences between splittings over VPC
groups and virtually abelian groups.

\begin{example}
\label{cannotdothreestageswithVPC}This is an example of a one-ended group $G$
with incommensurable polycyclic subgroups $H$ and $K$ of length $3$, and
$2$-canonical almost invariant sets $X$ and $Y$ over $H$ and $K$ respectively
which cross weakly. Thus Proposition \ref{commensurable2} cannot be
strengthened and there is no hope of enclosing $X$ and $Y$ in a $V_{0}$-vertex
group with stabiliser equal to the commensuriser of $H$ or $K$.

We start with an extension of $\mathbb{Z}\times\mathbb{Z}$ by $\mathbb{Z}$
which is given by an automorphism of $\mathbb{Z}\times\mathbb{Z}$ with no real
eigenvalues. This gives us a polycyclic group $H$ and we denote a lift of
$\mathbb{Z}$ into $H$ by $C_{1}$. Note that $H$ is the fundamental group of a
closed $3$-manifold $M$ which is a bundle over the circle with fibre the
torus. Our choice of $H$ implies that any polycyclic subgroup of length $2$ is
contained in the normal $\mathbb{Z}\times\mathbb{Z}$. We let $K$ denote a
second copy of $H$ and let $C_{2}$ denote the subgroup of $K$ corresponding to
$C_{1}$. Let $L$ denote the fundamental group of a hyperbolic surface $F$ with
one boundary component and denote the subgroup corresponding to $\partial F$
by $C_{3}$. Now we amalgamate $H$, $K$ and $L$ along the $C_{i}$'s to obtain
the desired group $G$, and denote by $C$ the identified copies of the $C_{i}%
$'s. Thus $G$ is the fundamental group of a space $Z$ which is the union of
two copies of $M$ and the surface $F$. Let $\Gamma$ denote the associated
graph of groups structure for $G$, which is a tree with four vertices carrying
the subgroups $C$, $H$, $K$ and $L$. Consider the subgroups $H\ast_{C}K$,
$H\ast_{C}L$, $K\ast_{C}L$ of $G$. Then $G$ can be obtained from the first two
groups by amalgamating over $H$. We let $X$ denote one of the standard
$H$-almost invariant subsets of $G$ associated to this splitting. Similarly,
the first and third groups give an amalgamated free product decomposition of
$G$ over $K$. We let $Y$ denote the corresponding $K$-almost invariant subset
of $G$.

Clearly $H$ and $K$ are not commensurable in $G$. Also it is clear that the
above splittings of $G$ over $H$ and $K$ are not compatible, so that $X$ and
$Y$ must cross. As $H\cap K=C$, and $e(H,C)=e(K,C)=1$, the splittings cannot
cross strongly, so that $X$ and $Y$ must cross weakly. It remains to show that
$X$ and $Y$ are $2$-canonical. We will do this by showing that they are
$1$-canonical and that $G$ has no nontrivial almost invariant subsets over any
VPC subgroups of length $2$.

We claim that if $W$ is a nontrivial almost invariant subset of $G$ over a
two-ended subgroup $A$, then $W$ is enclosed by the vertex of $\Gamma$ which
carries $L$. This can be seen by simply considering the covering space of $Z$
corresponding to a two-ended subgroup $A$. (It follows that the regular
neighbourhood of all the almost invariant subsets of $G$ over two-ended
subgroups has a single $V_{0}$-vertex of finite-by-Fuchsian type with
associated group $L$, has no $V_{0}$-vertices of commensuriser type and has
three isolated $V_{0}$-vertex groups which carry $C$. Collapsing the edges
which carry $C$ will yield $\Gamma$.)

Now any length two polycyclic subgroup of $G$ is conjugate into $H$ or $K$,
and any such subgroup of $H$ or $K$ is contained in the normal $\mathbb{Z}%
\times\mathbb{Z}$ subgroup of $H$ or of $K$. It is now easy to check, by
considering the covering space of $Z$ corresponding to the fibre torus of a
copy of $M$, that $G$ has no nontrivial almost invariant subsets over VPC
subgroups of length $2$. It follows that the splittings of $G$ over $H$ and
$K$ that we considered above are $2$-canonical, as required.
\end{example}

The above example shows that the process of refining our algebraic analogues
of the JSJ-decomposition is not possible over VPC groups for more than two
successive ranks even if we take $i$-canonical sets over VPC groups of length
$(i+1)$ at each stage. However, it is possible for virtually abelian groups
and we will indicate the necessary changes to the arguments.

The crucial properties we needed to obtain canonical decompositions in the
previous section were contained in the following two propositions which we
reproduce here for the reader's convenience.

\bigskip

\textbf{Proposition \ref{twocoends2}} \textit{Let }$G$\textit{ be a one-ended,
finitely generated group, and let }$X$\textit{ and }$Y$\textit{ be }%
$n$\textit{-canonical subsets of }$G$\textit{ over VPC subgroups }$H$\textit{
and }$K$\textit{ of length }$(n+1)$\textit{. If }$X$\textit{ crosses }%
$Y$\textit{ strongly, then }$Y$\textit{ crosses }$X$\textit{ strongly and the
number of coends in }$G$\textit{ of both }$H$\textit{ and }$K$\textit{ is }$2$\textit{.}

\bigskip

\textbf{Proposition \ref{commensurable2}} \textit{Let }$G$\textit{ be a
one-ended, finitely generated group without nontrivial almost invariant
subsets over VPC groups of length }$<n$\textit{. Let }$H$\textit{ and }%
$K$\textit{ be VPC\ subgroups of }$G$\textit{ of length }$(n+1)$\textit{, and
let }$X$\textit{ and }$Y$\textit{ be nontrivial }$n$\textit{-canonical subsets
of }$G$\textit{ over }$H$\textit{ and }$K$\textit{ respectively. Suppose that
}$X$\textit{ crosses }$Y$\textit{ weakly. Then }$H$\textit{ and }$K$\textit{
are commensurable.}

\bigskip

Note that in Proposition \textbf{\ref{twocoends2}} we only needed the almost
invariant sets to be $n$-canonical whereas in Proposition
\textbf{\ref{commensurable2}} we excluded the existence of nontrivial almost
invariant sets over VPC groups of length $<n$. The example at the end of the
previous section showed that the analogue of Proposition
\textbf{\ref{commensurable2}} is not true in general. However, the following
analogue holds when we restrict our attention to virtually abelian subgroups
of $G$.

\begin{proposition}
\textit{Let }$G$\textit{ be a one-ended, finitely generated group, and }let
$H$ and $K$ be virtually abelian\ subgroups of $G$ of rank $(n+1)$\textit{.
L}et $X$ and $Y$ be almost invariant subsets of $G$ over $H$ and $K$
respectively which are $n$-canonical with respect to abelian groups. Suppose
that $X$ crosses $Y$ weakly. Then $H$ and $K$ are commensurable.
\end{proposition}

\begin{proof}
Our argument is based on the proof of Proposition \ref{commensurable2}. As in
the first part of that proof, we know that $X$ and $Y$ cross each other weakly
and that $X\cap Y$ is $(H\cap K)$-almost invariant.

If $H$ and $K$ are not commensurable, then $H\cap K$ has infinite index in
both $H$ and $K$. In particular, $H\cap K$ has rank $\leq n$. As $H$ is
virtually abelian, there is $h$ in $H$ of infinite order which commutes with a
subgroup of $H\cap K$ of finite index. Thus $h$ commensurises $H\cap K$.
Further, by replacing $h$ by a suitable power, we can arrange that $h(X\cap
Y)$ is contained in $X\cap Y^{\ast}$, as in the proof of Proposition
\ref{twocoends2}. Thus $(X\cap Y)\cup h(X\cap Y)$ is an almost invariant
subset of $G$ over a subgroup of $H$ commensurable with $H\cap K$ which
crosses $Y$. This contradicts the assumption that $Y$ is $n$-canonical with
respect to abelian groups. Hence $H$ and $K$ are commensurable as required.
\end{proof}

Note that in the proof of Proposition \ref{commensurable2}, we proceeded
essentially as above in the case when $H\cap K$ had length $n$, but we
eliminated the possibility that $H\cap K$ had length $<n$ by using the
assumption that $G$ had no nontrivial almost invariant subsets over VPC
subgroups of length $<n$. In the case above, $G$ may have such subsets.
Instead we used the assumption that $H$ and $K\;$are virtually abelian, and
applied the same argument as when $H\cap K$ had length $n$.

With this proposition available, there is no difficulty in extending the main
decomposition theorems to almost invariant sets over virtually abelian groups
up to any length which are canonical with respect to abelian groups. We will
need the following definitions. In this section, it will be convenient to use
the notation VA for a virtually abelian group of finite rank.

\begin{definition}
Let $\Gamma$ be a minimal graph of groups decomposition of a group $G$. A
vertex $v$ of $\Gamma$ is of \textsl{VA-by-Fuchsian type} if $G(v)$ is an
extension of a VA group by a Fuchsian group, where the Fuchsian group is not
finite nor two-ended, and there is exactly one edge of $\Gamma$ which is
incident to $v$ for each peripheral subgroup $K$ of $G(v)$ and this edge
carries $K$. If the rank of the normal VA subgroup of $G(v)$ is $n$, we will
say that $G(v)$ is of \textsl{rank} $n$.
\end{definition}

Note that if $G=G(v)$, then the Fuchsian quotient group corresponds to a
closed orbifold.

\begin{definition}
Let $\Gamma$ be a minimal graph of groups decomposition of a group $G$. A
vertex $v$ of $\Gamma$ is $n$\textsl{-simple for abelian groups}, if whenever
$X$ is a nontrivial almost invariant subset of $G$ over a VA subgroup of rank
at most $n$ such that $X$ is enclosed by $v$, then $X$ is associated to an
edge splitting of $\Gamma$.
\end{definition}

The results we obtain follow. The proof consists of starting with the graph of
groups structure $\Gamma_{1,2}$ described in the previous section, and then
using the methods of that section to repeatedly refine it by splitting at
$V_{1}$-vertices. As usual, we state the existence result and then list the
properties of decomposition obtained.

\begin{theorem}
\label{regnbhdofaisetsoverabeliangpsoflengthlessthann}Let $G$ be a one-ended,
finitely presented group. Let $\mathcal{F}_{1,2,\ldots,n}$ denote the
collection of equivalence classes of all nontrivial almost invariant subsets
of $G$ which are over a virtually abelian subgroup of rank $i$, for $1\leq
i\leq n$, and are $(i-1)$-canonical with respect to abelian groups.

Then the regular neighbourhood construction of section
\ref{regnbhds:construction} works and yields a regular neighbourhood
$\Gamma_{1,2,...,n}=\Gamma(\mathcal{F}_{1,2,\ldots,n}:G)$.

Each $V_{0}$-vertex $v$ of $\Gamma_{1,2,...,n}$ satisfies one of the following conditions:

\begin{enumerate}
\item $v$ is isolated, so that $G(v)$ is VA of rank $\leq n$.

\item $G(v)$ is of VA-by-Fuchsian type of rank $k$, for some $k$ such that
$1\leq k\leq n-1$.

\item $G(v)$ is the full commensuriser $Comm_{G}(H)$ for some VA subgroup $H$
of rank at most $n$, such that $e(G,H)\geq2$.
\end{enumerate}

$\Gamma_{1,2,...,n}$ consists of a single vertex if and only if $\mathcal{F}%
_{1,2,\ldots,n}$ is empty, or $G$ itself satisfies one of the above three conditions.
\end{theorem}

Now we list the properties of $\Gamma_{1,2,...,n}$.

\begin{theorem}
Let $G$ be a one-ended, finitely presented group. Let $\mathcal{F}%
_{1,2,\ldots,n}$ denote the collection of equivalence classes of all
nontrivial almost invariant subsets of $G$ which are over a virtually abelian
subgroup of rank $i$, for $1\leq i\leq n$, and are $(i-1)$-canonical with
respect to abelian groups.

Then the regular neighbourhood $\Gamma_{1,2,...,n}=\Gamma(\mathcal{F}%
_{1,2,\ldots,n}:G)$ is a minimal bipartite graph of groups decomposition of
$G$ with the following properties:

\begin{enumerate}
\item each $V_{0}$-vertex $v$ of $\Gamma_{1,2,...,n}$ satisfies one of the
following conditions:

\begin{enumerate}
\item $v$ is isolated, so that $G(v)$ is VA of rank $\leq n$.

\item $G(v)$ is of VA-by-Fuchsian type of rank $k$, for some $k$ such that
$1\leq k\leq n-1$.

\item $G(v)$ is the full commensuriser $Comm_{G}(H)$ for some VA subgroup $H$
of rank at most $n$, such that $e(G,H)\geq2$.
\end{enumerate}

Further if $H$ is a VA subgroup of $G$ of rank $k\leq n$ such that
$e(G,H)\geq2$, if $H$ has large commensuriser and if there exists a nontrivial
$H$-almost invariant subset of $G$ which is $(k-1)$-canonical with respect to
abelian groups, then $\Gamma_{1,2,...,n}$ has a $V_{0}$-vertex $v$ such that
$G(v)=Comm_{G}(H)$.

\item If an edge of $\Gamma_{1,2,...,n}$ is incident to a $V_{0}$-vertex of
type a) or b) above, then it carries a VA group of some rank at most $n$.

\item any representative of an element of $\mathcal{F}_{1,2,\ldots,n}$ is
enclosed by some $V_{0}$-vertex of $\Gamma_{1,2,...,n}$, and each $V_{0}%
$-vertex of $\Gamma_{1,2,...,n}$ encloses such a subset of $G$. In particular,
if $1\leq i\leq n$, then any $(i-1)$-canonical splitting of $G$ over a VA
subgroup of rank $i$ is enclosed by some $V_{0}$-vertex of $\Gamma
_{1,2,...,n}$.

\item if $X$ is an almost invariant subset of $G$ over a finitely generated
subgroup $H$, and if $X$ does not cross any element of $\mathcal{F}%
_{1,2,\ldots,n}$, then $X$ is enclosed by a $V_{1}$-vertex of $\Gamma
_{1,2,...,n}$.

\item if $X$ is a $H$-almost invariant subset of $G$ associated to a splitting
of $G$ over $H$, and if $X$ does not cross any element of $\mathcal{F}%
_{1,2,\ldots,n}$, then $X$ is enclosed by a $V_{1}$-vertex of $\Gamma
_{1,2,...,n}$.

\item the $V_{1}$-vertex groups of $\Gamma_{1,2,...,n}$ are $n$-simple. In
particular, $\Gamma_{1,2,...,n}$ cannot be further refined by splitting at a
$V_{1}$-vertex along a virtually abelian subgroup of rank at most $n$.

\item If $\Gamma_{1}$ and $\Gamma_{2}$ are minimal bipartite graphs of groups
structures for $G$ which satisfy conditions 3 and 5 above, then they are
isomorphic provided there is a one-to-one correspondence between their
isolated $V_{0}$-vertices.

\item The graph of groups $\Gamma_{1,2,...,n}$ is invariant under the
automorphisms of $G$.

\item For $k\leq n$, the splittings of $G$ over a VA subgroup of rank $k$
which are $k$-canonical with respect to abelian groups, are precisely those
edge splittings of $\Gamma_{1,2,...,n}$ which are over such a subgroup. This
includes, but need not be limited to, all those edges of $\Gamma_{1,2,...,n}$
which are incident to $V_{0}$-vertices whose associated groups are VA of rank
$k$ or of VA-by-Fuchsian type of rank $k-1$.
\end{enumerate}
\end{theorem}

Recall from Example \ref{canonicaltorusmaynotbealgebraicallycanonical} that if
there is a regular neighbourhood of all the almost invariant subsets of $G$
which are over VA subgroups of length $1$ or $2$, then it cannot be a
refinement of $\Gamma_{1}$. Thus there may be almost invariant subsets of $G$
which are over VA subgroups of length $2$ but are not enclosed by any $V_{0}%
$-vertex of $\Gamma_{1,2}$. There are similar examples for higher rank groups.
However, we prove below that, as in Lemma \ref{groupsarealmostenclosed}, the
stabiliser of any almost invariant subset of $G$ which is over a VA subgroup
of rank at most $n$ must be `almost' conjugate into a $V_{0}$-vertex of
$\Gamma_{1,2,...,n}$. As for that lemma, it will be helpful to have an idea of
the level of a $V_{0}$-vertex. A $V_{0}$-vertex of $\Gamma_{1,2,...,n}$ has
level $n$ if it is sent to a $V_{1}$-vertex of $\Gamma_{1,2,...,n-1}$ by the
refinement projection, and has level $<n$ otherwise. This allows an inductive
definition of the level of any $V_{0}$-vertex of $\Gamma_{1,2,...,n}$.

\begin{proposition}
Let $G$ be a one-ended, finitely presented group and let $\Gamma_{1,2,...,n}$
denote the regular neighbourhood of the previous theorem. Let $X$ be any
almost invariant subset of $G$ over a virtually abelian subgroup $H$ of rank
$l+1\leq n$. Then either $X$ represents an element of $\mathcal{F}%
_{1,2,\ldots,n}$, and so is enclosed by some $V_{0}$-vertex of $\Gamma
_{1,2,...,n}$, or some subgroup of finite index in $H$ is conjugate into the
vertex group of a $V_{0}$-vertex of $\Gamma_{1,2,...,n}$ which is of large
commensuriser type and level $<n$. In the second case, there is an abelian
subgroup $A$ of $G$ of rank $<n$, such that $e(G,A)\geq2$, and a $V_{0}%
$-vertex $v$ of $\Gamma_{1,2,...,n}$ which is of large commensuriser type such
that $G(v)=Comm_{G}(A)$, and some subgroup of finite index in $H$ is conjugate
to a subgroup of $G(v)$ which contains $A$.
\end{proposition}

\begin{proof}
The proof is based on the proof of Proposition \ref{groupsarealmostenclosed}.
We will argue by induction on $l$. The induction starts when $l=1$, and this
is just the case $n=1$ of Proposition \ref{groupsarealmostenclosed}. It
suffices to consider only the case when $H\;$is abelian, so we will assume
this during our proof.

Now we will assume that the proposition holds for almost invariant subsets of
$G$ which are over a VA subgroup of rank at most $l$. If $X$ does not cross
any almost invariant subset of $G$ which is over a VA subgroup of rank at most
$l$, then $X$ is $l$-canonical and so lies in $\mathcal{F}_{1,2,\ldots,n}$.
Otherwise $X$ crosses some $K$-almost invariant subset $Y$ of $G$, where $K$
is VA of rank $\leq l$.

If $X$ crosses $Y$ strongly, then $H\cap K$ must have rank $l$, and hence $K$
has rank $l$ and $H\cap K$ has finite index in $K$. Now we apply our induction
hypothesis. If $Y$ lies in $\mathcal{F}_{1,2,\ldots,n}$, we use the fact that,
as $H$ is abelian, it centralises $H\cap K$. Thus $H\cap K$ has large
commensuriser, and there is a $V_{0}$-vertex $w$ of $\Gamma_{1,2,...,n}$ which
is of large commensuriser type such that $G(w)=Comm_{G}(H\cap K)$. As $H\cap
K$ has rank $l<n$, this proves the required result about $H$ in this case.
Otherwise, after simplifying by a conjugation, there is an abelian subgroup
$A$ of $G$ of rank $<l$, such that $e(G,A)\geq2$, and a $V_{0}$-vertex $v$ of
$\Gamma_{1,2,...,n}$ which is of large commensuriser type such that
$G(v)=Comm_{G}(A)$, and some subgroup of finite index in $K$ is a subgroup of
$G(v)$ which contains $A$. Now we simply note that $H$ centralises $H\cap K$
which has a subgroup of finite index which contains a subgroup of $A$ of
finite index. Thus $H$ commensurises $A$, and so lies in $G(v)$.

If $X$ and $Y$ cross each other weakly, then the first paragraph of the proof
of Proposition \ref{commensurable2} shows that one of $X^{(\ast)}\cap
Y^{(\ast)}$, call it $W$, is almost invariant over $H\cap K$. As $X$ and $Y$
cross, $W$ must be a nontrivial $(H\cap K)$-almost invariant subset of $G$. As
$G$ is one-ended, $H\cap K$ must have rank at least $1$. Now we apply the
argument of the preceding paragraph with $W$ in place of $Y$, to prove the
required result. Note that $X$ need not cross $W$. We simply need the fact
that $H$ contains $H\cap K$.

Finally suppose that $X$ crosses $Y$ weakly but $Y$ crosses $X$ strongly. As
in the proof of Proposition \ref{groupsarealmostenclosed}, we let $L=H\cap K$,
and replace $K$ by a subgroup of finite index to arrange that $L$ is normal in
$K$ and that $L\backslash K$ is infinite cyclic. Since $X$ crosses $Y$ weakly,
we can assume that $\delta X\cap Y$ is finite. By again replacing $K$ by a
subgroup of finite index, we may assume that for a generator $k$ of
$L\backslash K$, the translates of $\delta X\cap Y$ by the powers of $k$ do
not intersect. Consider the set $Z=X\cap Y\cap kX^{\ast}$. This set is almost
invariant over $L$. If $Z$ is $L$-finite, then it lies with in a finite
distance of $\delta Y$. As $\cup_{i\geq1}k^{i}Z=X\cap Y$, it follows that
$X\cap Y$ lies within a finite distance of $\delta Y$, contradicting the
hypothesis that $X$ and $Y$ cross. If $Z$ is not $L$-finite, then $X$ crosses
$Z\cup k^{-1}Z$, and as $L$ has lesser rank than $K$, this contradicts our
choice of $Y$ so as to minimise the rank of $K$.
\end{proof}

We end this section by discussing the behaviour of the sequence $\Gamma
_{1,2,...,n}$ of graphs of groups structures for a fixed group $G$ as $n$
increases. For brevity, we will denote $\Gamma_{1,2,...,n}$ by $\Gamma^{n}$ in
this paragraph only. Each $\Gamma^{n}$ is a refinement of $\Gamma^{n-1}$. We
would like to consider whether this sequence stabilises by applying Theorem
\ref{graphsstabilise2}, but this theorem does not apply because only those
edge groups incident to $V_{0}$-vertices of types a) or b) need to be VPC.
Instead we argue as follows. Each $V_{0}$-vertex of $\Gamma^{n}$ which is not
a vertex of $\Gamma^{n-1}$ encloses at least one splitting over a VA subgroup
of rank $n$. By picking one such splitting each time $\Gamma^{n}$ and
$\Gamma^{n-1}$ are distinct, we obtain a sequence $\sigma_{k}$ of compatible
splittings of $G$ over VA subgroups of rank at least $k$. This yields a
sequence of graphs of groups structures $\Delta_{k}$ for $G$ whose edge
splittings are precisely $\sigma_{1},\ldots,\sigma_{k}$. If the sequence
$\Gamma^{n}$ does not stabilise, the sequence $\Delta_{k}$ will be infinite.
This will not contradict Theorem \ref{graphsstabilise2}, but this can only
occur if there is a subsequence of the $\sigma_{i}$'s, say $\tau_{j}$, such
that $\tau_{j}$ is a splitting of $G$ over a VA subgroup $C_{j}$ of rank at
least $j$ such that $C_{j}\subset C_{j+1}$, for all $j$. Further the sequence
$\Gamma^{n}$ stabilises apart from such subsequences of edge splittings, so
there can only be finitely many such subsequences. In the case when there are
no such sequences, then the sequence $\Gamma_{1,2,...,n}$ of graphs of groups
decompositions of $G$ eventually stabilises yielding a decomposition which we
denote by $\Gamma_{\infty}$, whose $V_{0}$-vertices enclose all almost
invariant subsets of $G$ over any finitely generated subgroup which is
virtually abelian.

As in previous sections, one can also form a regular neighbourhood of only
those almost invariant subsets of $G$ which are associated to splittings. We
omit the formal statement.

\section{Previous decompositions over two-ended
subgroups\label{non-canonicaldecompositions}}

In this section, we will discuss the relationship between the
JSJ-decompositions of previous authors and the canonical decomposition which
we described in section \ref{JSJforlargecommensurisers}. Let $G$ be a
one-ended, finitely presented group. Our decomposition $\Gamma$ is a regular
neighbourhood of all equivalence classes of nontrivial almost invariant
subsets of $G$ which are over two-ended subgroups. It is clear that this is
not usually the same as any of the JSJ-decompositions of previous authors
because it may have edge groups which are not two-ended. Such edge groups can
only occur for edges which are incident to $V_{0}$-vertices which are of large
commensuriser type, and we will now describe how to alter $\Gamma$ so as to
obtain one of these other decompositions. Recall that each $V_{0}$-vertex $v$
of $\Gamma$ which is of large commensuriser type encloses at least one
splitting of $G$ over a two-ended subgroup. Such a splitting may not split the
vertex group $G(v)$, but in this case it must be an edge splitting for an edge
which is incident to $v$. For each such $V_{0}$-vertex $v$ of $\Gamma$, we
pick a maximal family of compatible splittings of $G$ which are each over a
two-ended group and are enclosed by $v$. This is possible by the accessibility
result in Theorem \ref{graphsstabilise}. We refine $\Gamma$ by splitting at
each such vertex using these splittings. The resulting graph of groups
structure is no longer canonical, as the splittings enclosed by $v$ are not
usually unique. Next we simply collapse each edge of $\Gamma$ which carries a
group which is not two-ended. The result is a graph of groups structure
$\Gamma^{\prime}$ for $G$ in which every edge group is two-ended. In
particular, it follows that every vertex group of $\Gamma^{\prime}$ is
finitely generated. Of course, $\Gamma^{\prime}$ is no longer bipartite.
Further, it is not true that any nontrivial almost invariant subset of $G$
over a two-ended subgroup is enclosed by a vertex of $\Gamma^{\prime}$.
However it follows from our construction of $\Gamma^{\prime}$ from $\Gamma$
that if $G$ possesses a nontrivial almost invariant subset over a two-ended
subgroup $H$, then $H$ has a subgroup of finite index which is conjugate into
some vertex group of $\Gamma^{\prime}$. The known JSJ-decompositions along
two-ended subgroups can all be refined to such a decomposition, but in
\cite{RS} there are some assumptions on unfoldedness which may somewhat
restrict the choice of splittings used to refine $\Gamma$. However these
decompositions are not canonical. We call any such decomposition of $G$ along
two-ended subgroups a non-canonical JSJ-decomposition of $G$.

Sela in \cite{S1} initiated the study of uniqueness of such decompositions up
to some moves which he called sliding, conjugation and modifying the boundary
homomorphism by a conjugation. In \cite{MF2}, Forester gave a complete
description of the uniqueness properties of these decompositions. In
\cite{MF}, Forester considered two moves on graphs of groups or equivalently
$G$-trees called a \textit{collapse move} and an \textit{expansion move}. The
first move involves selecting an edge $s$ in the graph which is not a loop and
such that the induced map from the edge group $G(s)$ to the initial vertex
group $G(v)$ is an isomorphism. One then collapses $s$ down to $v$. An
expansion move is the reverse of a collapse move. A move which factors as a
composition of expansions and collapses is called an \textit{elementary
deformation}. He showed that two cocompact $G$-trees are related by an
elementary deformation if and only if they have the same elliptic subgroups.

Consider one of the non-canonical JSJ-decompositions derived as above from our
canonical decomposition $\Gamma$. By construction none of the edge splittings
is crossed strongly by any almost invariant set over a two-ended subgroup.
Thus they are elliptic with respect to any splitting of $G$ over a two-ended
group. Next consider any two $G$-trees $T_{1}$ and $T_{2}$ corresponding to
two such decompositions. The Fuchsian vertex groups are the same in both and
are thus elliptic with respect to both the trees. The other vertex groups of
$T_{1}$ and $T_{2}$ do not admit any splittings over two-ended groups relative
to the edge groups. This is because in our refinement of $\Gamma$, we used a
maximal family of compatible splittings of $G$. It follows that the elliptic
subgroups of the $G$-trees $T_{1}$ and $T_{2}$ correspond. Thus, by Forester's
Theorem, we have:

\begin{theorem}
Let $G$ be a one-ended group and let $T_{1}$ and $T_{2}$ be two $G$-trees
corresponding to non-canonical JSJ-decompositions of $G$. Then $T_{1}$ and
$T_{2}$ are related by an elementary deformation.
\end{theorem}

\begin{corollary}
Let $G$ be a one-ended group and $\sigma$ a splitting of $G$ over a two-ended
subgroup. Let $\Gamma$ be any non-canonical JSJ-decomposition of $G$. Then
either $\sigma$ is enclosed by a $V_{0}$-vertex of $\Gamma$ of Fuchsian type,
or $\sigma$ can be obtained by collapses and expansions from the edge
splittings of $\Gamma$.
\end{corollary}

\section{Extensions\label{extensions}}

One might wonder how far the techniques above can be used to obtain canonical
decompositions enclosing almost invariant sets over other classes of groups.
In our arguments, the first step was to show that in any cross-connected
component, the crossings are either all weak or all strong. Then we used
different kinds of arguments in these two cases. For CCC's in which all
crossing is weak, we proved that such a CCC must enclose a splitting, then
proved a finiteness result for the CCC and finally we needed an accessibility
result. Here are two general results of this type. The first is a
reformulation of Theorem 3.13 of \cite{SS}, and asserts the existence of a
splitting very generally.

\begin{theorem}
Suppose that $H\subset G$ are finitely generated groups and that $G$ does not
contain any nontrivial almost invariant subsets over subgroups of infinite
index in $H$. Let $X$ be a nontrivial almost invariant subset of $G$ over $H$
and suppose that the translates of $X$ do not cross each other strongly. Then
$G$ splits over a subgroup commensurable with $H$.
\end{theorem}

Our arguments in section \ref{coendswhencommensuriserissmall} of this paper
extend to show the following accessibility result.

\begin{theorem}
\label{graphsstabilise3}Suppose $\mathcal{K}$ is a class of small groups
closed under commensurability. Suppose $G$ is a finitely presented group which
does not split over a subgroup of infinite index in an element of
$\mathcal{K}$. Let $\Gamma_{k}$ be a graph of groups decomposition of $G$
without trivial vertices and with all edge groups in $\mathcal{K}$, and
suppose that for each $k$, $\Gamma_{k+1}$ is a refinement of $\Gamma_{k}$.
Then, the sequence $\Gamma_{k}$ stabilises.
\end{theorem}

These two results can be used to handle more cases of CCC's in which all
crossing is weak. However, to handle CCC's in which all crossing is strong, we
used the results and arguments of Bowditch \cite{B2}, and of Dunwoody and
Swenson \cite{D-Swenson} which in turn depend on the special structure of VPC
groups. Some of these results can be summarised in the following theorems.

\begin{theorem}
\label{summaryofsplittingsinVAcase}Let $G$ be a finitely generated group and
let $X$ be an almost invariant subset of $G$ over a virtually abelian group
$H$ of rank $(n+1)$. Suppose that $X$ is $n$-canonical with respect to abelian
groups, i.e. it does not cross any almost invariant subset over a virtually
abelian group of rank $\leq n$. Then $G$ splits over a virtually abelian group
of rank $(n+1)$. If $X$ does not cross any translate of $X$ strongly, then $G$
splits over a subgroup commensurable with $H$.
\end{theorem}

\begin{theorem}
Let $G$ be a finitely generated group and let $X$ be an almost invariant
subset of $G$ over a VPC group $H$ of length $(n+1)$. Suppose $G$ does not
have any nontrivial almost invariant sets over VPC groups of length $<n$ and
that $X$ is $n$-canonical. Then $G$ splits over a VPC group of length $(n+1)$.
If $X$ does not cross any translate of $X$ strongly, then $G$ splits over a
subgroup commensurable with $H$.
\end{theorem}

In fact, we used relative versions of these theorems. Example
\ref{cannotdothreestageswithVPC} suggests that it may not be possible to
strengthen Theorem \ref{summaryofsplittingsinVAcase}. It is possible that
similar results may be provable for special classes of slender groups
considered in \cite{D-Sageev}, but these problems are still open. The
techniques of \cite{D-Sageev} and \cite{FP} enclose splittings rather than
almost invariant sets. The technique in \cite{FP} is particularly appealing.
Their enclosing group is an instance of our regular neighbourhood
construction. It is their construction that suggested to us regular
neighbourhoods and Bowditch's use of pretrees provided us with a crucial
technique. The crossing hypotheses used in their technique are weaker than
ours, provided of course that one starts with splittings. Thus there may be
further refinements of the decompositions that we obtained if one combines
their techniques with ours. We recall that even in the case of $3$-manifolds
the canonical decompositions obtained by enclosing splittings only are
different from the standard JSJ-decompositions (see \cite{NS}). So, our work
seems to suggest that there are several possible generalisations of
JSJ-decompositions to groups. Moreover, our theories of regular neighbourhoods
and canonical splittings are very general and these may apply to almost
invariant sets over groups more general than VPC groups.

\end{document}